\documentclass[a4paper,11pt,twoside,openany]{book}
\usepackage[utf8x]{inputenc}
\usepackage{cite} 
\usepackage{amssymb,amsmath}
\usepackage{graphicx}

\usepackage{verbatim} 
\usepackage{fancyhdr}
\usepackage{calc}

\usepackage{textcomp}
\usepackage{marvosym}
\usepackage{stmaryrd}
\usepackage{pifont}
\usepackage{array}
\usepackage{url}
\usepackage{subfigure}
\usepackage{txfonts}

\fancyheadoffset[LE,RO]{\marginparsep+\marginparwidth}

\fancypagestyle{plain}{%
\fancyhead{} 
}

\DeclareMathOperator\erf{erf}
\DeclareMathOperator\erfc{erfc}
\DeclareMathOperator\Ei{Ei}
\DeclareMathOperator\Ein{Ein}
\DeclareMathOperator\Li{Li}
\DeclareMathOperator\Si{Si}
\DeclareMathOperator\Ci{Ci}
\DeclareMathOperator\arccotan{arccotan}
\begin{document}

\thispagestyle{empty}

\begin{minipage}{1.0\textwidth}
\hspace{-3mm}\includegraphics[scale=0.36]{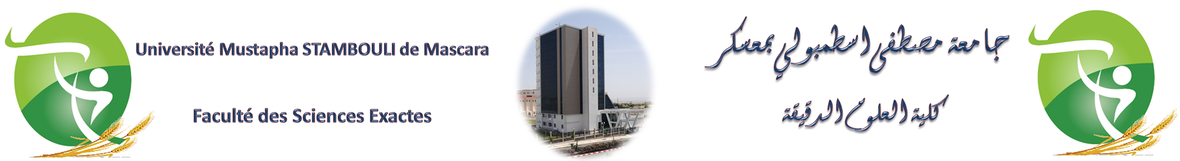}
\vspace{1cm}
\begin{center}
R\'epublique Alg\'erienne D\'emocratique et Populaire\\
Minist\`ere de l'Enseignement Sup\'erieur et de  la Recherche Scientifique\\

\vspace{1cm}
Universit\'e Mustapha Stambouli de Mascara\\
Facult\'e des sciences exactes

\vspace{1cm}
\textbf{Polycopi\'e de Cours}\\
\vspace{0.4cm}

\begin{Large}
\vspace{2cm}
\begin{huge}\textbf{Fonctions sp\'eciales et polyn\^omes orthogonaux}\end{huge}\\
\hspace{0.7cm}

\begin{normalsize}Pr\'esent\'e par\end{normalsize}
\vspace{0.5cm}

\textbf{Benaoumeur Bakhti}
\vspace{2cm}

\begin{normalsize}
 Cours destin\'e aux \'etudiants de la troisi\`eme ann\'ee licence physique
\end{normalsize}

\vspace{2cm}
\begin{Large}\textbf{Algerie 2020}\end{Large}

\end{Large}
\end{center}
\end{minipage}

\newpage
\thispagestyle{empty}

\begin{minipage}{1.0\textwidth}
\begin{center}
\vspace{8cm}
\begin{Large}
\textit{A mes parents, Meriem et Abdelkader.}
\end{Large}
\end{center}
\end{minipage}

\pagenumbering{roman}
\tableofcontents

\listoffigures
\addcontentsline{toc}{chapter}{Liste des figures}

\chapter*{Avant-propos}
\addcontentsline{toc}{chapter}{Avant-propos}\markboth{Avant-propos}{}
\label{chap:avant_propos}
Cet ouvrage est consacr\'e \`a l'\'etude des fonctions sp\'eciales les plus utilis\'ees en physique. Les fonctions sp\'eciales est une branche tr\`es vaste de math\'ematiques, de la physique th\'eorique et de la physique math\'ematique. Elles sont apparues au $xix^{eme}$ si\`ecle comme solutions d'\'equations de la physique math\'ematique, particuli\`erement les \'equations aux d\'eriv\'ees partielles d'ordre deux et quatre. Leur connaisance est indispensable \`a la bonne manipulation et compr\'ehension des probl\`emes actuels de la physique. Elles sont \'egalement li\'ees \`a l'art du calcul scientifique de la physique et des math\'ematiques. Les fonctions sp\'eciales sont incluses dans de nombreux logiciels de calcul formel tels que le Matlab, le Mathematica et le Maple et les \'etudiants sont fortement encourag\'es de prendre part \`a ce d\'eveloppement qui est devenu indispensable pour le traitement presque de tous les probl\`emes actuels de la physique. 

Conforme aux programmes LMD (Licence-Master-Doctorat), l'ouvrage contient six chapitres.

Dans le premier chapitre, nous d\'ecrivons les fonctions gamma et b\^eta qui sont importantes pour leur propre int\'er\^et math\'ematique et aussi parce que toutes les autres fonctions sp\'eciales d\'ependent essentiellement de ces deux fonctions. Ces deux fonctions ont de nombreuses applications. En physique et en particulier en th\'eorie des cordes, la fonction b\^eta (et la fonction gamma associ\'ee) est utilis\'ee pour calculer et reproduire l'amplitude de diffusion en fonction des trajectoires de Regge dans le "mod\`ele \`a double r\'esonance". En th\'eorie des probabilit\'es, elles sont utilis\'ees dans le processus d'attachement pr\'ef\'erentiel ou en g\'en\'eral dans le processus stochastique d'urne.

Le chapitre deux traite les fonctions de Bessel et leurs propri\'et\'es principales. Les fonctions de Bessel ont \'et\'e introduites et \'etudi\'ees d'abord par Euler, Lagrange et Bernoulli. Mais elles ont \'et\'e utilis\'es pour la permi\`ere fois par Friedrich Wilhelm Bessel pour expliquer le mouvement de trois corps, o\`u la fonction de Bessel a \'emerg\'ee dans le d\'eveloppement en s\'erie de la perturbation plan\'etaire. Les fonctions de Bessel sont extr\^emement utiles en physique et en ing\'enierie. Elles se sont av\'er\'es \^etre la solution de l'\'equation de Schr\"odinger dans une situation de sym\'etrie cylindrique. L'\'equation diff\'erentielle de Bessel d\'ecoule de la d\'etermination de solutions s\'eparables de l'\'equation de Laplace et de l'\'equation de Helmholtz en coordonn\'ees sph\'eriques et cylindriques. Les fonctions de Bessel sont \'egalement tr\`es importantes pour de nombreux probl\`emes de propagation des ondes, de potentiels statiques et dans la th\'eorie de la diffusion en m\'ecanique quantique. En ing\'enierie, elles sont utiles dans de nombreux probl\`emes tels que la conduction thermique, les ondes \'electromagn\'etiques dans un guide d'onde, le traitement du signal, les modes de vibration d'une membrane artificielle et dans l'acoustique. Nous allons pr\'esenter dans ce chapitre toutes les variantes des fonctions de Bessel qui sont les fonctions de Bessel de premi\`ere esp\`ece et de deuxi\`eme esp\`ece (appell\'ee aussi fonctions de Neumann), les fonctions de Bessel modifi\'ees, les fonction de Bessel sph\'eriques ainsi que les fonctions de Hankel et les fonctions de Hankel sph\'eriques. Nous donnerons des solutions d\'etaill\'ees des \'equations de Bessel en utilisant la m\'ethode de Frobenius ainsi que les d\'emonstrations de toutes leurs propri\'et\'es principales.

Les chapitres trois et quatre sont consacr\'es \`a d'autres fonctions d\'efinies par des int\'egrales
qui sont: la fonction erreur, int\'egrales de Fresnel, exponentielle int\'egrale, sinus int\'egrale, cosinus int\'egrale et logarithme int\'egrale.  Ses fonctions sont beaucoup utilis\'ees en physique, par exemple dans le domaine de l'optique et de l'\'el\'ecrtomagnetisme et dans le domaine des probabilit\'es et statistiques.

Le chapitre cinq est consacr\'e à l'\'etude des polyn\^omes orthogonaux. Nous \'etudions notamment les polyn\^omes de Legendre, d'harmoniques sph\'eriques, d'Her-mite, de Laguerre et de Chebyshev. Ces polyn\^omes sont des solutions des \'equations diff\'erentielles ordinaires d'ordre deux. Ces \'equations surviennent tr\`es souvent lorsqu'un probl\`eme poss\`ede une sym\'etrie sph\'erique. De tels probl\`emes peuvent survenir, par exemple, en m\'ecanique quantique, en th\'eorie de l'\'electromagnétisme, en hydrodynamique et en conduction thermique. En ing\'enierie, les polyn\^omes orthogonaux apparaissent dans de nombreuses applications, telles que la th\'eorie des lignes de transmission, la th\'eorie des circuits \'electriques, la physique des r\'eacteurs nucl\'eaires et en sismologie.

Dans le dernier chapitre, nous traitons en d\'etails les fonctions hyperg\'eom\'etriques. Ces dernieres ont \'et\'e introduites par Gauss comme une g\'en\'eralisation de la s\'erie g\'eom\'etrique. Nous traitons surtout les fonctions hyperg\'eom\'etriques les plus importantes qui sont les fonction hyperg\'eom\'etrique de Gauss et les fonctions hyperg\'eom\'etriques confluentes (appell\'ees aussi fonctions de Kummer). L'importance de ces fonctions est que toutes les fonctions pr\'esent\'ees pr\'ec\'edemment (Bessel, polyn\^omes orthogonaux,$\ldots$) peuvent \^etre exprim\'ees en termes de fonctions hypergom\'etriques.
Les fonctions hypergom\'etriques ont \'et\'e utilis\'ees dans une large gamme de probl\`emes en physique classique et quantique, en ing\'enierie et en math\'e-matiques appliqu\'ees. En physique, elles sont tr\`es utiles dans les probl\`emes de forces centrales, par exemple dans l'\'etude de l'atome d'hydrog\`ene et de l'oscillateur harmonique en m\'ecanique quantique. Leur int\'er\^et pour les math\'ematiques r\'eside dans le fait que de nombreuses \'equations (bien connues) aux d\'eriv\'ees partielles peuvent \^etre r\'eduites à l'\'equation hyperg\'eom\'etrique de Gauss par s\'eparation des variables.

Cet ouvrage s'adresse principalement aux \'etudiants de la troisi\`eme ann\'ee licence Sciences de la Mati\`ere (SM). Mais il sera utile \`a un cercle de lecteurs tr\`es \'etendue: \'etudiants en math\'ematiques et en sciences techniques. Il est conçu de façon \`a aplanir au mieux les difficult\'es inh\'erentes au discours scientifiques tout en conservant la rigueur n\'ecessaire. Cet ouvrage pr\'esente l'ensemble des notions de bases abord\'ees au cours "M\'ethodes Mathematiques" durant la troisi\'eme ann\'ee de Licence SM. Aussi, des exercices corrig\'es sont propos\'es \`a la fin de chaque chapitre permettant \`a l'\'etudiant de tester ses connaissances et de se pr\'eparer aux tests de controle et aux examens. Il est est le fruit de quelques ann\'ees d'enseignement du cours "M\'ethodes Math\'ematiques pour la Physique" dispens\'e au d\'epartement de Physique, facult\'e des sciences exactes de l'Universit\'e Mustapha Stambouli de Mascara.\\

Enfin, je tiens \`a remercier vivement mes amis Mohammed Elamine Sebih (Universit\'e Mustapha Stambouli de Mascara) et Mohamed Reda Chellali (Karlsruhe Institute of Technology) ansi que les examinateurs Prof. Boucif Abdesselam (Centre universitaire Ain Temouchent), Dr. Gherici Beldjilali et Dr. Abdelkader Segres (Universit\'e Mustapha Stambouli de Mascara) qui ont contribu\'es au perfectionnement de cet ouvrage par la lecture attentive du manuscrit et par leurs commentaires et propositions.\vspace{6mm}\\
\textit{Mascara, 2020}\hspace{9cm}B. Bakhti

\newpage
\pagenumbering{arabic}
\setcounter{page}{1}

\chapter{Les fonctions eul\'eriennes gamma et b\^eta}\label{chap: gamma_beta}
\section{Fonction gamma}
\subsection{D\'efinition}
La fonction gamma (not\'e $\Gamma$) a \'et\'e introduite par Euler en $1729$, elle est d\'efinie par l'int\'egrale
\begin{align}\label{eq:1_1}
\Gamma(x)=\int_0^{\infty}e^{-t}t^{x-1}dt, \hspace{10mm} x > 0
\end{align}
o\`u $x$ peut \^etre r\'eel ou complexe avec $Re(x)>0$. L'int\'egrale (\ref{eq:1_1}) appel\'e aussi l'int\'egrale d'Euler de premi\`ere esp\`ece
n'\'existe que si $x$ est strictement positif. Pour montrer cela, divisons l’intervalle d’int\'gration en deux parties: de $0$ \`a $\epsilon << 1$ et de $\epsilon$ à l'infini. On obtient
\begin{align*}
\Gamma(x)&=\int_0^{\epsilon}e^{-t}t^{x-1}dt+\int_{\epsilon}^{\infty}e^{-t}t^{x-1}dt\\
&=\frac{t^{x}}{x}\Big|_{0}^{\epsilon}+\int_{\epsilon}^{\infty}e^{-t}t^{x-1}dt\nonumber
\end{align*}
Dans le premier terme de droite, nous avons utilis\'e le fait que 
$e^{-t}t^{x-1}\simeq t^{x-1}$ pour $t<<1$. On remarque que le terme $t^{x}/x$ 
n'est fini que si $x$ est strictement positif, sinon l'int\'egrale diverge.
\subsection{Relation de r\'ecurrence}
\begin{figure}[hbt]
\centering
\includegraphics[scale=0.5]{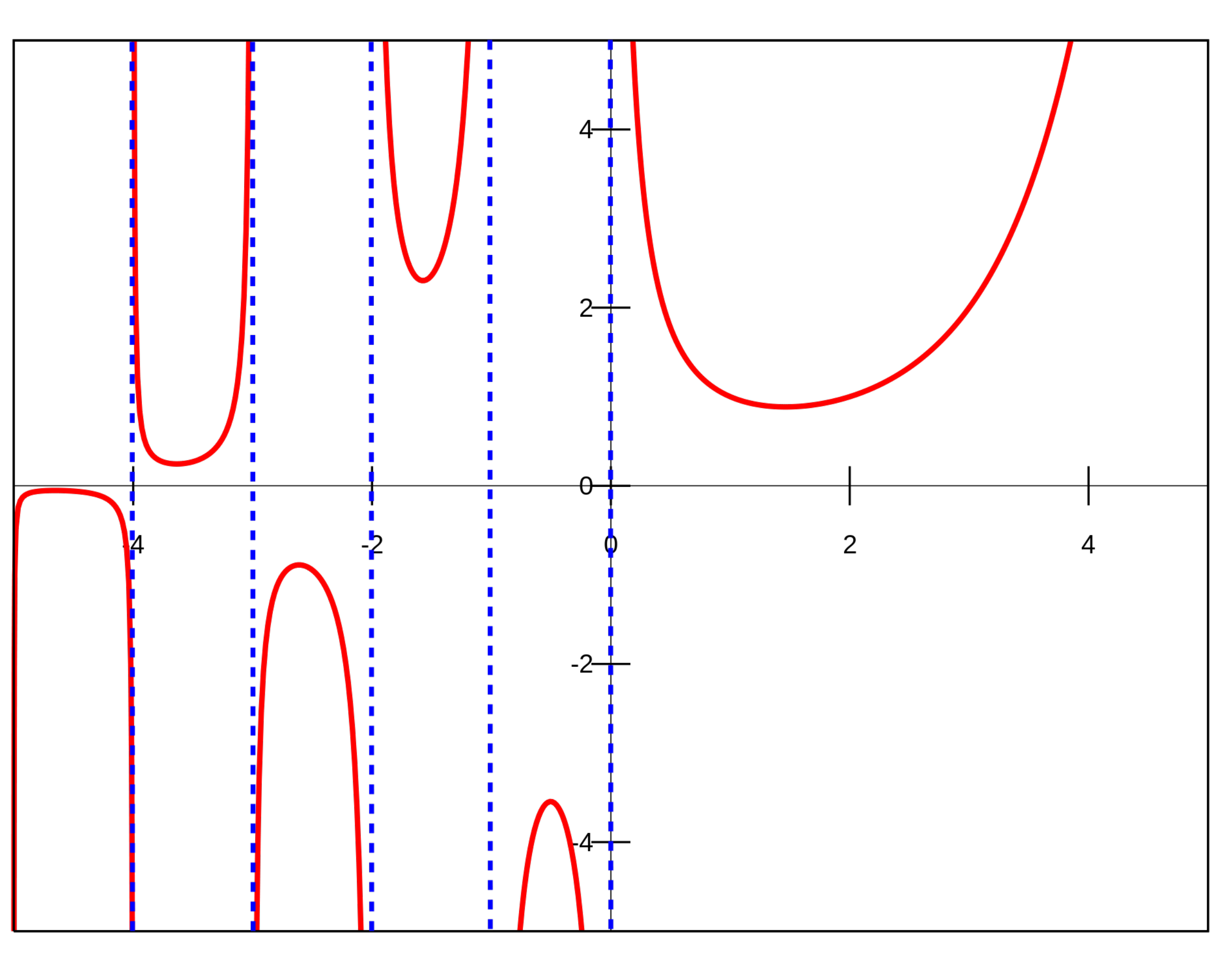}
\caption{\label{fig:gamma_plot} Trac\'e de la fonction gamma le long de l'axe des r\'eels.}
\end{figure}

La fonction gamma satisfait la relation
\begin{align}\label{eq:1_2}
\Gamma(x+1)=x\Gamma(x)
\end{align}
et si $x$ est un entier non n\'egatif, on d\'eduit que
\begin{align}\label{eq:1_3}
\Gamma(x+1)=x!,
\end{align}
\textit{D\'emonstration}\\
L'int\'gration par parties de l'int\'egrale d'Euler avec $u=t^{x}$ et $v'=e^{-t}$, donne
\begin{align*}
\Gamma(x+1)&=\int_0^{\infty}e^{-t}t^{x}dt\\
&=-e^{-t}t^{x}\Big|_{0}^{\infty}+x\int_0^{\infty}e^{-t}t^{x-1}dt\\
&=x\Gamma(x)
\end{align*}
Le terme $e^{-t}t^{x-1}$ donne $0$ \`a l'infini (\i.e. $t\rightarrow \infty$) parce que la d\'ecroissance de la fonction exponentielle $e^{-t}$ est plus rapide que la fonction $t^{x-1}$.\\
Si $x\in \mathsf{N}$, la relation (\ref{eq:1_2}) donne
\begin{align}
\Gamma(x+1)&=x\Gamma(x)\nonumber\\
 &= x\left\lbrace(x-1)\Gamma(x-1)\right\rbrace\nonumber\\
 &=x(x-1)\left\lbrace(x-2)\Gamma(x-2)\right\rbrace\nonumber\\
 &\vdots\nonumber\\
 &=x(x-1)(x-2)\ldots 3.2.1\Gamma(1)\nonumber\\
 &=x!\nonumber
\end{align}
La relation (\ref{eq:1_3}) montre que la fonction gamma prolonge la fonction factorielle \`a
l'ensemble des nombres complexes. Elle montre \'egalement que
\begin{align*}
\Gamma(1)&=0!\\
&=1
\end{align*}
qui peut \^etre d\'emontrer aussi facilement en calculant l'int\'egrale (\ref{eq:1_1}) pour $x=1$.\\
Bien que l'int\'egrale (\ref{eq:1_1}) n'\'existe que pour $Re(x)>0$, 
la fonction gamma peut \^etre prolong\'ee analytiquement sur l'ensemble des nombres complexes $x$ tel que $Re(x)<0$, except\'e pour $x = 0,  −1, −2, −3 \ldots$ qui sont des p\^oles.
Cela peut être fait en utilisant l'inverse de la relation de r\'ecurrence (\ref{eq:1_2})
\begin{align*}
\Gamma(x)=\frac{1}{x}\Gamma(x+1)
\end{align*}
par exemple
\begin{align*}
\Gamma\left(-\frac{1}{2}\right)&=\frac{1}{-1/2}\Gamma\left(-\frac{1}{2}+1\right)\\
&=-2\Gamma\left(\frac{1}{2}\right)
\end{align*}
et
\begin{align*}
\Gamma\left(-\frac{3}{2}\right)&=\frac{1}{-3/2}\Gamma\left(-\frac{1}{2}\right)\\
&=\frac{4}{3}\Gamma\left(\frac{1}{2}\right)
\end{align*}
et ainsi de suite. En plus, comme la fonction gamma diverge au point $x=0$ ($\Gamma(0)=+\infty$), elle diverge \'egalement aux points $x=-1,-2,-3,\ldots$. La fonction gamma est repr\'esent\'ee graphiquement sur la Fig.~(\ref{fig:gamma_plot})

\subsection{D'autre repr\'esentations de la fonction gamma}
La fonction gamma peut \^etre repr\'esent\'ee \'egalement par les int\'egrales suivantes
\begin{align}
&\Gamma(x)=2\int_0^{\infty}e^{-t^{2}}t^{2x-1}dt  \hspace{10mm} x > 0\label{eq:1_4}\\
&\Gamma(x)=\int_0^{1}\left[\ln\Big(\frac{1}{t}\Big)\right]^{x-1}dt  \hspace{10mm} x > 0\label{eq:1_5}
\end{align}
\textit{D\'emonstration}\\
Les formules (\ref{eq:1_4}) et (\ref{eq:1_5}) peuvent \^etre facilement d\'emontr\'ees en utilisant respectivement les changements de variables $u=t^{2}$ et $t=e^{-u}$ dans la d\'efinition (\ref{eq:1_1}).
\subsection{Relation de gamma avec fonctions trigonom\'etriques}
\begin{figure}[hbt]
\centering
\includegraphics[scale=0.2]{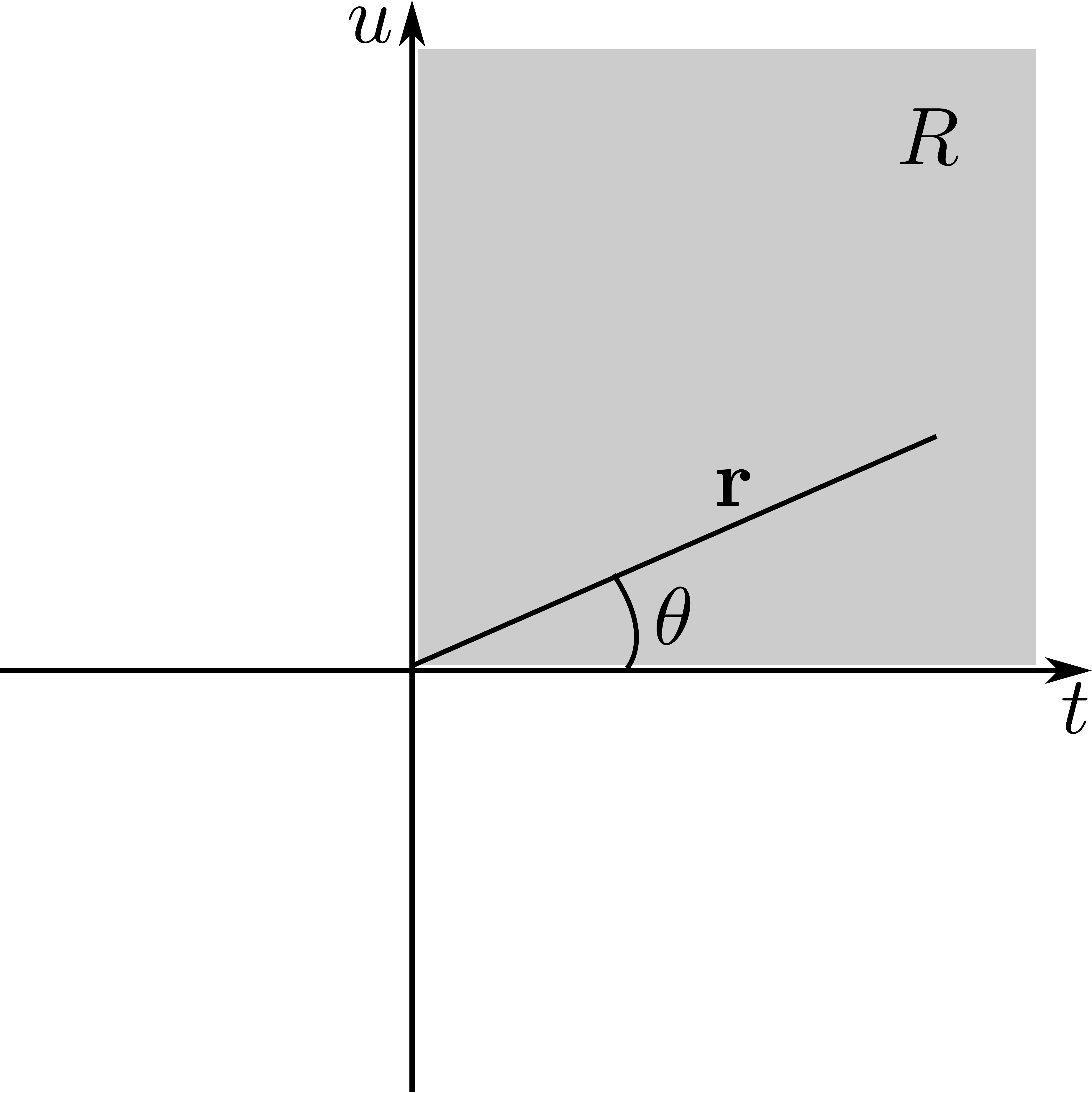}
\caption{\label{fig:tu_plan} Le plan $tu$}
\end{figure}

Une propri\'et\'es tr\`es importante de la fonction gamma est la suivante
\begin{align}\label{eq:1_6}
\int_0^{\pi/2}\cos^{2x-1}\theta\sin^{2y-1}\theta d \theta = \frac{\Gamma(x)\Gamma(y)}{2\Gamma(x+y)}
\end{align}
\textit{D\'emonstration}\\
L'Eq.~(\ref{eq:1_6}) peut \^etre prouv\'ee en \'evaluant l'int\'egrale suivante de deux mani\`eres diff\'erentes
\begin{align}\label{eq:1_7}
I=\int\int_Re^{-t^{2}-u^{2}}t^{2x-1}u^{2y-1}dtdu
\end{align}
o\`u $R$ est le premier quadrant du plan-$(t,u)$ (voir Fig.~(\ref{fig:tu_plan})).\\
D'abord, \'evaluons l'int\'egrale $I$ en utilisant l'expression Eq.~(\ref{eq:1_4}) dans les coordonn\'ees cart\'esiennes $(t,u)$, on a
\begin{align}
I&=\int_{0}^{\infty}\int_{0}^{\infty}e^{-t^{2}-u^{2}}t^{2x-1}u^{2y-1}dtdu\nonumber\\
&=\int_{0}^{\infty}e^{-t^{2}}t^{2x-1}dt \int_{0}^{\infty}e^{-u^{2}}u^{2y-1}du\nonumber\\
&=\frac{1}{2}\Gamma(x)\frac{1}{2}\Gamma(y)\nonumber\\
&=\frac{1}{4}\Gamma(x)\Gamma(y)\nonumber
\end{align}
o\`u nous avons utiliser l'Eq.~(\ref{eq:1_4}) dans la deuxi\`eme ligne.
Ensuite, \'evaluons l'int\'egrale $I$ dans les coordonn\'ees polaires en utilisant les changements de variables
\begin{align*}
t=r\cos\theta, \hspace{20mm}u=r\sin\theta
\end{align*}
nous obtenons
\begin{align}
I&=\int\int_Re^{-r^{2}\cos^{2}\theta-r^{2}\sin^{2}\theta}(r\cos\theta)^{2x-1}
(r\sin\theta)^{2y-1}rdrd\theta\nonumber\\
&=\int_{0}^{\infty}e^{-r^{2}}r^{2(x+y)-1}dr\int_0^{\pi/2}\cos^{2x-1}\theta\sin^{2y-1}\theta d\theta\nonumber\\
&=\frac{1}{2}\Gamma(x+y)\int_0^{\pi/2}\cos^{2x-1}\theta\sin^{2y-1}\theta d\theta\nonumber
\end{align}
En identifiant les deux expressions de l'int\'egrale $I$, on trouve l'\'equation d\'esir\'ee
\begin{align*}
\int_0^{\pi/2}\cos^{2x-1}\theta\sin^{2y-1}\theta d \theta = \frac{\Gamma(x)\Gamma(y)}{2\Gamma(x+y)}
\end{align*}
Si on pose $x=y=1/2$ dans l'Eq.~(\ref{eq:1_6}), il en r\'esulte que
\begin{align*}
\int_0^{\pi/2}d\theta &=\frac{\pi}{2} \\
&= \frac{\Gamma(1/2)\Gamma(1/2)}{2\Gamma(1)}
\end{align*}
donc
\begin{align*}
\Gamma\left(\frac{1}{2}\right)=\sqrt{\pi}
\end{align*}
\subsection{Formule de Weierstrass}
\begin{align}\label{eq:1_7a}
\Gamma(x)=\frac{1}{x}e^{-\gamma x}\prod_{n=1}^{\infty}e^{\frac{x}{n}}\left(1+\frac{x}{n}\right)^{-1}
\end{align}
o\`u $\gamma=0.5772$ est la constante d'Euler-Mascheroni.\\
\textit{D\'emonstration}\\
En commençant par la d\'efinition d'Euler de la fonction gamma (voir exercice 3 pour la d\'emonstration de cette relation)
\begin{align*}
\Gamma(x)=\lim_{m\rightarrow\infty}\frac{m!m^x}{x(x+1)(x+2)\ldots(x+m)}
\end{align*}
qui peut \^etre \'ecrit sous la forme
\begin{align*}
\Gamma(x)&=\lim_{m\rightarrow\infty}\left(\frac{m^x}{x}\prod_{n=1}^{m}\frac{n}{x+n}\right)\\
&=\lim_{m\rightarrow\infty}\left(\frac{m^x}{x}\prod_{n=1}^{m}\left(1+\frac{x}{n}\right)^{-1}\right)\\
&=\lim_{m\rightarrow\infty}\left(\frac{e^{x\ln m}}{x}\prod_{n=1}^{m}\left(1+\frac{x}{n}\right)^{-1}\right)\\
&=\lim_{m\rightarrow\infty}\left(e^{\sum_{n=1}^{m}\frac{x}{n}-\sum_{n=1}^{m}\frac{x}{n}}\frac{e^{x\ln m}}{x}\prod_{n=1}^{m}\left(1+\frac{x}{n}\right)^{-1}\right)\\
&=\lim_{m\rightarrow\infty}\left(e^{-\sum_{n=1}^{m}\frac{x}{n}+x\ln m}\frac{1}{x}\prod_{n=1}^{m} e^{\frac{x}{n}}\left(1+\frac{x}{n}\right)^{-1}\right)\\
&=\lim_{m\rightarrow\infty}\left(\frac{e^{-\gamma x}}{x}\prod_{n=1}^{m}e^{\frac{x}{n}}\left(1+\frac{x}{n}\right)^{-1}\right)\\
\end{align*}
avec 
\begin{align*}
\gamma=\lim_{m\rightarrow\infty}\left(\sum_{n=1}^{m}\frac{1}{n}-\ln m\right)
\end{align*}
donc
\begin{align*}
\Gamma(x)=\frac{e^{-\gamma x}}{x}\prod_{n=1}^{\infty}e^{\frac{x}{n}}\left(1+\frac{x}{n}\right)^{-1}
\end{align*}
\subsection{Formule de compl\'ement}
La formule de compl\'ement (appell\'ee aussi formule de reflexion d'Euler) est donn\'ee par
\begin{align}\label{eq:1_8}
\Gamma(x)\Gamma(1-x)=\frac{\pi}{\sin(\pi x)}
\end{align}
\textit{D\'emonstration}\\
D'apr\'es la formule de Weierstrass (\ref{eq:1_7a}), on a
\begin{align*}
\frac{1}{\Gamma(x)}=xe^{\gamma x}\prod_{n=1}^{\infty}e^{-x/n}\left(1+\frac{x}{n}\right)
\end{align*}
et
\begin{align*}
\frac{1}{\Gamma(-x)}=-xe^{-\gamma x}\prod_{n=1}^{\infty}e^{x/n}\left(1-\frac{x}{n}\right)
\end{align*}
alors
\begin{align*}
\frac{1}{\Gamma(x)\Gamma(-x)}=-x^2\prod_{n=1}^{\infty}\left(1-\frac{x^2}{n^2}\right)
\end{align*}
et d'apr\'es l'Eq.~(\ref{eq:1_2}), on a $\Gamma(-x)=-\Gamma(1-x)/x$, donc
\begin{align}\label{eq:1_7b}
\frac{1}{\Gamma(x)\Gamma(1-x)}=x\prod_{n=1}^{\infty}\left(1-\frac{x^2}{n^2}\right)
\end{align}
Pour montrer la formule de compl\'ement, il reste \`a montrer que
\begin{align*}
\prod_{n=1}^{\infty}\left(1-\frac{x^2}{n^2}\right)=\frac{\sin(\pi x)}{\pi x}
\end{align*}
On a
\begin{align}\label{eq:1_8a}
\sin\theta &= 2\sin\frac{\theta}{2}\cos\frac{\theta}{2}\\
&=2\sin\frac{\theta}{2}\sin\left(\frac{\theta}{2}+\frac{\pi}{2}\right)\nonumber
\end{align}
En appliquant ce r\'esultat aux deux termes du c\^ot\'e droit ($\sin(\theta/2)$ et $\sin((\theta+\pi)/2)$), nous obtenons
\begin{align}\label{eq:1_8b}
\sin\theta &=2\left[2\sin\frac{\theta}{4}\sin\left(\frac{\theta}{4}+\frac{\pi}{2}\right)\right]\left[2\sin\left(\frac{\theta}{4}+\frac{\pi}{4}\right)\sin\left(\frac{\theta}{4}+\frac{\pi}{4}+\frac{\pi}{2}\right)\right]\nonumber\\
&=2^{3}\sin\frac{\theta}{2^2}\sin\frac{\theta+\pi}{2^2}\sin\frac{\theta+2\pi}{2^2}\sin\frac{\theta+3\pi}{2^2}
\end{align}
En appliquant l'Eq.~(\ref{eq:1_8a}) encore une fois \`a (\ref{eq:1_8b}), on obtient
\begin{align}\label{eq:1_8c}
\sin\theta =&2^{7}\sin\frac{\theta}{2^3}\sin\frac{\theta+\pi}{2^3}\sin\frac{\theta+2\pi}{2^3}\sin\frac{\theta+3\pi}{2^3}\sin\frac{\theta+4\pi}{2^3}\nonumber\\
&\times\sin\frac{\theta+5\pi}{2^3}\sin\frac{\theta+6\pi}{2^3}\sin\frac{\theta+7\pi}{2^3}
\end{align}
Si en refaire ce processus $n-$fois, le r\'esultat g\'en\'eral peut \^etre d\'eduit des Eqs.~(\ref{eq:1_8a}), (\ref{eq:1_8b}) and (\ref{eq:1_8c}), et il s'\'ecrit
\begin{align}\label{eq:1_8d}
\sin\theta=2^{2^n-1}\sin\frac{\theta}{2^n}\sin\frac{\theta+\pi}{2^n}\sin\frac{\theta+2\pi}{2^n}\ldots\sin\frac{\theta+(2^n-1)\pi}{2^n}
\end{align}
Le dernier terme de l'Eq.~(\ref{eq:1_8d}) s'\'ecrit
\begin{align*}
\sin\frac{\theta+(2^n-1)\pi}{2^n}&=\sin\left(\pi - \frac{\pi-\theta}{2^{n}}\right)\\
&=\sin\frac{\pi-\theta}{2^{n}}
\end{align*}
le terme avant le dernier s'\'ecrit
\begin{align*}
\sin\frac{\theta+(2^n-2)\pi}{2^n}&=\sin\left(\pi - \frac{2\pi-\theta}{2^{n}}\right)\\
&=\sin\frac{2\pi-\theta}{2^{n}}
\end{align*}
en g\'en\'eral, le $r^{eme}$ terme avant le dernier s'\'ecrit
\begin{align*}
\sin\frac{\theta+(2^n-r)\pi}{2^n}=\sin\frac{r\pi-\theta}{2^{n}}
\end{align*}
Nous regroupons maintenant les termes apparaissant dans l'Eq.~(\ref{eq:1_8d}) en prenant ensemble la deuxi\`eme et le dernier, troisi\`eme et l'avant dernier,et ainsi de suite
\begin{align}\label{eq:1_8e}
\sin\theta =2^{2^n-1}&\sin\frac{\theta}{2^n}\left[\sin\frac{\pi+\theta}{2^n}\sin\frac{\pi-\theta}{2^n}\right]\left[\sin\frac{2\pi+\theta}{2^n}\sin\frac{2\pi-\theta}{2^n}\right]\ldots\nonumber\\
&\times\left[\sin\frac{(\frac{1}{2}2^n-1)\pi+\theta}{2^n}\sin\frac{(\frac{1}{2}2^n-1)\pi-\theta}{2^n}\right]\sin\frac{\frac{1}{2}2^{n}\pi+\theta}{2^{n}}
\end{align}
En utilisant la relation
\begin{align*}
\sin(a+b)\sin(a-b)&=\frac{1}{2}\left(\cos 2b - \cos 2a\right)\\
&=\sin^{2}a - \sin^2b
\end{align*}
l'Eq.~(\ref{eq:1_8e}) devient
\begin{align}\label{eq:1_8f}
\sin\theta =2^{2^n-1}\sin\frac{\theta}{2^n}&\left(\sin^{2}\frac{\pi}{2^n}-\sin^{2}\frac{\theta}{2^n}\right)\left(\sin^{2}\frac{2\pi}{2^n}-\sin^{2}\frac{\theta}{2^n}\right)\ldots\nonumber\\
&\times\left(\sin^{2}\frac{(\frac{1}{2}2^n-1)\pi}{2^n}-\sin^{2}\frac{\theta}{2^n}\right)\cos\frac{\theta}{2^n}
\end{align}
En devisant les deux c\^ot\'es de l'Eq.~(\ref{eq:1_8f}) par $\sin(\theta/2^n)$, en prenant la limit $\theta\rightarrow 0$ et en utilisant $\lim_{x\rightarrow 0}\sin x/x = 1$ afin que
\begin{align*}
\lim_{\theta\rightarrow 0}\frac{\sin\theta}{\sin(\theta/2^n)}&=\lim_{\theta\rightarrow 0}\frac{\sin\theta}{\theta}\frac{\theta/2^n}{\sin(\theta/2^n)}2^n\\
&=2^n
\end{align*}
On prend la limit pour le c\^ot\'e droit, on obtient
\begin{align}\label{eq:1_8g}
2^n=2^{2^n-1}\sin^{2}\frac{\pi}{2^n}\sin^{2}\frac{2\pi}{2^n}\ldots\sin^{2}\frac{(\frac{1}{2}2^n-1)\pi}{2^n}
\end{align}
En devisant l'Eq.~(\ref{eq:1_8f}) par l'Eq.~(\ref{eq:1_8g}), on trouve
\begin{align}
\frac{\sin\theta}{2^n}=\sin\frac{\theta}{2^{n}}\left[1-\frac{\sin^2(\theta/2^n)}{\sin^2(\pi/2^n)}\right]\left[1-\frac{\sin^2(\theta/2^n)}{\sin^2(2\pi/2^n)}\right]\ldots\nonumber\\
\times\left[1-\frac{\sin^2(\theta/2^n)}{\sin^2((\frac{1}{2}2^n-1)\pi/2^n)}\right]\cos\frac{\theta}{2^n}
\end{align}
On prend maintenant la limit $n\rightarrow\infty$ et en utilisant les limites
\begin{align*}
\lim_{x\rightarrow\infty}x\sin\frac{\theta}{x}&=\lim_{x\rightarrow\infty}\frac{\sin(\theta/x)}{\theta/x}x\\
&=\theta,
\end{align*}
\begin{align*}
\lim_{x\rightarrow\infty}\cos\frac{\theta}{x}=1,
\end{align*}
et
\begin{align*}
\lim_{x\rightarrow\infty}\frac{\sin^2(\theta/x)}{\sin^2(r\pi/x)}&=\lim_{x\rightarrow\infty}\left(\frac{\sin(\theta/x)}{\theta/x}\right)^2\left(\frac{\theta}{x}\right)^2\left(\frac{r\pi/x}{\sin(r\pi/x)}\right)^2\frac{1}{(r\pi/x)^2}\\
&=\frac{\theta^2}{r^2\pi^2}
\end{align*}
on trouve
\begin{align*}
\sin\theta &=\theta\left(1-\frac{\theta^2}{\pi^2}\right)\left(1-\frac{\theta^2}{2^2\pi^2}\right)\ldots\left(1-\frac{\theta^2}{r^2\pi^2}\right)\ldots\\
&=\theta\prod_{n=1}^\infty\left(1-\frac{\theta^2}{n^2\pi^2}\right)
\end{align*}
qui est le r\'esultat d\'esir\'e. En comparant le dernier r\'esultat avec (\ref{eq:1_7b}), il en r\'esulte la formule de compl\'ement
\begin{align*}
\Gamma(x)\Gamma(1-x)=\frac{\pi}{\sin\pi x}
\end{align*}
\subsection{Formule de duplication}
La fonction gamma v\'erifie la formule de duplication de Legendre
\begin{align}\label{eq:1_9}
\Gamma(2x)=\frac{2^{2x-1}}{\sqrt{\pi}}\Gamma(x)\Gamma\left(x+\frac{1}{2}\right)
\end{align}
Nous allons montrer cette formule plus tard dans ce chapitre lorsque nous d\'efinirons la fonction b\^eta.
\subsection{Formule de Stirling}
Lorsque $x$ est tr\`es grand, la fonction gamma se comporte comme
\begin{align}\label{eq:1_10}
\Gamma(x+1)\approx \sqrt{2\pi x}x^{x}e^{-x} \hspace{10mm} x\rightarrow +\infty
\end{align}
\textit{D\'emonstration}\\
On a
\begin{align*}
\Gamma(x+1) = \int_0^\infty e^{-t}t^{x}dt
\end{align*}
mais,
\begin{align*}
\ln e^{-t}t^{x} = x\ln t - t
\end{align*}
En faisant le changement de variable $u=t-x$, \i.e., $t=u+x$, on obtient
\begin{align*}
x\ln t - t &= x\ln(x+u) - x -u\\
&= x\ln x\left(1+\frac{u}{x}\right) -x-u\\
&=x\ln x + x\ln\left(1+\frac{u}{x}\right) -x-u
\end{align*}
Lorsque $x$ est tr\`es grand, $u/x << 1$, donc on peut d\'evelopper $\ln(1+u/x)$ en s\'erie de Taylor
\begin{align*}
\ln\left(1+\frac{u}{x}\right)=\sum_{n=1}^{\infty}\frac{(-1)^{n+1}}{n!}\frac{u^n}{x^n}
\end{align*}
il en r\'esulte que
\begin{align*}
x\ln(x+u) - x -u &= x\ln x + \sum_{n=1}^{\infty}\frac{(-1)^{n+1}}{n!}\frac{u^n}{x^{n-1}} -x-u\\
&=x\ln x -x  + \sum_{n=2}^{\infty}\frac{(-1)^{n+1}}{n!}\frac{u^n}{x^{n-1}}\\
&=x\ln x - x -\frac{u^2}{2x} + \frac{u^3}{3x^2} - \frac{u^4}{4x^3}\pm \ldots
\end{align*}
Lorsque $x$ est tr\`es grand on a
\begin{align*}
\ln\left(e^{-t}t^{x}\right) &=x\ln(x+u) - x -u\\
&\approx x\ln x - x -\frac{u^2}{2x}
\end{align*}
d'o\`u
\begin{align*}
e^{-t}t^{x}&\approx \frac{x^x}{e^{x}}e^{-\frac{u^2}{2x}}
\end{align*}
En remplaçant $e^{-t}t^{x}$ dans l'\'expression de $\Gamma(x+1)$ (avec le changement de variable $u=t-x$), on obtient
\begin{align*}
\Gamma(x+1)\approx\int_{-x}^\infty\frac{x^x}{e^{x}}e^{-\frac{u^2}{2x}}du
\end{align*}
En utilisant la formule de Gauss
\begin{align*}
\int_{-\infty}^{\infty}e^{-ax^2}dx=\sqrt{\frac{\pi}{a}}
\end{align*}
avec la condition que $x$ est tr\`es grand, on obtient
\begin{align*}
\Gamma(x+1)&\approx \frac{x^x}{e^{x}}\int_{-\infty}^\infty e^{-\frac{u^2}{2x}}du\\
&=\frac{x^x}{e^{x}}\sqrt{2x\pi}
\end{align*}
donc
\begin{align*}
\Gamma(x+1)\approx \sqrt{2\pi x}x^{x}e^{-x} \hspace{10mm} x\rightarrow +\infty
\end{align*}
Si $x$ est un nombre entier ($x\in\mathsf{N}$), on a
\begin{align*}
x!\approx \sqrt{2\pi x}x^{x}e^{-x} \hspace{10mm} x\rightarrow +\infty
\end{align*}
\subsection{Fonction gamma incompl\`ete}
Les fonction gamma incompl\`ete et gamma incompl\`ete compl\'ementaire sont d\'efinis respectivement par
\begin{align}
&\gamma(x,a)=\int_0^{a}e^{-t}t^{x-1}dt, \hspace{10mm} x> 0 \label{eq:1_11}\\
&\Gamma(x,a)=\int_a^{\infty}e^{-t}t^{x-1}dt, \hspace{10mm} x > 0\label{eq:1_12}
\end{align}
Elles remplissent les propri\'et\'es suivantes
\begin{align}
&(a)\hspace{1.5mm}\gamma(x,a)+\Gamma(x,a)=\Gamma(x),\nonumber\\
&(b)\hspace{1.5mm}\gamma(x+1,a)=x\gamma(x,a)-a^{x}e^{-a},\nonumber\\
&(c)\hspace{1.5mm}\Gamma(x+1,a)=x\Gamma(x,a)+a^{x}e^{-a},\label{eq:1_13}\\
&(d)\hspace{1.5mm}\frac{d}{da}\left[a^{-x}\gamma(x,a)\right]=-a^{-x-1}\gamma(x+1,a),\nonumber\\
&(e)\hspace{1.5mm}\frac{d}{da}\left[a^{-x}\Gamma(x,a)\right]=-a^{-x-1}\Gamma(x+1,a),\nonumber
\end{align}
\textit{D\'emonstration}\\
La preuve de la propri\'et\'e (\ref{eq:1_13})(a) est trivial.\\
Pour montrer l'Eq.~(\ref{eq:1_13})(b), on utilise l'int\'egration par parties avec $u'=e^{-t}$, et $v=t^{x}$, on a
\begin{align}
\gamma(x+1,a)&=\int_0^{a}e^{-t}t^{x}dt\nonumber\\
&=-e^{-t}t^{x}\Big|_0^{a} + x\int_0^{a}e^{-t}t^{x-1}dt\nonumber\\
&= x\gamma(x,a)-a^{x}e^{-a}\nonumber
\end{align}
La preuve de la propri\'et\'e (\ref{eq:1_13})(c) est identique \`a (b).\\
D\'emontrons (d), on a
\begin{align}
\frac{d}{da}\left[a^{-x}\gamma(x,a)\right]&= -xa^{-x-1}\gamma(x,a)+ a^{-x}\frac{d}{da}\gamma(x,a)\nonumber\\
&=-a^{-x-1}\Big[x\gamma(x,a)-a^{x}e^{-a}\Big]\nonumber\\
&=-a^{-x-1}\gamma(x+1,a)\nonumber
\end{align}
o\`u $d\gamma(x,a)/dx=e^{-a}a^{x-1}$ et dans la derni\`ere ligne nous avons utilis\'e la propri\'et\'e (b).\\
La preuve de la propri\'et\'e (e) est identique \`a (\ref{eq:1_13})(d).\\
\subsection{D\'eriv\'ee logarithmique}
La d\'eriv\'ee logarithmique de la fonction gamma d\'efinie une nouvelle fonction qui s'appelle la fonction psi (ou la fonction digamma)
\begin{align}\label{eq:1_14}
\psi(x)&=\frac{d}{dx}\left[\ln\Gamma(x)\right]\nonumber\\
&=\frac{\Gamma'(x)}{\Gamma(x)},\hspace{10mm}x\neq 0,-1,-2,-3,\ldots
\end{align}
qui a les propri\'et\'es
\begin{align}
&(a)\hspace{1.5mm}\psi(x+1)=\psi(x)+\frac{1}{x} ,\nonumber\\
&(b)\hspace{1.5mm}\psi(1-x)=\psi(x)+\pi\cot(\pi x) ,\label{eq:1_15}\\
&(c)\hspace{1.5mm}\psi(n+1)=-\gamma + \sum_{k=1}^{n}\frac{1}{k} \hspace{10mm} n=0,1,2,3,\ldots\nonumber
\end{align}
avec $\gamma=0.5772$ est la constante d'Euler.\\
\textit{D\'emonstration}\\
(a) Par d\'efinition, on a
\begin{align}
\psi(x+1)&=\frac{d}{dx}\ln\left[\Gamma(x+1)\right]\nonumber\\
&=\frac{d}{dx}\ln\left[x\Gamma(x)\right]\nonumber\\
&=\frac{\Gamma(x)+\Gamma'(x)}{x\Gamma(x)}\nonumber\\
&=\psi(x)+\frac{1}{x}\nonumber
\end{align}
(b) on a
\begin{align}
\psi(1-x)&=\frac{d}{d(1-x)}\ln\Gamma(1-x)\nonumber\\
&=-\frac{d}{dx}\ln\left(\frac{\pi}{\Gamma(x)\sin(\pi x)}\right)\nonumber\\
&=\frac{d}{dx}\ln\Gamma(x)+\frac{d}{dx}\ln\sin(\pi x)\nonumber\\
&=\psi(x)+\pi\cot(\pi x)\nonumber
\end{align}
Dans la deuxi\`eme ligne, nous avons utilis\'e la formule des compl\'ements Eq.~(\ref{eq:1_8}).\\
La $m^{\textit{\`eme}}$ d\'eriv\'e logarithmique de la fonction gamma d\'efinie la fonction polygamma, et elle est donn\'ee par
\begin{align}\label{eq:1_16}
\psi^{(m)}(x)=\frac{d^{m+1}}{dx^{m+1}}\left[\ln\Gamma(x)\right]\hspace{10mm} x\neq 0,-1,-2,-3,\ldots
\end{align}
\section{Fonction b\^eta}
\subsection{D\'efinition}
La fonction b\^eta (ou l'int\'egrale d'Euler de premi\`ere esp\`ece) est d\'efinie pour tous nombres complexes $x$ et $y$ de parties r\'eelles strictement positives par 
\begin{align}\label{eq:1_17}
B(x,y)=\int_0^{1}t^{x-1}(1-t)^{y-1}dt \hspace{10mm} x > 0, y > 0
\end{align} 
La fonction b\^eta peut \^etre d\'efinie \'egalement par l'int\'egrale suivante
\begin{align}\label{eq:1_18}
B(x,y)=\int_0^{\infty}\frac{t^{x-1}}{(1+t)^{x+y}}dt \hspace{10mm} x > 0, y > 0
\end{align}
\textit{D\'emonstration}\\
Pour d\'emontrer l'\'equivalence entre les deux d\'efinitions (\ref{eq:1_17}) et l'Eq.~(\ref{eq:1_18}), nous utilisons le changement de variable 
\begin{align*}
u=\frac{1}{1+t}
\end{align*}
Il en r\'esulte que
\begin{align*}
t=\frac{1-u}{u}, \hspace{20mm} dt=-\frac{du}{u^{2}}
\end{align*}
En remplaçant $t$ et $dt$ dans l'Eq.~(\ref{eq:1_17}), on obtient la formule (\ref{eq:1_18}).
\subsection{Relation entre les fonctions gamma et b\^eta}
La fonction b\^eta est li\'ee \`a la fonction gamma par la formule
\begin{align}\label{eq:1_19}
B(x,y) &= 2\int_0^{\pi/2}\cos^{2x-1}\theta\sin^{2y-1}\theta d\theta\nonumber\\
&=\frac{\Gamma(x)\Gamma(y)}{\Gamma(x+y)}
\end{align}
\textit{D\'emonstration}\\
En faisant le  changement de variable $t=\cos^{2}\theta$ dans l'Eq.~(\ref{eq:1_17}) nous obtenons la relation (\ref{eq:1_19}). Comme cons\'equence directe de cette relation, on a que
\begin{align*}
B(x,y)=B(y,x)
\end{align*}
\subsection{Propri\'et\'es de la fonction b\^eta}
Les propri\'et\'es les plus importantes de la fonction b\^eta sont
\begin{align}\label{eq:1_20}
&(a)\hspace{1.5mm} B(x+1,y)=\frac{x}{x+y}B(x,y),\\
&(b)\hspace{1.5mm}B(x,y+1)=\frac{y}{x+y}B(x,y),\nonumber\\
&(c)\hspace{1.5mm}B(x+1,y)+B(x,y+1)=B(x,y),\nonumber\\
&(d)\hspace{1.5mm}B(x,x)=2^{1-2x}B\left(x,\frac{1}{2}\right),\nonumber
\end{align}
\textit{D\'emonstration}\\
(a) En utilisant la d\'efinition de la fonction b\^eta, on a
\begin{align*}
B(x+1,y)&=\frac{\Gamma(x+1)\Gamma(y)}{\Gamma(x+y+1)}\\
&=\frac{x\Gamma(x)\Gamma(y)}{(x+y)\Gamma(x+y)}\\
&=\frac{x}{x+y}B(x,y)
\end{align*}
Pour (b) on a 
\begin{align*}
B(x,y+1)&=\frac{\Gamma(x)\Gamma(y+1)}{\Gamma(x+y+1)}\\
&=\frac{\Gamma(x)y\Gamma(y)}{(x+y)\Gamma(x+y)}\\
&=\frac{y}{x+y}B(x,y)
\end{align*}
La propri\'et\'e (c) d\'ecoule directement en faisant la somme de (a) et (b).\\
Pour montrer (d), on fait le  changement de variable $t=(1+u)/2$, $dt=du/2$. On a
\begin{align}
B(x,x)&=\int_0^1\left(\frac{1}{2}(1+u)\right)^{x-1}\left(\frac{1}{2}(1-u)\right)^{x-1}\frac{du}{2}\nonumber\\
&=\frac{1}{2^{2x-1}}\int_{-1}^1\left[(1+u)(1-u)\right]^{x-1}du\nonumber\\
&=\frac{1}{2^{2x-1}}\int_{-1}^1\left(1-u^2\right)^{x-1}du\nonumber\\
&=\frac{2}{2^{2x-1}}\int_{0}^1\left(1-u^2\right)^{x-1}du\nonumber
\end{align}
Dans la derni\`ere ligne nous avons utilis\'es le fait que l'int\'egrande est une fonction paire. Faisant un deuxi\`eme changement de variable $v=u^{2}$, on obtient
\begin{align*}
B(x,x)&=\frac{1}{2^{2x-1}}\int_{0}^1v^{-\frac{1}{2}}(1-v)^{x-1}dv\\
& = 2^{1-2x}B\left(x,\frac{1}{2}\right)
\end{align*}
Etant donn\'ee la propri\'et\'e (\ref{eq:1_20})(d), on peut d\'eduire maintenant la formule de duplication de Legendre Eq.~(\ref{eq:1_9}). En effet, en remplaçant la fonction b\^eta dans les deux c\^ot\'es de l'Eq.~(\ref{eq:1_20})(d) par la fonction gamma (en utilisant l'Eq.~(\ref{eq:1_19})), on obtient
\begin{align*}
\frac{\Gamma(x)\Gamma(x)}{\Gamma(2x)}=2^{1-2x}\frac{\Gamma(x)\Gamma(\frac{1}{2})}{\Gamma(x+\frac{1}{2})}
\end{align*}
Il en r\'esulte que
\begin{align*}
\Gamma(2x)=\frac{2^{2x-1}}{\sqrt{\pi}}\Gamma(x)\Gamma\left(x+\frac{1}{2}\right)
\end{align*}
Si $x\in \mathsf{N}$, on a $\Gamma(2x)=(2x-1)!$ et $\Gamma(x)=(x-1)!$. Dans ce cas la, la formule de duplication se r\'eduit \`a
\begin{align}\label{eq:1_9a}
\Gamma\left(x+\frac{1}{2}\right)=\frac{(2x)!}{x!}\frac{\sqrt{\pi}}{2^{2x}}
\end{align}
\section{Exercices}
$\mathbf{Exercice\hspace{1mm}1:}$ Calculer les int\'egrales suivantes en utilisant les propri\'et\'es  des fonctions gamma et b\^eta
\begin{align*}
&(a)\hspace{0.5mm}\int_{0}^{\infty}u^{4}e^{-u^{3}}du,\hspace{19mm}
(b)\hspace{0.5mm}\int_{a}^{\infty}e^{2au-u^{2}}du,\hspace{9mm}
(c)\hspace{0.5mm}\int_{0}^{\pi/2}\sqrt{\tan\theta}d\theta,\\
&(d)\hspace{0.5mm}\int_{0}^{\infty}u^{-3/4}e^{-\sqrt{u}}du,\hspace{10mm}
(e)\hspace{0.5mm}\int_{0}^{1}\frac{du}{\sqrt[3]{1-u^{4}}},\hspace{10mm}
(f)\hspace{0.5mm}\int_{-1}^{1}\left(\frac{1+u}{1-u}\right)^{1/2}du.
\end{align*}
\textbf{Solutions}\\
$(a)$ En faisant le changement de variable $t=u^3$, $du=dt/(3t^{2/3})$ on obtient
\begin{align*}
\int_{0}^{\infty}u^{4}e^{-u^{3}}du &=\int_{0}^{\infty}t^{\frac{4}{3}}e^{-t}\frac{dt}{3t^{\frac{2}{3}}}\\
&=\frac{1}{3}\int_{0}^{\infty}t^{\frac{2}{3}}e^{-t}dt
\end{align*}
En identifiant la derni\`ere integrale avec la d\'efinition de la fonction gamma ($\Gamma(x)=\int_{0}^{\infty}t^{x-1}e^{-t}dt$), on a $x-1=2/3$ ce qui signifie que $x=5/3$, donc
\begin{align*}
\int_{0}^{\infty}u^{4}e^{-u^{3}}du=\frac{1}{3}\Gamma\left(\frac{5}{3}\right)
\end{align*}
$(b)$ On a
\begin{align*}
\int_{a}^{\infty}e^{2au-u^{2}}du=\int_{a}^{\infty}e^{-(u-a)^{2}+a^2}du
\end{align*}
Posons $t=(u-a)^{2}$, \i.e., $dt=2(u-a)du=2\sqrt{t}du$, donc
\begin{align*}
\int_{a}^{\infty}e^{-(u-a)^{2}+a^2}du&=e^{a^2}\int_{0}^{\infty}e^{-t}\frac{dt}{2\sqrt{t}}\\
&=\frac{e^{a^2}}{2}\int_{0}^{\infty}t^{-\frac{1}{2}}e^{-t}dt\\
&=\frac{e^{a^2}}{2}\Gamma\left(\frac{1}{2}\right)\\
&=\frac{\sqrt{\pi} e^{a^2}}{2}
\end{align*}
$(c)$ On a
\begin{align*}
\int_{0}^{\pi/2}\sqrt{\tan\theta}d\theta=\int_{0}^{\pi/2}\sin^{1/2}\theta\cos^{-1/2}\theta d\theta
\end{align*}
En utilisant la relation Eq.~(\ref{eq:1_19}), avec $2x-1=-1/2$ et $2y-1=1/2$, \i.e., $x=1/4$ et $y=3/4$, on trouve
\begin{align*}
\int_{0}^{\pi/2}\sqrt{\tan\theta}d\theta=\frac{1}{2}\Gamma\left(\frac{1}{4}\right)\Gamma\left(\frac{3}{4}\right)
\end{align*}
$(d)$ En faisant le changement de variable $t=\sqrt{u}$, donc
\begin{align*}
\int_{0}^{\infty}u^{-3/4}e^{-\sqrt{u}}du &=\int_{0}^{\infty}t^{-3/2}e^{-t}2tdt\\
&=2\int_{0}^{\infty}t^{-1/2}e^{-t}dt\\
&=2\Gamma\left(\frac{1}{2}\right)\\
&=2\sqrt{\pi}
\end{align*}
$(e)$ Posons $t=u^4$, $dt=4t^{3/4}du$, on trouve
\begin{align*}
\int_{0}^{1}\frac{du}{\sqrt[3]{1-u^{4}}}&=\int_{0}^{1}\frac{1}{\sqrt[3]{1-t}}\frac{dt}{4t^{\frac{3}{4}}}\\
&=\frac{1}{4}\int_{0}^{1}t^{-\frac{3}{4}}(1-t)^{-\frac{1}{3}}dt
\end{align*}
En identifiant avec la d\'efinition de la fonction b\^eta on trouve que $x-1=-3/4$ et $y-1=-1/3$. D'o\`u
\begin{align*}
\int_{0}^{1}\frac{du}{\sqrt[3]{1-u^{3}}}=\frac{1}{4}B\left(\frac{1}{4},\frac{2}{3}\right)
\end{align*}
$(f)$ En utilisant le changement de variable $t=1+u$, on trouve
\begin{align*}
\int_{-1}^{1}\left(\frac{1+u}{1-u}\right)^{1/2}du&=\int_{-1}^{1}\left(1+u\right)^{\frac{1}{2}}\left(1-u\right)^{-\frac{1}{2}}du\\
&=\int_{0}^{2}t^{\frac{1}{2}}\left(2-t\right)^{-\frac{1}{2}}dt\\
&=\frac{1}{\sqrt{2}}\int_{0}^{2}t^{\frac{1}{2}}\left(1-\frac{t}{2}\right)^{-\frac{1}{2}}dt
\end{align*}
En faisant un deuxi\`eme changement de variable $v=t/2$, on trouve
\begin{align*}
\int_{-1}^{1}\left(\frac{1+u}{1-u}\right)^{1/2}du&=2\int_{0}^{1}v^{\frac{1}{2}}\left(1-v\right)^{-\frac{1}{2}}dv
\end{align*}
Par identification avec la d\'efinition de la fonction b\^eta, on obtient que $x-1=1/2$ et $y-1=-1/2$, donc
\begin{align*}
\int_{-1}^{1}\left(\frac{1+u}{1-u}\right)^{1/2}du&=2B\left(\frac{3}{2},\frac{1}{2}\right)
\end{align*}
$\textbf{Exercice\hspace{1mm}2:}$ Montrer que
\begin{align*}
&(a)\hspace{0.5mm}\int_{0}^{1}\sqrt{1-u^{2}}du=\frac{\pi}{4},\hspace{22mm}
(b)\hspace{0.5mm}\int_{0}^{1}\left(\frac{1}{u}-1\right)^{\frac{1}{4}}du=\Gamma\left(\frac{3}{4}\right)\Gamma\left(\frac{5}{4}\right),\\
&(c)\hspace{0.5mm}\int_{0}^{1}\left(\ln\frac{1}{u}\right)^{2}du=2,\hspace{25mm}
(d)\hspace{0.5mm}\int_{0}^{1}\frac{du}{\sqrt{1-u^4}}=\frac{\left[\Gamma\left(\frac{1}{4}\right)\right]^{2}}{4\sqrt{2\pi}},\\
&(e)\hspace{0.5mm}\int_{0}^{\infty}u^n e^{-au^{2}}du=\frac{\Gamma\left(\frac{n+1}{2}\right)}{2\hspace{0.4mm}a^{\frac{n+1}{2}}},\hspace{15mm}
(f)\hspace{0.5mm}B(n,n+1)=\frac{1}{2}\frac{\left[\Gamma(n)\right]^{2}}{\Gamma(2n)},\\
&(g)\hspace{0.5mm}\int_{0}^{\pi/2}\sin^{n}(\theta)d\theta =\!\! \int_{0}^{\pi/2}\cos^{n}(\theta)d\theta =\frac{\sqrt{\pi}}{2}\frac{\Gamma\left(\frac{1+n}{2}\right)}{\Gamma\left(\frac{2+n}{2}\right)},\\
&(h)\hspace{0.5mm}\int_{0}^{\pi/2}\tan^{n}\theta d\theta=\frac{1}{2}\Gamma\left(\frac{1+n}{2}\right)\Gamma\left(\frac{1-n}{2}\right),\\
&(i)\hspace{0.5mm}\int_{0}^{\infty}u^{m}e^{-u^{n}}du=\frac{1}{n}\Gamma\left(\frac{m+1}{n}\right), \hspace{15mm} \mbox{$m>-1, n>0$}.
\end{align*}
\textbf{Solutions}\\
$(a)$ Posons $t=u^2$, \i.e., $dt=2\sqrt{t}du$, donc
\begin{align*}
\int_{0}^{1}\sqrt{1-u^{2}}du=\int_{0}^{1}\sqrt{1-t}\frac{dt}{2\sqrt{t}}\
\end{align*}
Par identification avec la d\'efinition (\ref{eq:1_17}) de la fonction b\^eta, on obtient $x-1=-1/2$ et $y-1=1/2$, d'o\`u
\begin{align*}
\int_{0}^{1}\sqrt{1-u^{2}}du&=\frac{1}{2}B\left(\frac{1}{2},\frac{3}{2}\right)\\
&=\frac{1}{2}\Gamma\left(\frac{1}{2}\right)\Gamma\left(\frac{3}{2}\right)
\end{align*}
or $\Gamma(3/2)=(1/2)\Gamma(1/2)=\sqrt{\pi}/2$, donc
\begin{align*}
\int_{0}^{1}\sqrt{1-u^{2}}du=\frac{\pi}{4}
\end{align*}
$(b)$ On a
\begin{align*}
\int_{0}^{1}\left(\frac{1}{u}-1\right)^{\frac{1}{4}}du=\int_{0}^{1}u^{-\frac{1}{4}}\left(1-u\right)^{\frac{1}{4}}du
\end{align*}
D'apr\`es la d\'efinition de la fonction b\^eta on obtient $x-1=-1/4$ et $y-1=1/4$, il en r\'esulte que
\begin{align*}
\int_{0}^{1}\left(\frac{1}{u}-1\right)^{\frac{1}{4}}du&=B\left(\frac{3}{4},\frac{5}{4}\right)\\
&=\Gamma\left(\frac{3}{4}\right)\Gamma\left(\frac{5}{4}\right)
\end{align*}
$(c)$ Par identification avec la repr\'esentation (\ref{eq:1_5}) de la fonction gamma, on trouve que $x-1=2$ donc $x=3$, d'o\`u
\begin{align*}
\int_{0}^{1}\left(\ln\frac{1}{u}\right)^{2}du&=\Gamma(3)\\
&=2
\end{align*}
$(d)$ Posons $t=u^4$, il en r\'esulte que
\begin{align*}
\int_{0}^{1}\frac{du}{\sqrt{1-u^4}}&=\frac{1}{4}\int_{0}^{1}t^{-\frac{3}{4}}(1-t)^{-\frac{1}{2}}dt\\
&=\frac{1}{4}B\left(\frac{1}{4},\frac{1}{2}\right)\\
&=\frac{1}{4}\frac{\Gamma(1/4)\Gamma(1/2)}{\Gamma(3/4)}
\end{align*}
On peut calculer $\Gamma(3/4)$ \`a partir de la formule de duplication (pour $x=1/4$), on a 
\begin{align*}
\Gamma\left(2\frac{1}{4}\right)=\frac{2^{2\frac{1}{4}-1}}{\sqrt{\pi}}\Gamma\left(\frac{1}{4}\right)\Gamma\left(\frac{3}{4}\right)
\end{align*}
En remplaçant $\Gamma(3/4)$, on trouve
\begin{align*}
\int_{0}^{1}\frac{du}{\sqrt{1-u^4}}&=\frac{1}{4}\frac{\Gamma(1/4)\Gamma(1/2)\Gamma(1/4)}{\sqrt{2\pi}\Gamma(1/2)}\\
&=\frac{1}{4}\frac{\left[\Gamma(1/4)\right]^{2}}{\sqrt{2\pi}}
\end{align*}
$(e)$ En faisant le changement de variable $t^2=au^2$, \i.e., $dt=\sqrt{a}du$ et en utilisant la repr\'esentation (\ref{eq:1_4}) de la fonction gamma, on trouve
\begin{align*}
\int_{0}^{\infty}u^n e^{-au^{2}}du&=\int_{0}^{\infty}\frac{t^{n}}{a^{\frac{n}{2}}} e^{-t^{2}}\frac{dt}{\sqrt{a}}
\end{align*}
Par identification avec la formule (\ref{eq:1_4}), on trouve $2x-1=n$, \i.e., $x=(n+1)/2$, donc
\begin{align*}
\int_{0}^{\infty}u^n e^{-au^{2}}du=\frac{1}{2a^{\frac{n+1}{2}}}\Gamma\left(\frac{n+1}{2}\right)
\end{align*}
$(f)$ On a
\begin{align*}
B(n,n+1)&=\frac{\Gamma(n)\Gamma(n+1)}{\Gamma(2n+1)}\\
&=\frac{\Gamma(n)n\Gamma(n)}{2n\Gamma(2n)}\\
&=\frac{\left[\Gamma(n)\right]^{n}}{2\Gamma(2n)}
\end{align*}
$(g)$ En identifiant $\int_0^{\pi/2}\sin^{n}(\theta)d\theta$ avec la formule (\ref{eq:1_6}), on trouve $x=1/2$ et $y=(n+1)/2$, d'o\`u
\begin{align*}
\int_0^{\pi/2}\sin^{n}(\theta)d\theta &=\frac{\Gamma(1/2)\Gamma([n+1]/2)}{2\Gamma([n+2]/2)}\\
&=\frac{\sqrt{\pi}\Gamma([n+1]/2)}{2\Gamma([n+2]/2)}
\end{align*}
On fait le m\^eme raisonnement pour $\int_0^{\pi/2}\cos^{n}(\theta)d\theta$.\\
$(h)$ On a
\begin{align*}
\int_{0}^{\pi/2}\tan^{n}\theta d\theta=\int_{0}^{\pi/2}\sin^{n}\theta \cos^{-n}\theta d\theta
\end{align*}
Par identification avec la formule (\ref{eq:1_6}), on trouve $2x-1=-n$ et $2y-1=n$, c'est \`a dire $x=(1-n)/2$ et $y=(1+n)/2$, donc
\begin{align*}
\int_{0}^{\pi/2}\tan^{n}\theta d\theta=\frac{1}{2}\Gamma\left(\frac{1-n}{2}\right)\Gamma\left(\frac{1+n}{2}\right)
\end{align*}
$(i)$ Posons $t=u^n$, \i.e., $u=t^{1/n}$ et $du=dt/nt^{(n-1)/n}$, alors
\begin{align*}
\int_{0}^{\infty}u^{m}e^{-u^{n}}du&=\int_{0}^{\infty}t^{\frac{m}{n}}e^{-t}\frac{dt}{nt^{\frac{n-1}{n}}}\\
&=\frac{1}{n}\int_{0}^{\infty}t^{\frac{m+1}{n}-1}e^{-t}dt\\
&=\frac{1}{n}\Gamma\left(\frac{m+1}{n}\right)
\end{align*}
$\textbf{Exercice\hspace{1mm}3:}$ Montrer que
\begin{align*}
\Gamma(x)=\lim_{m\rightarrow \infty}\frac{m^xm!}{x(x+1)(x+2)\ldots(x+m)},\hspace{10mm}x>0
\end{align*}
\textbf{Solutions}\\
Posons $u=t/m$ dans l'int\'egral suivant, on obtient
\begin{align*}
\int_{0}^{m}\left(1-\frac{t}{m}\right)^{m}t^{x-1}dt=m^{x}\int_{0}^{1}\left(1-u\right)^{m}u^{x-1}du
\end{align*}
En faisant plusieurs int\'egrations par parties, on obtient
\begin{align*}
m^{x}\int_{0}^{1}\left(1-u\right)^{m}u^{x-1}du&=m^{x}\frac{m}{x}\int_{0}^{1}\left(1-u\right)^{m-1}u^{x}du\\
&=m^{x}\frac{m(m-1)}{x(x+1)}\int_{0}^{1}\left(1-u\right)^{m-2}u^{x+1}du\\
&=m^{x}\frac{m(m-1)(m-2)}{x(x+1)(x+2)}\int_{0}^{1}\left(1-u\right)^{m-3}u^{x+2}du\\
&\vdots\\
&=m^{x}\frac{m(m-1)(m-2)\ldots 1}{x(x+1)(x+2)\ldots(x+m-1)}\int_{0}^{1}u^{x+m-1}du\\
&=m^{x}\frac{m!}{x(x+1)(x+2)\ldots(x+m)}
\end{align*}
D'autre, d'apr\'es le th\'eor\`eme de convergence domin\'ee de Lebesgue, on a
\begin{align*}
\lim_{m\rightarrow\infty}\int_{0}^{m}\left(1-\frac{t}{m}\right)^{m}t^{x-1}dt&=\int_{0}^{m}\lim_{m\rightarrow\infty}\left(1-\frac{t}{m}\right)^{m}t^{x-1}dt\\
&=\int_{0}^{\infty}e^{-t}t^{x-1}dt\\
&=\Gamma(x)
\end{align*}
Dans la deuxi\`eme ligne on a utilis\'e le fait que $\lim_{t\rightarrow\infty}(1-t/m)^{m}=e^{-t}$. Donc on d\'eduit que
\begin{align*}
\Gamma(x)=\lim_{m\rightarrow\infty}\frac{m!m^{x}}{x(x+1)(x+2)\ldots(x+m)}
\end{align*}
\chapter{Les fonctions de Bessel}
\label{chap: bessel}
\section{D\'efinition}
\begin{figure}[!ht]
\centering
\includegraphics[scale=0.8]{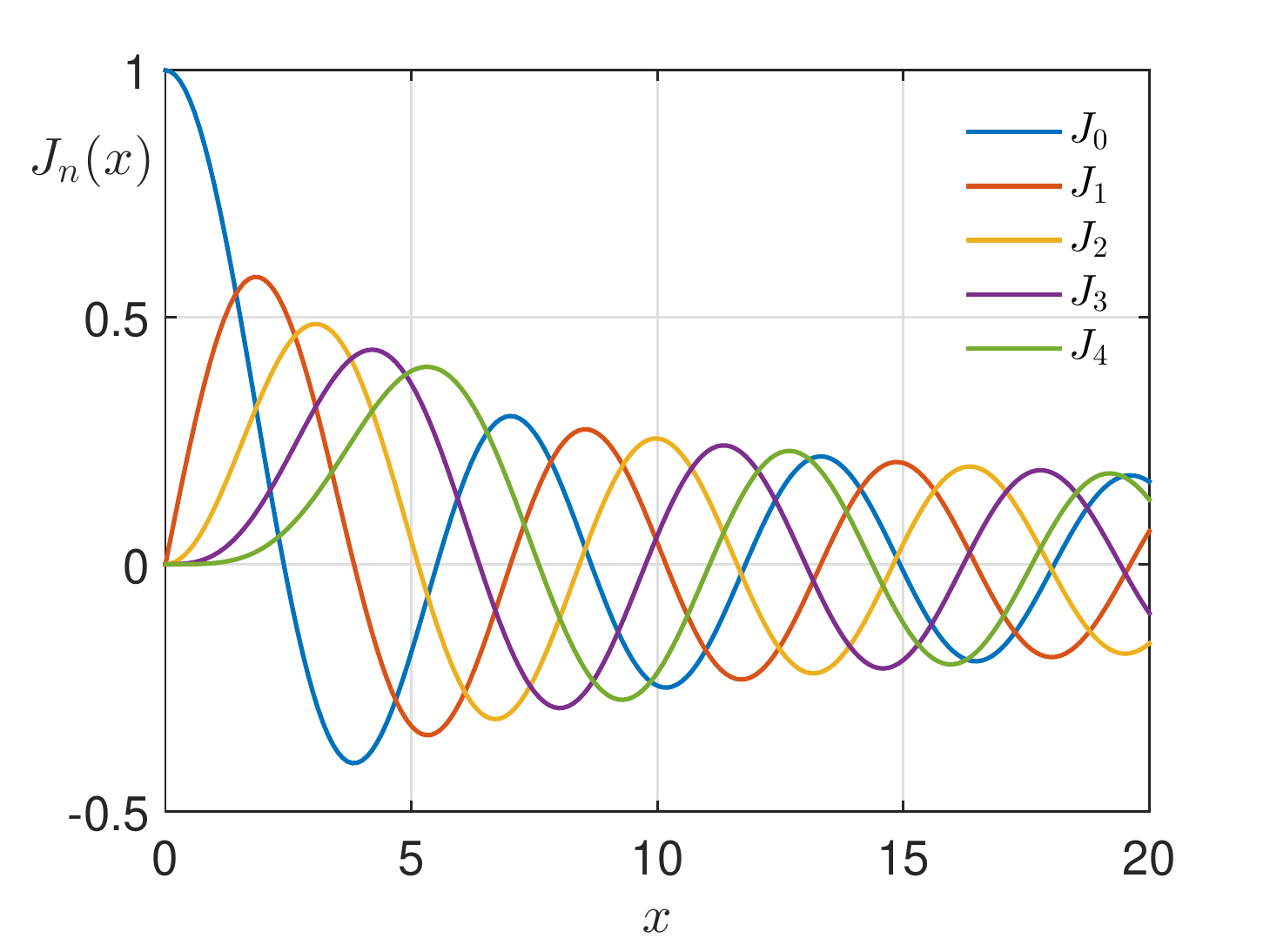}
\caption{\label{fig:bessel_j} Trac\'e de la fonction de Bessel $J_n(x)$ pour $n=0,1,2,3,4$.}
\end{figure}

Les fonctions de Bessel (appell\'ees aussi fonctions cylindriques ou d'harmoniques cylindriques) ont \'et\'e introduites par Bernouli et l'analyse de ses fonctions a \'et\'e 
d\'evelopp\'ee par Bessel en 1860. Les fonctions de Bessel n\^ot\'ees $J_n(x)$ sont les solutions de l'\'equation suivante dite de Bessel
\begin{align}\label{eq:2_1}
x^{2}\frac{d^{2}y}{dx^{2}} + x\frac{dy}{dx} + (x^{2}-n^{2})y=0,
\end{align}
o\`u $n$ est un nombre r\'eel positif (le plus souvent, $n$ est un entier naturel, ou un demi-entier). Cette \'equation admet deux solutions ind\'ependente $J_n(x)$ et $J_{-n}(x)$ o\`u $J_n(x)$
est donn\'ee par
\begin{align}\label{eq:2_2}
J_n(x)=\sum_{r=0}^{\infty}(-1)^{r}\frac{1}{r!\Gamma(n+r+1)}\left(\frac{x}{2}\right)^{2r+n}
\end{align}
\textit{D\'emonstration}\\
Pour r\'esoudre l'Eq.~(\ref{eq:2_1}), on utilise la m\'ethode de Frobenius, qui consiste \`a chercher des solutions d\'eveloppables en s\'eries enti\`eres. La m\'ethode de Frobenius permet alors de d\'eterminer une solution sous la forme
\begin{align*}
\mathit{y}(x,s)=\sum_{r=0}^{\infty}a_rx^{r+s}
\end{align*}
avec $a_0\neq 0$. On substitue dans l'\'equation diff\'erentielle ~(\ref{eq:2_1}), l'\'expression de $y(x,s)$, sa d\'eriv\'ee
\begin{align*}
\frac{d\mathit{y}}{dx}=\sum_{r=0}^{\infty}a_r(r+s)x^{r+s-1},
\end{align*}
et la d\'eriv\'ee seconde
\begin{align*}
\frac{d^{2}\mathit{y}}{dx^{2}}=\sum_{r=0}^{\infty}a_r(r+s)(r+s-1)x^{r+s-2}.
\end{align*}
Il en r\'esulte les \'equations indiciales suivantes
\begin{align*}
&a_0(s^{2}-n^{2})=0,\\
&a_1((s+1)^{2}-n^{2})=0,\\
&a_{r-2} +a_r\left\lbrace(r+s)^{2}-n^{2}\right\rbrace=0 \hspace{10mm} r\geq 2.
\end{align*}
Puisque $a_0\neq 0$, donc les solutions des \'equations indiciales devraient \^etre
\begin{align*}
&s=\pm n,\\
&a_1=0,\\
&a_r = -\frac{a_{r-2}}{(r+s)^{2}-n^{2}} \hspace{10mm} r\geq 2.
\end{align*}
On commence par $s=n$, dans ce cas le terme g\'en\'eral $a_r$ devient
\begin{align}\label{eq:2_3}
a_r = -\frac{a_{r-2}}{r(2n+r)} \hspace{10mm} r\geq 2
\end{align}
Pour trouver l'\'expression de $a_r$, $r=2,3,\ldots$ en fonction de $r$ et de $a_0$, consid\'erons quelques cas ($r=2,4,6$), puis \`a partir de ces cas, on d\'eduit la forme g\'en\'erale de $a_r$. Pour $n=2$, on a
\begin{align*}
a_2=-\frac{a_0}{2^{2}(n+1)}
\end{align*}
Pour $n=4$, on a
\begin{align*}
a_4&=-\frac{a_2}{2^{2}.2(n+2)}\\
&=\frac{a_0}{2^{4}.2(n+1)(n+2)}
\end{align*}
et pour $n=6$, on a
\begin{align*}
a_6&=-\frac{a_4}{6(2n+6)}\\
&=-\frac{a_0}{2^{6}3.2(n+1)(n+2)(n+3)}
\end{align*}
Donc lorsque l'indice de coefficient $a_r$ est pair, on peut \'ecrire 
\begin{align*}
a_{2r}=(-1)^{r}\frac{a_0}{2^{2r}r!(n+1)(n+2)\ldots(n+r)}
\end{align*}
En multipliant et en divisant par $\Gamma(n+1)$ et en utilisant la formule
\begin{align*}
(n+1)(n+2)\ldots(n+r)\Gamma(n+1)=\Gamma(n+r+1)
\end{align*}
on trouve
\begin{align*}
a_{2r}=(-1)^{r}\frac{a_0\Gamma(n+1)}{2^{2r}r!\Gamma(n+r+1)}
\end{align*}
Le fait que $a_1=0$ garantit que (d'apr\`es l'expression (\ref{eq:2_3})) tous les termes impairs sont nuls
\begin{align*}
a_1=a_3=a_5=\ldots=a_{2r+1}=\ldots=0
\end{align*}
Donc la solution pour $s=n$ est donn\'ee par
\begin{align*}
y(x,n)=\sum_{r=0}^{\infty}(-1)^{r}\frac{a_0\Gamma(n+1)}{2^{2r}r!\Gamma(n+r+1)}x^{2r+n}
\end{align*}
Si on choisit $a_0$ tel que
\begin{align*}
a_0=\frac{1}{2^{n}\Gamma(n+1)}
\end{align*}
on obtient le r\'esultat final qui s'appelle fonction de Bessel (ou fonction de Bessel de premi\`ere esp\`ece)
\begin{align*}
J_n(x)=\sum_{r=0}^{\infty}(-1)^{r}\frac{1}{r!\Gamma(n+r+1)}\left(\frac{x}{2}\right)^{2r+n}
\end{align*}
La Fig.~(\ref{fig:bessel_j})  repr\'esente les courbes des fonctions de Bessel $J_n(x)$
pour $n=0,1,2,3,4$.
La deuxi\`eme solution ind\'ependante not\'ee $J_{-n}(x)$ correspond \`a $s=-n$
\begin{align}\label{eq:2_4}
J_{-n}(x)=\sum_{r=0}^{\infty}(-1)^{r}\frac{1}{r!\Gamma(-n+r+1)}\left(\frac{x}{2}\right)^{2r-n}
\end{align}
\subsection{Fonction de Bessel de deuxi\`eme esp\`ece}
\begin{figure}[!ht]
\centering
\includegraphics[scale=0.7]{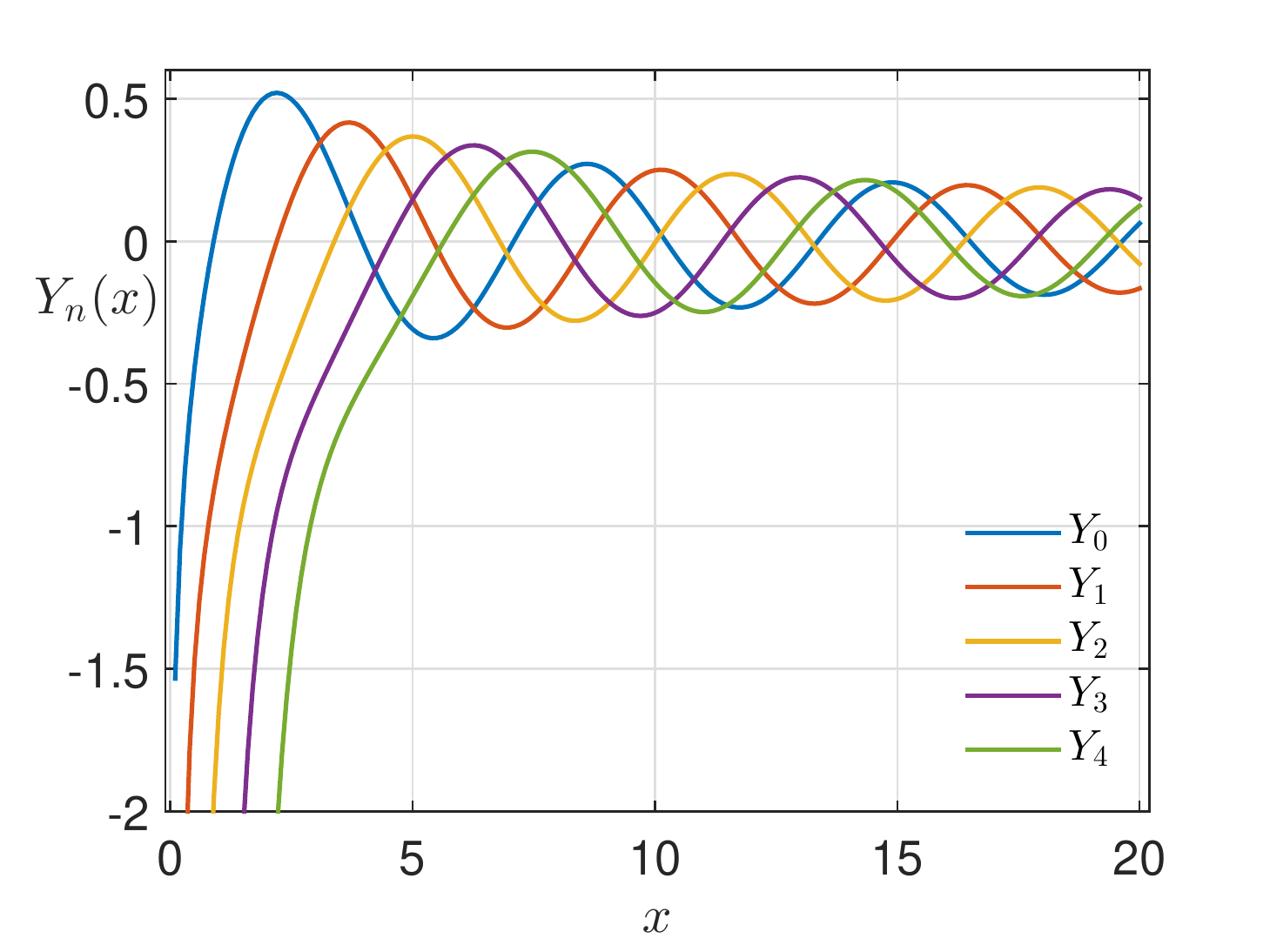}
\caption{\label{fig:bessel_y} Trac\'e de la fonction de Bessel de deuxi\`eme esp\`ece $Y_n(x)$ pour $n=0,1,2,3,4$.}
\end{figure}

Il est plus pratique dans beaucoup d'applications d'utiliser au lieu de la fonction $J_{-n}(x)$, une nouvelle fonction not\'ee $Y_n(x)$ et appell\'ee fonction de Bessel de de deuxi\'eme esp\`ece
(\'egalement appel\'ee fonction de Neumann ou encore fonction de Weber-Schl\"afli). Elle est d\'efinie par
\begin{align}\label{eq:2_5}
Y_n(x)=\frac{\cos(n\pi)J_n(x) - J_{-n}(x)}{\sin(n\pi)}.
\end{align}
Dans la suite, on prendra comme solutions ind\'ependantes de l'\'equation de Bessel (\ref{eq:2_1}) $J_n(x)$ et $Y_n(x)$.\\
\textit{D\'emonstration}\\
Dans le cas o\`u $n\notin\mathbb{N}$ on a $\sin(n\pi)\neq 0$ et on peut toujours \'ecrire $Y_n(x)$
comme une combinaison lineare de $J_n(x)$ et $J_{-n}(x)$
\begin{align*}
Y_n(x) &= \cot(n\pi)J_n(x)-\frac{1}{\sin(n\pi)}J_{-n}(x)\\
&=AJ_n(x) + BJ_{-n}(x)
\end{align*}
o\`u $A$ et $B$ sont des constantes. Cette combinaison lineare des deux solutions est aussi une solution.\\
Dans le cas o\`u $n$ est un entier naturel \i.e., $n\in\mathbb{N}$, le num\'erateur aussi bien que le d\'enominateur sont nuls, car $\sin(n\pi)= 0$ et $\cos(n\pi)J_n(x) - J_{-n}(x)=(-1)^{n}J_n(x) - J_{-n}(x)=0$ (voir Eq.~(\ref{eq:2_6}) qu'on va d\'emontrer dans la suite). Afin que nous puissions utiliser la
r\`egle de l'H\^opital, on definit $Y_n(x)$ comme suit
\begin{align*}
Y_n(x)&=\lim_{\nu\rightarrow n}Y_\nu (x)\\
&=\lim_{\nu\rightarrow n}\frac{\cos(\nu\pi)J_\nu(x) - J_{-\nu}(x)}{\sin(\nu\pi)}.
\end{align*}
En utilisant le r\`egle de l'H\^opital, on obtient
\begin{align*}
Y_n(x)&=\lim_{\nu\rightarrow n}\frac{\frac{\partial}{\partial\nu}\left[\cos(\nu\pi)J_\nu(x) - J_{-\nu}(x)\right]}{\frac{\partial}{\partial\nu}\sin(\nu\pi)}\\
&=\lim_{\nu\rightarrow n}\frac{-\pi\sin(\nu\pi)J_\nu(x) + \cos(\nu\pi)\frac{\partial}{\partial\nu}J_\nu(x) - \frac{\partial}{\partial\nu}J_{-\nu}(x)}{\pi\cos(\nu\pi)}.
\end{align*}
Apr\'es simplifications, on trouve
\begin{align*}
Y_n(x)=\frac{1}{\pi}\left[\frac{\partial}{\partial\nu}J_\nu(x) - (-1)^{n}\frac{\partial}{\partial\nu}J_{-\nu}(x)\right]_{\nu=n}.
\end{align*}
Comme $J_\nu(x)$ et $J_{-\nu}(x)$ sont deux solutions de l'\'equation de Bessel, donc elles 
verifient l'\'equation de Bessel
\begin{align*}
&x^{2}\frac{d^{2}}{dx^{2}}J_\nu(x) + x\frac{d}{dx}J_\nu(x) + (x^{2}-\nu^{2})J_\nu(x)=0\\
&x^{2}\frac{d^{2}}{dx^{2}}J_{-\nu}(x) + x\frac{d}{dx}J_{-\nu}(x) + (x^{2}-\nu^{2})J_{-\nu}(x)=0
\end{align*}
Si on d\'erive les deux \'equations pr\'ec\'edentes par rapport \`a $\nu$ on obtient
\begin{align*}
&x^{2}\frac{d^{2}}{dx^{2}}\left(\frac{\partial}{\partial\nu}J_\nu(x)\right) + x\frac{d}{dx}\left(\frac{\partial}{\partial\nu}J_\nu(x)\right) + (x^{2}-\nu^{2})\frac{\partial}{\partial\nu}J_\nu(x)-2\nu J_\nu(x)=0\\
&x^{2}\frac{d^{2}}{dx^{2}}\left(\frac{\partial}{\partial\nu}J_{-\nu}(x)\right) + x\frac{d}{dx}\left(\frac{\partial}{\partial\nu}J_{-\nu}(x)\right) + (x^{2}-\nu^{2})\frac{\partial}{\partial\nu}J_{-\nu}(x)-2\nu J_{-\nu}(x)=0,
\end{align*}
En multipliant la deuxi\`eme \'equation par $(-1)^{\nu}$ et on retranchant l'une de l'autre, on obtient
\begin{align*}
&x^{2}\frac{d^{2}}{dx^{2}}\left\lbrace\frac{\partial}{\partial\nu}J_\nu(x)-(-1)^{\nu}\frac{\partial}{\partial\nu}J_{-\nu}(x)\right\rbrace + x\frac{d}{dx}\left\lbrace\frac{\partial}{\partial\nu}J_\nu(x)-(-1)^{\nu}\frac{\partial}{\partial\nu}J_{-\nu}(x)\right\rbrace \nonumber \\
&+ (x^{2}-\nu^{2})\left\lbrace\frac{\partial}{\partial\nu}J_\nu(x)-(-1)^{\nu}\frac{\partial}{\partial\nu}J_{-\nu}(x)\right\rbrace - 2\nu\left\lbrace J_\nu(x) - (-1)^{\nu}J_{-\nu}(x)\right\rbrace=0.
\end{align*}
Quand $\nu\rightarrow n$, le dernier terme s'annule \`a cause de la formule (\ref{eq:2_6}), on obtient finalement
\begin{align*}
x^{2}\frac{d^{2}}{dx^{2}}Y_n(x) + x\frac{d}{dx}Y_n(x) + (x^{2} - n^{2})Y_n(x)=0.
\end{align*}
Ce qui prouve que $Y_n(x)$ est une solution de l'\'equation de Bessel. La repr\'esentation graphique de $Y_n(x)$ pour $n=0,1,2,3,4$ est donn\'ee dans la Fig.~(\ref{fig:bessel_y}) ci-dessous.
\subsection{Propri\'et\'es}
\textbf{1. } Si $n\in \mathbb{N}$, les deux solutions de l'\'equation de Bessel ne sont plus ind\'ependantes, elles sont li\'ees par la formule
\begin{align}\label{eq:2_6}
J_{-n}(x) = (-1)^{n}J_n(x)
\end{align}
\textit{D\'emonstration}\\
Si $n$ est un entier non n\'egatif ($n>0$), on a
\begin{align}\label{eq:2_7}
J_{-n}(x)=\sum_{r=0}^{\infty}(-1)^{r}\frac{1}{r!\Gamma(-n+r+1)}\left(\frac{x}{2}\right)^{2r-n}
\end{align}
mais comme $-n+r+1 \in \mathbb{N}$, $\Gamma(-n+r+1)=\infty$ pour tout $r$ telle que $-n+r+1<0$ (i.e pour $r=0,1,2,\ldots,n-1$). Donc
\begin{align*}
J_{-n}(x)=\sum_{r=n}^{\infty}(-1)^{r}\frac{1}{r!\Gamma(-n+r+1)}\left(\frac{x}{2}\right)^{2r-n}.
\end{align*}
En faisant le changement de variable $m=r-n$, on trouve
\begin{align*}
J_{-n}(x)&=\sum_{m=0}^{\infty}(-1)^{m+n}\frac{1}{(m+n)!\Gamma(m+1)}\left(\frac{x}{2}\right)^{2(m+n)-n}\nonumber\\
&=(-1)^{n}\sum_{m=0}^{\infty}(-1)^{m}\frac{1}{m!\Gamma(n+m+1)}\left(\frac{x}{2}\right)^{2m+n}\nonumber\\
&=(-1)^{n}J_n(x)
\end{align*}
o\`u dans la deuxi\`eme ligne de l'\'equation pr\'ec\'edente, on a utilis\'e
\begin{align*}
(m+n)!\Gamma(m+1)&=(m+n)(m+n-1)\ldots(m+1)m!\Gamma(m+1)\nonumber\\
&=m!\Gamma(n+m+1).
\end{align*}
Si $n<0$, on pose $n=-p$ avec $p>0$, on a alors
\begin{align*}
J_p(x)=(-1)^{-p}J_{-p}(x)
\end{align*}
ou de mani\`ere \'equivalente
\begin{align*}
(-1)^{p}J_p(x)=J_p(x).
\end{align*}
\textbf{2. } Si $n\in \mathbb{N}$ on a toujours
\begin{align}\label{eq:2_8}
Y_{-n}(x)=(-1)^{n}Y_n(x)
\end{align}
\textit{D\'emonstration}
\begin{align*}
Y_{-n}(x)&=\frac{1}{\pi}\left[\frac{\partial}{\partial (\nu)}J_{\nu}(x) - (-1)^{-n}\frac{\partial}{\partial\nu}J_{-\nu}(x)\right]_{\nu=-n}\nonumber\\
&=\frac{1}{\pi}\left[\frac{\partial}{\partial(-\nu)}J_{(-\nu)}(x) -
 (-1)^{-n}\frac{\partial}{\partial(-\nu)}J_{\nu}(x)\right]_{\nu=n}\nonumber\\
&=(-1)^{n}\frac{1}{\pi}\left[\frac{\partial}{\partial\nu}J_{\nu}(x) -
 (-1)^{n}\frac{\partial}{\partial\nu}J_{-\nu}(x)\right]_{\nu=n}\nonumber\\
&=(-1)^{n}Y_n(x)
\end{align*}
Ou on a utilis\'e le fait que $(-1)^{-n}=(-1)^{n}$.
\subsection{Fonction g\'en\'eratrice}
La fonction g\'en\'eratrice des fonctions de Bessel est donn\'ee par
\begin{align}\label{eq:2_9}
e^{\frac{x}{2}(t-\frac{1}{t})}=\sum_{n=-\infty}^{+\infty}t^{n}J_n(x)
\end{align}
\textit{D\'emonstration}
\begin{align*}
e^{\frac{x}{2}(t-\frac{1}{t})}&=e^{\frac{xt}{2}}e^{-\frac{x}{2t}}\nonumber\\
&=\sum_{s=0}^{\infty}\frac{1}{s!}\left(\frac{xt}{2}\right)^{s}\sum_{r=0}^{\infty}\frac{1}{r!}\left(-\frac{x}{2t}\right)^{r}\\
&=\sum_{r,s=0}^{\infty}\frac{x^{s}t^{s}(-1)^{r}x^{r}t^{-r}}{2^{s}2^{r}r!s!}\\
&=\sum_{r,s=0}^{\infty}(-1)^{r}\frac{x^{r+s}t^{s-r}}{2^{r+s}r!s!}
\end{align*}
En faisant le changement de variable $n=s-r$, on trouve
\begin{align*}
e^{\frac{x}{2}(t-\frac{1}{t})}&=\sum_{n=-\infty}^{\infty}t^{n}\sum_{r=0}^{\infty}(-1)^{r}\frac{1}{r!(r+n)!}\left(\frac{x}{2}\right)^{2r+n}\\
&=\sum_{n=-\infty}^{\infty}t^{n}\sum_{r=0}^{\infty}(-1)^{r}\frac{1}{r!\Gamma(n+r+1)}\left(\frac{x}{2}\right)^{2r+n}\\
&=\sum_{n=-\infty}^{\infty}t^{n}J_n(x)
\end{align*}
\subsection{Repr\'esentations int\'egrales}
\begin{align}
&(a)\hspace{1.5mm} J_n(x)=\frac{1}{\pi}\int_0^{\pi}\cos\left(n\phi - x\sin\phi\right)d\phi \hspace{30mm}n\in \mathbb{N}\label{eq:2_10}\\
&(b)\hspace{1.5mm}J_n(x)=\frac{\left(\frac{1}{2}x\right)^{n}}{\sqrt{\pi}\Gamma(n+\frac{1}{2})}\int_{-1}^{1}(1-t^{2})^{n-\frac{1}{2}}e^{ixt}dt\hspace{20mm} \left(n>-\frac{1}{2}\right)\label{eq:2_11}
\end{align}
\textit{D\'emonstration}\\
(a) Si $n\in \mathbb{N}$, on a $J_{-n}(x)=(-1)^{n}J_n(x)$, donc on peut \'ecrire
\begin{align*}
e^{\frac{1}{2}x\left(t-\frac{1}{t}\right)}&=J_0(x)+\sum_{n=-\infty}^{-1}t^{n}J_n(x)+\sum_{n=1}^{\infty}t^{n}J_n(x)\\
&=J_0(x)+\sum_{n=1}^{\infty}\left\lbrace t^{n}+(-1)^{n}t^{-n}\right\rbrace J_n(x)
\end{align*}
En faisant le changement de variable $t=e^{i\phi}$, donc
\begin{align*}
t-\frac{1}{t}=e^{i\phi}+e^{-i\phi}=2i\sin\phi
\end{align*}
l'\'equation pr\'ec\'edente devient
\begin{align*}
e^{ix\sin\phi}=J_0(x)+\sum_{n=1}^{\infty}\left\lbrace e^{in\phi}+(-1)^{n}e^{-in\phi}\right\rbrace J_n(x)
\end{align*}
mais on a
\begin{align*}
&e^{in\phi}+(-1)^{n}e^{-in\phi}=2\cos(n\phi)\hspace{20mm} si\hspace{1mm} n\hspace{1mm} est \hspace{1mm}pair\\
&e^{in\phi}+(-1)^{n}e^{-in\phi}=2i\sin(n\phi)\hspace{20mm} si\hspace{1mm} n\hspace{1mm} est \hspace{1mm}impair\\
\end{align*}
donc
\begin{align*}
e^{ix\sin\phi}&=J_0(x)+\sum_{n\hspace{1mm} pair}2\cos(n\phi)J_n(x)+\sum_{n\hspace{1mm} impair}2i\sin(n\phi)J_n(x)\\
&=J_0(x)+\sum_{r=1}^{\infty}2\cos(2r\phi)J_{2r}(x)+\sum_{r=1}^{\infty}2i\sin((2r-1)\phi)J_{2r-1}(x)
\end{align*}
En identifiant les parties r\'eelles et imaginaires des deux cot\'es de l'\'equation (avec $e^{ix\sin\phi}=\cos(x\sin\phi)+i\sin(x\sin\phi)$) on obtient
\begin{align}
&\cos(x\sin\phi)=J_0(x)+\sum_{r=1}^{\infty}2\cos(2r\phi)J_{2r}(x)\label{eq:2_10a}\\
&\sin(x\sin\phi)=\sum_{r=1}^{\infty}2\sin((2r-1)\phi)J_{2r-1}(x)\label{eq:2_10b}
\end{align}
En multipliant les deux cot\'es de l'Eq.~(\ref{eq:2_10a}) par $\cos(n\phi), (n>0)$ et les deux cot\'es de l'Eq.~(\ref{eq:2_10b}) par $\sin(n\phi), (n\geq 1)$, en int\'egrant de $0$ a $\pi$ et en utilisant les identiti\'es
\begin{align*}
\int_0^{\pi}\cos(m\phi)\cos(n\phi)d\phi=\left\lbrace\begin{array}{ll}
0 & m\neq n\\
\pi/2 & m=n\neq 0\\
\pi & m=n=0
\end{array}\right.
\end{align*}
et
\begin{align*}
\int_0^{\pi}\sin(m\phi)\sin(n\phi)d\phi=\left\lbrace\begin{array}{ll}
0 & m\neq n\\
\pi/2 & m=n\neq 0
\end{array}\right.
\end{align*}
on obtient
\begin{align*}
\int_0^{\pi}\cos(n\phi)\cos(x\sin\phi)d\phi=\left\lbrace\begin{array}{ll}
\pi J_n(x) & n\hspace{1mm} pair\\
0 & n\hspace{1mm} impair
\end{array}\right.
\end{align*}
et
\begin{align*}
\int_0^{\pi}\sin(n\phi)\sin(n\sin\phi)d\phi=\left\lbrace\begin{array}{ll}
0 &  n\hspace{1mm} pair\\
\pi J_n(x) & n\hspace{1mm} impair
\end{array}\right.
\end{align*}
En sommant ces deux \'equations on trouve
\begin{align*}
\int_0^{\pi}\left[\cos(n\phi)\cos(x\sin\phi)+\sin(n\phi)\sin(n\sin\phi)\right] d\phi=\pi J_n(x) \hspace{10mm} n\in \mathbb{N^{+}}
\end{align*}
donc
\begin{align*}
\int_0^{\pi}\cos\left(n\phi -x\sin\phi\right) d\phi=\pi J_n(x) \hspace{10mm} n\in \mathbb{N^{+}}
\end{align*}
Si $n$ est un entier n\'egatif ($n\in\mathbb{N^{-}}$), on fait le changement de variable $n=-m$ o\`u $m$ est positif, donc
\begin{align*}
\int_0^{\pi}\cos\left(-m\phi -x\sin\phi\right) d\phi=\pi J_{-m}(x)
\end{align*}
En faisant un deuxi\`eme changement de variable $\phi=\pi-\theta$, on trouve
\begin{align*}
\int_{\pi}^{0}&\cos\left[-m(\pi-\theta) -x\sin(\pi-\theta)\right](-d\theta)=\int_{0}^{\pi}\cos\left(-m\pi+\theta -x\sin\theta\right)d\theta\\
&=\int_{0}^{\pi}\left[\cos(m\theta-x\sin\theta)\cos(m\pi) +\sin(m\theta-x\sin\theta)\sin(m\pi)\right]d\theta\\
&=(-1)^{m}\int_{0}^{\pi}\cos(m\theta-x\sin\theta)d\theta\\
&=(-1)^{m}\pi J_m(x)\\
&=\pi J_{-m}(x)\\
&=\pi J_n(x)
\end{align*}
(b) Pour d\'emontrer la formule (\ref{eq:2_11}), consid\'erons l'int\'egrale $I$ d\'efinit par
\begin{align*}
I&=\int_{-1}^{1}\left(1-t^{2}\right)^{n-\frac{1}{2}}e^{ixt}dt\\
&=\int_{-1}^{1}\left(1-t^{2}\right)^{n-\frac{1}{2}}\sum_{r=0}^{\infty}\frac{(ixt)^{r}}{r!} dt\\
&=\sum_{r=0}^{\infty}\frac{i^{r}x^{r}}{r!}\int_{-1}^{1}\left(1-t^{2}\right)^{n-\frac{1}{2}}t^{r}dt
\end{align*}
Si $r$ est impair, la fonction \`a l'int\'erieur de l'int\'egrale est impair, donc l'int\'egrale est nulle (parce que le domaine d'int\'egration est sym\'etrique). Par contre si $r$ est pair (\i.e., $r=2s, s=0,1,2,\ldots$) on peut \'ecrire
\begin{align*}
\int_{-1}^{1}\left(1-t^{2}\right)^{n-\frac{1}{2}}t^{r}dt=2\int_{0}^{1}\left(1-t^{2}\right)^{n-\frac{1}{2}}t^{2s}dt
\end{align*}
En faisant le changement de variable $u=t^2$, on obtient
\begin{align*}
I&=\sum_{s=0}^{\infty}\frac{i^{2s}x^{2s}}{(2s)!}2\int_{0}^{1}\left(1-u\right)^{n-\frac{1}{2}}u^s\frac{du}{2\sqrt{u}}\\
&=\sum_{s=0}^{\infty}\frac{(-1)^{s}x^{2s}}{(2s)!}B\Big(n+\frac{1}{2},s+\frac{1}{2}\Big)\\
&=\sum_{s=0}^{\infty}\frac{(-1)^{s}x^{2s}}{(2s)!}\frac{\Gamma(n+\frac{1}{2})\Gamma(s+\frac{1}{2})}{\Gamma(n+s+1)}
\end{align*}
On applique la formule de duplication (\ref{eq:1_9a}) pour $\Gamma(s+1/2)$, on obtient
\begin{align*}
I&=\sum_{s=0}^{\infty}\frac{(-1)^{s}x^{2s}}{(2s)!}\frac{\Gamma(n+\frac{1}{2})}{\Gamma(n+s+1)}\frac{(2s)!\sqrt{\pi}}{s!2^{2s}}\\
&=\Gamma\left(n+\frac{1}{2}\right)\sqrt{\pi}\left(\frac{x}{2}\right)^{-n}\sum_{s=0}^{\infty}(-1)^{s}\frac{1}{s!\Gamma(n+s+1)}\left(\frac{x}{2}\right)^{2s+n}\\
&=\Gamma\left(n+\frac{1}{2}\right)\sqrt{\pi}\left(\frac{x}{2}\right)^{-n}J_n(x)
\end{align*}
donc
\begin{align*}
J_n(x)=\frac{1}{\sqrt{\pi}\Gamma\left(n+\frac{1}{2}\right)}\left(\frac{x}{2}\right)^{n}I
\end{align*}
ce qui est le r\'esultat d\'esir\'e Eq.~(\ref{eq:2_11}).
\subsection{Relations de r\'ecurrence}
\begin{align}\label{eq:2_12}
&(a)\hspace{1.5mm}\frac{d}{dx}\left[x^{n}J_n(x)\right] = x^{n}J_{n-1}(x),\nonumber\\
&(b)\hspace{1.5mm}\frac{d}{dx}\left[x^{-n}J_n(x)\right] = -x^{-n}J_{n+1}(x),\nonumber\\
&(c)\hspace{1.5mm}J_n'(x) = J_{n-1}(x) - \frac{n}{x}J_n(x),\\
&(d)\hspace{1.5mm}J_n'(x) =  \frac{n}{x}J_n(x) - J_{n+1}(x),\nonumber\\
&(e)\hspace{1.5mm}J_n'(x) =  \frac{1}{2}\left(J_{n-1}(x) - J_{n+1}(x)\right),\nonumber\\
&(f)\hspace{1.5mm}\frac{2n}{x}J_n(x)=J_{n-1}(x) + J_{n+1}(x).\nonumber
\end{align}
\textit{D\'emonstration}\\
$(a)$ En multipliyant $J_n(x)$ par $x$, on trouve
\begin{align*}
x^{n}J_n(x)=\sum_{r=0}^{\infty}(-1)^{r}\frac{1}{r!\Gamma(n+r+1)}\frac{x^{2r+2n}}{2^{2r+n}}
\end{align*}
donc
\begin{align*}
\frac{d}{dx}\left[x^{n}J_n(x)\right]&=\sum_{r=0}^{\infty}(-1)^{r}\frac{(2r+2n)}{r!(n+r)!}\frac{x^{2r+2n-1}}{2^{2r+n}}\\
&=x^{n}\sum_{r=0}^{\infty}(-1)^{r}\frac{1}{r!(r+n-1)!}\left(\frac{x}{2}\right)^{2r+n-1}\\
&=x^{n}\sum_{r=0}^{\infty}(-1)^{r}\frac{1}{r!\Gamma(r+n)}\left(\frac{x}{2}\right)^{2r+n-1}\\
&=x^{n}J_{n-1}(x)
\end{align*}
$(b)$ On a
\begin{align*}
x^{-n}J_n(x)=\sum_{r=0}^{\infty}(-1)^{r}\frac{1}{r!\Gamma(n+r+1)}\frac{x^{2r}}{2^{2r+n}}
\end{align*}
donc
\begin{align*}
\frac{d}{dx}\left[x^{-n}J_n(x)\right]&=\sum_{r=0}^{\infty}(-1)^{r}\frac{2r}{r!\Gamma(n+r+1)}\frac{x^{2r-1}}{2^{2r+n}}\\
&=\sum_{r=1}^{\infty}(-1)^{r}\frac{2r}{r!\Gamma(n+r+1)}\frac{x^{2r-1}}{2^{2r+n}}\\
\end{align*}
La derni\`ere \'egalit\'e vient du fait que le premier terme de la s\'erie ($r=0$) est nul. Par cons\'equent
\begin{align*}
\frac{d}{dx}\left[x^{-n}J_n(x)\right]&=\sum_{r=0}^{\infty}(-1)^{r+1}\frac{2(r+1)}{(r+1)!\Gamma(n+r+2)}\frac{x^{2(r+1)-1}}{2^{2(r+1)+n}}\\
&=-x^{-n}\sum_{r=0}^{\infty}(-1)^{r}\frac{1}{r!\Gamma(r+n+2)}\left(\frac{x}{2}\right)^{2r+n+1}\\
&=-x^{-n}J_{n+1}(x)
\end{align*}
$(c)$ En d\'eveloppant la d\'eriv\'ee de $x^{n}J_n(x)$ dans la formule (a), on trouve
\begin{align*}
x^{n}J_{n-1}(x)&=\frac{d}{dx}\left[x^{n}J_n(x)\right] \\
&= nx^{n-1}J_n(x) + x^{n}\frac{d}{dx}J_n(x)
\end{align*}
d'o\`u
\begin{align*}
\frac{d}{dx}J_n(x) = J_{n-1}(x) - \frac{n}{x}J_n(x)
\end{align*}
$(d)$ En d\'eveloppant la d\'eriv\'ee de $x^{-n}J_n(x)$ dans (b), on obtient
\begin{align*}
-x^{-n}J_{n+1}(x)&=\frac{d}{dx}\left[x^{-n}J_n(x)\right] \\
&= -nx^{-n-1}J_n(x)+x^{-n}\frac{d}{dx}J_n(x)
\end{align*}
donc
\begin{align*}
\frac{d}{dx}J_n(x) = \frac{n}{x}J_n(x) - J_{n+1}(x)
\end{align*}
$(e)$ Si on fait la somme de $(c)$ et $(d)$, il en r\'esulte $(e)$.\\
$(f)$ On obtient $(f)$ en faisant la soustraction entre $(c)$ et $(d)$.\\
Toutes les relations de  r\'ecurrence de $J_n(x)$ sont valables pour la fonction $Y_n(x)$ et on peut les d\'emontrer exactement de la m\^eme mani\`ere.
\section{Fonctions de Hankel}
\begin{figure}[hbt]
\centering
\includegraphics[scale=0.4]{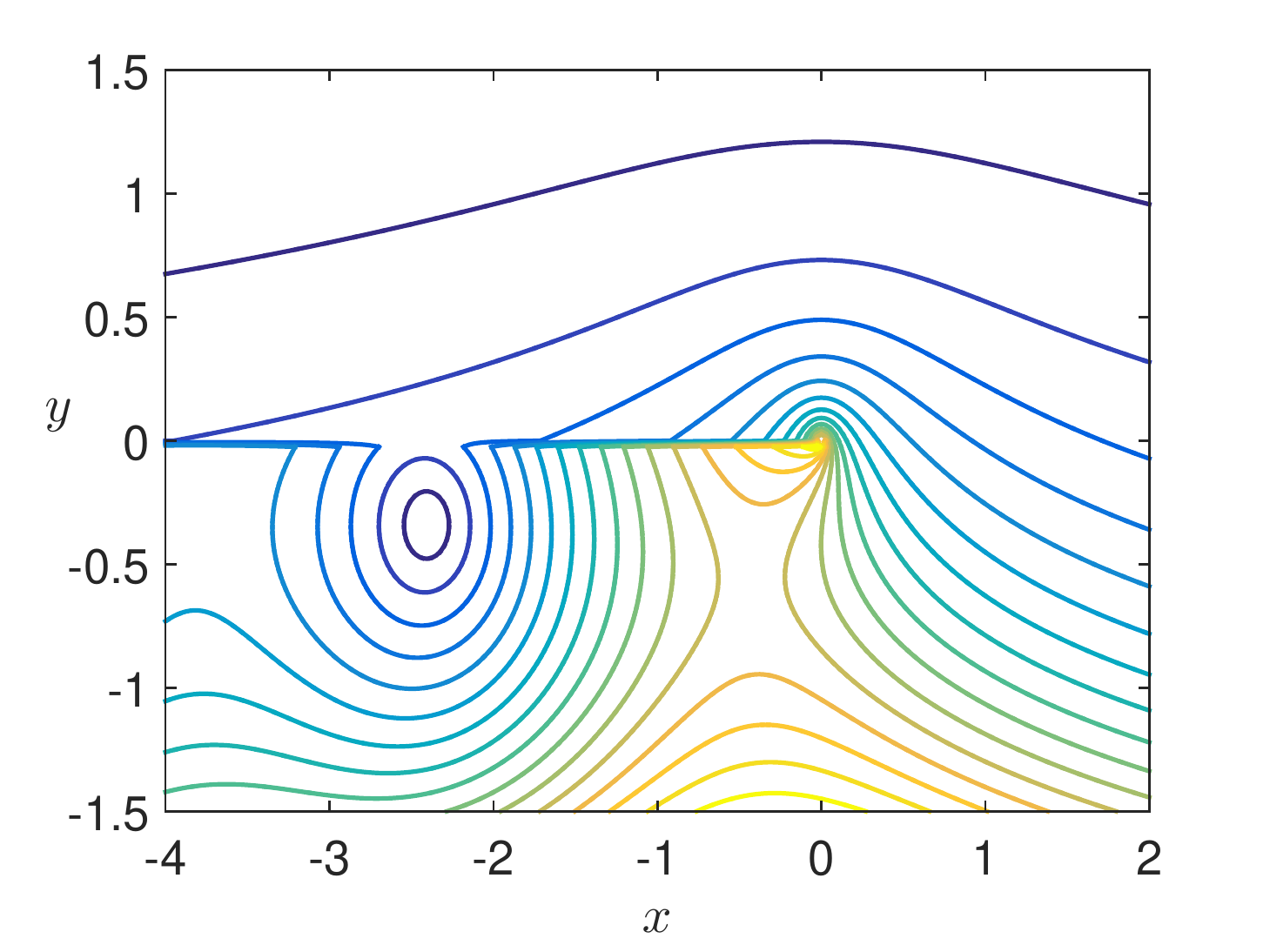}
\includegraphics[scale=0.4]{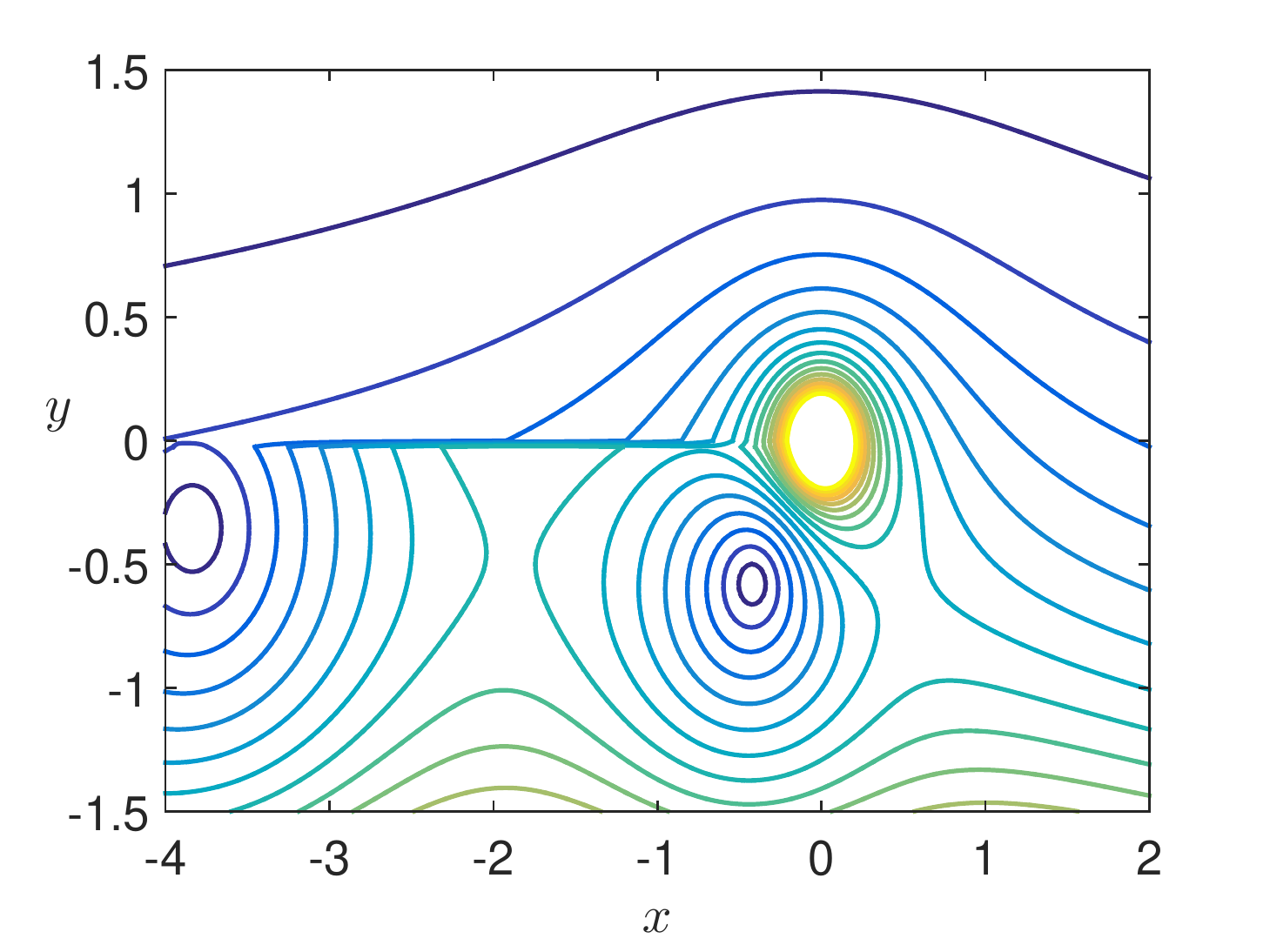}
\caption{\label{fig:hankel_h1} Trac\'e du module des fonctions de Hankel $H_n^{(1)}(z)$ dans le plan complexe pour $n=0$ (gauche) and $n=1$ (droite).}
\end{figure}

\begin{figure}[hbt]
\centering
\includegraphics[scale=0.4]{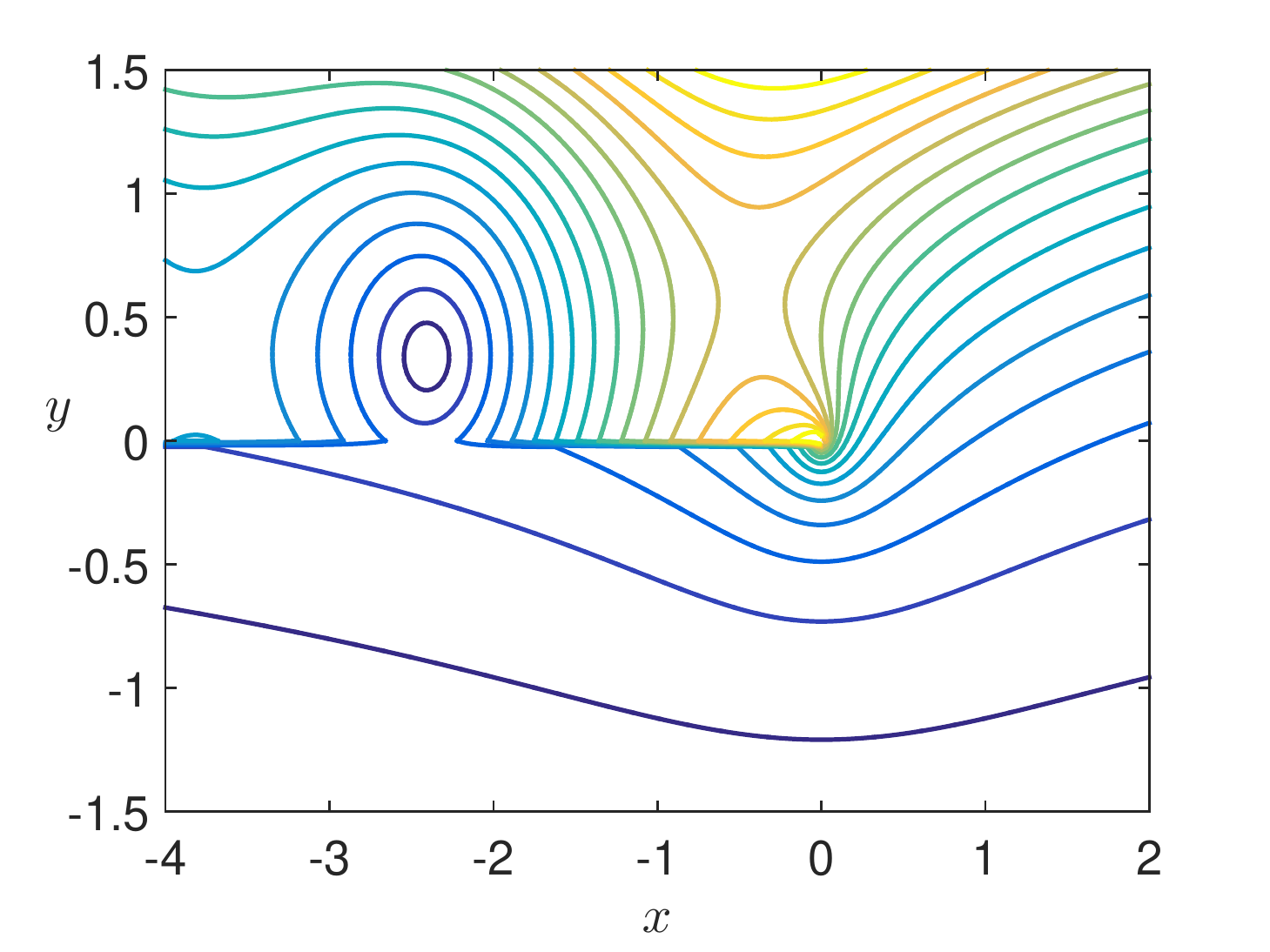}
\includegraphics[scale=0.4]{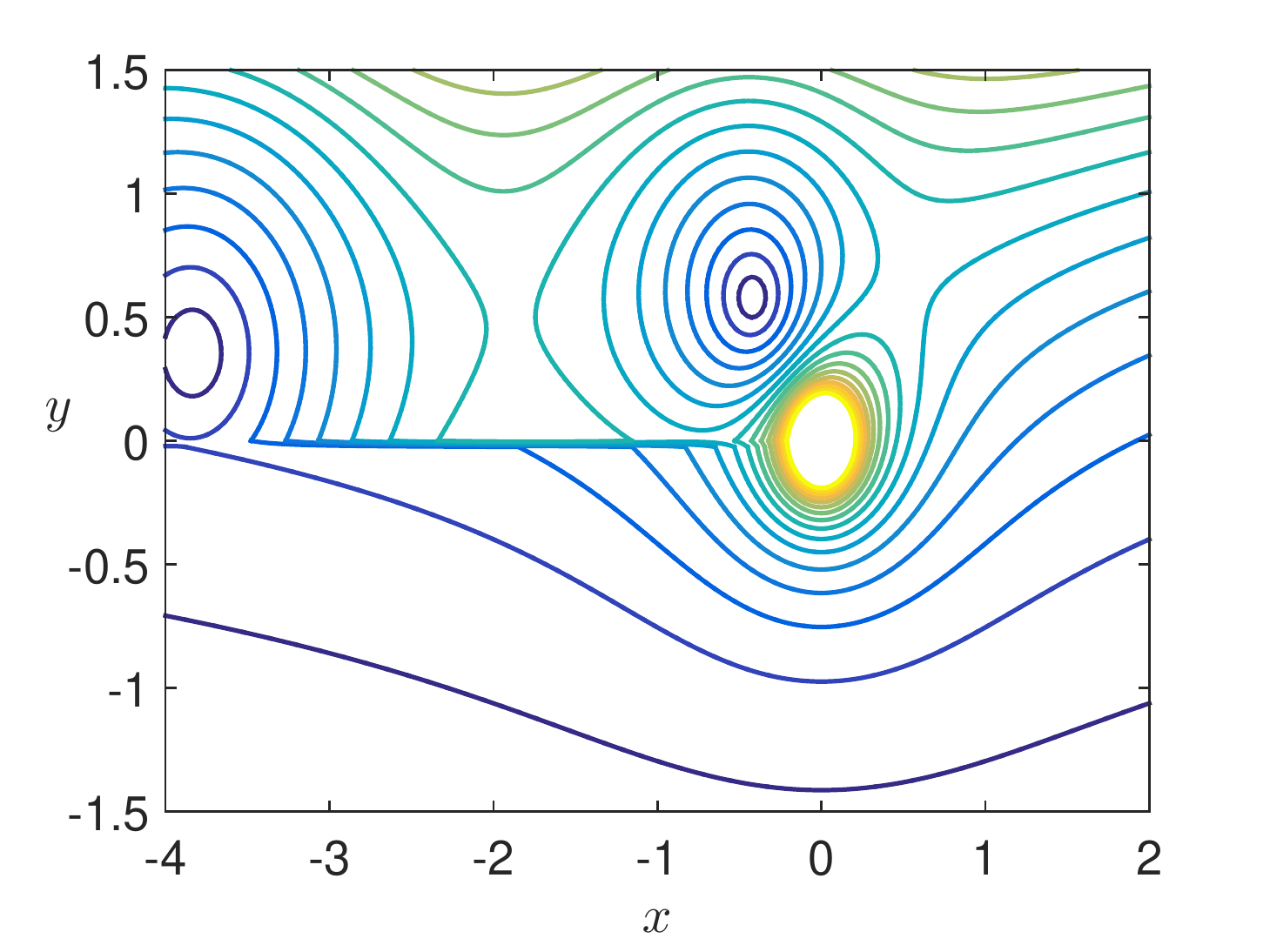}
\caption{\label{fig:hankel_h2} Trac\'e du module des fonctions de Hankel $H_n^{(2)}(z)$ dans le plan complexe pour $n=0$ (gauche) and $n=1$ (droite).}
\end{figure}

Les fonctions de Hankel de premi\`ere et de deuxi\`eme esp\`ece sont d\'efinies respectivement par
\begin{align}
&H_n^{(1)}(x) = J_n(x) + iY_n(x),\label{eq:2_13}\\
&H_n^{(2)}(x) = J_n(x) - iY_n(x)\label{eq:2_14}
\end{align}
et elles v\'erifient les m\^eme relations de r\'ecurrence que $J_n(x)$ et $Y_n(x)$. Les preuves de ces relations de r\'ecurrence sont triviaux en se basant sur les relations de r\'ecurrences de $J_n(x)$ et $Y_n(x)$. Les repr\'esentations graphiques des modules de $H_n^{(1)}$ et $H_n^{(2)}$ sont donn\'ees respectivement dans les Fig.~(\ref{fig:hankel_h1}) et (\ref{fig:hankel_h2}).
\section{Fonctions de Bessel modifi\'ees}
\begin{figure}[!ht]
\centering
\includegraphics[scale=0.8]{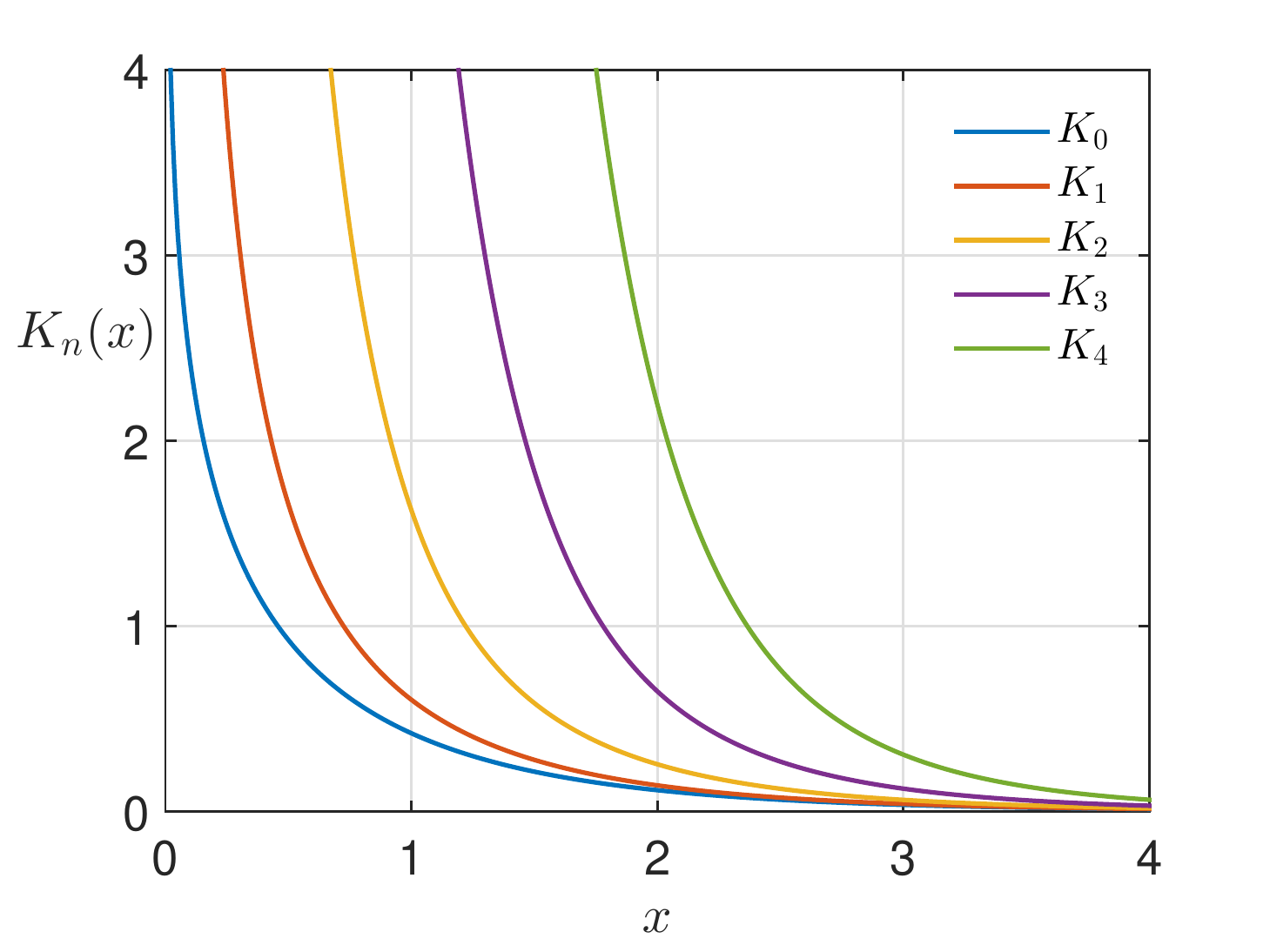}
\caption{\label{fig:bessel_k} Trac\'es des fonctions de Bessel modifi\'ees $K_n$ pour $n=0,1,2,3,4$.}
\end{figure}

L'\'equation de Bessel modifi\'ee est donn\'ee par
\begin{align}\label{eq:2_15}
x^{2}\frac{d^{2}y}{dx^{2}} + x\frac{dy}{dx} - (x^{2}+n^{2})y=0.
\end{align}
Si on utilise la m\'ethode de Frobenius \'etablit en detail pour l'\'equation de Bessel, on obtient deux solutions ind\'ependentes de l'Eq.~(\ref{eq:2_15})
\begin{align}
&I_n(x)=\sum_{r=0}^{\infty}\frac{1}{r!\Gamma(n+r+1)}\left(\frac{x}{2}\right)^{2r+n}\label{eq:2_16}\\
&I_{-n}(x)=\sum_{r=0}^{\infty}\frac{1}{r!\Gamma(-n+r+1)}\left(\frac{x}{2}\right)^{2r-n}\label{eq:2_17}
\end{align}
\begin{figure}[hbt]
\centering
\includegraphics[scale=0.8]{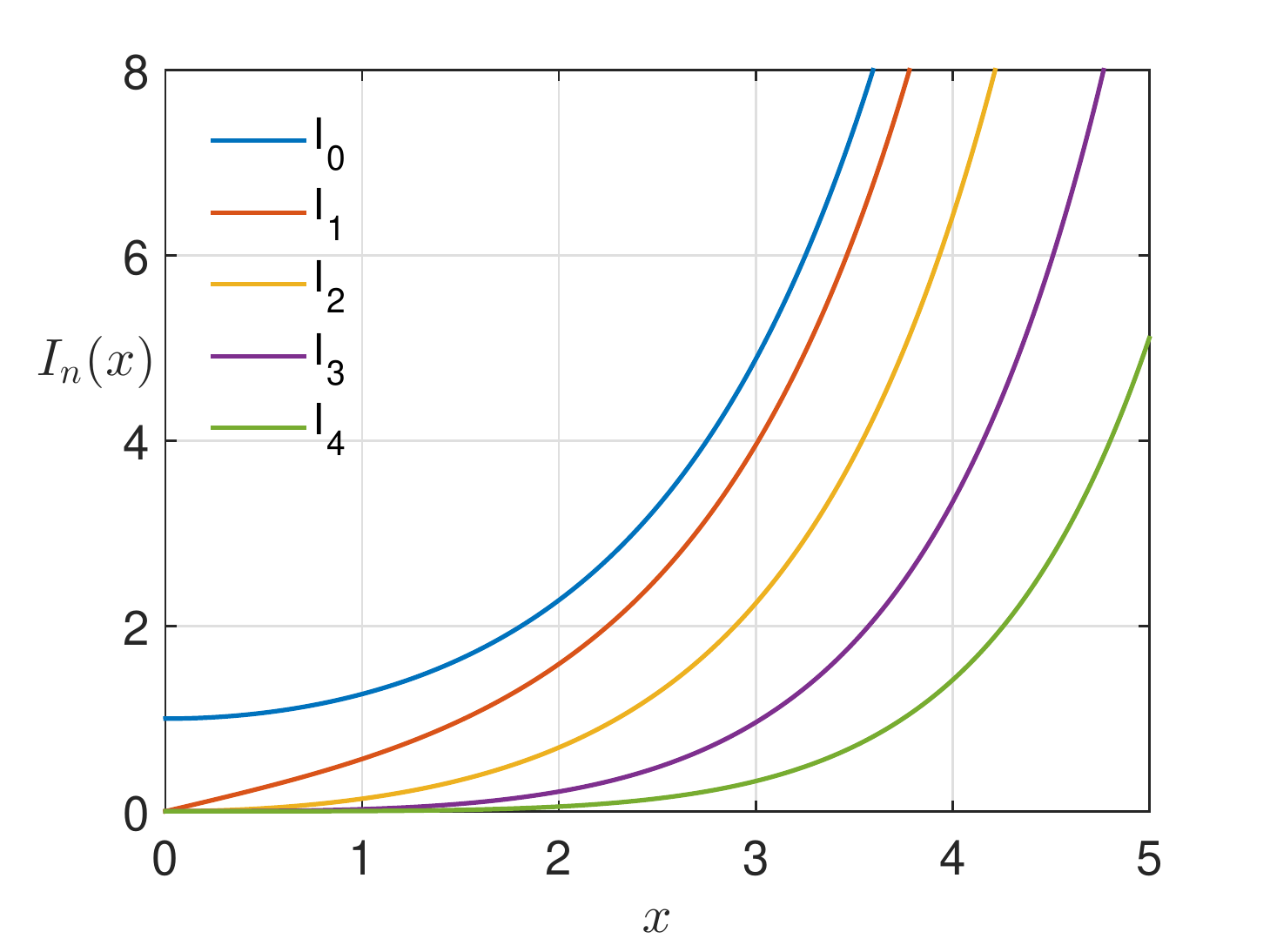}
\caption{\label{fig:bessel_i} Trac\'es des fonctions de Bessel modifi\'ees $I_n$ pour $n=0,1,2,3,4$.}
\end{figure}
Comme pour les fonctions de Bessel, il est plus pratique d'introduire au lieu de $I_{-n}(x)$, une nouvelle fonction not\'ee $K_n(x)$ et appell\'ee fonction de Bessel modifi\'ee de deuxi\'eme esp\`ece (\'egalement appel\'ee fonction de Neumann  modifi\'ee). Dans la suite, on considerera comme deux solutions ind\'ependantes de l'\'equation de Bessel  modifi\'ee (\ref{eq:2_1}), les fonctions $I_n(x)$ et $K_n(x)$. $K_n(x)$ est d\'efinit par
\begin{align}\label{eq:2_18}
K_n(x) = \frac{\pi}{2}\frac{I_{-n}(x) - I_n(x)}{\sin(n\pi)}
\end{align}
\textit{D\'emonstration}\\
La preuve est exactement la m\^eme que celle pour la formule (\ref{eq:2_5}). Les courbes 
des fonctions de Bessel  modifi\'ees $I_n(x)$ et $K_n(x)$ pour diff\'erentes valeurs de $n$ sont 
pr\'esent\'ees dans les Fig.~(\ref{fig:bessel_i}) et (\ref{fig:bessel_k}).
\subsection{Propri\'et\'es}
Il n'est pas difficile de montrer que les fonctions de Bessel  modifi\'ees $I_n(x)$ et $K_n(x)$ sont reli\'ees aux fonctions $J_n(x)$ et $Y_n(x)$
par les relations
\begin{align}
I_n(x)&=i^{-n}J_n(ix)\label{eq:2_19}\\
K_n(x) &= \frac{\pi}{2}i^{n+1}\left\lbrace J_n(ix) + iY_n(ix)\right\rbrace\label{eq:2_20}\\
&=\frac{\pi}{2}i^{n+1}H_n^{(1)}(ix)\nonumber
\end{align}
\subsection{Repr\'esentations int\'egrales}
\begin{align}
&(a)\hspace{1.5mm} I_n(x)=\frac{1}{\sqrt{\pi}\Gamma(n+\frac{1}{2})}\left(\frac{x}{2}\right)^{n}\int_{-1}^{1}e^{-xt}\left(1-t^{2}\right)^{n-\frac{1}{2}}dt \hspace{10mm} n> -\frac{1}{2}, \hspace{1mm} x> 0\label{eq:2_21}\\
&(b)\hspace{1.5mm}K_n(x)=\frac{\sqrt{\pi}}{\Gamma(n+\frac{1}{2})}\left(\frac{x}{2}\right)^{n}\int_{1}^{\infty}e^{-xt}\left(t^{2}-1\right)^{n-\frac{1}{2}}dt \hspace{10mm} n> -\frac{1}{2},\hspace{1mm} x> 0\label{eq:2_22}
\end{align}
\textit{D\'emonstration}\\
(a) D'apr\'es la formule (\ref{eq:2_19}), on a
\begin{align*}
I_n(x)=i^{-n}J_n(ix)
\end{align*}
Si on ins\`ere la repr\'esentation int\'egrale de la fonction de Bessel Eq.~(\ref{eq:2_11}) et en remplaçant $x$ par $ix$, on obtient
\begin{align*}
I_n(x)&=i^{-n}\frac{1}{\sqrt{\pi}\Gamma(n+\frac{1}{2})}i^{n}\left(\frac{x}{2}\right)^{n}\int_{-1}^1\left(1-t^{2}\right)^{n-\frac{1}{2}}e^{-xt}dt\\
&=\frac{1}{\sqrt{\pi}\Gamma(n+\frac{1}{2})}\left(\frac{x}{2}\right)^{n}\int_{-1}^1\left(1-t^{2}\right)^{n-\frac{1}{2}}e^{-xt}dt
\end{align*}
(b) Pour montrer la formule (\ref{eq:2_22}), on montre premi\`erement que l'int\'egrale
\begin{align*}
P=x^n\int_{1}^{\infty}e^{-xt}\left(t^{2}-1\right)^{n-\frac{1}{2}}dt
\end{align*}
est une solution de l'\'equation de Bessel modifi\'ee Eq.~(\ref{eq:2_15}). Donc on calcule la d\'eriv\'ee premi\`ere de $P$ ($P'=dP/dx$) et la d\'eriv\'ee seconde ($P''=d^{2}P/dx^{2}$), et on remplace $P$, $P'$ et $P''$ dans l'Eq.~(\ref{eq:2_15}). On trouve que $P$ satisfait l'\'equation de Bessel modifi\'ee
\begin{align*}
x^{2}\frac{d^{2}P}{dx^{2}} + x\frac{dP}{dx} - (x^{2}+n^{2})P=0.
\end{align*}
Donc $P$ s'\'ecrit sous la forme
\begin{align*}
P=AI_n(x)+BK_n(x)
\end{align*}
On va montrer maintenant que la constante $A$ est identiquement nulle. Pour cela, consid\'erons le cas o\`u
$x\rightarrow\infty$. Dans ce cas, $I_n(x)\rightarrow \infty$ puisque la s\'erie de $I_n(x)$ Eq.~(\ref{eq:2_16}) ne contient que des coefficients positifs. Par contre pour $P$, on utilise le fait que lorsque $x$ est tr\'es grand on a toujours $(t^2-1)^{n-\frac{1}{2}}<e^{xt/2}$, donc
\begin{align*}
P&<x^n\int_{1}^{\infty}e^{-xt}e^{xt/2}dt\\
&=2x^{n-1}e^{-x/2}
\end{align*}
Comme $P$ est toujours positive ($P>0$), l'in\'egalit\'e pr\'ec\'edente montre que lorsque $x\rightarrow\infty$ alors $P\rightarrow 0$ (parceque l'exponentielle domine la puissance). Alors, lorsque $x\rightarrow\infty$ on a $P\rightarrow 0$ et $I_n(x)\rightarrow \infty$, donc $P$ ne contient aucun multiple de $I_n(x)$ et $A$ doit \^etre nulle ($A=0$). Il en r\'esulte que
\begin{align*}
P=BK_n(x)
\end{align*}
il nous reste \`a d\'eterminer la constante $B$. pour cela, on \'etudie le comportement de l'int\'egrale $P$ et de la fonction $K_n(x)$ lorsque $x\rightarrow 0$. On consid\`ere le changement de variable $t=1+u/x$, \i.e., $dt=du/x$. Lorsque $t$ prend les valeurs $t=1$ et $t=\infty$, $u$ prend respectivement les valeurs $u=0$ et $u=\infty$, donc
\begin{align*}
P&=x^{n}\int_{0}^{\infty}e^{-x-u}\left(\frac{2u}{x}+\frac{u^{2}}{x^{2}}\right)^{n-\frac{1}{2}}\frac{du}{x}\\
&=e^{-x}x^{n}\int_{0}^{\infty}e^{-u}\left(\frac{u^2}{x^2}\right)^{n-\frac{1}{2}}\left(1+\frac{2x}{u}\right)^{n-\frac{1}{2}}\frac{1}{x}du
\end{align*}
Lorsque $x\rightarrow 0$, $(1+2u/x)^{n-1/2}\simeq 1$ et $e^{-x}\simeq 1$, donc
\begin{align*}
P\rightarrow &\frac{1}{x^n}\int_{0}^{\infty}e^{-u}u^{2n-1}du\\
&=\frac{1}{x^n}\Gamma(2n)
\end{align*}
Pour $K_n(x)$, on consid\`ere l'\'expression (\ref{eq:2_18})
\begin{align*}
K_n(x) = \frac{\pi}{2}\frac{I_{-n}(x) - I_n(x)}{\sin(n\pi)}
\end{align*}
o\`u $I_n(x)$ est donn\'ee par (\ref{eq:2_16})
\begin{align*}
I_n(x)=\sum_{r=0}^{\infty}\frac{1}{r!\Gamma(n+r+1)}\left(\frac{x}{2}\right)^{2r+n}
\end{align*}
Lorsque $x\rightarrow 0$, le terme le plus important $x$ est celui qui correspond \`a $r=0$, d'o\`u
\begin{align*}
I_n(x)\rightarrow \frac{1}{\Gamma(n+1)}\left(\frac{x}{2}\right)^{n}
\end{align*}
donc lorsque $x\rightarrow 0$, $K_n(x)$ tend vers
\begin{align*}
K_n(x)\rightarrow\frac{\pi}{2\sin(n\pi)}\frac{1}{\Gamma(-n+1)}\left(\frac{x}{2}\right)^{-n}
\end{align*}
En utilisant la formule de compl\'ement Eq.~(\ref{eq:1_8}), on obtient
\begin{align*}
K_n(x)\rightarrow\frac{\Gamma(n)2^{n-1}}{x^{n}}
\end{align*}
donc, la relation entre $P$ et $K_n(x)$ se r\'eduit \`a
\begin{align*}
\frac{1}{x^n}\Gamma(2n)=B\frac{\Gamma(n)2^{n-1}}{x^{n}}
\end{align*}
et il s'ensuit que
\begin{align*}
B=\frac{\Gamma(2n)}{\Gamma(n)2^{n-1}}
\end{align*}
En inserant la formule de duplication pour $\Gamma(2n)$, on obtient
\begin{align*}
B=\frac{2^{n}\Gamma(n+\frac{1}{2})}{\sqrt{\pi}}
\end{align*}
alors
\begin{align*}
P=\frac{2^{n}\Gamma(n+\frac{1}{2})}{\sqrt{\pi}}K_n(x)
\end{align*}
donc
\begin{align*}
K_n(x)&=\frac{\sqrt{\pi}}{2^{n}\Gamma(n+\frac{1}{2})}P\\
&=\frac{\sqrt{\pi}}{\Gamma(n+\frac{1}{2})}\left(\frac{x}{2}\right)^{n}\int_{1}^{\infty}e^{-xt}\left(t^{2}-1\right)^{n-\frac{1}{2}}dt
\end{align*}
\subsection{Relations de r\'ecurrence}
\begin{align}\label{eq:2_23}
&(a)\hspace{1.5mm}\frac{d}{dx}\left[x^{n}I_n(x)\right] = x^{n}I_{n-1}(x),\\
&(b)\hspace{1.5mm}\frac{d}{dx}\left[x^{-n}I_n(x)\right] = -x^{-n}I_{n+1}(x),\nonumber\\
&(c)\hspace{1.5mm}I_n'(x) = I_{n-1}(x) - \frac{n}{x}I_n(x),\nonumber\\
&(d)\hspace{1.5mm}I_n'(x) =  \frac{n}{x}I_n(x) - I_{n+1}(x),\nonumber\\
&(e)\hspace{1.5mm}I_n'(x) =  \frac{1}{2}\left(I_{n-1}(x) - I_{n+1}(x)\right),\nonumber\\
&(f)\hspace{1.5mm}\frac{2n}{x}I_n(x)=I_{n-1}(x) + I_{n+1}(x).\nonumber
\end{align}
\begin{align}\label{eq:2_24}
&(a)\hspace{1.5mm}\frac{d}{dx}\left[x^{n}K_n(x)\right] = -x^{n}K_{n-1}(x),\\
&(b)\hspace{1.5mm}\frac{d}{dx}\left[x^{-n}K_n(x)\right] = -x^{-n}K_{n+1}(x),\nonumber\\
&(c)\hspace{1.5mm}K_n'(x) = -K_{n-1}(x) - \frac{n}{x}K_n(x),\nonumber\\
&(d)\hspace{1.5mm}K_n'(x) =  \frac{n}{x}K_n(x) - K_{n+1}(x),\nonumber\\
&(e)\hspace{1.5mm}K_n'(x) =  -\frac{1}{2}\left(K_{n-1}(x) + K_{n+1}(x)\right),\nonumber\\
&(f)\hspace{1.5mm}-\frac{2n}{x}K_n(x)=K_{n-1}(x) - K_{n+1}(x).\nonumber
\end{align}
o\`u $I_n'(x)=dI_n(x)/dx$ et $K_n'(x)=dK_n(x)/dx$.\\
\textit{D\'emonstration}\\
Les Eqs.~(\ref{eq:2_23})(a)-(f) et Eqs.~(\ref{eq:2_24})(a)-(f) peuvent \^etre d\'emontr\'ees de la m\^eme mani\`ere que les relations de r\'ecurrence des fonctions de Bessel Eqs.~(\ref{eq:2_12})(a)-(f).
\section{Fonctions de Bessel sph\'eriques}
\begin{figure}[hbt]
\centering
\includegraphics[scale=1.1]{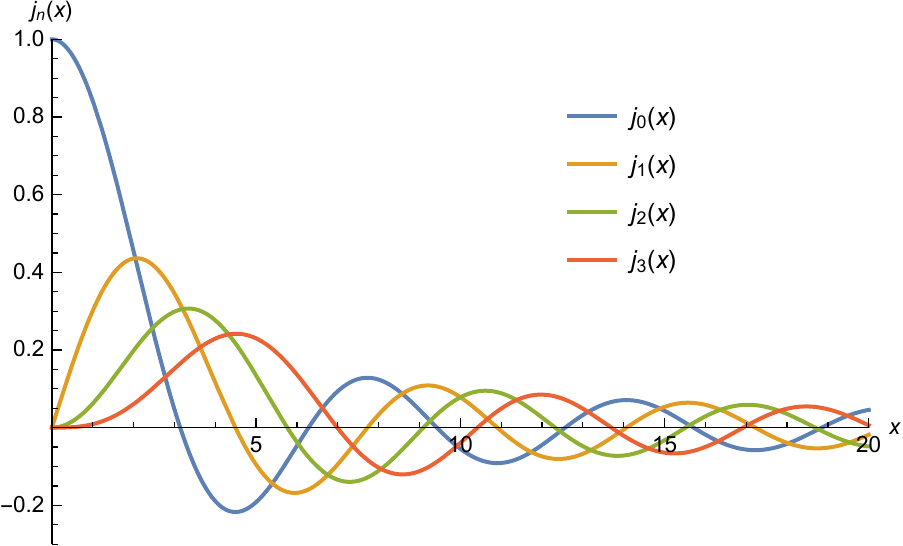}
\caption{\label{fig:sphericalbessel_j} Trac\'es des fonctions de Bessel sph\'eriques de premi\`eme esp\`ece  $j_n(x)$ pour $n=0,1,2,3$.}
\end{figure}

\begin{figure}[hbt]
\centering
\includegraphics[scale=1.1]{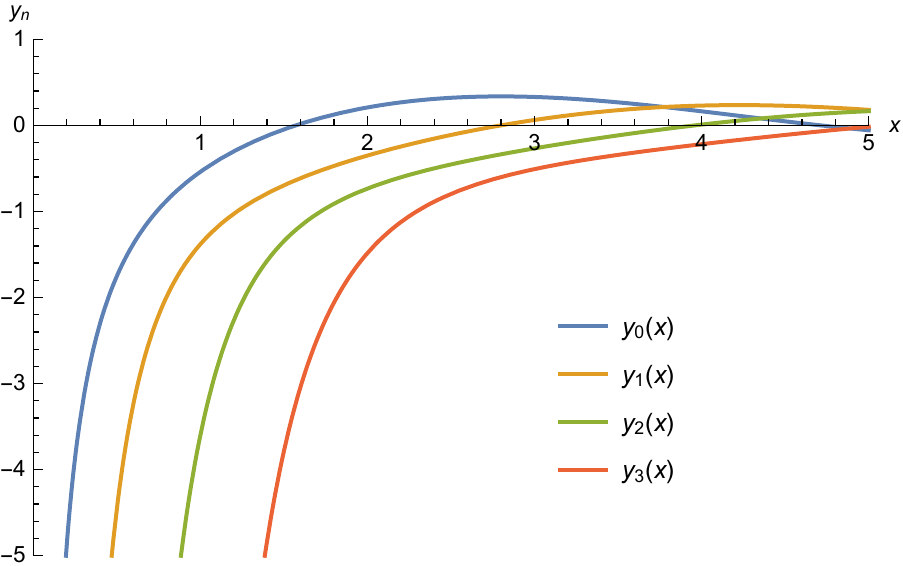}
\caption{\label{fig:sphericalbessel_y} Trac\'es des fonctions de Bessel sph\'eriques de deuxi\`eme esp\`ece $y_n(x)$ pour $n=0,1,2,3$.}
\end{figure}

L'\'equation de Bessel sph\'erique est donn\'ee par
\begin{align}\label{eq:2_25}
x^{2}\frac{d^{2}y}{dx^{2}} + x\frac{dy}{dx} + \left[x^{2}-n(n+1)\right] y=0, \hspace{10mm}{n\in \mathbb{N}}
\end{align}
Cette \'equation est solvable par la m\'ethode de Frobenius et admet deux solutions ind\'ependentes
not\'ees $j_n(x)$ et $y_n(x)$ et donn\'ees par
\begin{align}
&j_n(x)= \sqrt{\frac{\pi}{2x}}J_{n+\frac{1}{2}}(x)\label{eq:2_26}\\
&y_n(x)= \sqrt{\frac{\pi}{2x}}Y_{n+\frac{1}{2}}(x)\label{eq:2_27}
\end{align}
$j_n(x)$ sont appell\'ees fonctions de Bessel sph\'eriques de deuxi\`eme esp\`ece et $y_n(x)$
sont appell\'ees fonctions de Bessel sph\'eriques de premi\`ere esp\`ece (ou fonctions de Neuman sph\'eriques). Elles sont repr\'esent\'ees dans les Fig.~(\ref{fig:sphericalbessel_j}) et~(\ref{fig:sphericalbessel_y}).
\section{Fonctions de Hankel sph\'eriques}
Les fonctions de Hankel sph\'eriques not\'ees $h_n^{(1)}(x)$ et $h_n^{(2)}(x)$ sont d\'efinies par
\begin{align}
&h_n^{(1)}(x) = j_n(x) + iy_n(x)\label{eq:2_28}\\
&h_n^{(2)}(x) = j_n(x) - iy_n(x)\label{eq:2_29}
\end{align}
\subsection{Relations de r\'ecurrence}
Toutes les fonctions $j_n(x)$, $y_n(x)$, $h_n^{(1)}(x)$ et $h_n^{(2)}(x)$ v\'erifient
les m\^emes relations de r\'ecurrence. Si on  note par $f_n(x)$ l'une de ces fonctions
\begin{align*}
f_n(x) \equiv j_n(x), \hspace{1mm}y_n(x), \hspace{1mm} h_n^{(1)}(x), \hspace{1mm}h_n^{(2)}(x).
\end{align*}
on a
\begin{align}\label{eq:2_30}
&(a)\hspace{1.5mm}\frac{d}{dx}\left[x^{n+1}f_n(x)\right] = x^{n+1}f_{n-1}(x),\\
&(b)\hspace{1.5mm}\frac{d}{dx}\left[x^{-n}f_n(x)\right] = -x^{-n}f_{n+1}(x),\nonumber\\
&(c)\hspace{1.5mm}f_n'(x) = f_{n-1}(x) - \frac{n+1}{x}f_n(x),\nonumber\\
&(d)\hspace{1.5mm}f_n'(x) =  \frac{n}{x}f_n(x) - f_{n+1}(x),\nonumber\\
&(e)\hspace{1.5mm}(2n+1)\frac{d}{dx}f_n(x) =  nf_{n-1}(x) - f_{n+1}(x),\nonumber\\
&(f)\hspace{1.5mm}\frac{2n+1}{x}f_n(x)=f_{n-1}(x) + f_{n+1}(x).\nonumber
\end{align}
o\`u $f_n'(x)$ est la d\'eriv\'ee de $f_n(x)$ par rapport a $x$.\\
\textit{D\'emonstration}\\
Les Eqs.~(\ref{eq:2_30})(a)-(f) peuvent \^etre d\'emontr\'ees de la m\^eme
mani\`ere que les relations de r\'ecurrence des fonctions de Bessel Eqs.~(\ref{eq:2_12})(a)-(f).
\section{Comportement asymptotique}
Lorsque $x\rightarrow +\infty$
\begin{align*}
&J_n(x)\sim \left(\frac{2}{\pi x}\right)^{\frac{1}{2}}\cos\!\left(x-\Bigl(n+\frac{1}{2}\Bigr)\frac{\pi}{2}\right),\\
&Y_n(x)\sim \left(\frac{2}{\pi x}\right)^{\frac{1}{2}}\sin\!\left(x-\Bigl(n+\frac{1}{2}\Bigr)\frac{\pi}{2}\right),\\
&H_n^{(1)}(x)\sim \left(\frac{2}{\pi x}\right)^{\frac{1}{2}}e^{i(x-\left(n+\frac{1}{2})\frac{\pi}{2}\right)},\\
&H_n^{(2)}(x)\sim \left(\frac{2}{\pi x}\right)^{\frac{1}{2}}e^{-i\left(x-\left(n+\frac{1}{2}\right)\frac{\pi}{2}\right)},\\
&I_n(x) \sim \frac{1}{\sqrt{2\pi x}}e^{x},\\
&K_n(x) \sim \sqrt{\frac{1}{2 x}}e^{-x},\\
&j_n(x) \sim \frac{1}{x}\sin\!\left(x-\frac{n \pi}{2}\right),\\
&y_n(x) \sim \frac{1}{x}\cos\!\left(x-\frac{n \pi}{2}\right),\\
&h_n^{(1)}(x) \sim -\frac{i}{x}e^{i(x-\frac{n \pi}{2})},\\
&j_n(x) \sim \frac{1}{x}e^{-i(x-\frac{n \pi}{2})},\\
\end{align*}
Lorsque $x\rightarrow 0$ on a
\begin{align*}
&J_n(x)\sim \frac{1}{\Gamma(n+1)}\left(\frac{x}{2}\right)^{n},\\
&Y_n(x)\sim \left\lbrace\begin{array}{ll}
-\frac{1}{\pi}\Gamma(n)\left(\frac{2}{x}\right)^{n} &    n\neq 0\\
\frac{2}{\pi}\ln x                                  &    n=0
\end{array}\right.,\\
&j_n(x) \sim \frac{x^{n}}{(2n+1)!!} \hspace{10mm} n\in \mathbb{N},\\
&y_n(x) \sim -\frac{(2n-1)!!}{x^{n+1}}  \hspace{10mm} n\in \mathbb{N}.
\end{align*}
\section{Exercices}
$\mathbf{Exercice\hspace{1mm}1:}$ Montrer que
\begin{align*}
&1.\hspace{1.5mm} J_{\frac{1}{2}}(x)=Y_{-\frac{1}{2}}(x)=\sqrt{\frac{2}{\pi x}}\sin x\\
&2.\hspace{1.5mm}J_{-\frac{1}{2}}(x) = -Y_{\frac{1}{2}}(x)=\sqrt{\frac{2}{\pi x}}\cos x\\
&3.\hspace{1.5mm}H_{\frac{1}{2}}^{(1)}(x)=-iH_{-\frac{1}{2}}^{(1)}(x)=-i\sqrt{\frac{2}{\pi x}}e^{ix}\\
&4.\hspace{1.5mm}H_{\frac{1}{2}}^{(2)}(x)=iH_{-\frac{1}{2}}^{(2)}(x)=-i\sqrt{\frac{2}{\pi x}}e^{-ix}\\
&5.\hspace{1.5mm}I_{\frac{1}{2}}(x)=\sqrt{\frac{2}{\pi x}}\sinh x,\hspace{7mm} I_{-\frac{1}{2}}(x)=\sqrt{\frac{2}{\pi x}}\cosh x\\
&6.\hspace{1.5mm}K_{\frac{1}{2}}(x)=-K_{-\frac{1}{2}}(x)=\sqrt{\frac{2}{\pi x}}e^{-x}\\
\end{align*}
\textbf{Solutions}\\
$\textbf{1.} $ D'apr\'es l'Eq.~(\ref{eq:2_11}), on a 
\begin{align*}
J_{\frac{1}{2}}(x)&=\sqrt{\frac{x}{2\pi}}\frac{1}{\Gamma(1)}\int_{-1}^{1}e^{ixt}dt\\
&=\sqrt{\frac{x}{2\pi}}\frac{e^{ixt}}{ix}\Big|_{-1}^1\\
&=\sqrt{\frac{2}{\pi x}}\sin x
\end{align*}
A partir de la d\'efinition ~(\ref{eq:2_5}), on a
\begin{align*}
Y_{-\frac{1}{2}}(x)&=\frac{\cos(-\pi/2)J_{-\frac{1}{2}}(x)-J_{\frac{1}{2}}(x)}{\sin(-\pi/2)}\\
&=J_{\frac{1}{2}}(x)\\
&=\sqrt{\frac{2}{\pi x}}\sin x
\end{align*}
$\textbf{2.} $ D'apr\`es la relation (\ref{eq:2_4}) de $J_{-n}(x)$ on a
\begin{align*}
J_{-\frac{1}{2}}(x)=\sum_{r=0}^{\infty}(-1)^{r}\frac{1}{r!\Gamma(r+\frac{1}{2})}\left(\frac{x}{2}\right)^{2r-\frac{1}{2}}
\end{align*}
En utilisant la formule de duplication pour $\Gamma(r+\frac{1}{2})$, on obtient
\begin{align*}
J_{-\frac{1}{2}}(x)&=\sqrt{\frac{2}{x}}\sum_{r=0}^{\infty}(-1)^{r}\frac{2^{2r}r!}{r!\sqrt{\pi}(2r)!}\left(\frac{x}{2}\right)^{2r}\\
&=\sqrt{\frac{2}{x\pi}}\sum_{r=0}^{\infty}(-1)^{r}\frac{x^{2r}}{(2r)!}\\
&=\sqrt{\frac{2}{x\pi}}\cos x
\end{align*}
En utilisant la d\'efinition ~(\ref{eq:2_5})
\begin{align*}
Y_{\frac{1}{2}}(x)&=\frac{\cos(\pi/2)J_{\frac{1}{2}}(x)-J_{-\frac{1}{2}}(x)}{\sin(\pi/2)}\\
&=-J_{-\frac{1}{2}}(x)\\
&=-\sqrt{\frac{2}{x\pi}}\cos x
\end{align*}
$\textbf{3.} $ On a
\begin{align*}
H_{\frac{1}{2}}^{(1)}(x)&=J_{\frac{1}{2}}(x)+iY_{\frac{1}{2}}(x)\\
&=\sqrt{\frac{2}{x\pi}}\sin x - i\sqrt{\frac{2}{x\pi}}\cos x\\
&=-i\sqrt{\frac{2}{x\pi}}e^{ix}
\end{align*}
et
\begin{align*}
-iH_{-\frac{1}{2}}^{(1)}(x)&=-i\left(J_{-\frac{1}{2}}(x)+iY_{-\frac{1}{2}}(x)\right)\\
&=-i\sqrt{\frac{2}{x\pi}}\cos x + \sqrt{\frac{2}{x\pi}}\sin x\\
&=-i\sqrt{\frac{2}{x\pi}}e^{ix}
\end{align*}
$\textbf{4.} $ On a
\begin{align*}
H_{\frac{1}{2}}^{(2)}(x)&=J_{\frac{1}{2}}(x)-iY_{\frac{1}{2}}(x)\\
&=\sqrt{\frac{2}{x\pi}}\sin x + i\sqrt{\frac{2}{x\pi}}\cos x\\
&=i\sqrt{\frac{2}{x\pi}}e^{-ix}
\end{align*}
et
\begin{align*}
iH_{-\frac{1}{2}}^{(2)}(x)&=i\left(J_{-\frac{1}{2}}(x)-iY_{-\frac{1}{2}}(x)\right)\\
&=i\sqrt{\frac{2}{x\pi}}\cos x + \sqrt{\frac{2}{x\pi}}\sin x\\
&=i\sqrt{\frac{2}{x\pi}}e^{-ix}
\end{align*}
$\textbf{5.} $ En utilisant la repr\'esenation int\'egrale de fonction de Bessel modifi\'ee (\ref{eq:2_21}), on
\begin{align*}
I_{\frac{1}{2}}(x)&=\frac{1}{\sqrt{\pi}\Gamma(1)}\sqrt{\frac{x}{2}}\int_{-1}^{1}e^{-xt}dt\\
&=-\sqrt{\frac{x}{2\pi}}\frac{e^{-xt}}{x}\Big|_{-1}^1\\
&=\frac{1}{\sqrt{2\pi x}}2\sinh x\\
&=\sqrt{\frac{2}{\pi x}}\sinh x
\end{align*}
En utilisant l'\'expression (\ref{eq:2_17}) de $I_{-n}$, on a
\begin{align*}
I_{-\frac{1}{2}}(x)&=\sum_{r=0}^{\infty}\frac{1}{r!\Gamma\left(r+\frac{1}{2}\right)}\left(\frac{x}{2}\right)^{2r-\frac{1}{2}}\\
&=\sqrt{\frac{2}{x}}\sum_{r=0}^{\infty}\frac{2^{2r}r!}{r!(2r)!\sqrt{\pi}}\frac{x^{2r}}{2^{2r}}\\
&=\sqrt{\frac{2}{\pi x}}\sum_{r=0}^{\infty}\frac{x^{2r}}{(2r)!}\\
&=\sqrt{\frac{2}{\pi x}}\cosh x
\end{align*}
$\textbf{6.} $ En utilisant la d\'efinition (\ref{eq:2_18}) de $K_n$, on
\begin{align*}
K_{\frac{1}{2}}(x)&=\frac{\pi}{2}\frac{I_{-\frac{1}{2}}(x) -I_{\frac{1}{2}}(x)}{\sin(\pi/2)}\\
&=\sqrt{\frac{\pi}{2 x}}\left(\cosh x - \sinh x\right)\\
&=\sqrt{\frac{\pi}{2 x}}e^{-x}
\end{align*}
et
\begin{align*}
K_{-\frac{1}{2}}(x)&=\frac{\pi}{2}\frac{I_{\frac{1}{2}}(x) -I_{-\frac{1}{2}}(x)}{\sin(\pi/2)}\\
&=\sqrt{\frac{\pi}{2 x}}\left(\sinh x - \cosh x\right)\\
&=-\sqrt{\frac{\pi}{2 x}}e^{-x}
\end{align*}
$\mathbf{Exercice\hspace{1mm}2:}$ Montrer que
\begin{align*}
&1. \hspace{1.5mm} \int_{0}^{\infty}e^{-ax}J_{0}(bx)dx=\frac{1}{\sqrt{a^2+b^2}} \hspace{10mm}(a>0)\\
&2. \hspace{1.5mm}\int_{0}^{\infty}J_{n}(bx)dx=\frac{1}{b} \hspace{10mm}(\mbox{n est un entier non n\'egatif})\\
&3. \hspace{1.5mm}\int_{0}^{\infty}\frac{J_n(x)}{x}dx=\frac{1}{n} \\
\end{align*}
\textbf{Solutions}\\
$\textbf{1.} $ D'apr\`es l'Eq.~(\ref{eq:2_10}), on a 
\begin{align*}
J_0(x)=\frac{1}{\pi}\int_{0}^{\pi}\cos(x\sin\phi)d\phi
\end{align*}
donc
\begin{align*}
\int_{0}^{\infty}e^{-ax}J_{0}(bx)dx&=\int_{0}^{\infty}e^{-ax}\frac{1}{\pi}\int_{0}^{\pi}\cos(bx\sin\phi)d\phi dx\\
&=\frac{1}{\pi}\int_{0}^{\pi}\left[\int_{0}^{\infty}e^{-ax}\frac{1}{2}\left( e^{ibx\sin\phi}+ e^{-ibx\sin\phi}\right) dx \right]d\phi\\
&=\frac{1}{2\pi}\int_{0}^{\pi}\left[\frac{e^{-(a-ib\sin\phi)x}}{-a+ib\sin\phi}+\frac{e^{-(a+ib\sin\phi)x}}{-a-ib\sin\phi} \right]_{0}^{\infty}d\phi\\
&=\frac{1}{2\pi}\int_{0}^{\pi}\left[\frac{1}{a-ib\sin\phi}+\frac{1}{
a+ib\sin\phi} \right]d\phi\\
&=\frac{a}{\pi}\int_{0}^{\pi}\frac{1}{a^2+b^2\sin^2\phi}d\phi\\
\end{align*}
En faisant le changement de variable $u=\cot\phi = (\cos\phi/\sin\phi)$,  \i.e., $u^2=(1-\sin^2\phi)/\sin^2\phi$, donc $\sin^2\phi=1/(1+u^2)$ et $du=-d\phi/\sin^2\phi=(1+u^2)d\phi$, il en r\'esulte que
\begin{align*}
\int_{0}^{\infty}e^{-ax}J_{0}(bx)dx&=\frac{a}{\pi}\int_{+\infty}^{-\infty}\frac{du}{a^2+b^2\ + a^2u^2}\\
&=-\frac{a}{\pi}\frac{1}{a^2+b^2}\int_{-\infty}^{+\infty}\frac{du}{1+ \frac{a^2}{a^2+b^2}u^2}
\end{align*}
Posons $v^2=a^2u^2/(a^2+b^2)$, on obtient
\begin{align*}
\int_{0}^{+\infty}e^{-ax}J_{0}(bx)dx&=-\frac{a}{\pi}\frac{1}{a^2+b^2}\frac{\sqrt{a^2+b^2}}{a}\int_{-\infty}^{+\infty}\frac{dv}{1+ v^2}\\
&=-\frac{1}{\pi}\frac{1}{\sqrt{a^2+b^2}}\arccotan(v)\Big|_{-\infty}^{+\infty}\\
&=\frac{1}{\sqrt{a^2+b^2}}
\end{align*}
$\textbf{2.} $ On prouve d'abord le r\'esultat pour $n = 0$ et $n = 1$, et montrer ensuite que si le r\'esultat est vrai pour $n = N$, il est vrai aussi pour $n = N + 2$.
pour $n=0$, prenons la limite $a\rightarrow 0$ dans la premi\`re question de l'exercice $2.$, on obtient
\begin{align*}
\int_{0}^{\infty}J_{0}(bx)dx=\frac{1}{b}
\end{align*}
Pour $n=1$ en utiliant la relation (\ref{eq:2_12})(b) pour $n=0$, on a
\begin{align*}
\frac{d}{dx}J_0(x)=-J_1(x)
\end{align*}
En remplaçant $x$ par $bx$ et en int\'egrant par rapport \`a $x$ on obtient
\begin{align*}
\int_{0}^{\infty}\frac{d}{d(bx)}J_0(bx)dx = -\int_{0}^{\infty}J_1(bx)dx
\end{align*}
donc
\begin{align*}
\int_{0}^{\infty}J_1(bx)dx &= -\frac{1}{b}J_0(bx)\Big|_{0}^{\infty}\\
&=\frac{1}{b}
\end{align*}
parceque $J_0(\infty)=0$ et $J_0(0)=1$. Supposons maintenant que le r\'esultat est vrai pour $n = N$, et en int\'egrant l'Eq.~(\ref{eq:2_12})(e) entre $0$ et $\infty$, on obtient
\begin{align*}
J_n(x)\Big|_{0}^{\infty}=\frac{1}{2}\int_{0}^{\infty}\left[J_{n-1}(x)-J_{n+1}(x)\right]dx
\end{align*}
Pour $n>0$ on a $J_n(0)=J_n(\infty)=0$, donc
\begin{align*}
0=\frac{1}{2}\int_{0}^{\infty}\left[J_{n-1}(x)-J_{n+1}(x)\right]dx
\end{align*}
alors
\begin{align*}
\int_{0}^{\infty}J_{n+1}(x)dx=\int_{0}^{\infty}J_{n-1}(x)dx
\end{align*}
En remplaçant $x$ par $bx$ et $n-1$ par $n$ on trouve
\begin{align*}
\int_{0}^{\infty}J_{n+2}(bx)dx&=\int_{0}^{\infty}J_{n}(bx)dx\\
&=\frac{1}{b}
\end{align*}
Donc le r\'esultat $\int_{0}^{\infty}J_{n}(bx)dx=1/b$ est vrai quelque soit $n$.\\
$\textbf{3.} $ En int\'egrant l'Eq.~(\ref{eq:2_12})(f) de $0$ \`a $\infty$ on obtient
\begin{align*}
2n\int_{0}^{\infty}\frac{J_n(x)}{x}dx=\int_{0}^{\infty}J_{n-1}(x)dx + \int_{0}^{\infty}J_{n+1}(x)dx
\end{align*}
mais d'apr\'es la question pr\'ec\'edente, on a
\begin{align*}
\int_{0}^{\infty}J_{n}(x)dx=1\hspace{10mm} \forall n
\end{align*} 
donc
\begin{align*}
2n\int_{0}^{\infty}\frac{J_n(x)}{x}dx&=1 + 1\\
&=2
\end{align*}
alors
\begin{align*}
\int_{0}^{\infty}\frac{J_n(x)}{x}dx=\frac{1}{n}
\end{align*}
\chapter{Fonction erreur et int\'grales de Fresnel} 
\label{chap: erreur_fresnel}
\section{Fonction erreur}
La fonction erreur not\'ee $\erf$ a \'et\'e introduite en 1871 par Glaisher et elle est d\'efinit par l'int\'egrale
\begin{align}\label{eq:3_1}
\erf(x) = \frac{2}{\sqrt{\pi}}\int_0^{x}e^{-t^{2}}dt \hspace{10mm} -\infty < x < +\infty
\end{align}
Elle est beaucoup utilis\'ee dans le domaine des probabilit\'es et statistiques. Elle est li\'ee \`a la probabilit\'e pour qu'une variable normale centr\'ee r\'eduite prenne une valeure dans un intervalle $[-x,\hspace{1mm} x]$. Cette probabilit\'e est donn\'ee par $\erf(x/\!\!\sqrt{2})$. Autrement dit, lorsque les r\'esultats d'une s\'erie de mesures sont d\'ecrits par une distribution normale avec variance $\sigma$ et valeur moyenne $0$, alors $\erf(x/\!\sigma\sqrt{\pi})$ est la probabilit\'e que l'erreur d'une seule mesure se situe dans l'intervalle $[-x, \hspace{1mm}x]$. En physique, la fonction erreur est souvent utilis\'ee dans les probl\`emes de diffusion et notament dans la conduction en r\'egime transitoire. La fonction erreur compl\'ementaire est d\'efini par
\begin{align}\label{eq:3_2}
\erfc(x) =   \frac{2}{\sqrt{\pi}}\int_x^{\infty}e^{-t^{2}}dt \hspace{10mm} -\infty < x < +\infty
\end{align}
A partir des expressions (\ref{eq:3_1}) et (\ref{eq:3_2}), il est claire que $\erfc(x)$ est le compl\'ement \`a $1$ de la fonction erreur $\erf(x)$
\begin{align*}
\erfc(x) = 1- \erf(x)
\end{align*}
\subsection{Propri\'et\'es de la fonction erreur}
\begin{figure}[!ht]
\centering
\includegraphics[scale=0.7]{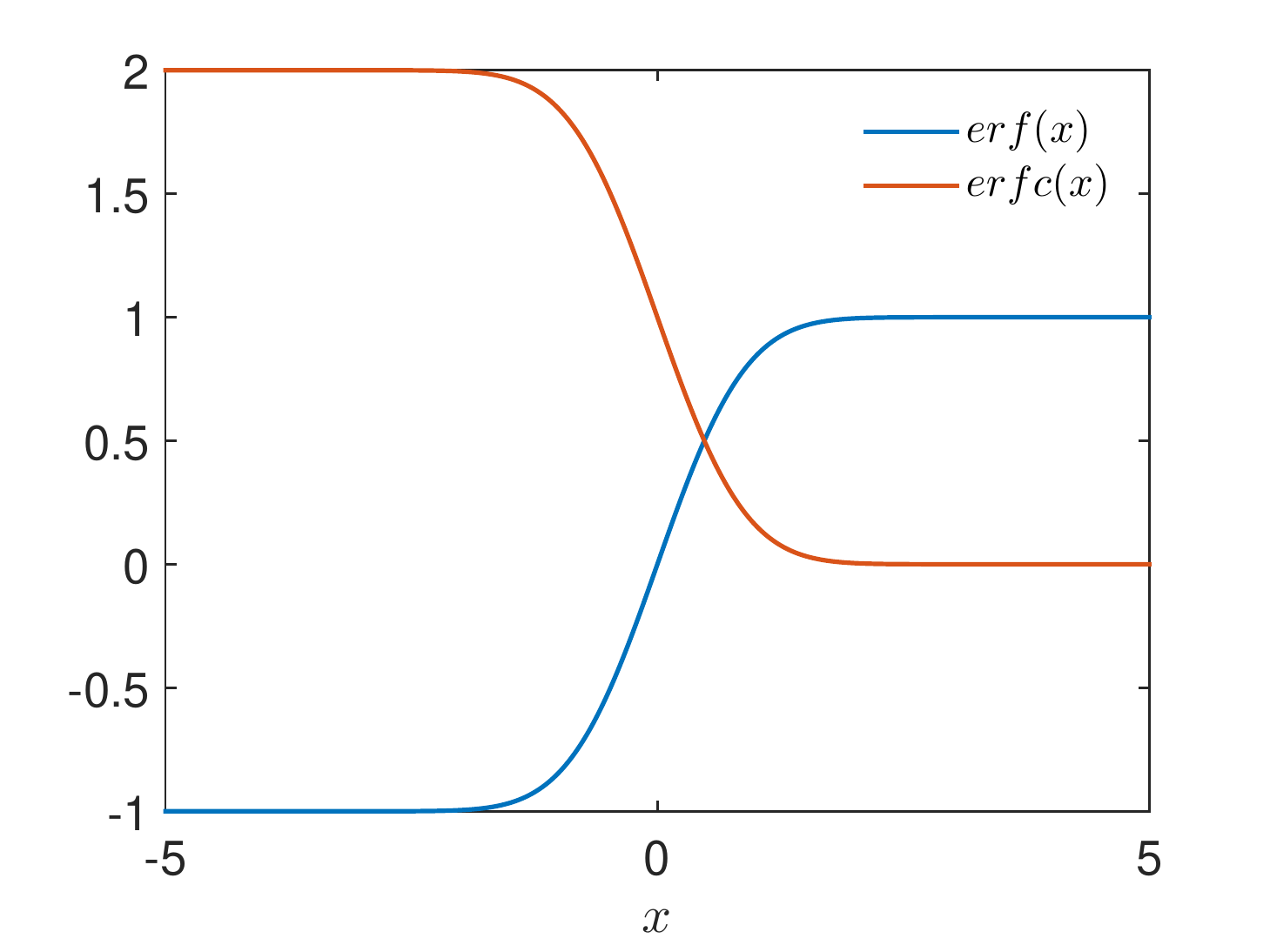}
\caption{\label{fig:ERF} Fonction erreur et fonction erreur compl\'ementaire.}
\end{figure}

\begin{align*}
&\erf(0) = 0, \hspace{5mm} \erf(\infty) = 1,\\
&\erfc(0) = 1, \hspace{5mm} \erfc(\infty) = 0,\\
&\erf(-x)=-\erf(x),\\
&\frac{d}{dx}\erf(x) = \frac{2}{\sqrt{\pi}}e^{-x^{2}},\\
&\erfc(x)\sim \frac{e^{-x^2}}{\sqrt{\pi}x} \hspace{10mm} x\rightarrow \infty,
\end{align*}
La fonction erreur et la fonction erreur compl\'ementaire sont illustr\'ees sur la Fig.~(\ref{fig:ERF}).
\subsection{Repr\'esentation en s\'erie de la fonction erreur}
\begin{align}\label{eq:3_3}
\erf(x) = \frac{2}{\sqrt{\pi}}\sum_{n=0}^{\infty}\frac{(-1)^{n}x^{2n+1}}{n!(2n+1)}\hspace{10mm} -\infty < x < +\infty
\end{align}
\textit{D\'emonstration}\\
On a
\begin{align*}
\erf(x)&=\frac{2}{\sqrt{\pi}}\int_{0}^{x}\sum_{n=0}^{\infty}\frac{(-t^{2})^{n}}{n!}dt\\
&=\frac{2}{\sqrt{\pi}}\sum_{n=0}^{\infty}\frac{(-1)^{n}}{n!}\int_0^xt^{2n}dt\\
&= \frac{2}{\sqrt{\pi}}\sum_{n=0}^{\infty}\frac{(-1)^{n}x^{2n+1}}{n!(2n+1)}
\end{align*}
\section{Int\'egrales de Fresnel}
\begin{figure}[!ht]
\centering
\includegraphics[scale=0.7]{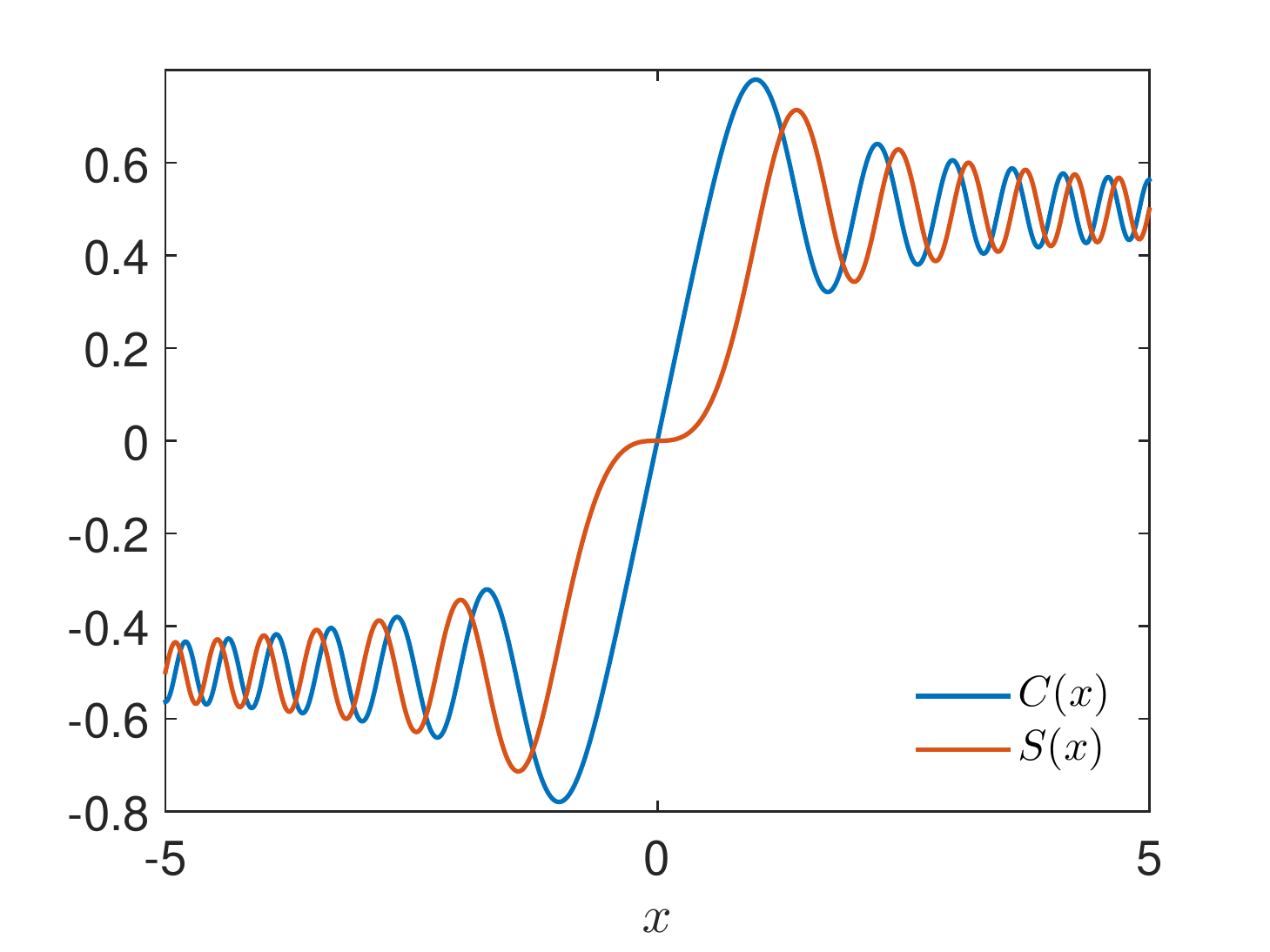}
\caption{\label{fig:C_S} Fonctions de Fresnel $C(x)$ et $S(x)$.}
\end{figure}

Les fonctions de Fresnel sont d\'efinies par les int\'egrales
\begin{align}
&C(x) = \int_0^{x}\cos\left(\frac{1}{2}\pi t^{2}\right)dt\label{eq:3_4}\\
&S(x) = \int_0^{x}\sin\left(\frac{1}{2}\pi t^{2}\right)dt\label{eq:3_5}
\end{align}
et elles sont repr\'esent\'ees dans la Fig.~(\ref{fig:C_S}).
\subsection{Properi\'et\'es des fonctions de Fresnel}
\begin{align}
&C(0) = S(0) = 0,\nonumber\\
&C(\infty) = S(\infty) = \frac{1}{2},\nonumber\\
&C(-x) = -C(x), \hspace{5mm} S(-x) = -S(x),\nonumber\\
&C'(x) = \cos\left(\frac{1}{2}\pi t^{2}\right),\nonumber\\
&S'(x) = \sin\left(\frac{1}{2}\pi t^{2}\right),\nonumber\\
&C(x) - iS(x) = \frac{1}{\sqrt{2i}}\erf\left(x\sqrt{\frac{i\pi}{2}}\right)\nonumber
\end{align}
\textit{D\'emonstration}\\
La preuve des cinq premi\`eres propri\'et\'es est simple. Pour la derni\`ere, on a
\begin{align*}
C(x) - iS(x) &=  \int_0^{x}\left\lbrace\cos\left(\frac{1}{2}\pi t^{2}\right) - i\sin\left(\frac{1}{2}\pi t^{2}\right)\right\rbrace dt\\
&=\int_0^{x}e^{-\frac{1}{2}i\pi t^{2}}dt\\
&=\int_0^{x\sqrt{\frac{i\pi}{2}}}e^{-u^2}\sqrt{\frac{2}{i\pi}} du\\
&=\frac{1}{\sqrt{2i}}\erf\left(x\sqrt{\frac{i\pi}{2}}\right)
\end{align*}
\subsection{Repr\'esentaions en s\'eries}
Les fonctions de Fresnel sont \'egalement d\'efinies par des d\'eveloppements en s\'eries enti\`eres
\begin{align}
&(a)\hspace{1.5mm} C(x) = \sum_{n=0}^{\infty}\frac{(-1)^{n}\left(\pi/2\right)^{2n}}{(2n)!(4n+1)}x^{4n+1}\label{eq:3_6}\\
&(b)\hspace{1.5mm} S(x) = \sum_{n=0}^{\infty}\frac{(-1)^{n}\left(\pi/2\right)^{2n+1}}{(2n+1)!(4n+3)}x^{4n+3}\label{eq:3_7}
\end{align}
\textit{D\'emonstration}\\
(a) Le d\'eveloppement en s\'erie de la fonction cosinus est donn\'ee par
\begin{align*}
\cos x=\sum_{n=0}^{\infty}\frac{(-1)^{n}}{(2n)!}x^{2n} 
\end{align*}
donc, \`a partir de la d\'efinition de $C(x)$, on a
\begin{align*}
C(x) &= \int_0^{x}\sum_{n=0}^{\infty}\frac{(-1)^{n}}{(2n)!}\left(\frac{1}{2}\pi t^{2}\right)^{2n}dt\\
&=\sum_{n=0}^{\infty}\frac{(-1)^{n}}{(2n)!}\left(\frac{\pi}{2}\right)^{2n}\int_0^{x}t^{4n}dt\\
&=\sum_{n=0}^{\infty}\frac{(-1)^{n}\left(\pi/2\right)^{2n}}{(2n)!(4n+1)}x^{4n+1}
\end{align*}
(b) Le d\'eveloppement en s\'erie de la fonction sinus est donn\'ee par
\begin{align*}
\sin x=\sum_{n=0}^{\infty}\frac{(-1)^{n}}{(2n+1)!}x^{2n+1}
\end{align*}
alors, on a
\begin{align*}
S(x) &= \int_0^{x}\sum_{n=0}^{\infty}\frac{(-1)^{n}}{(2n+1)!}\left(\frac{1}{2}\pi t^{2}\right)^{2n+1}dt\\
&=\sum_{n=0}^{\infty}\frac{(-1)^{n}}{(2n+1)!}\left(\frac{\pi}{2}\right)^{2n+1}\int_0^{x}t^{4n+2}dt\\
&=\sum_{n=0}^{\infty}\frac{(-1)^{n}\left(\pi/2\right)^{2n+1}}{(2n+1)!(4n+3)}x^{4n+3}
\end{align*}
\section{Exercices}
$\mathbf{Exercice\hspace{1mm}1:}$ Montrer que
\begin{align*}
&(a)\hspace{1.5mm}\erf(x)=\frac{1}{\sqrt{\pi}}\gamma\left(\frac{1}{2},x^{2}\right)\\
&(b)\hspace{1.5mm}\erfc(x)=\frac{1}{\sqrt{\pi}}\Gamma\left(\frac{1}{2},x^{2}\right)
\end{align*}
$\mathbf{Solutions}$\\
$(a)$ Posons $u=t^2$ dans la d\'efinition Eq.~(\ref{eq:3_1}) on obtient
\begin{align*}
\erf(x)&=\frac{2}{\sqrt{\pi}}\int_{0}^{x^2}e^{-u}\frac{du}{2\sqrt{u}}\\
&=\frac{1}{\sqrt{\pi}}\gamma\left(\frac{1}{2},x^{2}\right)
\end{align*}
$(b)$ D'apr\`es l'Eq.~(\ref{eq:1_13})(a) on a
\begin{align*}
\gamma\left(\frac{1}{2},x^{2}\right)+\Gamma\left(\frac{1}{2},x^{2}\right)=\Gamma\left(\frac{1}{2}\right)
\end{align*}
donc
\begin{align*}
\Gamma\left(\frac{1}{2},x^{2}\right)&=\sqrt{\pi}-\gamma\left(\frac{1}{2},x^{2}\right)\\
&=\sqrt{\pi}-\sqrt{\pi}\erf(x)\\
&=\sqrt{\pi}\erfc(x)
\end{align*}
$\mathbf{Exercice\hspace{1mm}2:}$ Montrer que
\begin{align*}
\erf\left(\frac{x}{\sqrt{2}}\right)=\frac{1}{\sqrt{\pi}}\int_x^{\infty}e^{-\frac{t^2}{2}}dt
\end{align*}
$\mathbf{Solution}$\\
On a
\begin{align*}
\erf\left(\frac{x}{\sqrt{2}}\right)=\sqrt{\frac{2}{\pi}}\int_{\frac{x}{\sqrt{2}}}^{\infty}e^{-t^2}dt
\end{align*}
En faisant le changement de variable $u=\sqrt{2}t$ on obtient le r\'esultat d\'edir\'e
\begin{align*}
\erf\left(\frac{x}{\sqrt{2}}\right)=\frac{1}{\sqrt{\pi}}\int_{x}^{\infty}e^{-\frac{u^2}{2}}du
\end{align*}
$\mathbf{Exercice\hspace{1mm}3:}$ Montrer que
\begin{align*}
&(a)\hspace{1.5mm}\int_{0}^{x}C(t)dt=xC(x)-\frac{1}{\pi}\sin\left(\frac{\pi}{2}x^2\right)\\
&(b)\hspace{1.5mm}\int_{0}^{x}S(t)dt=xS(x)+\frac{1}{\pi}\cos\left(\frac{\pi}{2}x^2\right)-\frac{1}{\pi}
\end{align*}
$\mathbf{Solutions}$\\
$(a)$ On a
\begin{align*}
\int_{0}^{x}C(t)dt&=\int_{0}^{x}\left[\int_{0}^{t}\cos\left(\frac{\pi}{2}u^2\right)du\right]dt\\
&=\int_{0}^{x}\int_{0}^{t}\sum_{n=0}^{\infty}\frac{(-1)^{n}}{(2n)!}\left(\frac{\pi}{2}u^2\right)^{2n}dudt\\
&=\int_{0}^{x}\sum_{n=0}^{\infty}\frac{(-1)^{n}}{(2n)!}\left(\frac{\pi}{2}\right)^{2n}\frac{t^{4n+1}}{4n+1}dt
\end{align*}
En int\'egrant par parties, on obtient
\begin{align*}
\int_{0}^{x}C(t)dt&=\left[ t\sum_{n=0}^{\infty}\frac{(-1)^{n}}{(2n)!}\left(\frac{\pi}{2}\right)^{2n}\frac{t^{4n+1}}{4n+1}\right]_{0}^{x}-\int_{0}^{x}t\sum_{n=0}^{\infty}\frac{(-1)^{n}}{(2n)!}\left(\frac{\pi}{2}\right)^{2n}t^{4n}dt\\
&=xC(x)-\sum_{n=0}^{\infty}\frac{(-1)^{n}}{(2n)!}\left(\frac{\pi}{2}\right)^{2n}\frac{x^{4n+2}}{2(2n+1)}\\
&=xC(x)-\frac{1}{\pi}\sum_{n=0}^{\infty}\frac{(-1)^{n}}{(2n+1)!}\left(\frac{\pi}{2}x^{2}\right)^{2n+1}\\
&=xC(x)-\frac{1}{\pi}\sin\left(\frac{\pi}{2}x^{2}\right)
\end{align*}
$(b)$ D'apr\`es l'Eq.~(\ref{eq:3_7}) on a
\begin{align*}
\int_{0}^{x}S(t)dt=\int_{0}^{x}\sum_{n=0}^{\infty}\frac{(-1)^{n}\left(\pi/2\right)^{2n+1}}{(2n+1)!(4n+3)}t^{4n+3}dt
\end{align*}
En int\'egrant par parties on obtient
\begin{align*}
\int_{0}^{x}S(t)dt&=\left[t\sum_{n=0}^{\infty}\frac{(-1)^{n}\left(\pi/2\right)^{2n+1}}{(2n+1)!(4n+3)}t^{4n+3}\right]_0^{x}-\int_{0}^{x}t\sum_{n=0}^{\infty}\frac{(-1)^{n}\left(\pi/2\right)^{2n+1}}{(2n+1)!}t^{4n+2}dt\\
&=xS(x)+ \frac{1}{\pi}\sum_{n=0}^{\infty}\frac{(-1)^{n+1}\left(\pi/2\right)^{2(n+1)}}{(2n+2)!}t^{(4n+4)}\Big|_{0}^{x}\\
&=xS(x)+ \frac{1}{\pi}\sum_{n=0}^{\infty}\frac{(-1)^{n+1}}{(2n+2)!}\left(\frac{\pi}{2}t^2\right)^{2(n+1)}\Big|_{0}^{x}\\
&=xS(x)+\frac{1}{\pi}\cos\left(\frac{\pi}{2}t^{2}\right)\Big|_{0}^{x}\\
&=xS(x)+\frac{1}{\pi}\cos\left(\frac{\pi}{2}x^{2}\right) - \frac{1}{\pi}
\end{align*}
\chapter{Exponentielle int\'egrale, sinus int\'egral et cosinus int\'egral}
\label{chap:exp_sin_cos_int}
\section{Exponentielle int\'egrale}
\begin{figure}[!ht]
\centering
\includegraphics[scale=0.42]{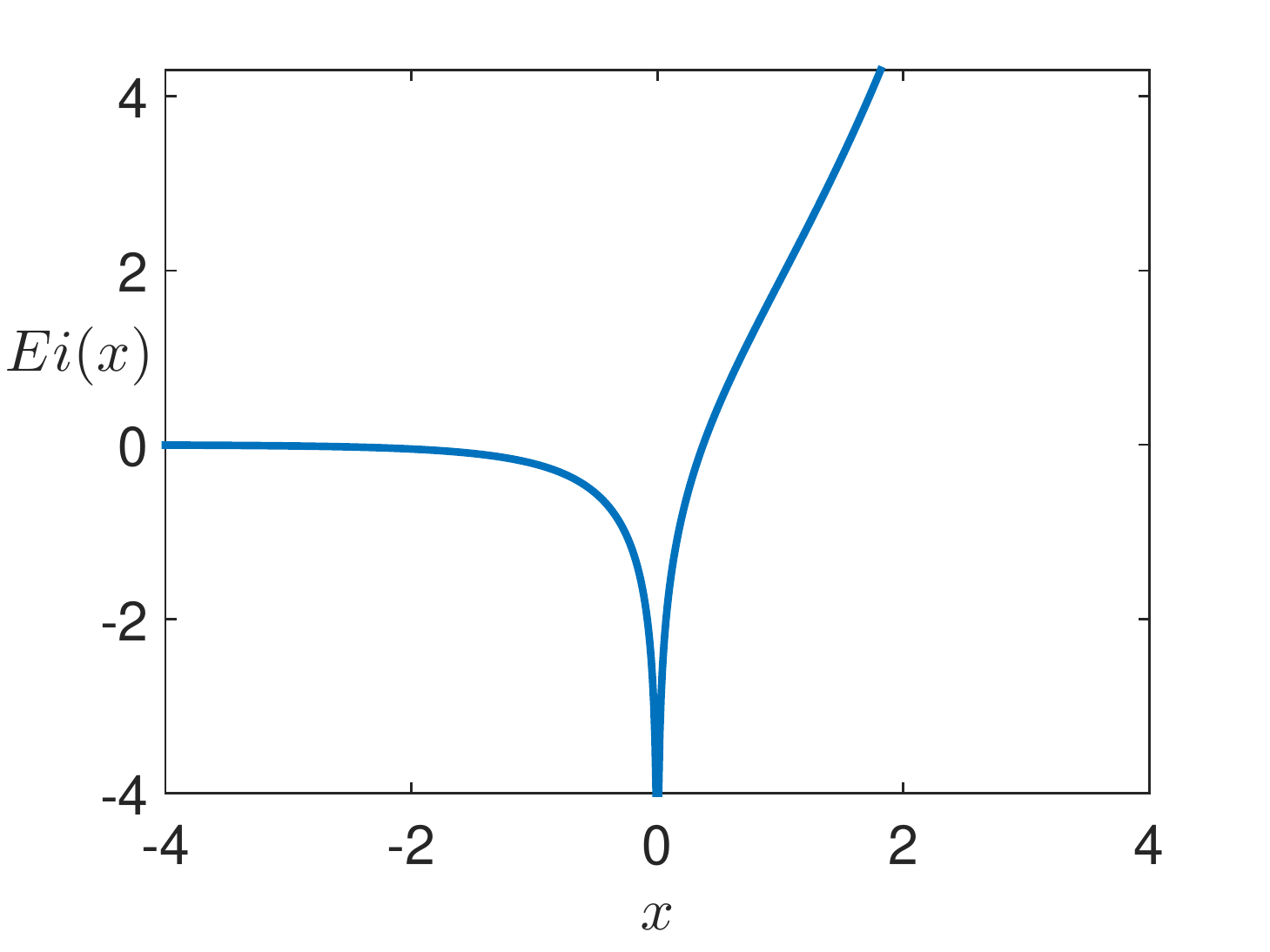}
\includegraphics[scale=0.42]{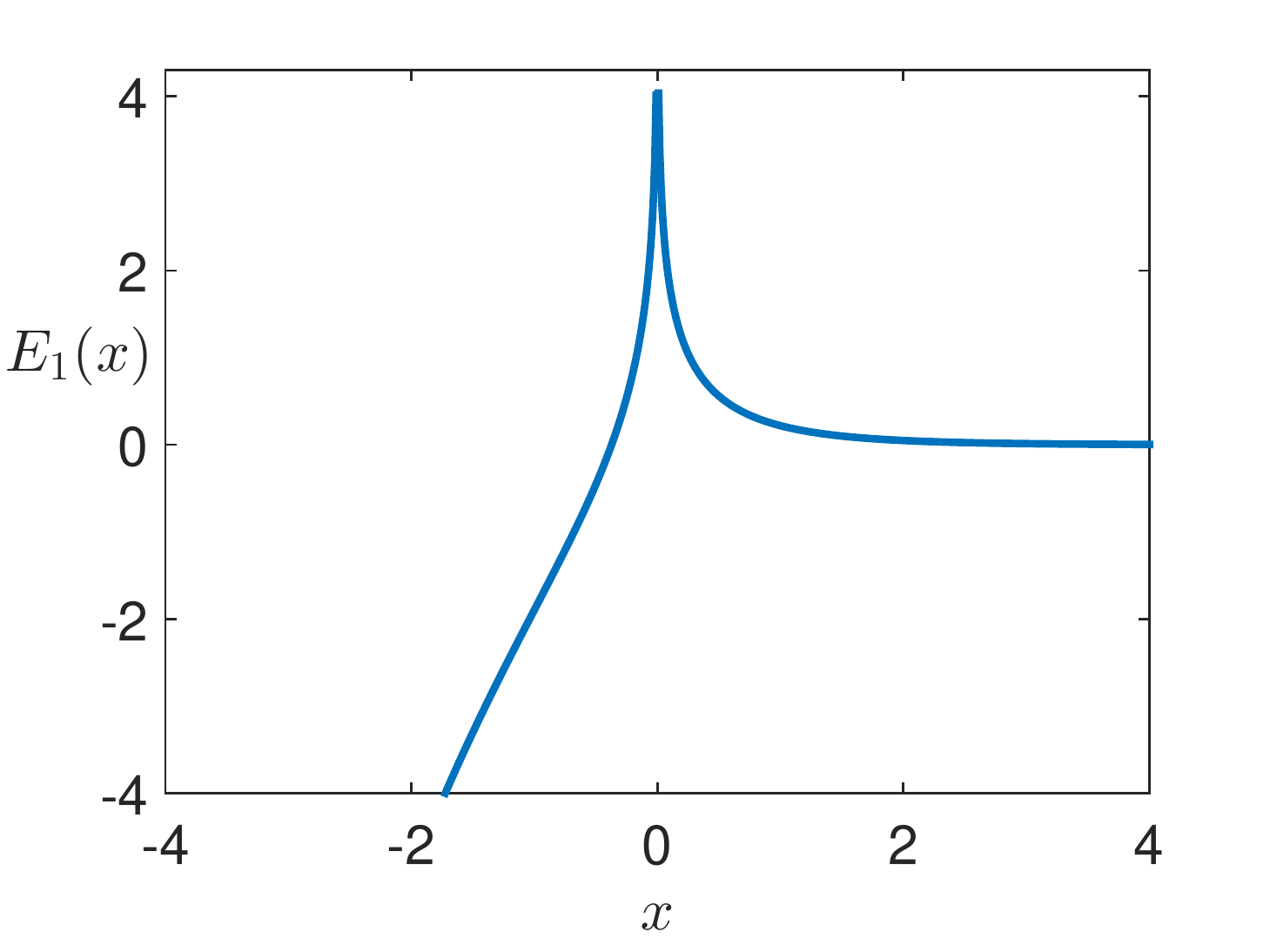}
\caption{\label{fig:E1_Ei}Exponentielle int\'egrale $\Ei(x)$ et $E_1(x)$.}
\end{figure}

La fonction exponentielle int\'egrale, not\'e $\Ei$, est d\'efinie par 
\begin{align}\label{eq:4_1}
\Ei(x) = \int_{-\infty}^{x}\frac{e^{t}}{t}dt \hspace{20mm} x\neq 0
\end{align}
L'exponentielle int\'egrale est reli\'ee \`a une autre fonction, not\'ee $E_1$ d\'efinie par
\begin{align}\label{eq:4_2}
-\Ei(-x) \equiv E_1(x) = \int_{x}^{\infty}\frac{e^{-t}}{t}dt \hspace{20mm} x > 0
\end{align}
Parfois on utilise l'appellation "exponentielle int\'egrale " pour $E_1(x)$.
Les repr\'esentations graphiques des fonctions $E_1$ et $\Ei$ sont donn\'ees dans la Fig.~(\ref{fig:E1_Ei}).

Une autre d\'efinition largement utilis\'ee est la suivante
\begin{align}\label{eq:4_21}
\Ein(x)=\int_{0}^{x}\frac{1-e^{-t}}{t}dt
\end{align}
$\Ein(x)$ est appell\'ee \'exponentielle int\'egrale modifi\'ee. Puisque on a
\begin{align*}
\frac{1-e^{-t}}{t}&=\frac{1}{t}\left(1-\sum_{n=0}^{\infty}\frac{(-1)^nt^n}{n!}\right)\\
&=\sum_{n=1}^{\infty}\frac{(-1)^{n+1}}{n!}t^{n-1}
\end{align*}
la fonction $\Ein$ peut \^etre d\'evelopp\'e en s\'erie enti\`ere comme
\begin{align}\label{eq:4_22}
\Ein(x)= \sum_{n=1}^{\infty}\frac{(-1)^{n+1}x^n}{n!n}
\end{align}
\subsection{Propri\'et\'es}
On peut montrer facilement que la fonction exponentielle int\'egrale $E_1(x)$ est reli\'ee \`a la fonction gamma par
\begin{align*}
E_1(x) = \Gamma(0,x)= \lim_{a\rightarrow 0}\left\lbrace\Gamma(a) - \gamma(a,x)\right\rbrace
\end{align*}
En utilisant la repr\'esentation en s\'erie de la fonction gamma, on peut \'ecrire
\begin{align}\label{eq:4_3}
E_1(x) &= -\gamma - \ln x - \sum_{n=1}^{\infty}\frac{(-1)^{n}x^{n}}{n!n} \hspace{10mm} x > 0\\
&= -\gamma - \ln x + \Ein(x)\nonumber
\end{align}
o\`u $\gamma=0.577$ est la constante d'Euler-Mascheroni. Les limites asymptotiques de $E_1(x)$ lorsque $x\rightarrow 0^{+}$ et $x\rightarrow +\infty$ sont donn\'ees par
\begin{align}
&E_1(x)\sim -\ln x \hspace{10mm} x\rightarrow 0^{+},\label{eq:4_4}\\
&E_1(x)\sim \frac{e^{-x}}{x} \hspace{10mm} x\rightarrow +\infty\label{eq:4_5}
\end{align}
\textit{D\'emonstration}\\
Montrons la formule (\ref{eq:4_3}), on a
\begin{align}\label{eq:4_3a}
E_1(x) &=\int_{x}^{\infty}\frac{1}{t}\sum_{n=0}^{\infty}\frac{(-t)^{n}}{n!} dt\nonumber\\
&=\sum_{n=0}^{\infty}\frac{(-1)^{n}}{n!}\int_{x}^{\infty}t^{n-1}dt\nonumber\\
&=\int_{x}^{\infty}t^{-1}dt + \sum_{n=1}^{\infty}\frac{(-1)^{n}}{n!}\int_{x}^{\infty}t^{n-1}dt\nonumber\\
&=-\ln x - \sum_{n=1}^{\infty}\frac{(-1)^{n}}{n!}\frac{x^n}{n}+ \ln t\Big|_{t\rightarrow\infty} + \sum_{n=1}^{\infty}\frac{(-1)^{n}}{n!}\frac{t^n}{n}\Big|_{t\rightarrow\infty}
\end{align}
Le dernier terme peut \^etre \'ecrit comme
\begin{align*}
\sum_{n=1}^{\infty}\frac{(-1)^{n}}{n!}\frac{t^n}{n}\Big|_{t\rightarrow\infty}&=\int_{0}^{\infty}\sum_{n=1}^{\infty}\frac{(-1)^{n}}{n!}t^{n-1}dt\\
&=\int_{0}^{\infty}\frac{1}{t}\left[\sum_{n=0}^{\infty}\frac{(-1)^{n}}{n!}t^{n} -1\right] dt\\
&=\int_{0}^{\infty}\frac{1}{t}\left[e^{-t} -1\right] dt
\end{align*}
En utilisant la repr\'esentation suivante de l'int\'egrale 
\begin{align*}
e^{-t}=\lim_{n\rightarrow\infty}\left(1-\frac{t}{n}\right)^{n}
\end{align*}
on obtient
\begin{align*}
\int_{0}^{\infty}\frac{1}{t}\left[e^{-t} -1\right] dt=-\lim_{n\rightarrow\infty}\int_{0}^{n}\frac{1}{t}\left[1-\left(1-\frac{t}{n}\right)^{n}\right] dt
\end{align*}
Posons $u=t/n$ on trouve
\begin{align*}
\int_{0}^{\infty}\frac{1}{t}\left[e^{-t} -1\right] dt=-\lim_{n\rightarrow\infty}\int_{0}^{1}\frac{1}{u}\left[1-\left(1-u\right)^{n}\right] du
\end{align*}
En faisant le changement de variable $v=1-u$, il en r\'esulte que
\begin{align*}
\int_{0}^{\infty}\frac{1}{t}\left[e^{-t} -1\right] dt=-\lim_{n\rightarrow\infty}\int_{0}^{1}\frac{1-v^n}{1-v}dv
\end{align*}
Il n'est pas difficile de montrer que
\begin{align*}
\frac{1-v^n}{1-v}=1+v+v^2+v^3+\ldots +v^{n-1}
\end{align*}
donc
\begin{align*}
\int_{0}^{\infty}\frac{1}{t}\left[e^{-t} -1\right] dt&=-\lim_{n\rightarrow\infty}\left[v+\frac{v^2}{2}+\frac{v^3}{3}+\ldots +\frac{v^n}{n}\right]_{0}^{1}\\
&=-\lim_{n\rightarrow\infty}\sum_{k=1}^{n} \frac{1}{k}
\end{align*}
Enfin les deux derniers termes de l'Eq.~(\ref{eq:4_3a}) peuvent \^etre \'ecrits comme
\begin{align*}
\ln t\Big|_{t\rightarrow\infty} + \sum_{n=1}^{\infty}\frac{(-1)^{n}}{n!}\frac{t^n}{n}\Big|_{t\rightarrow\infty}&=\lim_{n\rightarrow\infty}\left[\ln n - \sum_{k=1}^{n} \frac{1}{k}\right]\\
&=-\gamma
\end{align*}
o\`u $\gamma$ est la constante d'Euler-Mascheroni $\gamma= 0.5772156649$. Donc il en r\'esulte l'\'equation d\'esir\'ee Eq.~(\ref{eq:4_3}).
\section{Logarithme int\'egral}
Le logarithme int\'egral, not\'e $\Li$, est d\'efini pour tout nombre r\'eel strictement positif $x\neq 1$ par l'int\'egrale 
\begin{align}\label{eq:4_6}
\Li(x) = \int_0^{x}\frac{dt}{\ln t}
\end{align}
Il est li\'e \`a la fonction exponentielle int\'egrale $E_1$ par
\begin{align*}
\Li(x)= -E_1(\ln x) \hspace{20mm} 0<x<1
\end{align*}
\section{Sinus int\'egral et cosinus int\'egral}
\begin{figure}[!ht]
\centering
\includegraphics[scale=0.42]{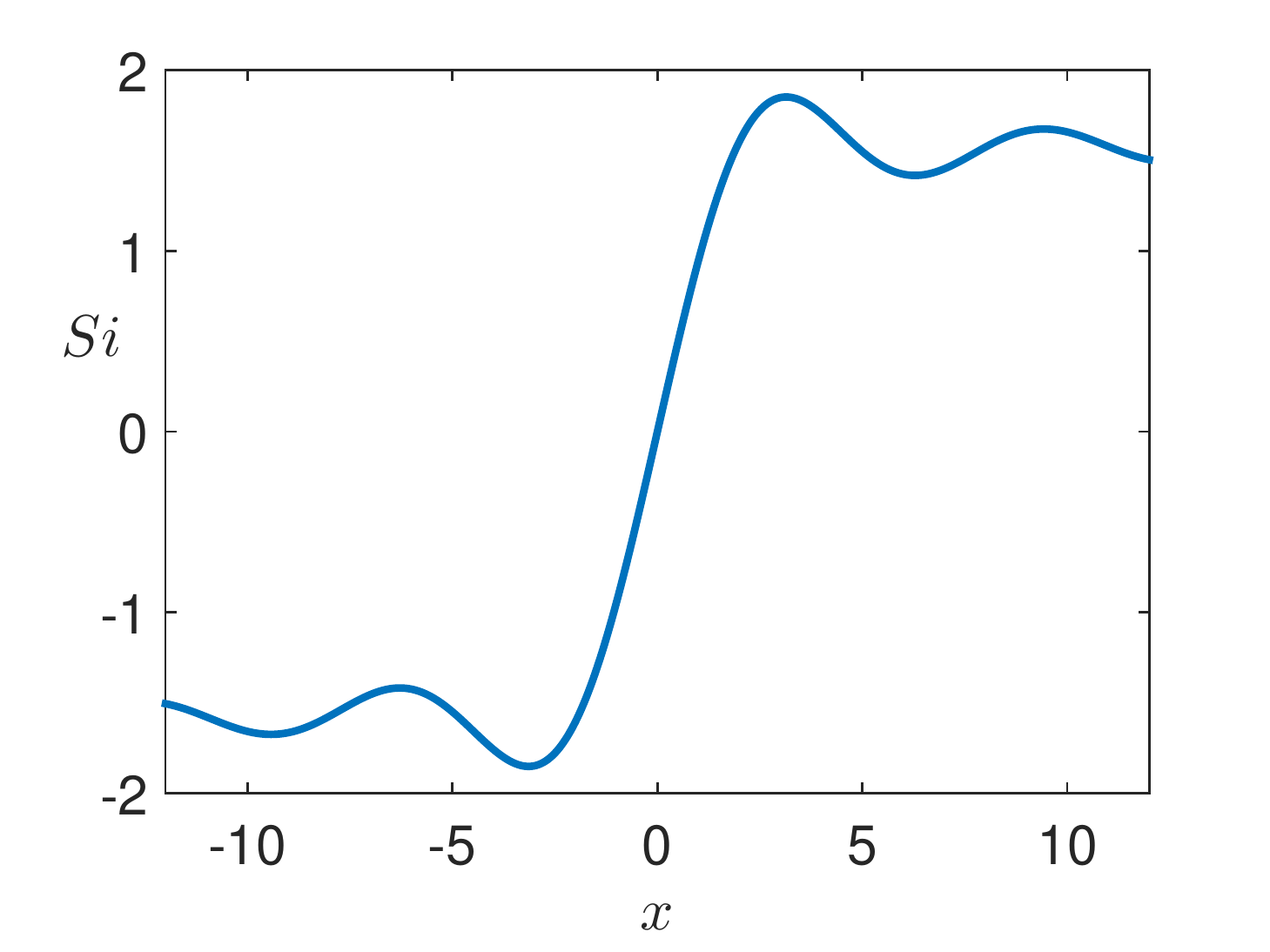}
\includegraphics[scale=0.42]{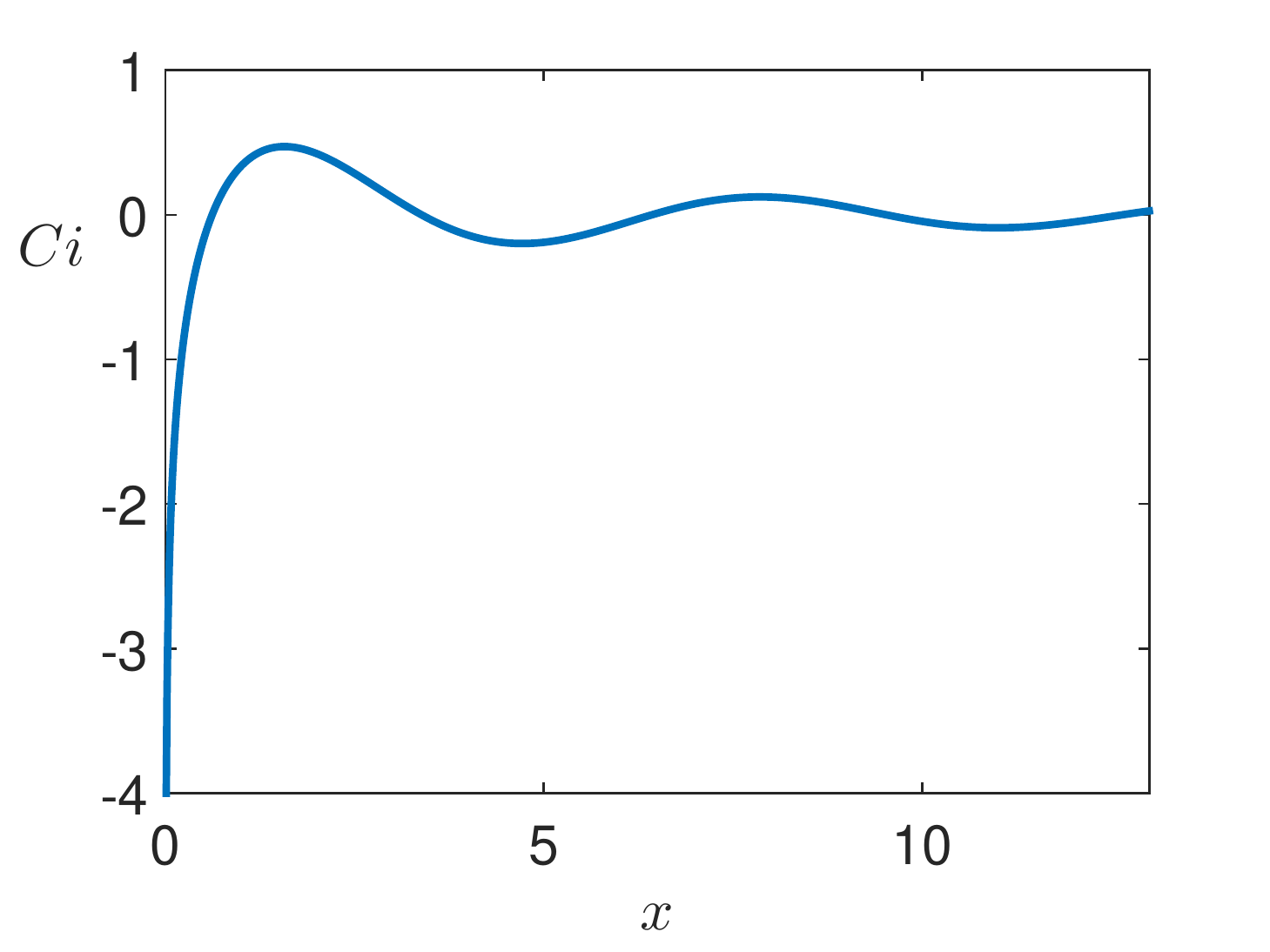}
\caption{\label{fig:Si_Ci}Sinus int\'egral $\Si(x)$ et cosinus int\'egral $Ci(x)$.}
\end{figure}

Les fonctions sinus et cosinus int\'egral, not\'ees respectivement $\Si$ et $\Ci$ sont d\'efinies par les int\'egrales
\begin{align}
&\Si(x) = \int_0^{x}\frac{\sin t}{t}dt \hspace{20mm} x>0\label{eq:4_7}\\
&\Ci(x) = -\int_x^{\infty}\frac{\cos t}{t}dt \hspace{20mm} x>0\label{eq:4_8}
\end{align}
elle sont repr\'esent\'ees sur la Fig.~(\ref{fig:Si_Ci}).
\subsection{Propri\'et\'es}
\begin{align*}
&\Si(0) = 0, \hspace{10mm} \Si(\infty)=\frac{\pi}{2},\\
&\Ci(0) = -\infty, \hspace{10mm} \Ci(\infty)=0,\\
&\Si'(x) = \frac{\sin x}{x}, \hspace{10mm} \Ci'(x) = \frac{\cos x}{x}
\end{align*}
\textit{D\'emonstration}\\
La d\'emonstration des expressions pr\'ec\'edentes est triviale sauf celui pour $\Si(\infty)$ et $\Ci(0)$. Pour les d\'emontrer consid\'erons le fait que
\begin{align*}
\frac{1}{t}=\int_{0}^{\infty}e^{-xt}dx
\end{align*}
donc
\begin{align*}
\Si(\infty)&=\int_0^{\infty}\left[\int_{0}^{\infty}e^{-xt}dx\right]\sin t dt\\
&=\int_0^{\infty}\left[\int_{0}^{\infty}e^{-xt}\sin t dt\right]dx
\end{align*}
En faisant une int\'egration par parties, on obtient
\begin{align*}
\int_{0}^{\infty}e^{-xt}\sin t dt= \left[-\frac{1}{x}e^{-xt}\sin t\right]_0^{\infty}+\frac{1}{x}\int_{0}^{\infty}e^{-xt}\cos t dt
\end{align*}
Le premier terme de droite est nul. En faisant une deuxi\`eme int\'egration par parties on trouve
\begin{align*}
\int_{0}^{\infty}e^{-xt}\sin t dt= \frac{1}{x}\left[-\frac{1}{x}e^{-xt}\cos t\right]_0^{\infty}-\frac{1}{x^2}\int_{0}^{\infty}e^{-xt}\sin t dt
\end{align*}
et il s’ensuit que
\begin{align*}
\left(1+\frac{1}{x^2}\right)\int_{0}^{\infty}e^{-xt}\sin t dt = \frac{1}{x^2}
\end{align*}
d'o\`u
\begin{align*}
\int_{0}^{\infty}e^{-xt}\sin t dt=\frac{1}{1+x^2}
\end{align*}
En remplaçant dans l'expression de $\Si(\infty)$ on obtient
\begin{align*}
\Si(\infty)&=\int_{0}^{\infty}\frac{dx}{1+x^2}\\
&=\arctan(x)\Big|_0^{\infty}\\
&=\frac{\pi}{2}
\end{align*}
En utilisant exactement le m\^eme raisonnement pour $\Ci(0)$ on obtient
\begin{align*}
\Ci(0)&=-\ln\left(1+t^2\right)\Big|_0^{\infty}\\
&=-\infty
\end{align*}
\subsection{Repr\'esentations en s\'eries}
Les fonctions $\Si$ et $\Ci$ sont d\'eveloppables en s\'eries enti\`eres sur $\mathbb{R}$, et leurs repr\'esentations sont donn\'ees par
\begin{align}
&(a)\hspace{1.5mm}\Si(x) = \sum_{n=0}^{\infty}\frac{(-1)^{n}}{(2n+1)(2n+1)!}x^{2n+1}\label{eq:4_9}\\
&(b)\hspace{1.5mm}\Ci(x)=\gamma + \ln x + \sum_{n=1}^{\infty}\frac{(-1)^{n}}{2n(2n)!}x^{2n}\label{eq:4_10}
\end{align}
o\`u $\gamma =0.577$ est la constante d'Euler-Mascheroni.\\
\textit{D\'emonstration}\\
(a) Le d\'eveloppement en s\'erie de la fonction sinus est donn\'e par
\begin{align*}
\sin x=\sum_{n=0}^{\infty}\frac{(-1)^{n}}{(2n+1)!}x^{2n+1}
\end{align*}
alors
\begin{align*}
\Si(x)&=\int_0^{x}\sum_{n=0}^{\infty}\frac{1}{t}\frac{(-1)^{n}}{(2n+1)!}t^{2n+1}dt \\
&=\left.\sum_{n=0}^{\infty}\frac{(-1)^{n}}{(2n+1)!}\frac{t^{2n+1}}{2n+1}\right]_0^{x}\\
&=\sum_{n=0}^{\infty}\frac{(-1)^{n}}{(2n+1)(2n+1)!}x^{2n+1}
\end{align*}
(b) Le d\'eveloppement en s\'erie de la fonction cosinus est donn\'e par
\begin{align*}
\cos x=\sum_{n=0}^{\infty}\frac{(-1)^{n}}{(2n)!}x^{2n} 
\end{align*}
alors
\begin{align*}
\Ci(x)&=-\int_x^{\infty}\sum_{n=0}^{\infty}\frac{1}{t}\frac{(-1)^{n}}{(2n)!}t^{2n}dt \\
&=-\int_x^{\infty}\frac{dt}{t}-\sum_{n=1}^{\infty}\frac{(-1)^{n}}{(2n)!}\int_x^{\infty}t^{2n-1}dt\\
&=-\left.\ln t\right]_x^{\infty} - \left.\sum_{n=1}^{\infty}\frac{(-1)^{n}}{(2n)!}\frac{t^{2n}}{2n}\right]_x^{\infty}\\
&=\ln x +\sum_{n=1}^{\infty}\frac{(-1)^{n}}{(2n)!}\frac{x^{2n}}{2n}+\left[-\sum_{n=1}^{\infty}\frac{(-1)^{n}}{(2n)!}\frac{t^{2n}}{2n}-\ln t\right]_{t\rightarrow\infty}\\
&=\gamma + \ln x + \sum_{n=1}^{\infty}\frac{(-1)^{n}}{(2n)!}\frac{x^{2n}}{2n}
\end{align*}
o\`u on peur montrer (de la m\^eme mani\'ere que pour l'Eq.~(\ref{eq:4_3})) que le dernier terme converge vers la constante $\gamma$
\begin{align*}
\left[-\sum_{n=1}^{\infty}\frac{(-1)^{n}}{(2n)!}\frac{t^{2n}}{2n}-\ln t\right]_{t\rightarrow\infty}&=\lim_{m\rightarrow\infty}\left[\sum_{n=1}^{m}\frac{1}{n}-\ln m\right]\\
&= \gamma
\end{align*}
\section{Exercices}
$\mathbf{Exercice\hspace{1mm}1:}$ Montrer que
\begin{align*}
&(a)\hspace{1.5mm}\int_{x}^{\infty}\frac{\sin t}{t}dt = \frac{1}{2}\pi -\Si(x)\\
&(b)\hspace{1.5mm}\int_{0}^{x}\frac{\cos t}{t}dt = \Ci(x)
\end{align*}
$\mathbf{Solutions}$\\
$(a)$ On a
\begin{align*}
\int_{0}^{\infty}\frac{\sin t}{t}dt= \int_{0}^{x}\frac{\sin t}{t}dt +\int_{x}^{\infty}\frac{\sin t}{t}dt
\end{align*}
donc
\begin{align*}
\Si(\infty)=\Si(x) + \int_{x}^{\infty}\frac{\sin t}{t}dt
\end{align*}
et il en r\'esulte que
\begin{align*}
\int_{x}^{\infty}\frac{\sin t}{t}dt=\frac{\pi}{2}-\Si(x)
\end{align*}
$(b)$ On a
\begin{align*}
\int_{0}^{\infty}\frac{\cos t}{t}dt= \int_{0}^{x}\frac{\cos t}{t}dt+ \int_{x}^{\infty}\frac{\cos t}{t}dt
\end{align*}
donc
\begin{align*}
\Ci(0)=\int_{0}^{x}\frac{\cos t}{t}dt-\Ci(x)
\end{align*}
Puisque $\Ci(0)=0$, il en r\'esulte l'\'equation d\'esir\'e.\\
$\mathbf{Exercice\hspace{1mm}2:}$ Montrer que
\begin{align*}
&(a)\hspace{1.5mm}E_1\left(xe^{\frac{1}{2}i\pi}\right)=-\Ci(x)+i\left[-\frac{1}{2}\pi+\Si(x)\right]\\
&(b)\hspace{1.5mm}E_1\left(xe^{\frac{1}{2}i\pi}\right)=\int_{0}^{\infty}\frac{\cos(x+t)}{x+t}dt-i\int_{0}^{\infty}\frac{\sin(x+ t)}{x+t}dt
\end{align*}
et d\'eduire que
\begin{align*}
&(c)\hspace{1.5mm}\Si(x)=-f(x)\cos x -g(x)\sin x +\frac{1}{2}\pi\\
&(d)\hspace{1.5mm}\Ci(x)=f(x)\sin x -g(x)\cos x
\end{align*}
o\`u 
\begin{align*}
f(x)=\int_{0}^{\infty}\frac{\sin t}{x+t}dt,\hspace{10mm}g(x)=\int_{0}^{\infty}\frac{\cos t}{x+t}dt,\hspace{5mm}x\neq 0,\hspace{5mm}|arg x|<\pi
\end{align*}
$\mathbf{Solution}$\\
$(a)$ On a
\begin{align*}
xe^{i\frac{\pi}{2}}=ix
\end{align*}
donc
\begin{align*}
E_1\left(xe^{\frac{1}{2}i\pi}\right)&=\int_{ix}^{\infty}\frac{e^{-t}}{t}dt\\
&=\int_{x}^{\infty}\frac{e^{-it}}{t}dt\\
&=\int_{x}^{\infty}\frac{\cos t}{t}dt - i \int_{x}^{\infty}\frac{\sin t}{t}dt
\end{align*}
En utilisant les r\'esultats de l'exercice pr\'ec\'edent, on obtient
\begin{align*}
E_1\left(xe^{\frac{1}{2}i\pi}\right)=-\Ci(x)-i\left[\frac{\pi}{2}-\Si(x)\right]
\end{align*}
$(b)$ On a 
\begin{align*}
E_1\left(xe^{\frac{1}{2}i\pi}\right)=\int_{x}^{\infty}\frac{\cos t}{t}dt - i \int_{x}^{\infty}\frac{\sin t}{t}dt
\end{align*}
En faisant le changemet de variable $u=t-x$, il s'ensuit
\begin{align*}
E_1\left(xe^{\frac{1}{2}i\pi}\right)=\int_{0}^{\infty}\frac{\cos(u+x)}{u+x}du - i \int_{0}^{\infty}\frac{\sin(u+x)}{u+x}du
\end{align*}
$(c)$ D'apr\`es la derni\`ere relation on a
\begin{align*}
E_1\left(xe^{\frac{1}{2}i\pi}\right)&=\cos(x)\int_{0}^{\infty}\frac{\cos(u)}{u+x}du - \sin(x)\int_{0}^{\infty}\frac{\sin(u)}{u+x}du  \\
&- i\cos(x) \int_{0}^{\infty}\frac{\sin(u)}{u+x}du - i\sin(x) \int_{0}^{\infty}\frac{\cos(u)}{u+x}du
\end{align*}
Par identification avec $(a)$ on obtient
\begin{align*}
&-\Ci(x)=\cos(x)\int_{0}^{\infty}\frac{\cos(u)}{u+x}du - \sin(x)\int_{0}^{\infty}\frac{\sin(u)}{u+x}du  \\
&-\frac{\pi}{2}+\Si(x)=- \cos(x) \int_{0}^{\infty}\frac{\sin(u)}{u+x}du - \sin(x) \int_{0}^{\infty}\frac{\cos(u)}{u+x}du
\end{align*}
donc il en r\'esulte
\begin{align*}
Ci(x)&=\left(\int_{0}^{\infty}\frac{\sin(u)}{u+x}du\right)\sin(x) -\left(\int_{0}^{\infty}\frac{\cos(u)}{u+x}du\right)\cos(x)\\
&=f(x)\sin(x)-g(x)\cos(x)
\end{align*}
et
\begin{align*}
\Si(x)&=- \left(\int_{0}^{\infty}\frac{\sin(u)}{u+x}du\right)\cos(x) - \left(\int_{0}^{\infty}\frac{\cos(u)}{u+x}du\right)\sin(x) +\frac{\pi}{2}\\
&=-f(x)\cos(x) -g(x)\sin(x)+\frac{\pi}{2}
\end{align*}
$\mathbf{Exercice\hspace{1mm}3:}$ Montrer que
\begin{align*}
\int_{0}^{\infty}e^{-t}\ln(t)dt=-\gamma
\end{align*}
o\`u $\gamma$ est la constante d'Euler.\\
$\mathbf{Solution}$\\
D'apr\`es l'Eq.~(\ref{eq:4_4}) on a
\begin{align*}
-\gamma &= E_1(1)- \Ein(1)\\
&= \int_{1}^{\infty}\frac{e^{-t}}{t}dt - \int_{0}^{1}\frac{1-e^{-t}}{t}dt
\end{align*}
En faisant une int\'egration par parties on obtient
\begin{align*}
-\gamma &= \left[e^{-t}\ln t\right]_{1}^{\infty}+\int_{1}^{\infty}e^{-t}\ln t dt - \left[(1-e^{-t})\ln t\right]_{0}^{1} + \int_{0}^{1}e^{-t}\ln t dt\\
&=\int_{0}^{\infty}e^{-t}\ln t dt
\end{align*}
\vspace{1.5mm}\\
$\mathbf{Exercice\hspace{1mm}4:}$ Montrer que
\begin{align*}
\int_{0}^{1}\frac{1-e^{-t}-e^{-1/t}}{t}dt=\gamma
\end{align*}
$\mathbf{Solution}$\\
On a
\begin{align*}
\int_{0}^{1}\frac{1-e^{-t}-e^{-1/t}}{t}dt&=\int_{0}^{1}\frac{1-e^{-t}}{t}dt-\int_{0}^{1}\frac{e^{-1/t}}{t}dt
\end{align*}
En remplaçant $1/t$ par $t$ dans la deuxi\`eme int\'egrale \`a droite, on obtient
\begin{align*}
\int_{0}^{1}\frac{1-e^{-t}-e^{-1/t}}{t}dt&=\int_{0}^{1}\frac{1-e^{-t}}{t}dt-\int_{1}^{\infty}te^{-t}\frac{dt}{t^2}\\
&=\Ein(1)-E_1(1)\\
&=\gamma
\end{align*}
$\mathbf{Exercice\hspace{1mm}5:}$ Montrer que
\begin{align*}
\int_{-\infty}^{x}\frac{e^{t}}{t^{2}(1-t)}dt=\frac{e^{x}}{x}-2\Ei(x)+e\Ei(x-1),\hspace{10mm}x<0
\end{align*}
$\mathbf{Solution}$\\
L'int\'egrant peut \^etre factoris\'e comme
\begin{align*}
\frac{e^{t}}{t^{2}(1-t)}=\frac{(t-1)e^{t}}{t^2}+\frac{2e^{t}}{t} +\frac{e^{t}}{1-t}
\end{align*}
Donc
\begin{align*}
\int_{-\infty}^{x}\frac{e^{t}}{t^{2}(1-t)}dt=\int_{-\infty}^{x}\frac{(t-1)e^{t}}{t^2}dt+\int_{-\infty}^{x}\frac{2e^{t}}{t}dt + \int_{-\infty}^{x}\frac{e^{t}}{1-t}dt
\end{align*}
En faisant un changement de variable $u=t-1$ dans la troisi\`eme int\'egrale de c\^ot\'e droite, on obtient
\begin{align*}
\int_{-\infty}^{x}\frac{e^{t}}{t^{2}(1-t)}dt&=\frac{e^{x}}{x}-2\Ei(x) + \int_{-\infty}^{x-1}\frac{e^{1+u}}{u}du\\
&=\frac{e^{x}}{x}-2\Ei(x)+e\Ei(x-1)
\end{align*}
\chapter{Les polyn\^omes orthogonaux}
\label{chap: poly_orthog}
Les polyn\^omes orthogonaux sont tr\`es importants en physique. Ce sont les solutions d'\'equations qui surviennent tr\`es souvent lorsqu'un probl\`eme poss\`ede une sym\'etrie sph\'erique. Les polyn\^omes orthogonaux sont profond\'ement li\'es \`a la m\'ecanique quantique et ils sont pr\'esents dans de nombreuses applications, telles que la th\'eorie \'electromagnétique, l'hydrodynamique et la conduction thermique. En ing\'enierie, les polyn\^omes orthogonaux apparaissent dans de nombreuses applications, par exemple dans la th\'eorie des lignes de transmission, la th\'eorie des circuits \'electriques, la physique des r\'eacteurs nucl\'eaires ainsi que la sismologie.
\section{Polyn\^omes de Legendre}
\begin{figure}[hbt]
\centering
\includegraphics[scale=0.7]{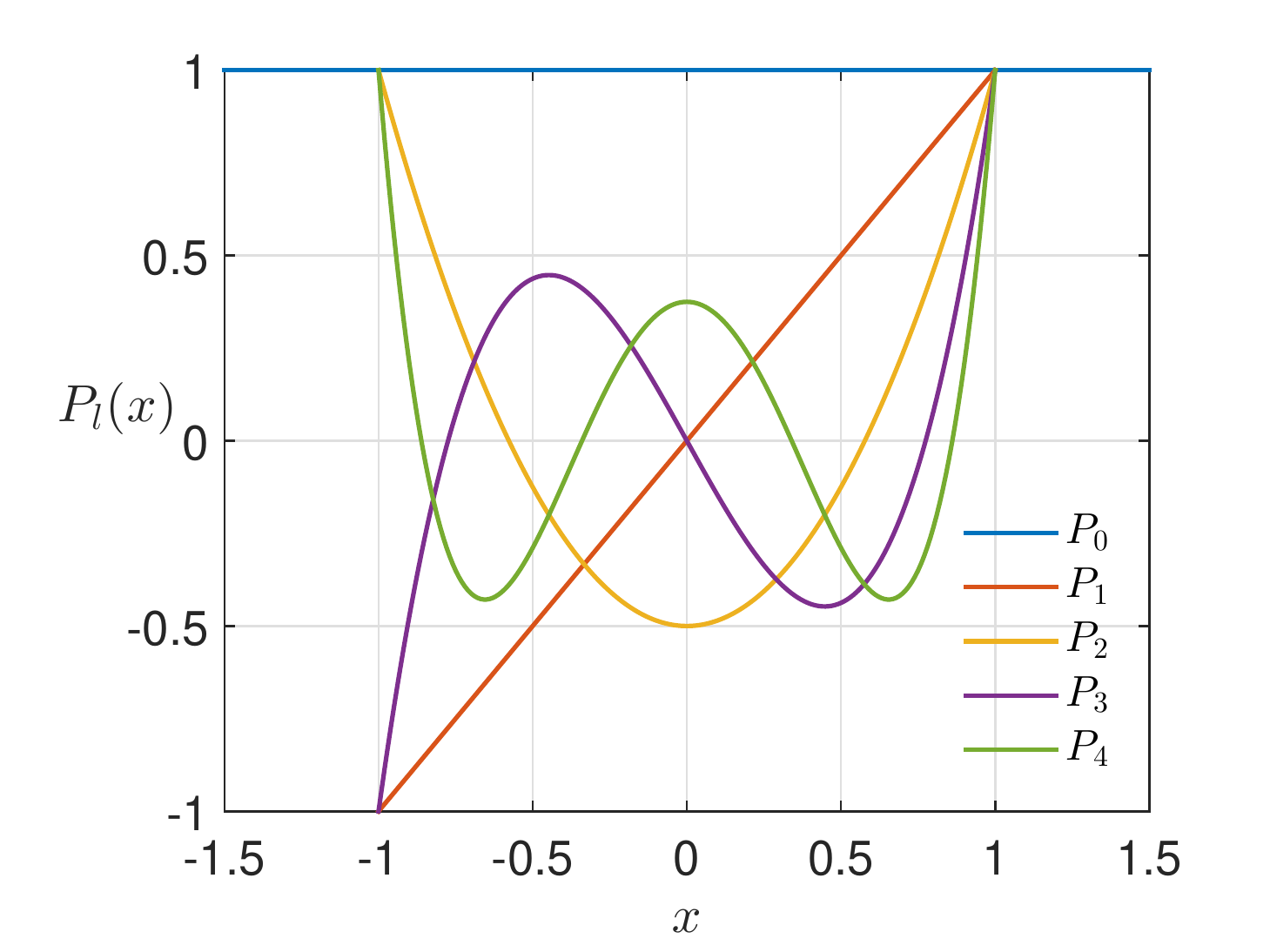}
\caption{\label{fig:legendre} Polyn\^omes de Legendre pour $n=0,1,2,3,4$.}
\end{figure}

Adrien-Marie Legendre a introduit, en $1784$, les polyn\^omes de Legendre, tout en \'etudiant l'attraction des sph\'ero\"ides et des ellipso\"ides. Ces polyn\^omes sont les solutions d'une \'equation diff\'erentielle ordinaire appel\'ee \'equation diff\'erentielle de Legendre. Cette \'equation est fr\'equemment rencontr\'ee en physique et en ing\'enierie. En particulier, cela se produit lors de la r\'esolution de l’\'equation de Laplace en coordonn\'ees sph\'eriques.
L'\'equation de Legendre est donn\'ee par 
\begin{align}\label{eq:5_1}
(1-x^{2})\frac{d^{2}y}{dx^{2}} -2x\frac{dy}{dx} +l(l+1)y=0
\end{align}
dans le cas g\'en\'erale $l\in\mathbb{R}$. L'\'equation (\ref{eq:5_1}) est solvable par la m\'ethode de Frobenius (qui a \'et\'e utilis\'ee pour r\'esoudre les \'equations de Bessel). Les solutions ont la forme de s\'eries enti\`eres et elles sont convergentes uniquement dans le domaine $-1<x<1$. Dans le cas o\`u $l\in \mathbb{N}$, il est possible d'obtenir des solutions qui sont r\'eguli\`eres aux points $x\pm 1$, et pour lesquelles la s\'erie s'arr\^ete au terme de degr\'e $l$. Dans le derni\`er cas, la solution de l'Eq.~(\ref{eq:5_1}) est appell\'ee polyn\^one de Legendre. Elle est donn\'ee par
\begin{align}\label{eq:5_2}
P_l(x)=\sum_{r=0}^{[l/2]}(-1)^{r}\frac{(2l-2r)!}{2^{l}r!(l-r)!(l-2r)!}x^{l-2r}\hspace{10mm}-1\leq x \leq 1
\end{align}
Notons que la solution~(\ref{eq:5_2}) est convergente seulement dans le domaine $[-1,\hspace{0.5mm}1]$. Les points $x\pm 1$ sont des points singuliers r\'eguliers de l'\'equation différentielle (\ref{eq:5_1}).\\
\textit{D\'emonstration}\\
Avant de chercher la solution de l'Eq.~(\ref{eq:5_1}), cherchons d'abord le domaine de convergence. Pour cela, on utilise le th\'eor\`eme suivant; si on a une \'equation diff\'erentielle de la forme
\begin{align}\label{eq:5_3}
x^{2}\frac{d^{2}y}{dx^{2}} -xq(x)\frac{dy}{dx} +p(x)y=0
\end{align}
et si $q(x)$ et $p(x)$ peuvent \^etre developp\'ees sous les formes
\begin{align*}
&q(x)=\sum_{n=0}^{\infty}q_nx^{n},\\
&p(x)=\sum_{n=0}^{\infty}p_nx^{n}
\end{align*}
donc, le domaine de convergence de la solution de l'Eq.~(\ref{eq:5_3}) est la m\^eme que le domaine de convergence de $q(x)$ et $p(x)$.\\
Pour d\'eterminer les fonction $p(x)$ et $q(x)$, on \'ecrit l'Eq.~(\ref{eq:5_1}) sous la forme
\begin{align*}
x^{2}\frac{d^{2}y}{dx^{2}} -x\frac{2x^{2}}{1-x^2}\frac{dy}{dx} +\frac{l(l+1)x^2}{1-x^2}y=0
\end{align*}
En comparant avec l'Eq.~(\ref{eq:5_3}), on d\'eduit que
\begin{align*}
&q(x)=\frac{2x^{2}}{1-x^2},\\
&p(x)=\frac{l(l+1)x^2}{1-x^2}
\end{align*}
Le d\'eveloppement en s\'erie enti\`ere des fonctions $q(x)$ et $p(x)$ est
\begin{align*}
&q(x)=2x^{2}\sum_{m=0}^{\infty}x^{2m},\\
&p(x)=l(l+1)x^2\sum_{m=0}^{\infty}x^{2m}
\end{align*}
donc les fonctions $q(x)$ et $p(x)$ sont convergentes seulement dans l'intervalle $-1<x<1$. Pour r\'esoudre l'Eq.~(\ref{eq:5_1}), on utilise la m\'ethode de Frobenius qui consiste \`a chercher des solutions sous forme de s\'eries enti\`eres
\begin{align*}
y(x,s)=\sum_{r=0}^{\infty}a_rx^{r+s}
\end{align*}
Si en ins\`ere $y$, sa d\'eriv\'ee $y'$ et la d\'eriv\'ee seconde $y''$ dans l'Eq.~(\ref{eq:5_1}),
on trouve trois \'equations indiciales
\begin{align*}
&a_0s(s-1)=0,\\
&a_1s(s+1)=0,\\
&a_{n+2}(s+n+2)(s+n+1)-a_n\left\lbrace(s+n)(s+n+1)-l(l+1)\right\rbrace\hspace{5mm} n\geq 0
\end{align*}
Pour $a_0$ et $a_1$ arbitraires et diff\'erents de $0$, les \'equations indiciales admettent comme solutions
\begin{align*}
s=0, \hspace{10mm} s=1
\end{align*}
et
\begin{align*}
a_{n+2}&=a_n\frac{n(n+1)- l(l+1)}{(n+1)(n+2)}\\
&=a_n\frac{(n-l)(l+n+1)}{(n+1)(n+2)}
\end{align*}
Pour d\'eduire la forme du terme g\'en\'eral $a_n$ en fonction de $a_0$(ou $a_1$) et $n$, on consid\'ere quelques cas particuliers
\begin{align*}
a_2 &= -a_0\frac{l(l+1)}{2}\\
a_3 &=a_1\frac{(1-l)(l+2)}{2.3}\\
a_4 &=-a_2\frac{(l-2)(l+3)}{3.4}\\
&=a_0\frac{l(l-2)(l+1)(l+3)}{1.2.3.4}\\
a_5&=-a_3\frac{(l-3)(l+4)}{4.5}\\
&=a_1\frac{(l-1)(l-3)(l+2)(l+4)}{1.2.3.4.5}
\end{align*}
A partir des \'equations pr\'ec\'edentes on peut d\'eduire le terme g\'en\'eral qui est 
donn\'e par
\begin{align*}
&a_{2n}=\frac{(-1)^{n}a_0}{(2n)!}l(l-2)(l-4)\ldots(l-2n+2)(l+1)(l+3)\ldots(l+2n-1)\\
&a_{2n+1}=\frac{(-1)^{n}a_1}{(2n+1)!}(l-1)(l-3)\ldots(l-2n+1)(l+2)(l+4)\ldots(l+2n)
\end{align*}
Donc, la solution g\'en\'erale peut \^etre \'ecrite sous la forme
\begin{align*}
y(x,0)&=a_0\left\lbrace1+\sum_{n=1}^{\infty}(-1)^{n}\frac{l(l-2)(l-4)\ldots(l-2n+2)(l+1)(l+3)\ldots(l+2n-1)}{(2n)!}x^{2n}\right\rbrace\nonumber\\
&+a_1\left\lbrace x+\sum_{n=1}^{\infty}(-1)^{n}\frac{(l-1)(l-3)\ldots(l-2n+1)(l+2)(l+4)\ldots(l+2n)}{(2n+1)!}x^{2n+1}\right\rbrace\nonumber\\
&=a_0y_1(x)+a_1y_2(x)
\end{align*}
Nous soulignons que cette solution est valable dans l'intervalle $-1<x<1$. Afin d'\'etendre la solution aux points $x=-1$ and $x=1$, nous utilisons l'observation suivante:
si $l$ est pair, \i.e., $l=2n$, on a
\begin{align*}
&a_2\neq 0, \hspace{1mm}a_4\neq 0, \ldots,a_{2n}\neq 0,\\
&a_{2n+2}=a_{2n+4}=\ldots= 0
\end{align*}
et si $l$ est impair, \i.e., $l=2n+1$
\begin{align*}
&a_3\neq 0, \hspace{1mm}a_5\neq 0, \ldots,a_{2n+1}\neq 0,\\
&a_{2n+3}=a_{2n+5}=\ldots= 0
\end{align*}
Donc pour $l=2n$, $y_1(x)$ devient fini pour tout $x$ \`a l'int\'erieur de l'intervalle $[-1,1]$ et on a l'inverse pour $l=2n+1$ o\`u le polyn\^ome $y_2(x)$ devient fini pour tout  $x$ \`a l'int\'erieur de l'intervalle $[-1,1]$. La solution de l'\'equation de Legendre dans l'intervalle $[-1,1]$ est donc donn\'ee soit par $y_1(x)$ (dans le cas o\`u $l$ est pair) ou par $y_2(x)$ (si $l$ est impair). $y(x)$ peut \^etre \'ecrit sous la forme g\'en\'erale
\begin{align}
y(x) &= a_lx^{l}+a_{l-2}x^{l-2}+a_{l-4}x^{l-4}+\ldots+\left\lbrace\begin{array}{ll}
a_0 & \mbox{si $l$ est pair}\\
a_1& \mbox{si $l$ est impair}
\end{array}\right.\nonumber\\
&=\sum_{r=0}^{[l/2]}a_{l-2r}x^{l-2r}\label{eq:5_4}
\end{align}
$[l/2]$ est la fonction de plafond ("ceiling function") et elle d\'efinie par
\begin{align*}
\left[\frac{l}{2}\right]=\left\lbrace\begin{array}{ll}
\vspace{2mm}
\frac{l}{2} & \mbox{si $l$ est pair}\\
\frac{l-1}{2} & \mbox{si $l$ est impair}
\end{array}\right.
\end{align*}
$a_{l-2r}$ peut \^etre d\'eduit de la forme g\'en\'erale de $a_{n}$. On a
\begin{align*}
a_n=-a_{n+2}\frac{(n+2)(n+1)}{(l-n)(l+n+1)}
\end{align*}
En remplaçant $n$ par $l-2$ dans l'equation pr\'ec\'edente, on trouve 
\begin{align*}
a_{l-2}= -a_l\frac{l(l-1)}{2(2l-1)}
\end{align*}
pour $n=l-4$ on obtient
\begin{align*}
a_{l-4}&=a_{l-2}\frac{(l-2)(l-3)}{4(2l-3)}\\
&=a_l\frac{l(l-1)(l-2)(l-3)}{2.4(2l-1)(2l-3)}
\end{align*}
A partir des \'equations de $a_{l-2}$ et $a_{l-4}$, on peut d\'eduire la forme g\'en\'erale
\begin{align*}
a_{l-2r} = (-1)^{r}\frac{l(l-1)(l-2)\ldots(l-2r+1)}{2.4\ldots(2r)(2l-1)(2l-3)\ldots(2l-2r+1)}a_l
\end{align*}
En remplaçant $a_{l-2r}$ dans l'Eq.~(\ref{eq:5_4}), on obtient
\begin{align*}
y(x)=a_l\sum_{r=0}^{[l/2]}(-1)^{r}\frac{l(l-1)(l-2)\ldots(l-2r+1)}{2.4\ldots(2r)(2l-1)(2l-3)\ldots(2l-2r+1)}x^{l-2r}
\end{align*}
On peut simplifier l'expression de $y(x)$ en \'ecrivant autrement le num\'erateur et le d\'enominateur
\begin{align*}
l(l-1)\ldots(l-2r+1)&=l(l-1)\ldots(l-2r+1)\frac{(l-2r)!}{(l-2r)!}\nonumber\\
&=\frac{l!}{(l-2r)!},
\end{align*}
\begin{align*}
2.4.6\ldots2r=(2.1)(2.2)(2.3)\ldots(2.r)=2^{r}r!
\end{align*} 
et
\begin{align*}
(2l-1)(2l-3)\ldots(2l-2r+1)&=\frac{2l(2l-1)(2l-2)(2l-3)\ldots(2l-2r+1)(2l-2r)!}{2l(2l-2)\ldots(2l-2r+2)(2l-2r)!}\nonumber\\
&=\frac{(2l)!}{2^{r}l(l-1)\ldots(l-r+1)(2l-2r)!}\nonumber\\
&=\frac{(2l)!(l-r)!}{2^{r}l!(2l-2r)!}
\end{align*} 
donc
\begin{align*}
y(x)&=a_l\sum_{r=0}^{[l/2]}(-1)^{r}\frac{l!}{(l-2r)!}\frac{1}{2^{r}r!}\frac{2^{r}l!(2l-2r)!}{(2l)!(l-r)!}x^{l-2r}\nonumber\\
&=a_l\sum_{r=0}^{[l/2]}(-1)^{r}\frac{(l!)^{2}(2l-2r)!}{r!(l-2r)!(2l)!}x^{l-2r}
\end{align*}
Choisissant $a_l$ d'\^etre
\begin{align*}
a_l=\frac{(2l)!}{2^{l}(l!)^{2}}
\end{align*}
On obtient le r\'esultat final pour la solution de l'\'equation de Legendre Eq.~(\ref{eq:5_1})
\begin{align*}
P_l(x)=\sum_{r=0}^{[l/2]}(-1)^{r}\frac{(2l-2r)!}{2^{l}r!(l-r)!(l-2r)!}x^{l-2r}\hspace{10mm}-1\leq x \leq 1
\end{align*}
Les cinq premiers polyn\^omes sont 
\begin{align*}
&P_0(x)=1\\
&P_1(x)=x\\
&P_2(x)=\frac{1}{2}(3x^2-1)\\
&P_3(x)=\frac{1}{2}(5x^3-3x)\\
&P_4(x)=\frac{1}{8}(35x^4-30x^2+3)
\end{align*}
et ils sont repr\'esent\'es dans la Fig.~(\ref{fig:legendre}).
\subsection{Fonction g\'en\'eratrice}
\begin{align}\label{eq:5_5}
\frac{1}{\sqrt{1-2tx+t^{2}}}=\sum_{l=0}^{\infty}t^{l}P_l(x) \hspace{5mm} \mbox{if $|t|<1, |x| \leq 1$}
\end{align}
\textit{D\'emonstration}\\
On a
\begin{align*}
\frac{1}{\sqrt{1-2tx+t^{2}}}&=\left[1-t(2x-t)\right]^{-\frac{1}{2}}\\
&=1+\left(-\frac{1}{2}\right)\left[-t(2x-t)\right]+\frac{(-1/2)(-3/2)}{2!}\left[-t(2x-t)\right]^{2}\\
&+\ldots+\frac{(-1/2)(-3/2)\ldots(-(2r-1)/2)}{r!}\left[-t(2x-t)\right]^{r}+\ldots\\
&=\sum_{r=0}^{\infty}\frac{1.3.5\ldots(2r-1)}{2^{r}r!}(-1)^{r}t^{r}(2x-t)^{r}\\
&=\sum_{r=0}^{\infty}\frac{(2r)!}{2^{2r}(r!)^{2}}t^{r}(2x-t)^{r}
\end{align*}
D'apr\'es la formule du bin\^ome de Newton, on a
\begin{align*}
(2x-t)^{r}=\sum_{s=0}^{r}\frac{r!}{s!(r-s)!}(2x)^{r-s}(-t)^{s}
\end{align*}
En remplaçant dans l'\'equation pr\'ec\'edente, on trouve
\begin{align*}
\frac{1}{\sqrt{1-2tx+t^{2}}}&=\sum_{r=0}^{\infty}\frac{(2r)!}{2^{2r}(r!)^{2}}\sum_{s=0}^{r}(-1)^{s}\frac{r!}{s!(r-s)!}(2x)^{r-s}t^{r+s}
\end{align*}
Si on pose $n=r+s$, on a $n$ qui varie entre $0$ et $\infty$,  et puisque $s$ varie entre $0$ et $r$, donc $r$ varie entre $[n/2]$ et $n$, alors
\begin{align*}
\frac{1}{\sqrt{1-2tx+t^{2}}}&=\sum_{n=0}^{\infty}\sum_{r=[n/2]}^{n}t^{n}(-1)^{n-r}\frac{(2r)!}{2^{n}r!(n-r)!(2r-n)!}x^{2r-n}
\end{align*}
En faisant le changement de variable $l=n-r$, on obtient finalement
\begin{align*}
\frac{1}{\sqrt{1-2tx+t^{2}}}&=\sum_{n=0}^{\infty}t^{n}\sum_{l=[n/2]}^{0}(-1)^{l}\frac{(2n-2l)!}{2^{n}l!(n-l)!(n-2l)!}x^{n-2l}\\
&=\sum_{n=0}^{\infty}t^{n}P_n(x)
\end{align*}
\subsection{Formule de Rodrigues}
\begin{align}\label{eq:5_6}
P_l(x)=\frac{1}{2^{l}l!}\frac{d^{l}}{dx^{l}}(x^{2}-1)^{l}
\end{align}
\textit{D\'emonstration}\\
D'apr\`es la formule du bin\^ome de Newton, on a
\begin{align}
\left(x^2-1\right)^{l}=\sum_{r=0}^l\frac{l!}{r!(l-r)!}(-1)^{r}x^{2(l-r)}
\end{align}
donc
\begin{align*}
\frac{1}{2^{l}l!}\frac{d^{l}}{dx^{l}}(x^{2}-1)^{l}=\frac{1}{2^{l}l!}\sum_{r=0}^l\frac{l!}{r!(l-r)!}(-1)^{r}\frac{d^{l}}{dx^{l}}x^{2(l-r)}
\end{align*}
mais
\begin{align*}
\frac{d^{l}}{dx^{l}}x^{2(l-r)}=0 \hspace{10mm} \mbox{ si\hspace{2mm} $2(l-r) <l$, \hspace{1mm}i.e.,\hspace{1mm} si \hspace{1mm}$r> l/2$}
\end{align*}
On peut donc remplacer $\sum_{r=0}^l$ par $\sum_{r=0}^{[l/2]}$.\\
Si $r\leq l/2$, on a
\begin{align*}
\frac{d^{l}}{dx^{l}}x^{2(l-r)}&=\{2(l-r)\}\{2(l-r)-1\}\{2(l-r)-2\}\ldots\{2(l-r)-l+1\}x^{2(l-r)-l}\nonumber\\
&=\frac{(2l-2r)!}{(l-2r)!}x^{l-2r}
\end{align*}
et on obtient alors
\begin{align*}
\frac{1}{2^{l}l!}\frac{d^{l}}{dx^{l}}(x^{2}-1)^{l}&=\sum_{r=0}^{[l/2]}(-1)^{r}\frac{(2l-2r)!}{2^{l}r!(l-r)!(l-2r)!}x^{l-2r}\\
&=P_l(x)
\end{align*}
\subsection{Repr\'esentation int\'egrale de Laplace}
\begin{align}\label{eq:5_7}
P_l(x)=\frac{1}{\pi}\int_0^{\pi}\left(x+\sqrt{x^{2}-1}\cos\theta\right)^{l}d\theta
\end{align}
\textit{D\'emonstration}\\
Consid\'erons la relation (qui peut \^etre prouv\'e en utilisant le changement de variable $t=\tan\theta$)
\begin{align}\label{eq:5_7a}
\int_{0}^{\pi}\frac{d\theta}{1+\lambda\cos\theta}=\frac{\pi}{\sqrt{1-\lambda^2}}
\end{align}
Posons $\lambda=-u\sqrt{x^2-1}/(1-ux)$, En d\'eveloppant les deux c\^ot\'es de l'Eq.~(\ref{eq:5_7a}) en puissance de $u$, on obtient
\begin{align*}
\frac{1}{1+\lambda\cos\theta}&=\frac{1}{1+\frac{u\sqrt{x^2-1}}{1-ux}\cos\theta}\\
&=(1-ux)\left[1-u\left(x+\sqrt{x^2-1}\cos\theta\right) \right]^{-1}\\
&=(1-ux)\sum_{l=0}^{\infty}u^{l}\left( x+\sqrt{x^2-1}\cos\theta \right)^{l}
\end{align*}
o\`u on autilis\'e le th\'eor\`erme de bin\^ome de Newton $(1-a)^{-1}=\sum_{n=0}^{\infty}a^n$.
\begin{align*}
\frac{1}{\sqrt{1-\lambda^2}}&=\frac{1}{\sqrt{1-\frac{u^2(x^2-1)}{(1-ux)^2}}}\\
&=\frac{1-ux}{\sqrt{(1-ux)^2-u^2(x^2-1)}}\\
&=\frac{1-ux}{\sqrt{1-2ux+u^2}}
\end{align*}
Substitution dans l'Eq.~(\ref{eq:5_7a}) donne
\begin{align*}
\int_{0}^{\pi}\sum_{l=0}^{\infty}u^{l}\left( x+\sqrt{x^2-1}\cos\theta \right)^{l}d\theta=\frac{\pi}{\sqrt{1-2ux+u^2}}
\end{align*}
On ins\`ere la formule de fonction g\'en\'eratrice (\ref{eq:5_7a}) dans l'\'equation pr\'ec\'edente, on obtient
\begin{align*}
\sum_{l=0}^{\infty}u^{l}\int_{0}^{\pi}\left( x+\sqrt{x^2-1}\cos\theta \right)^{l}d\theta=\pi\sum_{l=0}^{\infty}u^{l}P_l(x)
\end{align*}
Par identification des coefficients de $u^l$ on obtient
\begin{align*}
\pi P_l(x) = \int_{0}^{\pi}\left( x+\sqrt{x^2-1}\cos\theta \right)^{l}d\theta
\end{align*}
\subsection{Propri\'et\'es des polyn\^omes de Legendre}
\begin{align}\label{eq:5_8}
&(a)\hspace{0.5mm} P_l(1)=1,\nonumber\\
&(b)\hspace{0.5mm}P_l(-1)=(-1)^{l},\nonumber\\
&(c)\hspace{0.5mm}P_l'(1)=\frac{1}{2}l(l+1),\\
&(d)\hspace{0.5mm}P_l'(-1)=(-1)^{l-1}\frac{1}{2}l(l+1),\nonumber\\
&(e)\hspace{0.5mm}P_{2l}(0)=(-1)^{l}\frac{2l}{2^{2l}(l!)^{2}},\nonumber\\
&(f)\hspace{0.5mm}P_{2l+1}(0)=0\nonumber
\end{align}
\textit{D\'emonstration}\\
(a) On met $x=1$ dans l'Eq.~(\ref{eq:5_5}), on obtient
\begin{align*}
\frac{1}{1-t}=\sum_{l=0}^{\infty}t^{l}P_l(1)
\end{align*}
Le d\'eveloppement en s\'erie de $1/(1-t)$ donne
\begin{align*}
\frac{1}{1-t}=\sum_{l=0}^{\infty}t^{l}
\end{align*}
En identifiant les deux expressions pr\'ec\'edentes, il en r\'esulte que $P_l(1) = 1$, $\forall\hspace{1mm} l$.\\
(b) la preuve est exactement similaire \`a (a) avec $x=-1$.\\
(c) En posant $x=1$ dans l'\'equation de Legendre (\ref{eq:5_1}), on obtient
\begin{align*}
-2P_l'(1)+l(l+1)P_l(1)=0
\end{align*}
et puisque $P_l(1)=1$, il en r\'esulte que $P_l'(1)=l(l+1)/2$.\\
(d) M\^eme preuve que (c) en posant ici $x=-1$.\\
Pour montrer (e) et (f), on pose $x=0$ dans l'Eq.~(\ref{eq:5_6}). En faisant le d\'eveloppement en s\'erie, on trouve
\begin{align*}
\frac{1}{\sqrt{1+t^{2}}}&=\sum_{l=0}^{\infty}t^{l}P_l(0)\nonumber\\
&=\sum_{l=0}^{\infty}(-1)^{l}\frac{(2l)!}{2^{2l}(l!)^{2}}t^{2l}
\end{align*}
En identifiant les coefficients de la m\^eme puissance de $t$, on obtient
\begin{align*}
&P_{2l}(0)=(-1)^{l}\frac{(2l)!}{2^{2l}(l!)^{2}}\\
&P_{2l+1}(0)=0
\end{align*}
\subsection{Relation d'orthogonalit\'e}
\begin{align}\label{eq:5_9}
\int_{-1}^{1}P_l(x)P_m(x)dx=\left\lbrace\begin{array}{ll}
\vspace{2mm}
0  & \mbox{si $l \neq m$}\\
\frac{2}{2l+1} & \mbox{si $l= m$}
\end{array}\right.
\end{align}
\textit{D\'emonstration}\\
$P_l(x)$ et $P_m(x)$ satisfont l'\'equation de Legendre qui peut \^etre \'ecrite sous la forme
\begin{align}
&\frac{d}{dx}\left\lbrace(1-x^{2})\frac{dP_l(x)}{dx}\right\rbrace+l(l+1)P_l(x)=0\label{eq:5_10}\\
&\frac{d}{dx}\left\lbrace(1-x^{2})\frac{dP_m(x)}{dx}\right\rbrace+m(m+1)P_m(x)=0\label{eq:5_11}
\end{align}
En multipliant (\ref{eq:5_10}) par $P_m(x)$ et (\ref{eq:5_11}) par $P_l(x)$ et en faisant la soustraction des deux \'equations et en int\'egrant par rapport \`a $x$ de $-1$ \`a $1$, on trouve
\begin{align*}
&\int_{-1}^{1}\left[P_m(x)\frac{d}{dx}\left\lbrace(1-x^{2})\frac{dP_l(x)}{dx}\right\rbrace-P_l(x)\frac{d}{dx}\left\lbrace(1-x^{2})\frac{dP_m(x)}{dx}\right\rbrace\right]dx\\
&\hspace{30mm} +\left\lbrace l(l+1)-m(m+1)\right\rbrace\int_{-1}^{1}P_l(x)P_m(x) dx=0
\end{align*}
Maintenant, en remplaçant les deux termes de la premi\`ere int\'egrale par
\begin{align*}
&P_m(x)\frac{d}{dx}\left\lbrace(1-x^{2})\frac{dP_l(x)}{dx}\right\rbrace=\frac{d}{dx}\left\lbrace P_m(x)(1-x^{2})\frac{dP_l(x)}{dx}\right\rbrace-\frac{dP_m(x)}{dx}(1-x^{2})\frac{dP_l(x)}{dx}\\
&P_l(x)\frac{d}{dx}\left\lbrace(1-x^{2})\frac{dP_m(x)}{dx}\right\rbrace=\frac{d}{dx}\left\lbrace P_l(x)(1-x^{2})\frac{dP_m(x)}{dx}\right\rbrace-\frac{dP_l(x)}{dx}(1-x^{2})\frac{dP_m(x)}{dx}
\end{align*}
on obtient
\begin{align*}
&\left[P_m(x)(1-x^{2})\frac{dP_l(x)}{dx}-P_l(x)(1-x^{2})\frac{dP_m(x)}{dx}\right]_{-1}^{1}\\
&\hspace{30mm}+(l-m)(l+m+1)\int_{-1}^{1}P_l(x)P_m(x) dx=0
\end{align*}
ce qui montre que si $l\neq m$, on devrait avoir
\begin{align*}
\int_{-1}^{1}P_l(x)P_m(x)dx=0
\end{align*}
Pour montrer que $\int_{-1}^{1}(P_l)^{2} dx=2/(2l+1)$, on utilise la fonction g\'en\'eratrice Eq.~(\ref{eq:5_6}),
\begin{align*}
\frac{1}{1-2tx+t^{2}}&=\left[\sum_{l=0}^{\infty}t^{l}P_l(x)\right]^{2}\\
&=\sum_{l=0}^{\infty}t^{l}P_l(x)\sum_{m=0}^{\infty}t^{m}P_m(x)
\end{align*}
En int\'egrant les deux c\^ot\'es par rapport \`a $x$, on obtient
\begin{align*}
\sum_{l,m=0}^{\infty}t^{l+m}\int_{-1}^{1}P_l(x)P_m(x)dx&=\left[-\frac{1}{2t}\ln\left(1+t^{2}-2t\right)\right]_{-1}^{1}\\
&=\frac{1}{t}\left[\ln(1+t)-\ln(1-t)\right]\\
&=2\sum_{l=0}^{\infty}\frac{t^{2l}}{2l+1}
\end{align*}
La derni\'ere ligne r\'esulte du d\'eveloppement en s\'erie de $\ln(1+t)$ et $\ln(1-t)$. En identifiant les coefficients de m\^eme puissance de $t$, il en r\'esulte que si $l=m$
\begin{align*}
\int_{-1}^{1}\left[P_l(x)\right]^{2} dx=\frac{2}{2l+1}
\end{align*}
\subsection{S\'eries de Legendre}
Si $f(x)$ est un polyn\^ome de degr\'e $n$, alors il peut \^etre \'ecrit sous la forme
\begin{align}\label{eq:5_12}
f(x)=\sum_{r=0}^{n}c_rP_r(x)
\end{align}
avec les coefficients $c_r$ donn\'es par
\begin{align}\label{eq:5_13}
c_r=\left(r+\frac{1}{2}\right)\int_{-1}^1f(x)P_r(x)dx
\end{align}
\textit{D\'emonstration}\\
Si $f(x)$ est un polyn\^ome de degr\'e $n$, on peut \'ecrire
\begin{align*}
f(x)=b_nx^{n} + b_{n-1}x^{n-1}+\ldots+b_1x+b_0
\end{align*}
D'apr\'es l'Eq.~(\ref{eq:5_4}), le polyn\^ome de Legendre s'\'ecrit sous la forme
\begin{align*}
P_n(x)=a_nx^{n}+a_{n-2}x^{n-2}+\ldots
\end{align*}
Si on multiplie la derni\`ere expression par $b_n/a_n$ et en le soustrayant de $f(x)$, on trouve
que la diff\'erence est un polyn\^ome de degr\'e $(n-1)$
\begin{align*}
f(x)-c_nP_n(x)=g_{n-1}(x)
\end{align*}
o\`u $c_n=b_n/a_n$ et $g_{n-1}(x)$ est un polyn\^ome de degr\'e $n-1$.
En faisant la m\^eme chose pour $g_{n-1}(x)$, on peut d\'emontrer facilement que
$g_{n-1}(x)$ peut s'\'ecrire sous la forme
\begin{align*}
g_{n-1}(x)=c_{n-1}P_{n-1}(x) + g_{n-2}(x)
\end{align*}
donc
\begin{align*}
f(x)=c_nP_n(x)+c_{n-1}P_{n-1}(x) + g_{n-2}(x)
\end{align*}
On fait la m\^eme chose pour $g_{n-2}(x)$, $g_{n-3}(x)$ et ainsi de suite, on obtient le r\'esultat d\'esir\'e
\begin{align*}
f(x)&=c_nP_n(x)+c_{n-1}P_{n-1}(x) + c_{n-2}P_{n-2}(x)+\ldots+c_0P_0(x)\\
&=\sum_{r=0}^{n}c_rP_r(x)
\end{align*}
Les coefficients $c_n$ peuvent \^etre calculer de la mani\`ere suivante
\begin{align*}
\int_{-1}^{1}f(x)P_l(x)dx&=\sum_{r=0}^{n}\int_{-1}^{1}c_nP_n(x)P_l(x)dx\\
&=\frac{2c_l}{2l+1}
\end{align*}
donc
\begin{align*}
c_l=\left(l+\frac{1}{2}\right)\int_{-1}^{1}f(x)P_l(x)dx
\end{align*}
Si $f(x)$ est un polyn\^ome de degr\'e inf\'erieur \`a $l$, alors
\begin{align*}
\int_{-1}^1f(x)P_l(x)dx=0
\end{align*}
Si $f(x)$ est de degr\'e $n$ tel que $n<l$, d'apr\'es la relation d'orthogonalit\'e
on a
\begin{align*}
c_l&=\left(n+\frac{1}{2}\right)\int_{-1}^{1}f(x)P_l(x)dx\\
&=\left(n+\frac{1}{2}\right)\sum_{r=0}^{n}\int_{-1}^{1}P_r(x)P_l(x)dx\\
&=0
\end{align*}
car $r$ est toujours inf\'erieur \`a $l$ ($r<l$).
\subsection{Relations de r\'ecurrence}
\begin{align}\label{eq:5_14}
&(a)\hspace{1.5mm} xP_l(x)=\frac{l+1}{2l+1}P_{l+1}(x)+\frac{l}{2l+1}P_{l-1}(x)\nonumber\\
&(b)\hspace{1.5mm} P_l'(x)=\sum_{r=0}^{[\frac{1}{2}(l+1)]}(2l-4r-1)P_{l-2r-1}(x)\nonumber\\
&(c)\hspace{1.5mm}P_{l+1}'(x)-P_{l-1}'(x)=(2l+1)P_l(x)\\
&(d)\hspace{1.5mm}xP_l'(x)-P_{l-1}'(x)=lP_{l}(x)\nonumber\\
&(e)\hspace{1.5mm}P_l'(x)-xP_{l-1}'(x)=lP_{l-1}(x)\nonumber
\end{align}
\textit{D\'emonstration}:\\
(a) $P_l(x)$ est un polyn\^ome d'ordre $l$, donc $xP_l(x)$ est un polyn\^ome d'ordre $l+1$. D'apr\'es (\ref{eq:5_12}), on a
\begin{align*}
xP_l(x)=c_{l+1}P_{l+1}(x)+c_{l-1}P_{l-1}(x)+\ldots+\left\lbrace\begin{array}{ll}
c_1P_1(x) & \mbox{si $l$ est pair}\\
c_0P_0(x) & \mbox{si $l$ est impair}
\end{array}\right.
\end{align*}
o\`u les coefficients $c_r$ sont donn\'es par
\begin{align*}
c_r&=\left(r+\frac{1}{2}\right)\int_{-1}^{1}xP_l(x)P_r(x)dx\\
&=\left(r+\frac{1}{2}\right)\int_{-1}^{1}P_l(x)\left\lbrace xP_r(x)\right\rbrace dx
\end{align*}
$xP_r(x)$ est un polyn\^ome d'ordre $r+1$, donc $c_r=0$ si $r+1<l$, \i.e. si $r<l-1$. alors
\begin{align}\label{eq:5_15}
xP_l(x)=c_{l+1}P_{l+1}(x)+c_{l-1}P_{l-1}(x)
\end{align}
Pour calculer $c_{l+1}$ et $c_{l-1}$, on pose $x=1$ dans l'expression (\ref{eq:5_15}) et dans sa d\'eriv\'ee, on obtient
\begin{align*}
&P_l(1)=c_{l+1}P_{l+1}(1)+c_{l-1}P_{l-1}(1)\\
&P_l(1)+P_l'(1)=c_{l+1}P_{l+1}'(1)+c_{l-1}P_{l-1}'(1)
\end{align*}
En utilisant les propri\'et\'es (\ref{eq:5_8})(a) et (c), on obtient les \'equations suivantes
\begin{align*}
&1=c_{l+1}+c_{l-1}\\
&1+\frac{1}{2}l(l+1)=c_{l+1}\left[\frac{1}{2}(l+1)(l+2)\right]+c_{l-1}\left[\frac{1}{2}(l-1)l\right]
\end{align*}
qui admettent comme solutions
\begin{align*}
c_{l+1}=\frac{l+1}{2l+1},\hspace{20mm}c_{l-1}=\frac{l}{2l+1}
\end{align*}
On obtient finalement
\begin{align*}
xP_l(x)=\frac{l+1}{2l+1}P_{l+1}(x)+\frac{l}{2l+1}P_{l-1}(x)
\end{align*}
(b) $P_l'(x)$ est un polyn\^ome de degr\'e $l-1$, donc on peut l'\'ecrire en s\'erie de Legendre
\begin{align*}
P_l'(x)=c_{l-1}P_{l-1}(x)+c_{l-3}P_{l-3}(x)+\ldots+&c_{l-2r+1}P_{l-2r+1}(x)\\
&+\ldots+\left\lbrace\!\!\begin{array}{ll}
c_1P_1(x) & \mbox{si $l$ est pair}\\
c_0P_0(x) & \mbox{si $l$ est impair}
\end{array}\right.
\end{align*}
avec 
\begin{align*}
c_n=\left(n+\frac{1}{2}\right)\int_{-1}^{1}P_l'(x)P_n(x)dx
\end{align*}
En int\'egrant par parties, on trouve
\begin{align*}
c_n=\left(n+\frac{1}{2}\right)\left\lbrace\left[P_l(x)P_n(x)\right]_{-1}^{1}-\int_{-1}^{1}P_l(x)P_n'(x)dx\right\rbrace
\end{align*}
et puisque $P_n'(x)$ est un polyn\^ome de degr\'e $n-1<l$, donc l'int\'egrale dans l'\'equation pr\'ec\'edente est nulle. En utilisant les propri\'et\'es (\ref{eq:5_8})(a) et (b) des polyn\^omes de Legendre on obtient
\begin{align*}
c_n=\left(n+\frac{1}{2}\right)\left[1-(-1)^{n+l}\right]
\end{align*}
$n$ prends les valeurs $l-1$, $l-3$, $l-5$,\ldots, donc $n+l$ prends les valeurs $2l-1$, $2l-3$, $2l-5$,\ldots. Puisque toutes les valeurs sont impaires, $(-1)^{s+l}=-1$, donc
\begin{align*}
c_n=2n+1
\end{align*} 
et
\begin{align*}
c_{l-2r-1}=2l-4r-1
\end{align*}
Il en r\'esulte que
\begin{align*}
P_l'(x)&=(2l-1)P_{l-1}(x)+(2l-5)P_{l-5}(x)+\ldots+(2l-4r-1)P_{2l-4r-1}(x)\\
&\hspace{20mm}+\ldots+\left\lbrace\begin{array}{ll}
3P_1(x) & \mbox{si $l$ est pair}\\
P_0(x) & \mbox{si $l$ est impair}
\end{array}\right.\\
&=\sum_{r=0}^{[\frac{1}{2}(l-1)]}(2l-4r-1)P_{2l-4r-1}(x)
\end{align*}
(c) On applique la relation (b) pour $P_{l+1}'(x)$ et $P_{l-1}'(x)$, on obtient
\begin{align*}
&P_{l+1}'(x)=(2l+1)P_l(x)+(2l-3)P_{l-2}(x)+(2l-7)P_{l-4}(x)+\ldots\\
&P_{l-1}'(x)=(2l-3)P_{l-2}(x)+(2l-7)P_{l-4}(x)+\ldots
\end{align*}
la soustraction des deux expressions donne
\begin{align*}
P_{l+1}'(x)-P_{l-1}'(x)=(2l+1)P_l(x)
\end{align*}
(d) En multipliant (a) par $2l+1$ et d\'erivant par rapport \`a $x$, on obtient
\begin{align*}
(l+1)P_{l+1}'(x)-(2l+1)\left\lbrace P_l(x)+xP_l'(x)\right\rbrace+lP_{l-1}'(x)=0
\end{align*}
On remplace $P_{l+1}'(x)$ par $(2l+1)P_l(x)-P_{l-1}'(x)$ (Eq.~(\ref{eq:5_14})(c)), et on r\'earrange les termes de l'\'equation on trouve que
\begin{align*}
xP_l'(x)-P_{l-1}'(x)=lP_l(x)
\end{align*}
(e) En multipliant (c) par $x$ et en substituant le r\'esultat pour $(2l+1)xP_l(x)$ dans (\ref{eq:5_14})(d), on obtient
\begin{align}\label{eq:5_16}
(l+1)P_{l+1}(x)+lP_{l-1}(x)=xP_{l+1}'(x)-xP_{l-1}'(x)
\end{align}
Si on \'ecrit (d) avec $l$ remplac\'e par $l+1$, on obtient
\begin{align*}
xP_{l+1}'(x)-P_{l}'(x)=(l+1)P_{l+1}(x)
\end{align*}
Maintenant, en substituant $(l+1)P_{l+1}(x)$ dans l'Eq.~(\ref{eq:5_16}) et en r\'earrangeant les termes, on obtient le r\'esultat d\'esir\'e
\begin{align*}
P_{l}'(x)-xP_{l-1}'(x)=lP_{l-1}(x)
\end{align*}
\section{Polyn\^ome associ\'e de Legendre}
\begin{figure}[htb]
\centering
\includegraphics[height=9cm,width=13cm]{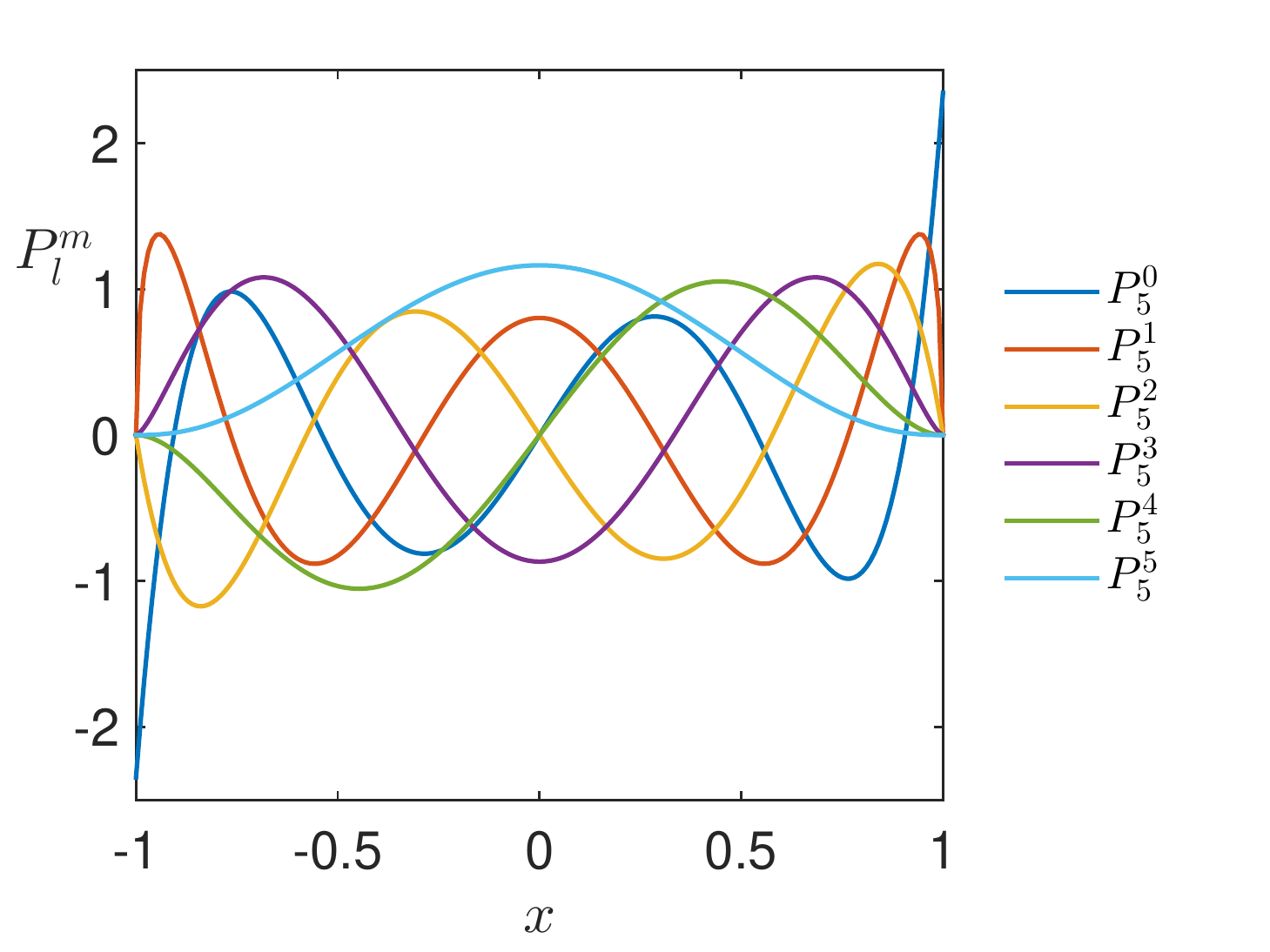}
\caption{\label{fig:Alegendre5} Polyn\^omes associ\'es de Legendre $P_l^{m}(x)$ pour $l=5$ et $0<m<l$.}
\end{figure}

L'\'equation associ\'ee de Legendre est donn\'ee par
\begin{align}\label{eq:5_17}
(1-x^{2})\frac{d^{2}y}{dx^{2}}-2x\frac{dy}{dx}+\left[ l(l+1)-\frac{m^{2}}{1-x^{2}}\right]y=0
\end{align}
L'Eq.~(\ref{eq:5_17}) admet une solution r\'eguli\`ere uniquement dans l'intervalle $[-1,\hspace{1mm}1]$ avec $0\leq m \leq l$. Pour $m=0$, l'Eq.~(\ref{eq:5_17}) se r\'eduit \`a l'\'equation de Legendre. La solution de cette \'equation not\'e $P_l^{m}(x)$ est donn\'ee par
\begin{align}\label{eq:5_18}
P_l^{m}(x)=(1-x^{2})^{m/2}\frac{d^{m}}{dx^{m}}P_l(x)
\end{align}
La solution (\ref{eq:5_18}) est appell\'ee polyn\^ome associ\'e de Legendre. Pour des valeurs n\'egatives de $m$, on utilise
\begin{align}\label{eq:5_19}
P_l^{-m}(x)=(-1)^{m}\frac{(l-m)!}{(l+m)!}P_l^{m}(x)
\end{align}
o\`u encore $0\leq m \leq l$.\\
\textit{D\'emonstration}\\
Pour montrer que (\ref{eq:5_18}) est une solution de l'Eq.~(\ref{eq:5_17}), on doit montrer que
\begin{align*}
(1-x^{2})\frac{d^{2}}{dx^{2}}P_l^{m}(x) -2x\frac{d}{dx}P_l^{m}(x)+\left[l(l+1)-\frac{m^{2}}{1-x^{2}}\right]P_l^{m}(x)=0
\end{align*}
Puisque $P_l(x)$ est une solution de l'\'equation de Legendre, on a
\begin{align*}
\left(1-x^{2}\right)\frac{d^{2}}{dx^{2}}P_l(x) -2x\frac{d}{dx}P_l(x)+ l(l+1)P_l(x)=0
\end{align*}
En d\'erivant l'\'equation pr\'ec\'edente  $m$ fois par rapport a $x$, elle devient
\begin{align*}
\frac{d^{m}}{dx^{m}}\left[\left(1-x^{2}\right)\frac{d^{2}}{dx^{2}}P_l(x)\right] -2\frac{d^{m}}{dx^{m}}\left[ x\frac{d}{dx}P_l(x)\right]+ l(l+1)\frac{d^{m}}{dx^{m}}P_l(x)=0
\end{align*}
En utilisant la r\`egle de Leibniz pour le $m-$i\`eme d\'erivatives 
\begin{align*}
\frac{d^{m}}{dx^{m}}\left\lbrace fg\right\rbrace=\sum_{n=0}^{m}\frac{m!}{n!(n-m)!}\frac{d^{n}f}{dx^{n}}\frac{d^{n-m}g}{dx^{n-m}}
\end{align*} 
l'\'equation pr\'ec\'edente devient
\begin{align*}
\left(1-x^{2}\right)\frac{d^{m+2}}{dx^{m+2}}P_l(x)&+m\frac{d}{dx}\left(1-x^{2}\right)\frac{d^{m+1}}{dx^{m+1}}P_l(x)+\frac{m(m-1)}{2}\frac{d^{2}}{dx^{2}}\left(1-x^{2}\right)\frac{d^{m}}{dx^{m}}P_l(x)\\
&-2\left[x\frac{d^{m+1}}{dx^{m+1}}P_l(x)+m\frac{d}{dx}x\frac{d^{m}}{dx^{m}}P_l(x)\right] +l(l+1)\frac{d^{m}}{dx^{m}}P_l(x)=0
\end{align*}
En regroupant les coefficients de $d^{m+2}P_l/dx^{m+2}$, $d^{m+1}P_l/dx^{m+1}$ et $d^{m}P_l/dx^{m}$, on obtient
\begin{align}\label{eq:5_19a}
\left(1-x^{2}\right)\frac{d^{m+2}}{dx^{m+2}}&P_l(x)-2x(m+1)\frac{d^{m+1}}{dx^{m+1}}P_l(x)\\
&+\left[l(l+1)-m(m+1)-2m\right]\frac{d^{m}}{dx^{m}}P_l(x)=0
\end{align}
qui devient, on d\'enote $y_{l}^{m}(x)=d^{m}P_l(x)/dx^{m}$
\begin{align*}
\left(1-x^{2}\right)\frac{d^{2}}{dx^{2}}y_l^{m}(x)-2x(m+1)\frac{d}{dx}y_l^m(x)+\left[l(l+1)-m(m+1)-2m\right]y_l^m(x)=0
\end{align*}
Si on remplace maintenant
\begin{align*}
P_l^m(x)=\left(1-x^{2}\right)^{m/2}y_l^m(x)
\end{align*}
dans l'\'equation pr\'ec\'edente, on obtient
\begin{align}\label{eq:5_19b}
\left(1-x^{2}\right)\frac{d^{2}}{dx^{2}}&\left[\left(1-x^{2}\right)^{-m/2}P_l^{m}(x)\right]-2(m+1)x\frac{d}{dx}\left[\left(1-x^{2}\right)^{-m/2}P_l^{m}(x)\right]\nonumber\\
&+\left[l(l+1)-m(m+1)\right]\left(1-x^{2}\right)^{-m/2}P_l^{m}(x)=0
\end{align}
mais
\begin{align*}
\frac{d}{dx}\left[\left(1-x^{2}\right)^{-m/2}P_l^{m}(x)\right]=\frac{dP_l^{m}(x)}{dx}\left(1-x^{2}\right)^{-m/2}+mx\left(1-x^{2}\right)^{-(m/2)-1}P_l^{m}(x)
\end{align*}
et
\begin{align*}
\frac{d^{2}}{dx^{2}}&\left[\left(1-x^{2}\right)^{-m/2}P_l^{m}(x)\right]=\\
&\frac{d^{2}P_l^{m}(x)}{dx^{2}}\left(1-x^{2}\right)^{-m/2}+mx\left(1-x^{2}\right)^{-(m/2)-1}\frac{dP_l^{m}(x)}{dx}+mx\left(1-x^{2}\right)^{-(m/2)-1}P_l^m(x)
\end{align*}
donc, l'Eq.~(\ref{eq:5_19b}) devient
\begin{align*}
\left(1-x^{2}\right)^{-(m/2)+1}&\frac{d^{2}P_l^{m}(x)}{dx^{2}}+2mx\left(1-x^{2}\right)^{-m/2}\frac{dP_l^{m}(x)}{dx}+m\left(1-x^{2}\right)^{-m/2}P_l^m(x)\\
&+m(m+2)\left(1-x^{2}\right)^{-(m/2)-1}x^2P_l^m(x)\\
&-2(m+1)x\left[\left(1-x^{2}\right)^{-m/2}\frac{dP_l^{m}(x)}{dx} +mx\left(1-x^{2}\right)^{-(m/2)-1}P_l^m(x)\right]\\
&+\left[l(l+1)-m(m+1)\right]\left(1-x^{2}\right)^{-m/2}P_l^m(x)=0
\end{align*}
\newpage
On annule le facteur commun $\left(1-x^{2}\right)^{-m/2}$, on obtient
\begin{align*}
\left(1-x^{2}\right)&\frac{d^{2}P_l^{m}(x)}{dx^{2}} +\left[2mx-2(m+1)x\right]\frac{dP_l^{m}(x)}{dx}\\
&+\left[m+\frac{m(m+1)}{1-x^2}-\frac{2(m+1)mx^2}{1-x^2}+l(l+1)-m(m+1) \right]P_l^m(x)=0
\end{align*}
qui peut \^etre simplifi\'e \`a
\begin{align*}
\left(1-x^{2}\right)\frac{d^{2}P_l^{m}(x)}{dx^{2}}-2x\frac{dP_l^{m}(x)}{dx}+\left[l(l+1)-\frac{m^2}{1-x^2}\right]P_l^m(x)=0
\end{align*}
qui est l'\'equation d\'esir\'e.\\
Les polyn\^omes associ\'e de Legendre $P_l^{m}(x)$ pour $l=5$ et $0<m<l$ sont repr\'esent\'es dans la Fig.~(\ref{fig:Alegendre5}).
Les premiers polyn\^omes associ\'es de Legendre, y compris ceux des valeurs n\'egatives de $m$, sont les suivants
\begin{align*}
&P_0^0(x)=1\\
&P_1^{1}(x)=-\sqrt{1-x^2}\\
&P_1^0(x)=x\\
&P_1^{-1}=-\frac{1}{2}P_1^{1}(x)=\frac{1}{2}\sqrt{1-x^2}\\
&P_2^2(x)=3(1-x^2)\\
&P_2^1(x)=-3x\sqrt{1-x^2}\\
&P_2^0(x)=\frac{1}{2}(3x^2-1)\\
&P_2^{-1}(x)=-\frac{1}{6}P_2^1(x)=\frac{1}{2}x\sqrt{1-x^2}\\
&P_2^{-2}(x)=\frac{1}{24}P_2^2(x)=\frac{1}{8}(1-x^2)\\
&P_3^3(x)=-15(1-x^2)^{3/2}\\
&P_3^2(x)=15x(1-x^2)\\
&P_3^1(x)=-\frac{3}{2}(5x^2-1)\sqrt{1-x^2}\\
&P_3^0(x)=\frac{1}{2}(5x^3-3x)\\
&P_3^{-1}(x)=-\frac{1}{12}P_3^1(x)=-\frac{1}{8}(5x^2-1)\sqrt{1-x^2}\\
&P_3^{-2}(x)=\frac{1}{120}P_3^2(x)=\frac{1}{8}x(1-x^2)\\
&P_3^{-3}(x)=-\frac{1}{720}P_3^3(x)=\frac{1}{48}(1-x^2)^{3/2}\\
&P_4^4(x)=105(1-x^2)^2\\
&P_4^3(x)=-105x(1-x^2)^{3/2}\\
&P_4^2(x)=\frac{15}{2}(7x^2-1)(1-x^2)\\
&P_4^1(x)=-\frac{5}{2}(7x^3-3x)\sqrt{1-x^2}\\
&P_4^0(x)=\frac{1}{8}(35x^4-30x^2+3)\\
&P_4^{-1}(x)=-\frac{1}{20}P_4^1(x)=\frac{1}{8}(7x^3-3x)\sqrt{1-x^2}\\
&P_4^{-2}(x)=\frac{1}{360}P_4^2(x)=\frac{1}{48}(7x^2-1)(1-x^2)\\
&P_4^{-3}(x)=-\frac{1}{5040}P_4^3(x)=\frac{x}{48}(1-x^2)^{3/2}\\
&P_4^{-4}(x)=\frac{1}{40320}P_4^4(x)=\frac{1}{384}(1-x^2)^2\\
\end{align*}
\subsection{Relation d'orthogonalit\'e}
\begin{align}\label{eq:5_20}
\int_{-1}^{1}P_l^{m}(x)P_{l'}^{m}(x)dx=\frac{2(l+m)!}{(2l+1)(l-m)}\delta_{l,l'}
\end{align}
\textit{D\'emonstration}\\
Pour d\'emontrer que si $l\neq l'$, l'int\'egral donne $0$, nous suivons les m\^emes \'etapes que celle pour la relation d'orthogonalit\'e de polyn\^ome de Legendre Eq.~(\ref{eq:5_9}).
$P_l^{m}(x)$ et $P_{l'}^{m}(x)$ satisfy l'\'equation de Legendre associ\'ee qui peut \^etre \'ecrire sous la form
\begin{align}
&\frac{d}{dx}\left\lbrace(1-x^{2})\frac{dP_l^{m}}{dx}\right\rbrace+\left[l(l+1)-\frac{m^{2}}{1-x^2}\right] P_l^{m}=0\label{eq:5_20a}\\
&\frac{d}{dx}\left\lbrace(1-x^{2})\frac{dP_{l'}^{m}}{dx}\right\rbrace+\left[l'(l'+1)-\frac{m^{2}}{1-x^2}\right] P_{l'}^{m}=0\label{eq:5_20b}
\end{align}
En multipliant Eq.~(\ref{eq:5_20a}) par $P_{l'}^{m}(x)$ et l'Eq.~(\ref{eq:5_20b}) par $P_l^{m}(x)$, en fiasant la substraction des deux \'equations et int\'egrant par rapport \`a $x$ de $-1$ \`a $1$, on trouve
\begin{align*}
&\int_{-1}^{1}\left[P_{l'}^{m}\frac{d}{dx}\left\lbrace(1-x^{2})\frac{dP_l^{m}}{dx}\right\rbrace-P_l^{m}\frac{d}{dx}\left\lbrace(1-x^{2})\frac{dP_{l'}^{m}}{dx}\right\rbrace\right]dx\\
&\hspace{20mm} +\left\lbrace l(l+1)-l'(l'+1)\right\rbrace\int_{-1}^{1}P_l^{m}P_{l'}^{m} dx=0
\end{align*}
En remplaçant les deux termes du premi\`ere int\'egrale par
\begin{align*}
&P_{l'}^{m}\frac{d}{dx}\left\lbrace(1-x^{2})\frac{dP_l^{m}}{dx}\right\rbrace=\frac{d}{dx}\left\lbrace P_{l'}^{m}(1-x^{2})\frac{dP_l^{m}}{dx}\right\rbrace-\frac{dP_{l'}^{m}}{dx}(1-x^{2})\frac{dP_l^{m}}{dx}\\
&P_l^{m}\frac{d}{dx}\left\lbrace(1-x^{2})\frac{dP_{l'}^{m}}{dx}\right\rbrace=\frac{d}{dx}\left\lbrace P_l^{m}(1-x^{2})\frac{dP_{l'}^{m}}{dx}\right\rbrace-\frac{dP_l^{m}}{dx}(1-x^{2})\frac{dP_{l'}^{m}}{dx}
\end{align*}
il s'ensuit
\begin{align*}
\left[P_{l'}^{m}(1-x^{2})\frac{dP_l^m}{dx}-P_l^m(1-x^{2})\frac{dP_{l'}^{m}}{dx}\right]_{-1}^{1}+(l-l')(l+l'+1)\int_{-1}^{1}P_l^mP_{l'}^{m} dx=0
\end{align*}
donc
\begin{align*}
(l-l')(l+l'+1)\int_{-1}^{1}P_l^m(x)P_{l'}^{m}(x) dx=0
\end{align*}
ce qui montre que si $l\neq l'$, on devrait avoir
\begin{align*}
\int_{-1}^{1}P_l^m(x)P_{l'}^{m}(x)dx=0
\end{align*}
Montrons maintenant que si $l=l'$, on a
\begin{align*}
\int_{-1}^{1}\left[ P_l^m(x)\right]^2dx=\frac{2(l+m)!}{(2l+1)(l-m)!}
\end{align*}
Supposons tout d'abord que $m>0$. D'apr\`es l'Eq.~(\ref{eq:5_18}), on a
\begin{align}
\int_{-1}^{1}\left[ P_l^m(x)\right]^2dx=&\int_{-1}^{1}\left(1-x^2\right)^{m}\left[\frac{d^m}{dx^m}P_l(x)\right]\left[\frac{d^m}{dx^m}P_l(x)\right]dx\nonumber\\
&=\left[\left(\frac{d^{m-1}}{dx^{m-1}}P_l(x)\right)\left(\left(1-x^2\right)^{m}\frac{d^m}{dx^m}P_l(x)\right)\right]_{-1}^{1}\nonumber\\
&-\int_{-1}^{1}\left[\frac{d^{m-1}}{dx^{m-1}}P_l(x)\right]\frac{d}{dx}\left[\left(1-x^2\right)^{m}\frac{d^m}{dx^m}P_l(x)\right]dx\nonumber\\
&=-\int_{-1}^{1}\left[\frac{d^{m-1}}{dx^{m-1}}P_l(x)\right]\frac{d}{dx}\left[\left(1-x^2\right)^{m}\frac{d^m}{dx^m}P_l(x)\right]dx\label{eq:5_20c}
\end{align}
o\`u nous avons int\'egr\'e par parties. A partir de l'Eq.~(\ref{eq:5_19a}) avec $m$ est remplac\'e par $m-1$, on a
\begin{align*}
\left(1-x^{2}\right)\frac{d^{m+1}}{dx^{m+1}}P_l(x)-2xm\frac{d^m}{dx^m}P_l(x)+\left[l(l+1)-m(m-1)\right]\frac{d^{m-1}}{dx^{m-1}}P_l(x)=0
\end{align*}
En multipliant par $\left(1-x^{2}\right)^{m-1}$, l'\'equation devient
\begin{align*}
\left(1-x^{2}\right)^{m}&\frac{d^{m+1}}{dx^{m+1}}P_l(x)-2xm\left(1-x^{2}\right)^{m-1}\frac{d^m}{dx^m}P_l(x)\\
&+\left[l(l+1)-m(m-1)\right]\left(1-x^{2}\right)^{m-1}\frac{d^{m-1}}{dx^{m-1}}P_l(x)=0
\end{align*}
qui peut \^etre \'ecrit comme
\begin{align*}
\frac{d}{dx}\left[\left(1-x^{2}\right)^{m}\frac{d^{m}}{dx^{m}}P_l(x)\right]=-(l+m)(l-m+1)\left(1-x^{2}\right)^{m-1}\frac{d^{m-1}}{dx^{m-1}}P_l(x)
\end{align*}
On ins\`ere ce r\'esultat dans l'Eq.~(\ref{eq:5_20c}), on obtient
\begin{align*}
\int_{-1}^{1}&\left[ P_l^m(x)\right]^2dx\\
&=\int_{-1}^{1}\left[\frac{d^{m-1}}{dx^{m-1}}P_l(x)\right](l+m)(l-m+1)\left(1-x^{2}\right)^{m-1}\left[\frac{d^{m-1}}{dx^{m-1}}P_l(x)\right]dx\\
&=(l+m)(l-m+1)\int_{-1}^{1}\left(1-x^{2}\right)^{m-1}\left[\frac{d^{m-1}}{dx^{m-1}}P_l(x)\right]^{2}dx\\
&=(l+m)(l-m+1)\int_{-1}^{1}\left[P_l^{m-1}(x)\right]^{2}dx
\end{align*}
Dans la derni\`ere ligne, nous avons utilis\'e la d\'efinition de polyn\^one de Legendre associ\'e Eq.~(\ref{eq:5_18}). En appliquant \`a nouveau ce r\'esultat, on trouve
\begin{align*}
\int_{-1}^{1}\left[ P_l^m(x)\right]^2dx&=(l+m)(l-m+1)(l+m-1)(l-m+2)\int_{-1}^{1}\left[P_l^{m-2}(x)\right]^{2}dx\\
&=(l+m)(l+m-1)(l-m+1)(l-m+2)\int_{-1}^{1}\left[P_l^{m-2}(x)\right]^{2}dx
\end{align*}
En r\'ep\'etant le processus $m$ fois, on trouve
\begin{align*}
\int_{-1}^{1}&\left[ P_l^m(x)\right]^2dx\\
&=(l+m)(l+m-1)\ldots(l+1)(l-m+1)(l-m+2)\ldots l\int_{-1}^{1}\left[P_l^{0}(x)\right]^{2}dx\\
&=(l+m)(l+m-1)\ldots(l+1)l(l-1)\ldots(l-m+2)(l-m+1)\frac{2}{2l+1}\\
&=\frac{(l+m)!}{(l-m)!}\frac{2}{2l+1}
\end{align*}
qui est le r\'esultat d\'esir\'e. Supposons maintenant que $m<0$, \i.e., $m=-n$ avec $n>0$, donc
\begin{align*}
\int_{-1}^{1}\left[ P_l^m(x)\right]^2dx&=\int_{-1}^{1}\left[ P_l^{-n}(x)\right]^2dx\\
&=\int_{-1}^{1}\left[ (-1)^nP\frac{(l-n)!}{(l+n)!}P_l^n(x)\right]^2dx\\
&=\left[P\frac{(l-n)!}{(l+n)!}\right]^2\int_{-1}^{1}\left[P_l^n(x)\right]^2dx\\
&=\left[P\frac{(l-n)!}{(l+n)!}\right]^2\frac{(l+n)!}{(l-n)!}\frac{2}{2l+1}\\
&=\frac{(l-n)!}{(l+n)!}\frac{2}{2l+1}\\
&=\frac{(l+m)!}{(l-m)!}\frac{2}{2l+1}
\end{align*}
\section{Harmoniques sph\'eriques}
\begin{figure}[!ht]
\centering
\includegraphics[scale=0.34]{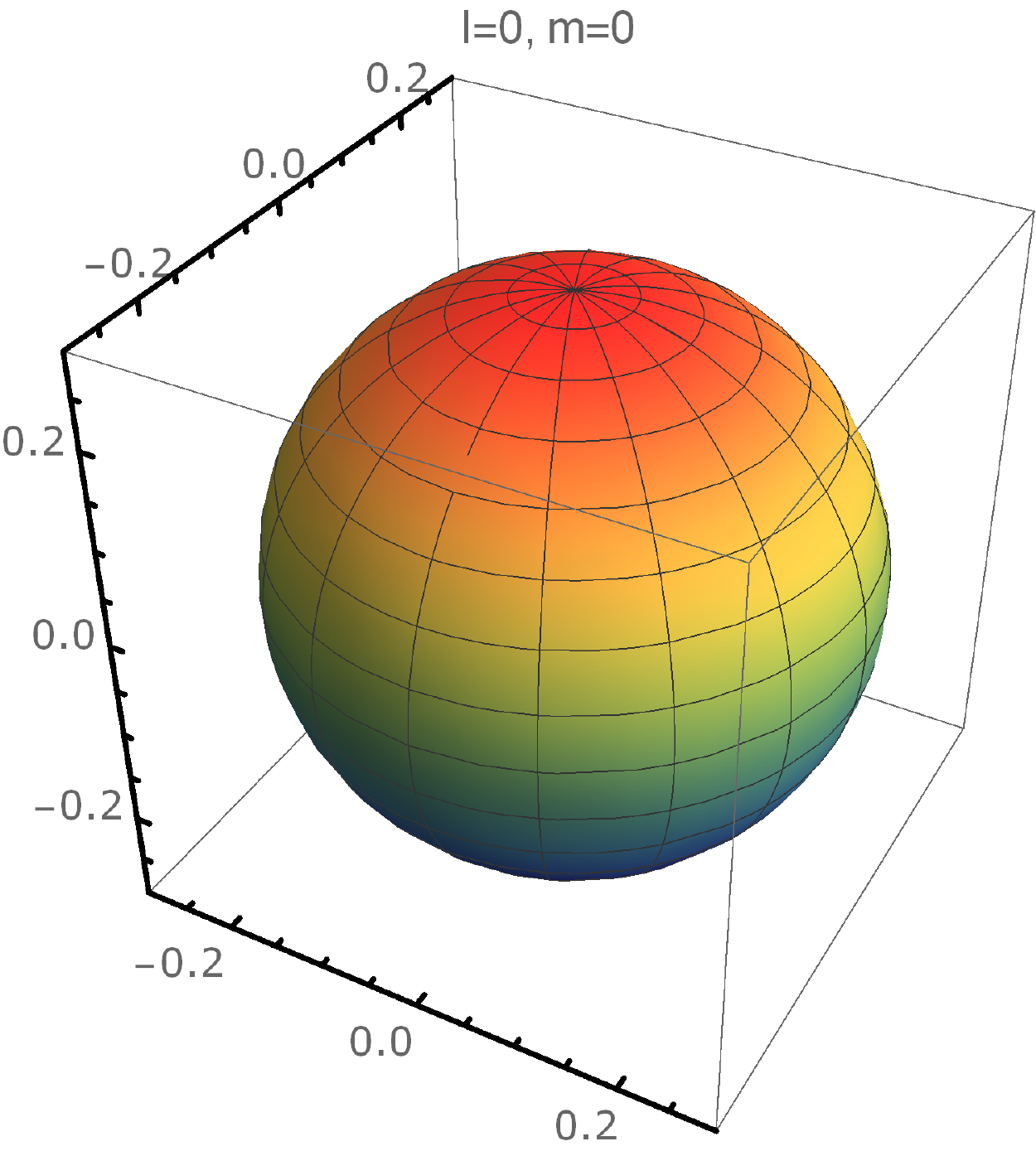}
\includegraphics[scale=0.34]{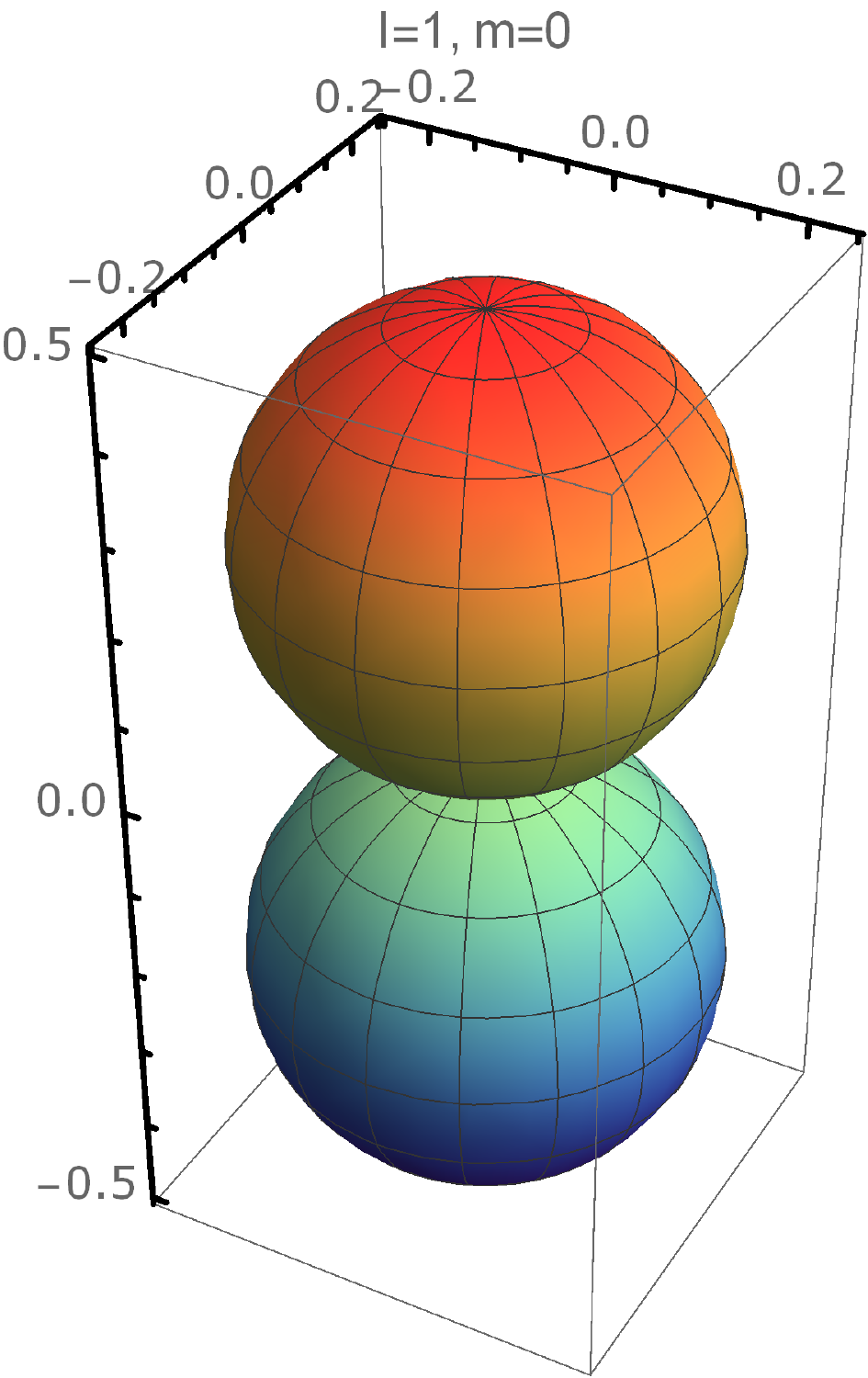}
\includegraphics[scale=0.34]{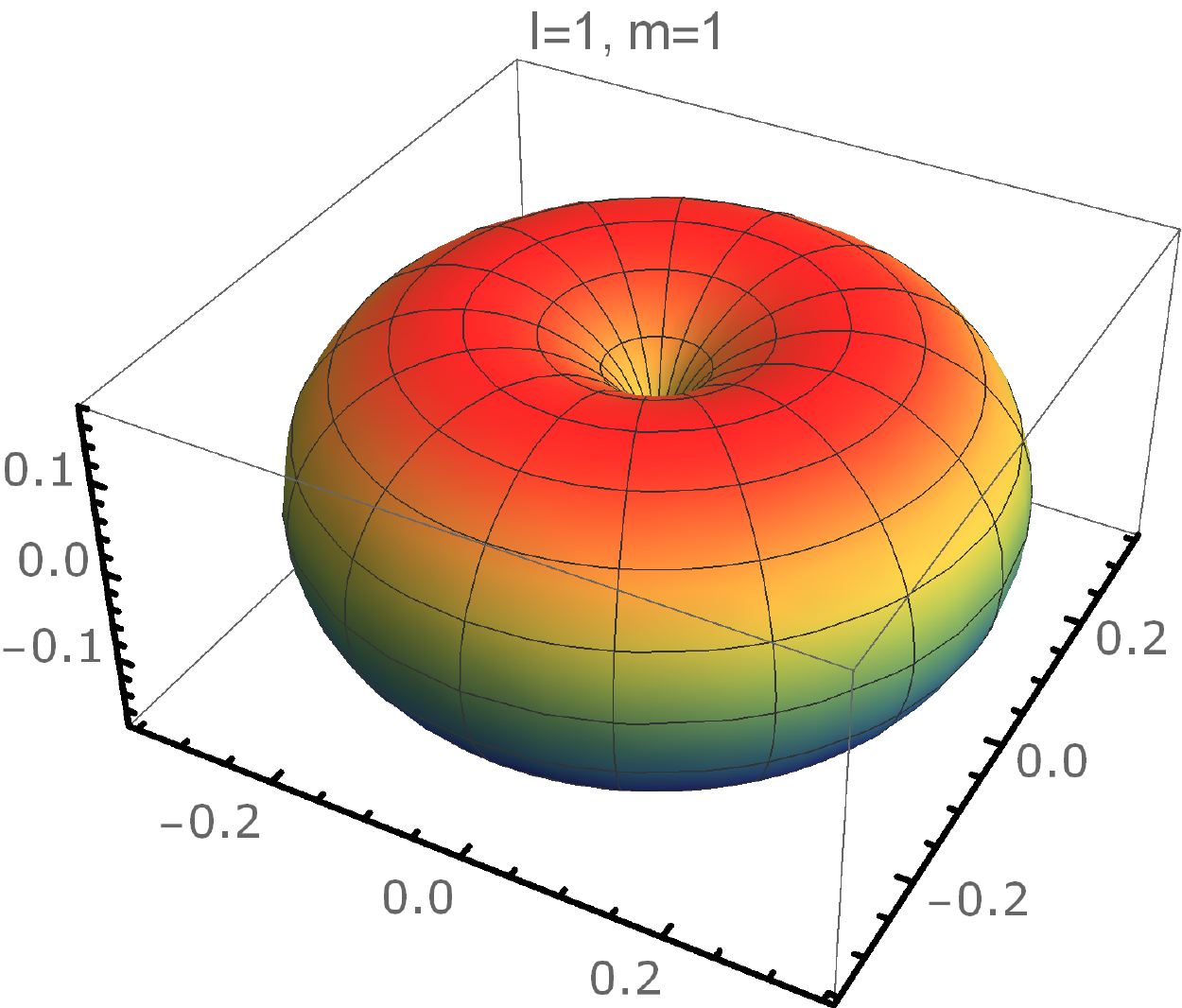}\\
\includegraphics[scale=0.34]{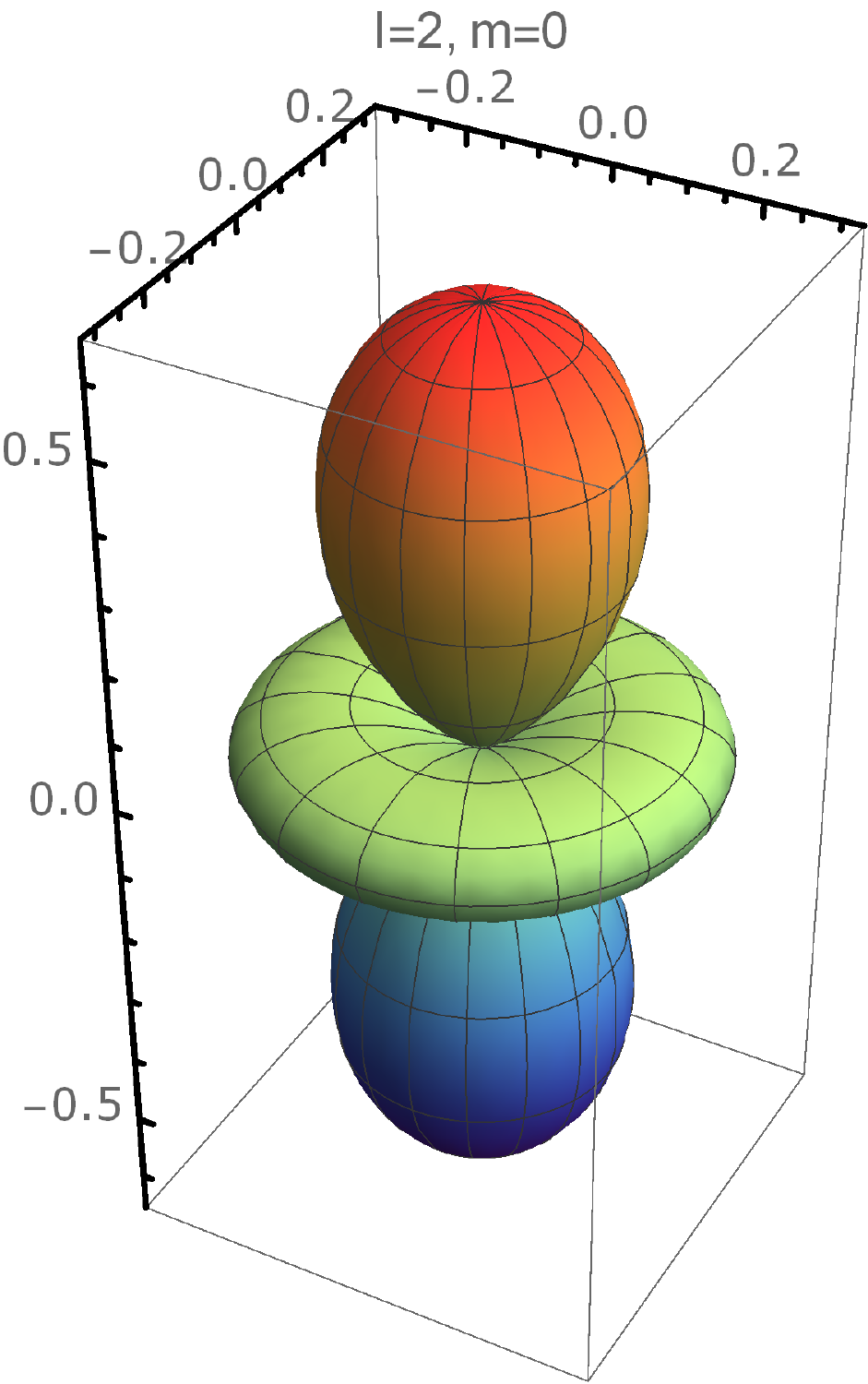}
\includegraphics[scale=0.34]{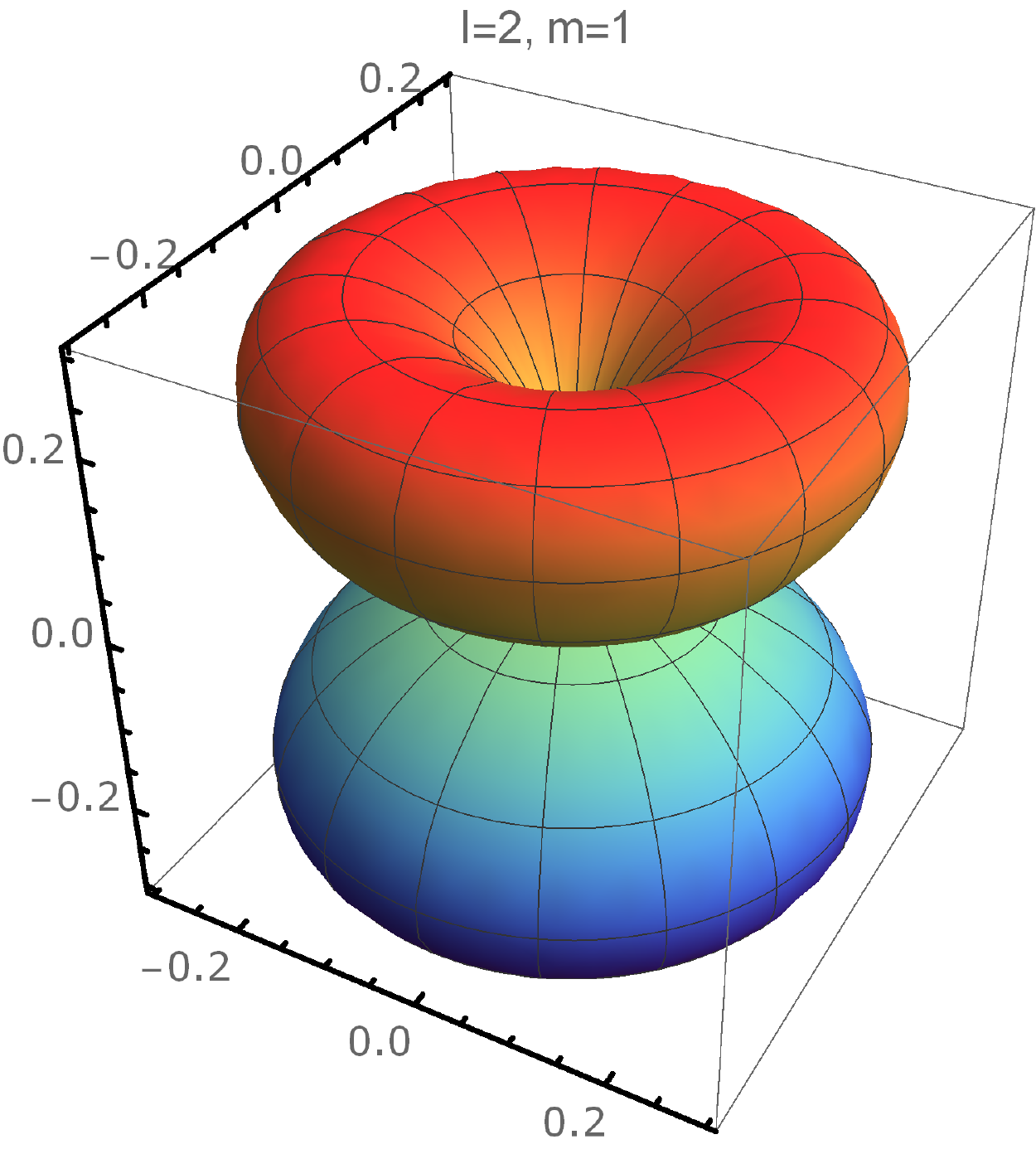}
\includegraphics[scale=0.34]{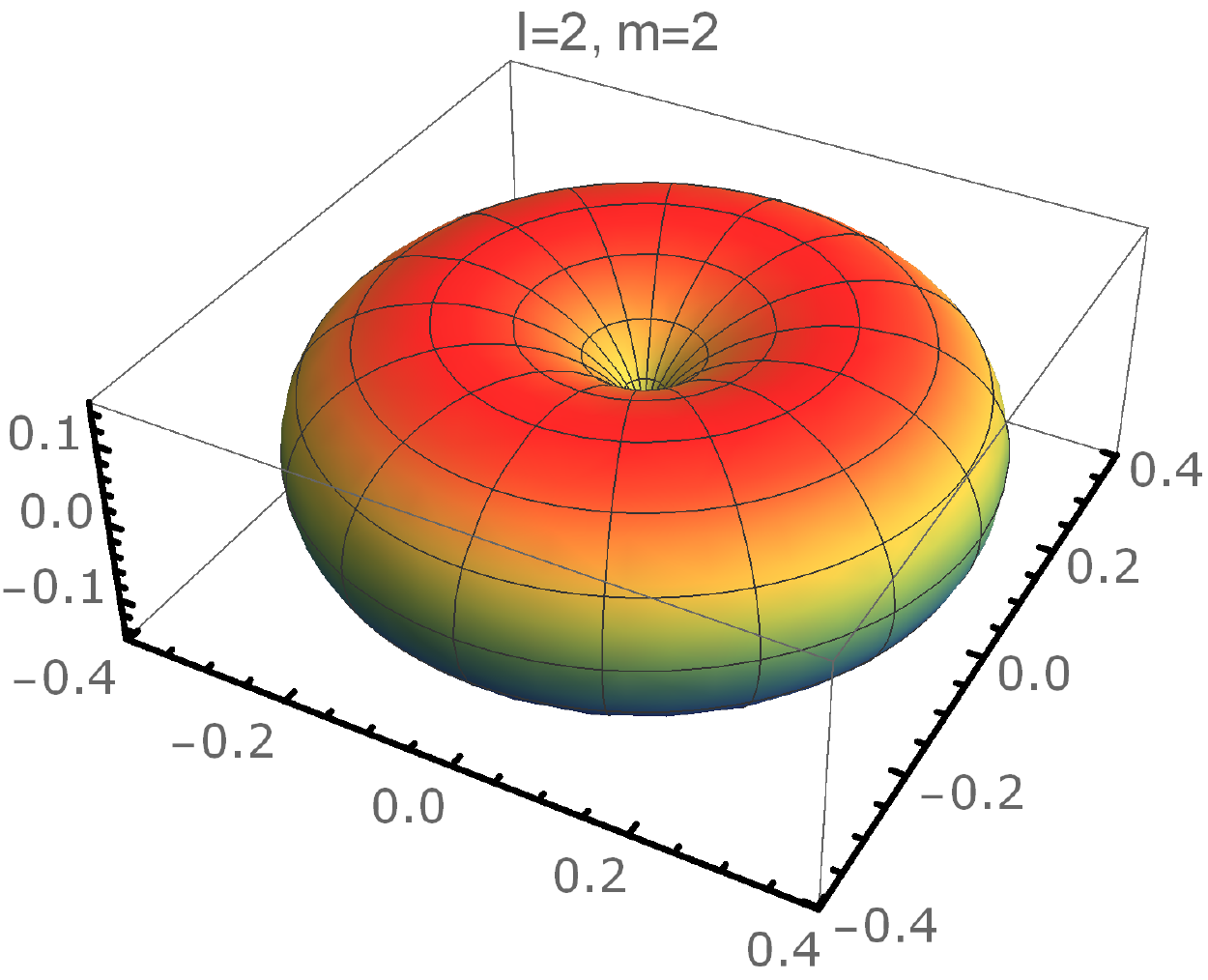}\\
\includegraphics[scale=0.34]{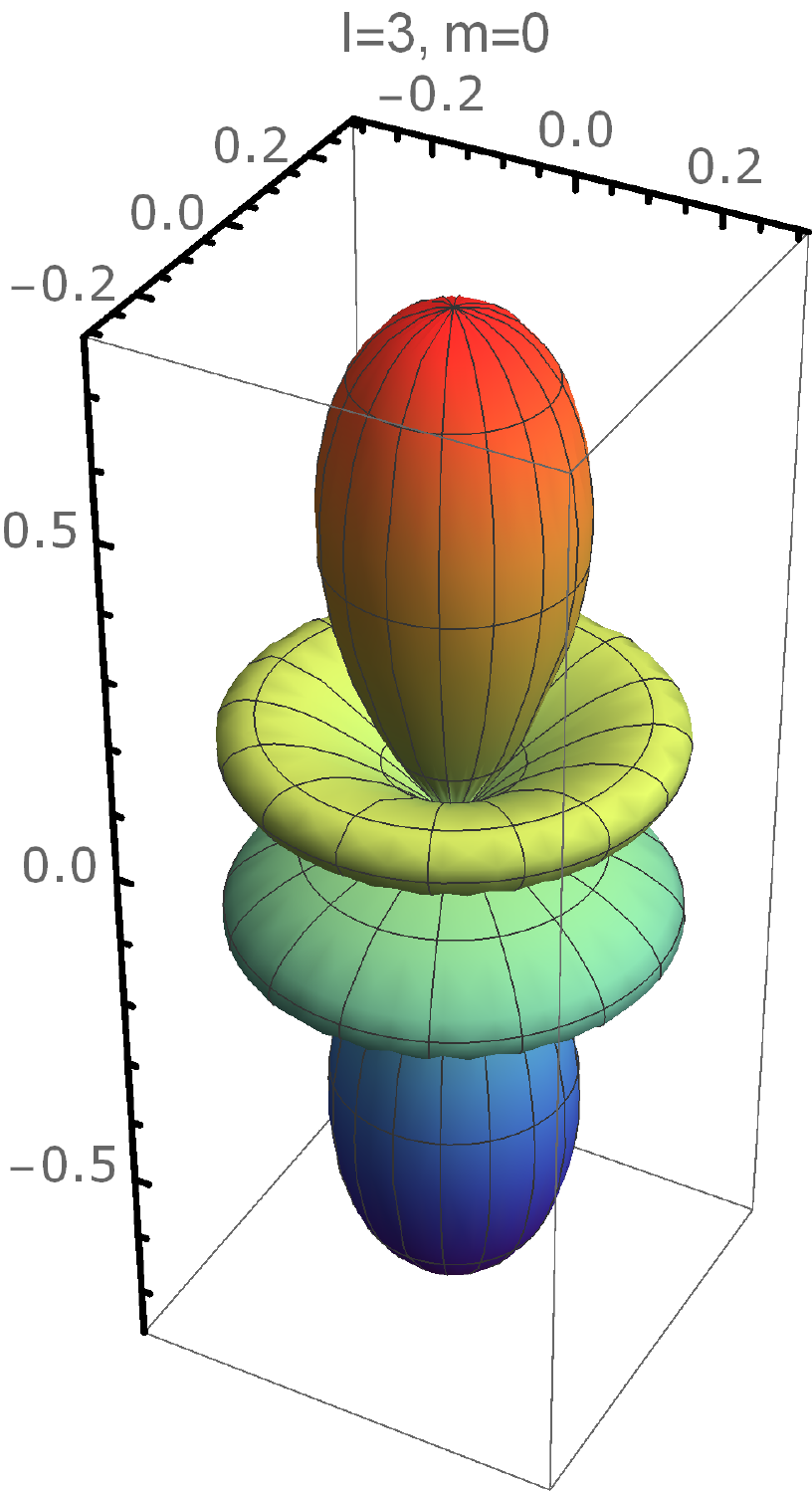}
\includegraphics[scale=0.34]{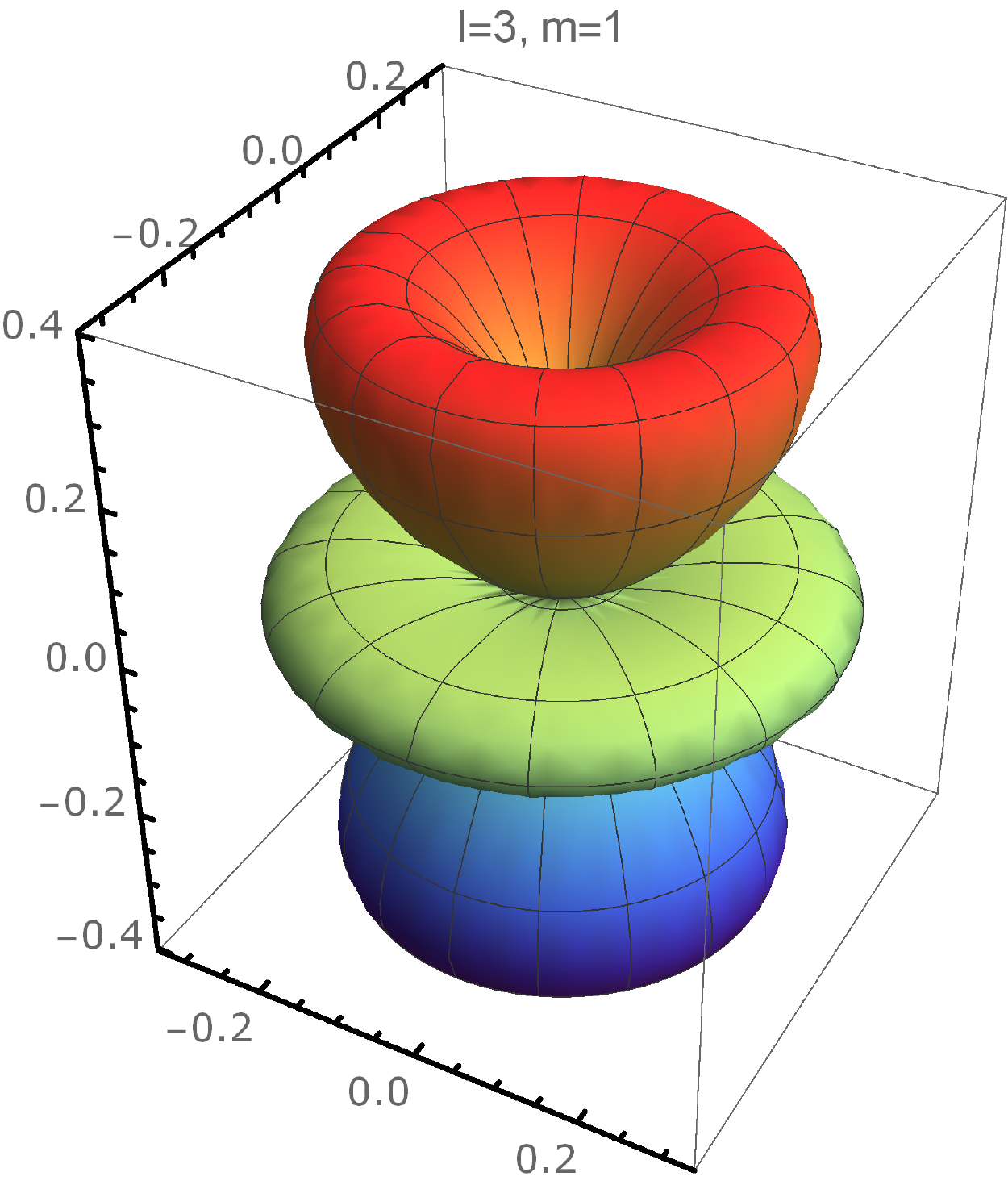}
\includegraphics[scale=0.34]{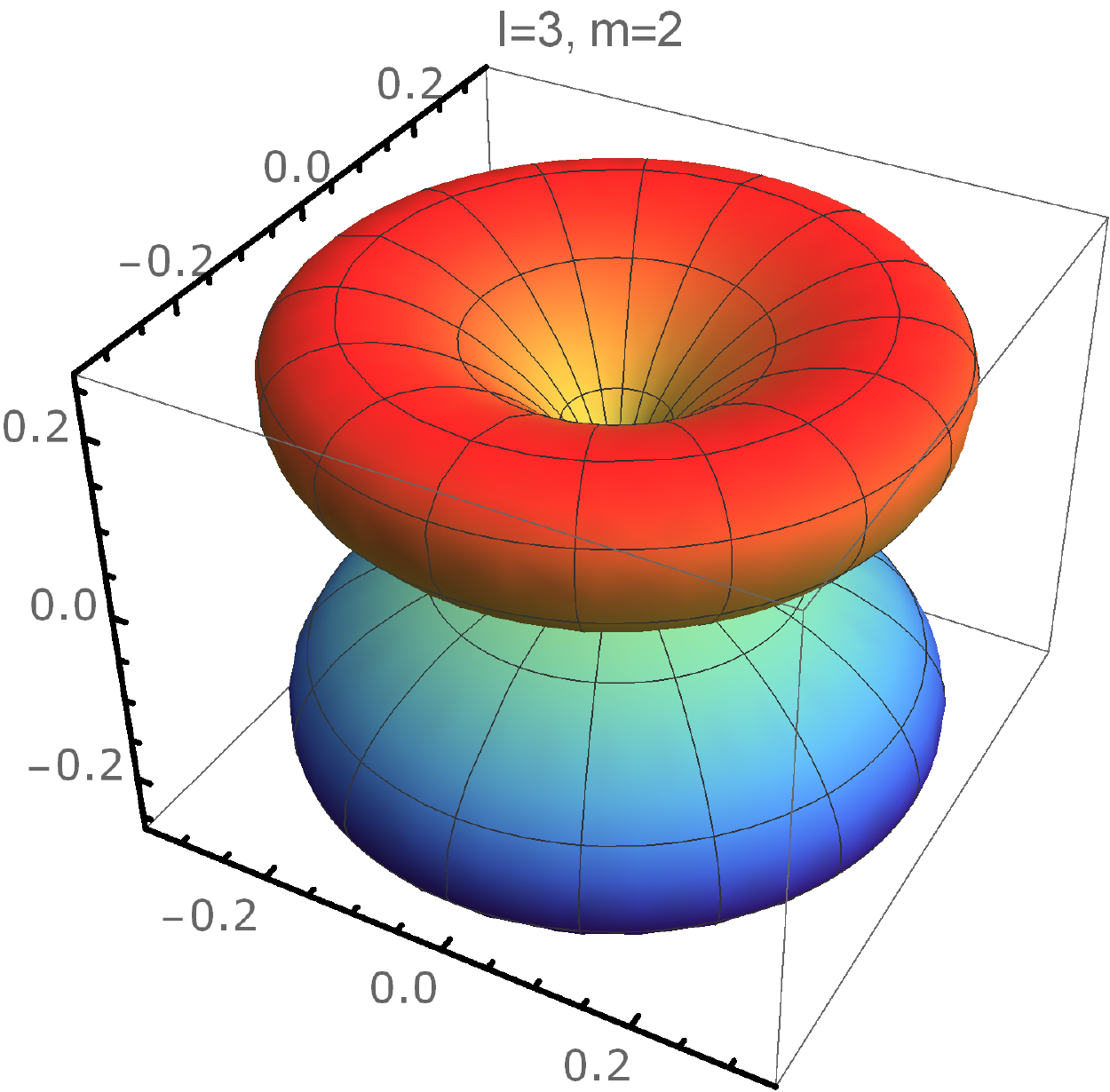}
\includegraphics[scale=0.34]{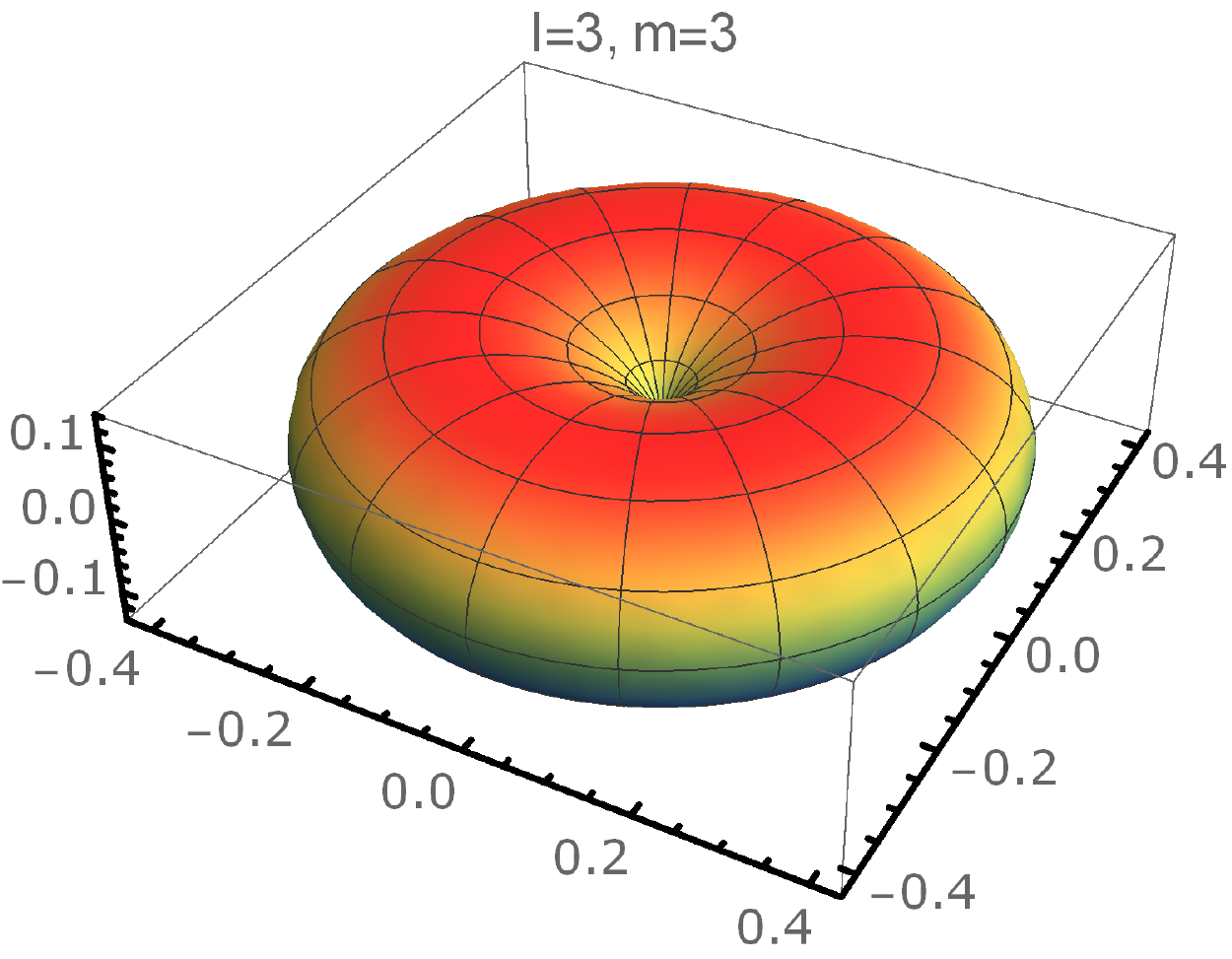}
\caption{\label{fig:harmonique_spherique} Harmoniques sph\'eriques $Y_l^{m}(\theta,\phi)$ pour $l=0,1,2$ et $0\leq m\leq l$.}
\end{figure}

On appelle harmoniques spheriques, les solutions de l'\'equation de type
\begin{align}\label{eq:5_21}
\frac{1}{\sin\theta}\left(\frac{\partial}{\partial\theta}\sin\theta\frac{\partial\Psi}{\partial\theta}\right)+\frac{1}{\sin^{2}\theta}\frac{\partial^{2}\Psi}{\partial\phi^{2}}+l(l+1)\Psi=0
\end{align}
Les harmoniques sph\'eriques sont particuli\`erement utiles pour les probl\'emes qui ont une sym\'etrie de rotation, car elles sont les vecteurs propres de certains op\'erateurs li\'es aux rotations. En physique, elles sont souvent utilis\'ees en acoustique, en cristallographie, en mechanique quantique, en cosmologie, etc. Cette \'equation peut \^etre resolue par la m\'ethode de s\'eparation des variables, \i.e., on cherche les solutions sous la forme $\Psi(\theta,\phi)=\Theta(\theta)\Phi(\phi)$. En remplaçant dans l'Eq.~(\ref{eq:5_21}), cela nous donne  
\begin{align*}
\frac{\Phi}{\sin\theta}\left\lbrace\frac{d}{d\theta}\left(\sin\theta\frac{d\Theta}{d\theta}\right)\right\rbrace+\frac{\Theta}{\sin^{2}\theta}\frac{d^{2}\Phi}{d\phi^{2}}+l(l+1)\Theta\Phi=0
\end{align*}
En devisant par $\Theta(\theta)\Phi(\phi)$ et multipliant  par $\sin^{2}\theta$, on trouve (apr\'es r\'earrangement)
\begin{align*}
\frac{\sin\theta}{\Theta}\frac{d}{d\theta}\left(\sin\theta\frac{d\Theta}{d\theta}\right)+l(l+1)\sin^{2}\theta=-\frac{1}{\Phi}\frac{d^{2}\Phi}{d\phi^{2}}
\end{align*}
On remarque que le c\^ot\'e gauche contient toute et seulement la d\'ependance en $\theta$, 
tandis que le c\^ot\'e droit contient toute et seulement la d\'ependance en $\phi$. Le deux c\^ot\'es doivent donc \^etre \'egaux \`a une constante. Nous choisissons de noter cette constante $m^2$, on a
\begin{align}
&\frac{\sin\theta}{\Theta}\frac{d}{d\theta}\left(\sin\theta\frac{d\Theta}{d\theta}\right)+l(l+1)\sin^{2}\theta=m^{2}\label{eq:5_22}\\
&-\frac{1}{\Phi}\frac{d^{2}\Phi}{d\phi^{2}} = m^{2}\label{eq:5_23}
\end{align}
La solution de l'Eq.~(\ref{eq:5_23}) est donn\'ee par
\begin{align*}
\Phi(\phi)=Ae^{im\phi}+Be^{-im\phi}
\end{align*}
Pour que la solution soit continue, il faut que $\Phi(2\pi)=\Phi(0)$, cela impose que $m$ soit un entier positif, n\'egatif ou nul. 
Eq.~(\ref{eq:5_22}) peut \^etre \'ecrite sous la forme
\begin{align*}
\frac{1}{\sin\theta}\frac{d}{d\theta}\left(\sin\theta\frac{d\Theta}{d\theta}\right)+\left\lbrace l(l+1)-\frac{m^{2}}{\sin^{2}\theta}\right\rbrace\Theta=0
\end{align*}
En faisant le changement de variable $\cos\theta=x$, on trouve
\begin{align*}
\frac{d}{dx}\left\lbrace(1-x^{2})\frac{d\Theta}{dx}\right\rbrace+\left\lbrace l(l+1)-\frac{m^{2}}{1-x^{2}}\right\rbrace\Theta=0
\end{align*}
qui est l'\'equation associ\'ee de Legendre. Il n'y aura de solutions finies aux points $\theta=0$ et $\pi$ ($x=+1$ et $x=-1$) que si $l$ est un entier. Donc cette \'equation a pour solutions les
polyn\^omes de Legendre associ\'ees d\'efinis par
\begin{align*}
\Theta(\theta)=P_l^{m}(x)=P_l^{m}(\cos\theta)
\end{align*}
avec $-l<m<l$. Donc, la solution g\'en\'erale de l'Eq.~(\ref{eq:5_21}) est donn\'ee par
\begin{align*}
\Psi(\theta,\phi)=\left(Ae^{im\phi}+Be^{-im\phi}\right)P_l^{m}(\cos\theta)
\end{align*}
En utilisant la formule (\ref{eq:5_19}) pour $P_{l}^{-m}(x)$, on aura
\begin{align*}
\Psi(\theta,\phi)=A_1e^{im\phi}P_l^{m}(\cos\theta)+A_2e^{-im\phi}P_l^{-m}(\cos\theta)
\end{align*}
o\`u 
\begin{align*}
&A_1=A\\
&A_2=(-1)^{m}\frac{(l+m)!}{(l-m)!}B
\end{align*}
Posons
\begin{align*}
y_l^{m}(\theta,\phi)=e^{im\phi}P_l^{m}(\cos\theta)
\end{align*}
$\Psi$ devient
\begin{align*}
\Psi(\theta,\phi)=A_1y_l^{m}(\theta,\phi)+A_2y_l^{-m}(\theta,\phi)
\end{align*}
qui est la solution de l'Eq.~(\ref{eq:5_21}) pour toute valeur de $m$. Puisque l'Eq.~(\ref{eq:5_21}) est homog\`ene, la solution g\'en\'erale est donn\'ee par
\begin{align}
\Psi(\theta,\phi)=\sum_{m=0}^{l}\left\lbrace A_1^{(m)}y_l^{m}(\theta,\phi)+A_2^{(m)}y_l^{-m}(\theta,\phi)\right\rbrace
\end{align}
Dans de nombreuses applications en physique et en ing\'enierie, il est plus pratique de d\'efinir la solution de l'Eq.~(\ref{eq:5_21}) comme un multiple des $y_l^{m}$ qui s'appellent les harmoniques sph\'eriques not\'ee $Y_l^{m}$ et d\'efinis par
\begin{align}
Y_l^{m}(\theta,\phi)&=(-1)^{m}\frac{1}{\sqrt{2\pi}}\sqrt{\frac{(2l+1)(l-m)!}{2(l+m)!}}y_l^{m}(\theta,\phi)\nonumber\\
&=(-1)^{m}\frac{1}{\sqrt{2\pi}}\sqrt{\frac{(2l+1)(l-m)!}{2(l+m)!}}e^{im\phi}P_l^{m}(\cos\theta)\label{eq:5_23a}
\end{align}
Les expressions explicites pour $Y_0^{0},Y_0^{0},Y_{1}^{\pm 1},Y_2^{0},Y_2^{\pm 1},Y_{2}^{\pm 2}$ sont donn\'ees ci-dessous
\begin{align*}
&Y_0^{0}=\frac{1}{2\sqrt{\pi}}\\
&Y_1^{0}=\sqrt{\frac{3}{4\pi}}\cos\theta, \hspace{20mm}Y_1^{\pm 1}=\sqrt{\frac{3}{8\pi}}\sin\theta e^{\pm i\phi}\\
&Y_2^{0}=\sqrt{\frac{15}{16\pi}}\left(3\cos^{2}\theta -1\right),\hspace{10mm}Y_2^{\pm 1}=\sqrt{\frac{15}{8\pi}}\cos\theta\sin\theta e^{\pm i\phi}\\
&Y_2^{\pm 2}=\sqrt{\frac{15}{32\pi}}\sin^{2}\theta e^{\pm 2i\phi}
\end{align*}
Elles sont pr\'esent\'ees dans la figure Fig.~(\ref{fig:harmonique_spherique}).
\subsection{Relation d'orthogonalit\'e}
\begin{align}\label{eq:5_24}
\int_{0}^{2\pi}d\phi\int_{0}^{\pi}d\theta\sin\theta\left\lbrace Y_l^{m}(\theta,\phi)\right\rbrace^{*}Y_{l'}^{m'}(\theta,\phi)=\delta_{ll'}\delta_{mm'}
\end{align}
\textit{D\'emonstration}\\
On a
\begin{align*}
\int_{0}^{2\pi}d\phi\int_{0}^{\pi}d\theta\sin\theta &\left\lbrace Y_l^{m}(\theta,\phi)\right\rbrace^{*}Y_{l'}^{m'}(\theta,\phi)\\
&=(-1)^{m+m'}\frac{1}{2\pi}\sqrt{\frac{(2l+1)(2l'+1)(l-m)!(l'-m')!}{4(l+m)!(l'+m')!}}\\
&\times\int_{0}^{2\pi}e^{i(m-m')\phi}d\phi\int_{0}^{\pi}P_l^{m}(\cos\theta)P_{l'}^{m'}(\cos\theta)\sin\theta d\theta\\
&=(-1)^{m+m'}\frac{1}{2\pi}\sqrt{\frac{(2l+1)(2l'+1)(l-m)!(l'-m')!}{4(l+m)!(l'+m')!}}\\
&\times 2\pi\delta_{mm'}\int_{-1}^{1}P_l^{m}(x)P_{l'}^{m'}(x)dx\\
\end{align*}
dans la derni\`ere ligne, nous avons fait le changement de variable $x=\cos\theta$. 
Si on utilise l'Eq.~(\ref{eq:5_20}), on obtient
\begin{align*}
\int_{0}^{2\pi}d\phi\int_{0}^{\pi}d\theta &\sin\theta\left\lbrace Y_l^{m}(\theta,\phi)\right\rbrace^{*}Y_{l'}^{m'}(\theta,\phi)\\
&=(-1)^{2m}\sqrt{\frac{(2l+1)(2l'+1)(l-m)!(l'-m)!}{4(l+m)!(l'+m)!}}
\delta_{mm'}\int_{-1}^{1}P_l^{m}(x)P_{l'}^{m}(x)dx\\
&=\sqrt{\frac{(2l+1)(2l'+1)(l-m)!(l'-m)!}{4(l+m)!(l'+m)!}}\delta_{mm'}\frac{2(l+m)!}{(2l+1)(l-m)!}\delta_{ll'}\\
&=\delta_{mm'}\delta_{ll'}
\end{align*}
\subsection{Properi\'et\'es}
\begin{align}\label{eq:5_25}
\left\lbrace Y_l^{m}(\theta,\phi)\right\rbrace^{*}=(-1)^{m}Y_l^{-m}(\theta,\phi)
\end{align}
\textit{D\'emonstration}\\
En utilisant l'expression (\ref{eq:5_23a}), on a
\begin{align*}
\left\lbrace Y_l^{m}(\theta,\phi)\right\rbrace^{*}=(-1)^{m}\frac{1}{\sqrt{2\pi}}\sqrt{\frac{(2l+1)(l-m)!}{2(l+m)!}}e^{-im\phi}P_l^{m}(\cos\theta)
\end{align*}
En remplaçant $P_l^{m}(x)$ par $P_l^{-m}(x)$ (Eq.~(\ref{eq:5_19}), on obtient
\begin{align*}
\left\lbrace Y_l^{m}(\theta,\phi)\right\rbrace^{*}&=(-1)^{m}\frac{1}{\sqrt{2\pi}}\sqrt{\frac{(2l+1)(l-m)!}{2(l+m)!}}e^{-im\phi}(-1)^{m}\frac{(l+m)!}{(l-m)!}P_l^{-m}(\cos\theta)\\
&=(-1)^{m}(-1)^{m}\frac{1}{\sqrt{2\pi}}\sqrt{\frac{(2l+1)(l+m)!}{2(l-m)!}}e^{-im\phi}P_l^{-m}(\cos\theta)\\
&=(-1)^{m}Y_l^{-m}(\theta,\phi)
\end{align*}
\section{Polyn\^omes d'Hermite}
\begin{figure}[!ht]
\centering
\includegraphics[scale=0.7]{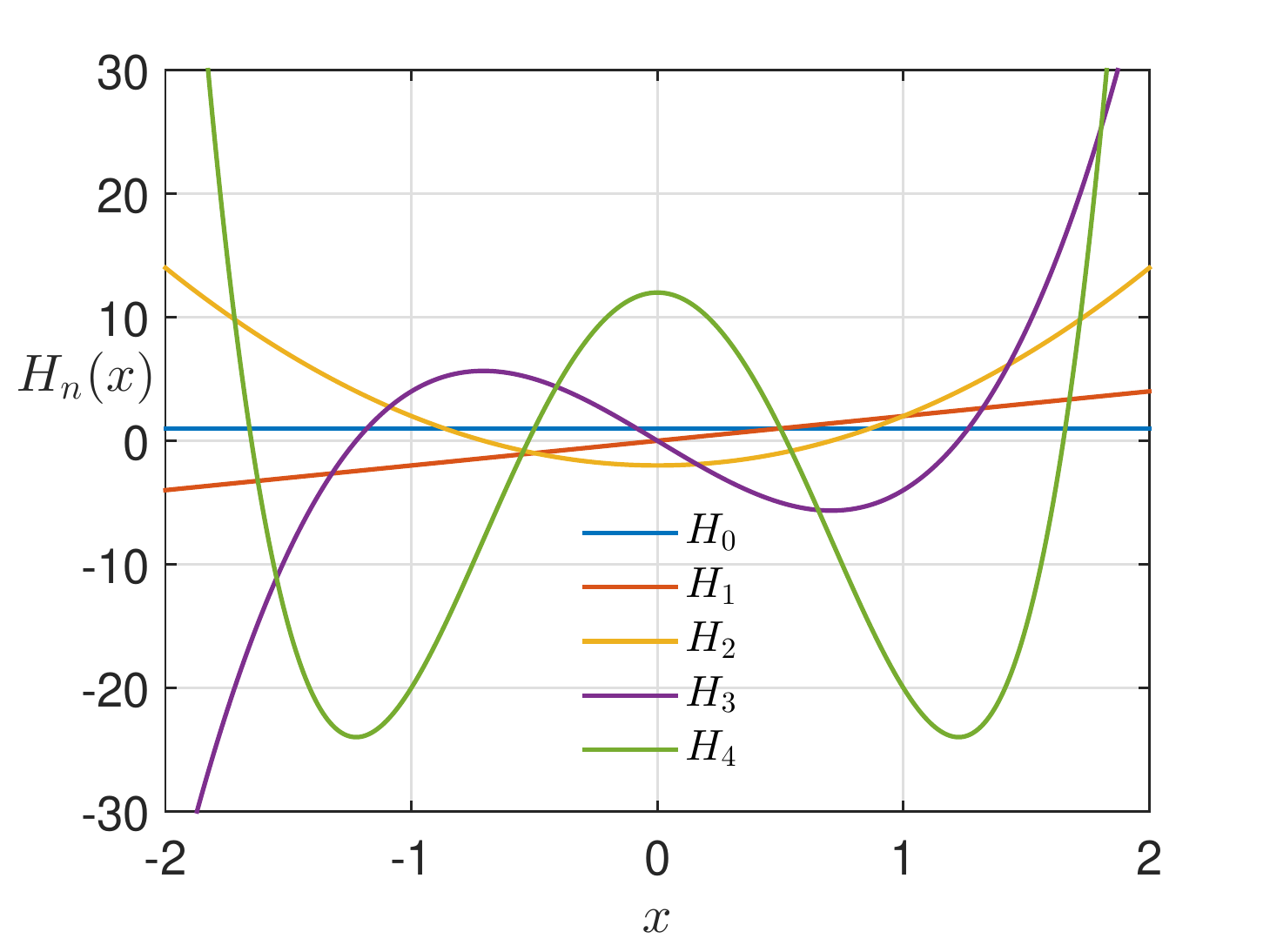}
\caption{\label{fig:hermite} Polyn\^omes d'Hermite pour $n=0,1,2,3,4$.}
\end{figure}

Les polyn\^omes d'Hermite ont \'et\'e d\'efinis par Laplace en $1810$ bien que sous une forme difficilement reconnaissable, et ont \'et\'e \'etudi\'es en d\'etail par Chebyshev en $1859$. Les travaux de Chebyshev ont \'et\'e oubli\'es et ils ont \'et\'e nomm\'es plus tard en l'honneur de Charles Hermite qui \'ecrivait sur les polyn\^mes en $1864$ et les d\'ecrivait comme nouveaux. Par cons\'equent, ils n'\'etaient pas nouveaux, bien que plus tard dans des papiers de $1865$, Hermite fut le premier \`a d\'efinir les polyn\^omes multidimensionnels. Les polyn\^omes d'Hermite sont les solution de l'\'equation d'Hermite qui est donn\'ee par
\begin{align}\label{eq:5_26}
\frac{d^{2}y}{dx^{2}}-2x\frac{dy}{dx}+2ny=0
\end{align}
En utilisant la m\'ethode de Frobenius d\'ecrite en d\'etail au Chap.~(\ref{chap: bessel}), on trouve comme solutions
\begin{align}\label{eq:5_27}
H_n(x)=\sum_{r=0}^{[\frac{n}{2}]}(-1)^{r}\frac{n!}{r!(n-2r)!}(2x)^{n-2r}
\end{align}
qui sont les polyn\^omes d'Hermite.\\
Quelques expressions explicites pour $H_n(x)$ pour $n=0,1,2,3,4$ sont donn\'ees ci-dessous
\begin{align*}
&H_0(x)=1\\
&H_1(x)=2x\\
&H_2(x)=4x^{2}-2\\
&H_3(x)=8x^{3}-12x\\
&H_4(x)=16x^{4}-48x^{2}+12
\end{align*}
et elles sont repr\'esent\'ees sur la Fig.~(\ref{fig:hermite}).
\subsection{Fonction g\'en\'eratrice}
\begin{align}\label{eq:5_28}
e^{2tx-t^{2}}=\sum_{n=0}^{\infty}\frac{t^{n}}{n!}H_n(x)
\end{align}
\textit{D\'emonstration}
\begin{align*}
e^{2tx-t^{2}}&=e^{2tx}e^{-t^{2}}\\
&=\sum_{r=0}^{\infty}\frac{(2tx)^{r}}{r!}\sum_{s=0}^{\infty}\frac{(-t^{2})^{s}}{s!}\\
&=\sum_{r,s=0}^{\infty}(-1)^{s}\frac{2^{r}x^{r}}{r!s!}t^{r+2s}\\
\end{align*}
On prend $n=r+2s$, \i.e. $r=n-2s$, et puisque $r=n-2s\geq 0$ alors $s\leq n/2$ et par suite
\begin{align*}
y&=\sum_{n=0}^{\infty}\sum_{s=0}^{[\frac{n}{2}]}(-1)^{s}\frac{2^{n-2s}x^{n-2s}}{(n-2s)!s!}t^{n}\\
&=\sum_{n=0}^{\infty}\frac{t^{n}}{n!}H_n(x)
\end{align*}
\subsection{Relation d'orthogonalit\'e}
\begin{align}\label{eq:5_29}
\int_{-\infty}^{\infty}e^{-x^{2}}H_n(x)H_m(x)dx=2^{n}n!\sqrt{\pi}\delta_{nm}
\end{align}
\textit{D\'emonstration}\\
On a
\begin{align*}
&e^{-t^{2}+2tx}=\sum_{n=0}^{\infty}\frac{t^{n}}{n!}H_n(x)\\
&e^{-s^{2}+2sx}=\sum_{m=0}^{\infty}\frac{s^{m}}{m!}H_m(x)\\
\end{align*}
donc
\begin{align}\label{eq:5_30}
\sum_{n,m=0}^{\infty}\frac{t^{n+m}}{n!m!}\int_{-\infty}^{\infty}e^{-x^{2}}H_n(x)H_m(x)dx=\int_{-\infty}^{\infty}e^{-x^{2}}e^{-t^{2}+2tx}e^{-s^{2}+2sx}dx
\end{align}
mais
\begin{align}
\int_{-\infty}^{\infty}e^{-x^{2}}e^{-t^{2}+2tx}e^{-s^{2}+2sx}dx &=e^{-t^{2}-s^{2}}\int_{-\infty}^{\infty}e^{-x^{2}+2(s+t)x}dx\nonumber\\
&=e^{-t^{2}-s^{2}}\int_{-\infty}^{\infty}e^{-(x-(s+t))^{2}+(s+t)^{2}}dx\nonumber\\
&=e^{2st}\int_{-\infty}^{\infty}e^{-(x-(s+t))^{2}}dx\nonumber
\end{align}
En faisant le changement de variable $u=x-(s+t)$, on obtient
\begin{align*}
\int_{-\infty}^{\infty}e^{-x^{2}}e^{-t^{2}+2tx}e^{-s^{2}+2sx}dx&=e^{2st}\int_{-\infty}^{\infty}e^{-u^{2}}dx\\
&=e^{2st}\sqrt{\pi}\\
&=\sqrt{\pi}\sum_{n=0}^{\infty}\frac{2^{n}s^{n}t^{n}}{n!}
\end{align*}
En identifiant les coefficients de $t^{n}$ du c\^ot\'e gauche de l'expression (\ref{eq:5_30}) avec ceux du c\^ot\'e droit de  la derni\`ere expression, on d\'eduit que
\begin{align}
\int_{-\infty}^{\infty}e^{-x^{2}}H_n(x)H_m(x)dx=\left\lbrace\begin{array}{ll}
0 & \mbox{si $n\neq m$}\\
\sqrt{\pi}2^{n}n! & \mbox{si $n= m$}
\end{array}\right.
\end{align}
\subsection{Relations de r\'ecurrece}
\begin{align}
&(a)\hspace{1.5mm}H_n'(x)=2nH_{n-1}(x) \hspace{5mm}(n\geq 1); \hspace{5mm}H_0'(x)=0,\label{eq:5_31}\\
&(b)\hspace{1.5mm}H_{n+1}(x)=2xH_n(x)-2nH_{n-1}(x)\hspace{3mm}(n\geq 1);\hspace{3mm}H_1(x)=2xH_0(x)\nonumber
\end{align}
\textit{D\'emonstration}\\
(a) Si on d\'erive les deux c\^ot\'es de l'expression de la fonction g\'en\'eratrice (\ref{eq:5_28}) par rapport \`a $x$, on obtient
\begin{align*}
\sum_{n=0}^{\infty}\frac{t^{n}}{n!}H_n'(x)&=2te^{2tx-t^{2}}\\
&=2t\sum_{n=0}^{\infty}\frac{t^{n}}{n!}H_n(x)\\
&=2\sum_{n=0}^{\infty}\frac{t^{n+1}}{n!}H_n'(x)\\
&=2\sum_{n=1}^{\infty}\frac{t^{n}}{(n-1)!}H_{n-1}'(x)
\end{align*}
En comparant les coefficients de $t^{n}$ des deux c\^ot\'es, on trouve pour $n=0$
\begin{align*}
H_0'(x)=0
\end{align*}
et pour $n\geq 1$
\begin{align*}
\frac{H_n'(x)}{n!}=\frac{2H_{n-1}(x)}{(n-1)!}
\end{align*}
qui peut \^etre simplifier \`a
\begin{align*}
H_n'(x)=2nH_{n-1}(x)
\end{align*}
(b) D\'erivons les deux c\^ot\'es de l'Eq.~(\ref{eq:5_28}) par rapport \`a $t$, on obtient
\begin{align*}
(2x-2t)e^{2tx-t^{2}}=\sum_{n=0}^{\infty}\frac{nt^{n-1}}{n!}H_n(x)
\end{align*}
pour $n=0$, le terme de c\^ot\'e droite est nul, donc la somme commence pour $n=1$,
\begin{align*}
(2x-2t)\sum_{n=0}^{\infty}\frac{t^{n}}{n!}H_n(x)=\sum_{n=1}^{\infty}\frac{nt^{n-1}}{n!}H_n(x)
\end{align*}
En d\'eveloppant le terme du c\^ot\'e gauche, on trouve
\begin{align*}
2x\sum_{n=0}^{\infty}\frac{t^{n}}{n!}H_n(x)-2\sum_{n=0}^{\infty}\frac{t^{n+1}}{n!}H_n(x)=\sum_{n=1}^{\infty}\frac{nt^{n-1}}{n!}H_n(x)
\end{align*}
qui peut \^etre simplifier \`a
\begin{align*}
2x\sum_{n=0}^{\infty}\frac{t^{n}}{n!}H_n(x)-2\sum_{n=1}^{\infty}\frac{t^{n}}{(n-1)!}H_n(x)=\sum_{n=1}^{\infty}\frac{t^{n}}{n!}H_{n+1}(x)
\end{align*}
En comparant les coefficients de $t^{n}$ pour $n\geq 1$, on obtient
\begin{align*}
\frac{2x}{n!}H_n(x)-\frac{2}{(n-1)!}H_{n-1}(x)=\frac{1}{n!}H_{n+1}(x)
\end{align*}
En multipliant par $n!$, on trouve
\begin{align*}
2xH_n(x)-2nH_{n-1}(x)=H_{n+1}(x)
\end{align*}
et pour $n=0$ on obtient
\begin{align*}
2xH_{0}(x)=H_{1}(x)
\end{align*}
\section{Polyn\^ome de Laguerre}
\begin{figure}[!ht]
\centering
\includegraphics[scale=0.7]{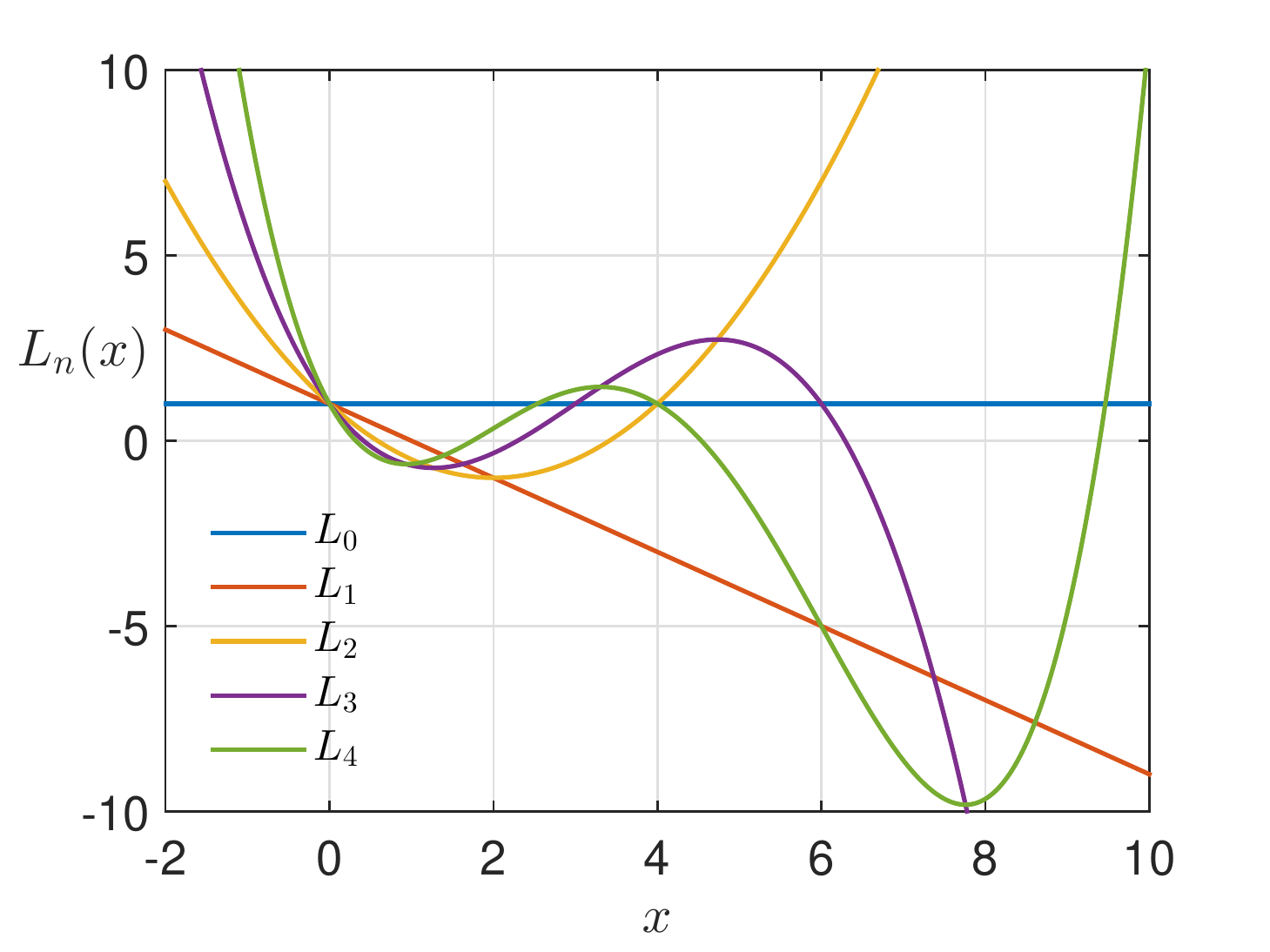}
\caption{\label{fig:laguerre} Polyn\^omes de Laguerre pour $n=0,1,2,3,4$.}
\end{figure}

L'\'equation de Laguerre est donn\'ee par
\begin{align}\label{eq:5_33}
x\frac{d^{2}y}{dx^{2}}+(1-x)\frac{dy}{dx}+ny=0
\end{align}
On applique la m\^eme m\'ethode de Frobenius du Chap.~(\ref{chap: bessel}), c'est \`a dire, on cherche une solution sous la forme
\begin{align*}
y(x,s)=\sum_{r=0}^{\infty}a_rx^{r+s}
\end{align*}
on obtient les \'equations indiciales suivantes
\begin{align*}
&a_0s^2=0\\
&a_1(s+1)^2=0\\
&a_r(s+r-n) - a_{r+1}(s+r+1)^2=0,\hspace{10mm} r\geq 2
\end{align*}
Les deux premi\`eres \'equations admettent une double racine $s=0$. La derni\`ere \'equation donne une relation de r\'ecurrence pour les coefficients $a_r$ pour $r\geq 2$
\begin{align*}
a_{r+1}=a_r\frac{s+r-n}{(s+r+1)^{2}}
\end{align*}
En remplaçant $s=0$ on obtient
\begin{align*}
a_{r+1}=a_r\frac{r-n}{(r+1)^{2}}
\end{align*}
avec cette forme de coefficients $a_r$, on remarque que la s\'erie infinie $y(x,0)$ qui en r\'esulte se comporte comme un exponentielle pour $x$ tr\'es grand, donc la s\'erie diverge. Pour \'eviter ce probleme, la s\'erie devrait \^etre tronquer: on prend $a_r=0$ si $r>n$ \.i.e., pour $r=n+1,n+2,\ldots$. Ceci est possible si $n$ est un entier positif. On r\'e\'ecrit $a_r$  sous la forme
\begin{align*}
a_{r+1}=-a_r\frac{n-r}{(r+1)^{2}}
\end{align*}
la solution devient
\begin{align*}
y&=a_0\left\lbrace 1-\frac{n}{1^{2}}x+\frac{n(n-1)}{(2!)^{2}}x^{2}+\ldots
+(-1)^{r}\frac{n(n-1)\ldots(n-r+1)}{(r!)^{2}}x^{r}+\ldots\right\rbrace\\
&=a_0\sum_{r=0}^{n}(-1)^{r}\frac{n(n-1)\ldots(n-r+1)}{(r!)^{2}}x^{r}\\
&=a_0\sum_{r=0}^{n}(-1)^{r}\frac{n!}{(n-r)!(r!)^{2}}x^{r}
\end{align*}
En choisissant $a_0=1$, on obtient la solution finale de l'Eq.~(\ref{eq:5_22}) not\'ee $L_n(x)$ et qui s'appelle le polyn\^ome de Laguerre
\begin{align}\label{eq:5_34}
L_n(x)=\sum_{r=0}^{n}(-1)^{r}\frac{n!}{(n-r)!(r!)^{2}}x^{r}
\end{align}
Quelques expressions explicites pour $L_n(x)$ pour $n=0,1,2,3,4$ sont donn\'ees ci-dessous
\begin{align*}
&L_0(x)=1\\
&L_1(x)=-x+1\\
&L_2(x)=\frac{1}{2}\left(x^{2}-4x+2\right)\\
&L_3(x)=\frac{1}{3!}\left(-x^{3}+9x^{2}-18x+6\right)\\
&L_4(x)=\frac{1}{4!}\left(x^{4}-16x^{3}+72x^{2}-96x+24\right)
\end{align*}
et elles sont repr\'esent\'ees dans la Fig.~(\ref{fig:laguerre}).
\subsection{Fonction g\'en\'eratrice}
\begin{align}\label{eq:5_35}
\frac{1}{1-t}\exp\left\lbrace -xt/(1-t)\right\rbrace=\sum_{n=0}^{\infty}t^{n}L_n(x)
\end{align}
\textit{D\'emonstration}\\
On a
\begin{align*}
\frac{1}{1-t}\exp\left\lbrace -xt/(1-t)\right\rbrace &=\frac{1}{1-t}\sum_{r=0}^{\infty}\frac{1}{r!}\left(-\frac{xt}{1-t}\right)^{r}\\
&=\sum_{r=0}^{\infty}\frac{(-1)^{r}}{r!}\frac{x^{r}t^{r}}{(1-t)^{r+1}}
\end{align*}
Le d\'eveloppement en s\'erie de $1/(1-t)^{r+1}$ donne
\begin{align*}
\frac{1}{(1-t)^{r+1}}=\sum_{s=0}^{\infty}\frac{(r+s)!}{r!s!}t^{s}
\end{align*}
donc
\begin{align*}
\frac{1}{1-t}\exp\left\lbrace -xt/(1-t)\right\rbrace &=\sum_{r,s=0}^{\infty}(-1)^{r}\frac{(r+s)!}{(r!)^{2}s!}x^{r}t^{r+s}
\end{align*}
Posons $n=r+s$, on obtient
\begin{align*}
\frac{1}{1-t}\exp\left\lbrace -xt/(1-t)\right\rbrace &=\sum_{n=0}^{\infty}t^{n}\sum_{r=0}^{\infty}(-1)^{r}\frac{n!}{(r!)^{2}(n-r)!}x^{r}\\
&=\sum_{n=0}^{\infty}t^{n}L_n(x)
\end{align*}
\subsection{Relation d'orthogonalit\'e}
\begin{align}\label{eq:5_36}
\int_0^{\infty}e^{-x}L_n(x)L_m(x)dx=\delta_{nm}
\end{align}
\textit{D\'emonstration}:\\
A partir de la fonction g\'en\'eratrice, on a
\begin{align}
&\frac{1}{1-t}\exp\left\lbrace -xt/(1-t)\right\rbrace =\sum_{n=0}^{\infty}t^{n}L_n(x)\nonumber\\
&\frac{1}{1-s}\exp\left\lbrace -xs/(1-s)\right\rbrace =\sum_{m=0}^{\infty}s^{m}L_m(x)\nonumber
\end{align}
donc
\begin{align}
\sum_{n,m=0}^{\infty}t^{n}s^{m}\int_0^{\infty}e^{-x}L_n(x)L_m(x)dx &=\frac{1}{1-t}\frac{1}{1-s}\int_0^{\infty}e^{-x}e^{-xt/(1-t)}e^{-xs/(1-s)}dx\nonumber\\
&=\frac{1}{1-t}\frac{1}{1-s}\left[-\frac{1}{1+\left[t/(1-t)\right]+\left[s/(1-s)\right]}\right.\nonumber\\
&\left.\times\exp\left\lbrace -x\left(1+\frac{t}{1-t}+\frac{s}{1-s}\right) \right\rbrace \right]_0^{\infty}\nonumber\\
&=\frac{1}{1-st}\nonumber\\
&=\sum_{n=0}^{\infty}s^{n}t^{n}\nonumber
\end{align}
En comparant les coefficients de $t^{n+m}$ des deux c\^ot\'es, on trouve que
\begin{align}
\int_0^{\infty}e^{-x}L_n(x)L_m(x)dx=\delta_{nm}\nonumber
\end{align}
\subsection{Relations de r\'ecurrence}
\begin{align}\label{eq:5_37}
&(a)\hspace{1.5mm} (n+1)L_{n+1}(x)=(2n+1-x)L_n(x)-nL_{n-1}(x)\nonumber\\
&(b)\hspace{1.5mm}xL_n'(x)=nL_n(x)-nL_{n-1}(x)\\
&(c)\hspace{1.5mm}L_n'(x)=-\sum_{r=0}^{n-1}L_r(x)\nonumber
\end{align}
\textit{D\'emonstration}\\
En d\'erivant les deux c\^ot\'es de l'expression (\ref{eq:5_35}) par rapport a $t$, on obtient
\begin{align*}
\sum_{n=0}^{\infty}nt^{n-1}L_n(x)=\frac{1}{(1-t)^{2}}e^{-xt/(1-t)}-\frac{x}{(1-t)^{2}}\frac{1}{1-t}e^{-xt/(1-t)}
\end{align*}
En remplaçant les deux termes du c\^ot\'e droit par la fonction g\'en\'eratrice (\ref{eq:5_35}), on obtient
\begin{align*}
\sum_{n=0}^{\infty}nt^{n-1}L_n(x)=\frac{1}{(1-t)}\sum_{n=0}^{\infty}t^{n}L_n(x)-\frac{x}{(1-t)^{2}}\sum_{n=0}^{\infty}t^{n}L_n(x)
\end{align*}
En multipliant par $(1-t)^2$, il en r\'esulte que
\begin{align*}
(1-t)^{2}\sum_{n=0}^{\infty}nt^{n-1}L_n(x)=(1-t)\sum_{n=0}^{\infty}t^{n}L_n(x)-x\sum_{n=0}^{\infty}t^{n}L_n(x)
\end{align*}
qui peut \^etre simplifier \'a
\begin{align*}
\sum_{n=0}^{\infty}nt^{n-1}L_n(x)-2\sum_{n=0}^{\infty}&nt^{n}L_n(x) + \sum_{n=0}^{\infty}nt^{n+1}L_n(x)=\nonumber\\
&\sum_{n=0}^{\infty}t^{n}L_n(x)-\sum_{n=0}^{\infty}t^{n+1}L_n(x)-x\sum_{n=0}^{\infty}t^{n}L_n(x)
\end{align*}
En faisant glisser les indices de $n-1$ et $n+1$ \`a $n$ pour que $t$ a une puissance $n$ dans tous les  termes, on obtient
\begin{align}
\sum_{n=-1}^{\infty}(n+1)t^{n}L_{n+1}(x)-2\sum_{n=0}^{\infty}&nt^{n}L_n(x) + \sum_{n=1}^{\infty}(n-1)t^{n}L_{n-1}(x)=\nonumber\\
&\sum_{n=0}^{\infty}t^{n}L_n(x)-\sum_{n=1}^{\infty}t^{n}L_{n-1}(x)-x\sum_{n=0}^{\infty}t^{n}L_n(x)\nonumber
\end{align}
En identifiant les coefficients de $t^{n}$, on obtient la relation d\'esir\'e
\begin{align}
(n+1)L_{n+1}(x)=(2n+1-x)L_n(x)-nL_{n-1}(x)\nonumber
\end{align}
(b) D\'erivons maintenant l'\'equation de la fonction g\'en\'eratrice (\ref{eq:5_35}) par rapport \`a $x$, nous trouvons
\begin{align}
\sum_{n=0}^{\infty}t^{n}L_n'(x)&=-\frac{t}{1-t}\frac{1}{1-t}e^{-xt/(1-t)}\nonumber\\
&=-\frac{t}{1-t}\sum_{n=0}^{\infty}t^{n}L_n(x)\label{eq:5_38}
\end{align}
En multipliant par $(1-t)$, on trouve
\begin{align}
\sum_{n=0}^{\infty}t^{n}L_n'(x)-\sum_{n=0}^{\infty}t^{n+1}L_n'(x)=-\sum_{n=0}^{\infty}t^{n+1}L_n(x)\nonumber
\end{align}
qui peut \^etre \'ecrite sous la forme
\begin{align*}
\sum_{n=0}^{\infty}t^{n}L_n'(x)-\sum_{n=1}^{\infty}t^{n}L_{n-1}'(x)=-\sum_{n=1}^{\infty}t^{n}L_{n-1}(x)\nonumber
\end{align*}
Par identification des coefficients de $t^{n}$ des deux c\^ot\'es pour $n\geq 1$, on trouve
\begin{align*}
L_n'(x)-L_{n-1}'(x)=-L_{n-1}(x)\hspace{10mm}(n\geq 1)
\end{align*}
cette \'equation peut \^etre \'ecrite de deux mani\`eres differentes (pour $n$ et $n+1$)
\begin{align}
&L_{n-1}'(x)=L_n'(x)+L_{n-1}(x)\nonumber\\
&L_{n+1}'(x)=L_n'(x)-L_n(x)\nonumber
\end{align}
La  d\'eriv\'ee de la formule (\ref{eq:5_37} a) par rapport \`a $x$ s'\'ecrit
\begin{align*}
(n+1)L_{n+1}'(x)=(2n+1-x)L_n'(x)-L_n(x)-nL_{n-1}'(x)\nonumber
\end{align*}
En remplaçant l'expressions de $L_{n-1}'(x)$ et $L_{n+1}'(x)$, on trouve apr\'es simplifications
\begin{align*}
xL_n'(x)=nL_n(x)-nL_{n-1}(x)
\end{align*}
(c) En d\'eveloppant $1/(1-t)$ en s\'erie enti\`ere et en remplaçant dans (\ref{eq:5_38}), on obtient
\begin{align}
\sum_{n=0}^{\infty}t^{n}L_n'(x)&=-t\sum_{r=0}^{\infty}\sum_{s=0}^{\infty}t^rt^{s}L_s(x)\nonumber\\
&=-\sum_{r,s=0}^{\infty}t^{r+s+1}L_s(x)\nonumber
\end{align}
Posons $n=r+s+1$, \i.e. $r=n-s-1$, et comme $r\geq 0$ on a $n-s-1 \geq 0$, cela imlique que $s\leq n-1$, donc il en r\'esulte que
\begin{align}
\sum_{n=0}^{\infty}t^{n}L_n'(x)=-\sum_{n=0}^{\infty}t^{n}\sum_{s=0}^{n-1}L_{s}(x)\nonumber
\end{align}
Par identification des coefficients de $t^{n}$ des deux c\^ot\'es, on obtient
\begin{align*}
L_n'(x)=-\sum_{s=0}^{n-1}L_{s}(x)\nonumber
\end{align*}
\section{Poly\^nome de Laguerre associ\'e}
L'\'equation de Laguerre associ\'ee est donn\'ee par
\begin{align}\label{eq:5_39}
x\frac{d^{2}y}{dx^{2}}+(k+1-x)\frac{dy}{dx}+ny=0
\end{align}
la solution de cette \'equation  est donn\'ee par 
\begin{align*}
y_n^{k}(x)=\frac{d^{k}}{dx^{k}}L_{n+k}(x)
\end{align*}
cette solution multipli\'ee par la constante $(-1)^{k}$ est appell\'ee polyn\^ome associ\'e de Laguerre not\'e $L_n^{k}(x)$ 
\begin{align}\label{eq:5_40}
L_n^{k}(x)=(-1)^{k}\frac{d^{k}}{dx^{k}}L_{n+k}(x)
\end{align}
et il est donn\'e explicitement par
\begin{align}\label{eq:5_41}
L_n^{k}(x)=\sum_{r=0}^{n}(-1)^{k}\frac{(n+k)!}{(n-r)!(k+r)!}x^{r}
\end{align}
\subsection{Fonction g\'en\'eratrice}
\begin{align}\label{eq:5_42}
\frac{e^{-xt/(1-t)}}{(1-t)^{k+1}}=\sum_{n=0}^{\infty}t^{n}L_n^{k}(x)
\end{align}
\textit{D\'emonstration}\\
On d\'erive l'\'equation de la fonction g\'en\'eratrice du polyn\^ome de Laguerre (Eq.~(\ref{eq:5_35})) $k$ fois par rapport \`a $x$, on trouve
\begin{align*}
\frac{1}{1-t}\frac{d^{k}}{dx^{k}}e^{-xt/(1-t)}=\frac{d^{k}}{dx^{k}}\sum_{n=k}^{\infty}L_{n}(x)t^{n}\nonumber
\end{align*}
la somme du c\^ot\'e droit commence par $n=k$ parceque $L_n(x)$ est un polyn\^ome de degr\'e $n$ et la d\'eriv\'ee d'ordre $k$ de $L_n(x)$ donne $0$ si $n<k$. La d\'eriv\'ee d'ordre $k$ de l'exponentielle donne
\begin{align*}
\left(-\frac{t}{1-t}\right)^{k}\frac{1}{1-t}e^{-xt/(1-t)}=\frac{d^{k}}{dx^{k}}\sum_{n=0}^{\infty}L_{n+k}(x)t^{n+k}\nonumber
\end{align*}
si on ins\`ere l'expression (\ref{eq:5_35})) de la fonction g\'en\'eratrice dans l'\'equation pr\'ec\'edente on obtient
\begin{align*}
(-1)^{k}\frac{t^{k}}{(1-t)^{k+1}}e^{-xt/(1-t)}=\sum_{n=0}^{\infty}(-1)^{k}L_{n}^{k}(x)t^{n+k}\nonumber
\end{align*}
Apr\`es simplifications on trouve le r\'esultat d\'esir\'e
\begin{align*}
\frac{1}{(1-t)^{k+1}}e^{-xt/(1-t)}=\sum_{n=0}^{\infty}t^{n}L_{n}^{k}(x)\nonumber
\end{align*}
\subsection{Relation d'orthogonalit\'e}
\begin{align}\label{eq:5_43}
\int_{0}^{\infty}e^{-x}x^{k}L_n^{k}(x)L_m^{k}(x)dx=\frac{(n+k)!}{n!}\delta_{nm}
\end{align}
\textit{D\'emonstration}\\
D'apr\'es la formule (\ref{eq:5_42}), on a 
\begin{align}
&\frac{1}{(1-t)^{k+1}}e^{-xt/(1-t)}=\sum_{n=0}^{\infty}t^{n}L_{n}^{k}(x)\nonumber\\
&\frac{1}{(1-s)^{k+1}}e^{-xs/(1-s)}=\sum_{m=0}^{\infty}s^{m}L_{m}^{k}(x)\nonumber
\end{align}
donc
\begin{align}\label{eq:5_44}
\sum_{n,m=0}^{\infty}t^{n+m}\int_{0}^{\infty}e^{-x}x^{k}L_{n}^{k}(x)L_{m}^{k}(x)dx=&\frac{1}{(1-t)^{k+1}}\frac{1}{(1-s)^{k+1}}\\
&\times\int_{0}^{\infty}e^{-x}x^{k}e^{-xt/(1-t)}e^{-xs/(1-s)}dx\nonumber
\end{align}
l'int\'egrale peut \^etre calculer en posant: $\left\lbrace -1-t/(1-t)-s/(1-s)\right\rbrace x = -ax = u$ o\`u $a=\left\lbrace -1-t/(1-t)-s/(1-s)\right\rbrace$, \i.e. $dx=-du/a$, on obtient
\begin{align}
\int_{0}^{\infty}x^{k}\exp\left\lbrace -x-\frac{xt}{(1-t)} -\frac{xs}{(1-s)}\right\rbrace dx &=\int_{0}^{\infty}\left(\frac{u}{a}\right)^{k}e^{-u}\frac{du}{a}\nonumber\\
&=\frac{1}{a^{k+1}}\int_{0}^{\infty}u^{k}e^{-u}du\nonumber\\
&=\frac{1}{a^{k+1}}\Gamma(k+1)\nonumber\\
&=\frac{k!(1-t)^{k+1}(1-s)^{k+1}}{(1-st)^{k+1}}\nonumber
\end{align}
En remplaçant le r\'esultat obtenu de l'int\'erale dans l'expression (\ref{eq:5_44}), on trouve
\begin{align}
\sum_{n,m=0}^{\infty}t^{n+m}\int_{0}^{\infty}e^{-x}x^{k}L_{n}^{k}(x)L_{m}^{k}(x)dx=\frac{k!}{(1-st)^{k+1}}\nonumber
\end{align}
Le d\'eveloppement en s\'erie de $1/(1-st)^{k+1}$ donne
\begin{align*}
\frac{1}{(1-st)^{k+1}}=\sum_{n=0}^{\infty}\frac{(n+k)!}{k!}\frac{s^{n}t^{n}}{n!}
\end{align*}
donc, on otient finalement
\begin{align}
\sum_{n,m=0}^{\infty}t^{n+m}\int_{0}^{\infty}e^{-x}x^{k}L_{n}^{k}(x)L_{m}^{k}(x)dx=\sum_{n=0}^{\infty}\frac{(n+k)!}{n!}s^{n}t^{n}\nonumber
\end{align}
En identifiant les coefficients de $s^{m}t^{n}$ des deux c\^ot\'es, on trouve
\begin{align*}
\int_{0}^{\infty}e^{-x}x^{k}L_n^{k}(x)L_m^{k}(x)dx=\frac{(n+k)!}{n!}\delta_{nm}
\end{align*}
\subsection{Relations de r\'ecurrence}
\begin{align}\label{eq:5_45}
&(a)\hspace{1.5mm} L_{n-1}^{k}(x)+L_n^{k-1}(x)=L_n^{k}(x)\nonumber\\
&(b)\hspace{1.5mm}(n+1)L_{n+1}^{k}(x)=(2n+k+1-x)L_n^{k}(x)-(n+k)L_{n-1}^{k}(x)\nonumber\\
&(c)\hspace{1.5mm} x {L_n^{k}}'(x)=nL_n^{k}(x)-(n+k)L_{n-1}^{k}(x)\nonumber\\
&(d)\hspace{1.5mm} {L_n^{k}}'(x)=-\sum_{r=0}^{n-1}L_r^{k}(x)\\
&(e)\hspace{1.5mm} {L_n^{k}}'(x)=-L_{n-1}^{k+1}(x)\nonumber\\
&(f)\hspace{1.5mm} L_n^{k+1}(x)=\sum_{r=0}^{n}L_r^{k}(x)\nonumber
\end{align}
\textit{D\'emonstration}\\
(a) On a
\begin{align*}
&L_{n-1}^{k}(x)+L_n^{k-1}(x)=\sum_{r=0}^{n-1}\frac{(-1)^{r}(n-1+k)!}{(n-1-r)!(k+r)!r!}x^{r}+\sum_{r=0}^{n}\frac{(-1)^{r}(n+k-1)!}{(n-r)!(k-1+r)!r!}x^{r}\\
&=\sum_{r=0}^{n-1}\frac{(-1)^{r}(n-1+k)!}{(n-1-r)!(k+r)!r!}x^{r}+\sum_{r=0}^{n-1}\frac{(-1)^{r}(n+k-1)!}{(n-r)!(k-r+1)!r!}x^{r}+\frac{(-1)^{n}(n-1+k)!}{(n-n)!(k-1+n)!n!}x^{r}\\
&=\sum_{r=0}^{n-1}(-1)^{r}\frac{(n+k-1)!}{(n-r-1)!(k+r-1)!r!}\left\lbrace \frac{1}{k+r} +\frac{1}{n-r}\right\rbrace x^{r} + (-1)^{n}\frac{x^{n}}{n!}\\
&=\sum_{r=0}^{n-1}(-1)^{r}\frac{(n+k-1)!}{(n-r-1)!(k+r-1)!r!} \frac{n+k}{(k+r)(n-r)} x^{r} + (-1)^{n}\frac{x^{n}}{n!}\\
&=\sum_{r=0}^{n-1}(-1)^{r}\frac{(n+k)!}{(n-r)!(k+r)!r!} + (-1)^{n}\frac{x^{n}}{n!}\\
&=\sum_{r=0}^{n}(-1)^{r}\frac{(n+k)!}{(n-r)!(k+r)!r!}\\
&=L_n^{k}(x)
\end{align*}
(b) D\'erivons $k$ fois l'expression (\ref{eq:5_37})(a) avec $n$ remplac\'e par $n+k$. On trouve
\begin{align}
(n+k+1)\frac{d^{k}}{dx^{k}}L_{n+k+1}(x)=&(2n+2k+1)\frac{d^{k}}{dx^{k}}L_{n+k}(x)-\frac{d^{k}}{dx^{k}}\left\lbrace xL_{n+k}(x)\right\rbrace\nonumber\\
&-(n+k)\frac{d^{k}}{dx^{k}}L_{n+k-1}(x)\nonumber
\end{align}
En utilisant la formule de d\'erivation de Leibniz, on a
\begin{align}\label{eq:5_46}
\frac{d^{k}}{dx^{k}}\left\lbrace fg\right\rbrace=\sum_{n=0}^{k}\frac{k!}{n!(n-k)!}\frac{d^{n}f}{dx^{n}}\frac{d^{n-k}g}{dx^{n-k}}
\end{align}
alors
\begin{align*}
\frac{d^k}{dx^k} \left\lbrace xL_{n+k}(x)\right\rbrace &=\sum_{r=0}^{k}\frac{k!}{r!(k-r)!}\frac{d^{r}}{dx^{r}}(x) \frac{d^{k-r}}{dx^{k-r}}L_{n+k}(x)\\
&=x\frac{d^{k}}{dx^{k}}L_{n+k}(x)+k\frac{d^{k-1}}{dx^{k-1}}L_{n+k}(x)
\end{align*}
donc on obtient
\begin{align*}
(n+k+1)\frac{d^{k}}{dx^{k}}L_{n+k+1}(x)&=(2n+2k+1)\frac{d^{k}}{dx^{k}}L_{n+k}(x)-x\frac{d^{k}}{dx^{k}}L_{n+k}(x)-k\frac{d^{k-1}}{dx^{k-1}}L_{n+k}(x)\\
&-(n+k)\frac{d^{k}}{dx^{k}}L_{n+k-1}(x)
\end{align*}
En utilisant la d\'efinition du polyn\^ome de Laguerre associ\'e (\ref{eq:5_41}), l'\'equation pr\'ec\'edente devient
\begin{align*}
(n+k+1)(-1)^{k}L_{n+1}^{k}(x)&=(2n+2k+1)(-1)^{k}L_{n}^{k}(x)-x(-1)^{k}L_{n}^{k}(x)-k(-1)^{k-1}L_{n+1}^{k-1}(x)\\
&-(n+k)(-1)^{k}L_{n-1}^{k}(x)
\end{align*}
et si on prend en consid\'eration la relation (\ref{eq:5_45})(a) avec $n$ remplac\'e par $n+1$, 
on arrive au r\'esultat d\'esir\'e
\begin{align*}
(n+1)L_{n+1}^{k}(x)=(2n+k+1-x)L_n^{k}(x)-(n+k)L_{n-1}^{k}(x)
\end{align*}
(c) En d\'erivant $k$ fois l'expression (\ref{eq:5_45})(b) par rapport \`a $x$, on obtient
\begin{align*}
\frac{d^{k}}{dx^{k}}\left\lbrace xL_{n+k}'(x)\right\rbrace =(n+k)\frac{d^{k}}{dx^{k}}L_{n+k}(x)-(n+k)\frac{d^{k}}{dx^{k}}L_{n+k-1}(x)
\end{align*}
En utilisant la r\'egle de Leibniz Eq.~(\ref{eq:5_46}), on trouve
\begin{align*}
x\frac{d^{k}}{dx^{k}}{L_{n+k}'}(x)+k\frac{d^{k}}{dx^{k}}L_{n+k}(x)=(n+k)\frac{d^{k}}{dx^{k}}L_{n+k}(x)-(n+k)\frac{d^{k}}{dx^{k}}L_{n+k-1}(x)
\end{align*}
\`a partir de la d\'efinition du polyn\^ome de Laguerre associ\'e, on obtient (apr\'es simplifications)
\begin{align*}
x{L_n^{k}}'(x)=nL_n^{k}(x)-(n+k)L_{n-1}^{k}(x)
\end{align*}
(d) D\'erivons $k$ fois l'expression (\ref{eq:5_45})(c) par rapport \`a $x$. Il en r\'esulte
\begin{align*}
\frac{d^{k}}{dx^{k}}L_{n+k}'(x)=-\sum_{r=0}^{n+k-1}\frac{d^{k}}{dx^{k}}L_{r}(x)
\end{align*}
ou $d^{k}L_{n+k}'(x)/dx^{k}=(-1)^{k}{L_{n}^{k}}'(x)$. Comme $L_{r}(x)$ est un polyn\^ome de degr\'e $r$, la d\'eriv\'ee d'ordre $k$ de $L_{r}(x)$ 
vaut $0$ si $k>r$ donc 
\begin{align*}
(-1)^{k}{L_{n}^{k}}'(x)=-\sum_{r=0}^{n+k-1}\frac{d^{k}}{dx^{k}}L_{r}(x)
\end{align*}
En faisant le changement de variable $s=r-k$, on trouve 
\begin{align*}
{L_n^{k}}'(x)=-\sum_{r=0}^{n-1}L_r^{k}(x)
\end{align*}
(e) Par d\'efinition du polyn\^ome de Laguerre associ\'e, on a
\begin{align*}
{L_{n}^{k}}'(x)&=\frac{d}{dx}\sum_{r=0}^{n}(-1)^{r}\frac{(n+k)!}{(n-r)!(k+r)!r!}x^r\\
&=\sum_{r=1}^{n}(-1)^{r}\frac{(n+k)!}{(n-r)!(k+r)!(r-1)!!}x^{r-1}\\
&=\sum_{s=0}^{n}(-1)^{s+1}\frac{(n+k)!}{(n-s-1)!(k+s+1)!s!}x^s\\
&=-\sum_{s=0}^{n}(-1)^{s+1}\frac{(n-1+k+1)!}{(n-s-1)!(k+s+1)!s!}x^s\\
&=-L_{n-1}^{k+1}(x)
\end{align*}
(f) En comparant les r\'esultats pr\'ec\'edents (d) et (e), on a
\begin{align*}
-\sum_{r=0}^{n-1}L_r^{k}(x)=-L_{n-1}^{k+1}(x)
\end{align*}
En remplaçant $n$ par $n+1$, on trouve
\begin{align*}
L_n^{k+1}(x)=-\sum_{r=0}^{n}L_r^{k}(x)
\end{align*}
\section{Polyn\^ome de Chebyshev}
L'\'equation de Chebyshev est donn\'ee par
\begin{align}\label{eq:5_47}
(1-x^{2})\frac{d^{2}y}{dx^{2}}-x\frac{dy}{dx}+n^{2}y=0
\end{align}
Les deux solutions ind\'ependentes de cette \'equation sont appell\'ees polyn\^omes de Chebyshev de premi\`ere esp\`ece (not\'e $T_n(x)$) et de deuxi\`eme esp\`ece 
(not\'e $U_n(x)$) et sont d\'efinis par
\begin{align}
&T_n(x)=\cos\left(n\cos^{-1}x\right)\label{eq:5_48}\\
&U_n(x)=\sin\left(n\cos^{-1}x\right)\label{eq:5_49}
\end{align}
o\`u $n$ est un entier non-n\'egatif.
Les polyn\^omes de Chebyshev peuvent \^etre \'ecrit sous les formes
\begin{align}
&T_n(x)=\frac{1}{2}\left\lbrace\left[x+i\sqrt{1-x^2}\right]^{n}+\left[x-i\sqrt{1-x^2}\right]^{n}\right\rbrace\label{eq:5_50}\\
&U_n(x)=-\frac{i}{2}\left\lbrace\left[x+i\sqrt{1-x^2}\right]^{n}-\left[x-i\sqrt{1-x^2}\right]^{n}\right\rbrace\label{eq:5_51}
\end{align}
\textit{D\'emonstration}\\
En faisant le changement de variable $x=\cos\theta$ dans l'expression (\ref{eq:5_48}), on obtient
\begin{align*}
T_n(x)&=\cos\left(n\cos^{-1}\cos\theta\right)\\
&=\cos\left(n\theta\right)\\
&=\frac{1}{2}\left\lbrace e^{in\theta}+e^{-in\theta}\right\rbrace\\
&=\frac{1}{2}\left\lbrace\left(e^{i\theta}\right)^{n}+\left(e^{-i\theta}\right)^{n}\right\rbrace\\
&=\frac{1}{2}\left\lbrace \left[\cos\theta + i\sin\theta\right]^{n} +\left[\cos\theta - i\sin\theta\right]^{n}\right\rbrace\\
&=\frac{1}{2}\left\lbrace \left[x + i\sqrt{1-x^2}\right]^{n} +\left[x - i\sqrt{1-x^2}\right]^{n}\right\rbrace
\end{align*}
La preuve de la deuxi\`eme relation (\ref{eq:5_49}) pour $U_n(x)$ est similaire \`a celle de $T_n(x)$.
\subsection{Repr\'esentations en s\'eries}
\begin{align}
&T_n(x)=\sum_{r=0}^{[\frac{n}{2}]}(-1)^{r}\frac{n!}{(2r)!(n-2r)!}\left(1-x^{2}\right)^{r}x^{n-2r}\label{eq:5_52}\\
&U_n(x)=\sum_{r=0}^{[\frac{n-1}{2}]}(-1)^{r}\frac{n!}{(2r+1)!(n-2r-1)!}\left(1-x^{2}\right)^{r+\frac{1}{2}}x^{n-2r-1}\label{eq:5_53}
\end{align}
\textit{D\'emonstration}\\
En utilsant la formule du bin\^ome de Newton pour les deux termes de la formule (\ref{eq:5_50}), on obtient
\begin{align*}
T_n(x)&=\frac{1}{2}\left\lbrace\left[x+i\sqrt{1-x^2}\right]^{n}+\left[x-i\sqrt{1-x^2}\right]^{n}\right\rbrace\\
&=\frac{1}{2}\left\lbrace\sum_{r=0}^{n}\frac{n!}{r!(n-r)!}x^{n-r}\left[i\sqrt{1-x^{2}}\right]^{r}+\sum_{r=0}^{n}\frac{n!}{r!(n-r)!}x^{n-r}\left[-i\sqrt{1-x^{2}}\right]^{r}\right\rbrace\\
&=\frac{1}{2}\sum_{r=0}^{n}\frac{n!}{r!(n-r)!}x^{n-r}\left(1-x^{2}\right)^{r/2}i^{r}\left(1+(-1)^{r}\right)
\end{align*}
si $r$ est impair, on a $1+(-1)^{r}=0$ et si $r$ est pair, on a $1+(-1)^{r}=2$, posant $r=2s$, on obtient
\begin{align*}
T_n(x)=\frac{1}{2}\sum_{s=0}^{[\frac{n}{s}]}\frac{n!}{2!(n-2s)!}x^{n-2s}\left(1-x^{2}\right)^{s}2(-1)^{s}
\end{align*}
La preuve de la deuxi\`eme expression (\ref{eq:5_53}) est similaire \`a celle de (\ref{eq:5_52}).\\
Quelques expressions explicites pour $T_n(x)$ et $U_n(x)$ (calcul\'ees a partir des formules (\ref{eq:5_52}) et (\ref{eq:5_53})) sont donn\'ees ci-dessous
\begin{align*}
&T_0(x)=1,\\
&T_1(x)=x,\\
&T_2(x)=2x^2-1,\\
&T_3(x)=4x^3-3x,\\
&T_4(x)=8x^4-8x^2+1,\\
&T_5(x)=16x^5-20x^3+5x.
\end{align*}
\begin{figure}[!ht]
\centering
\includegraphics[scale=0.7]{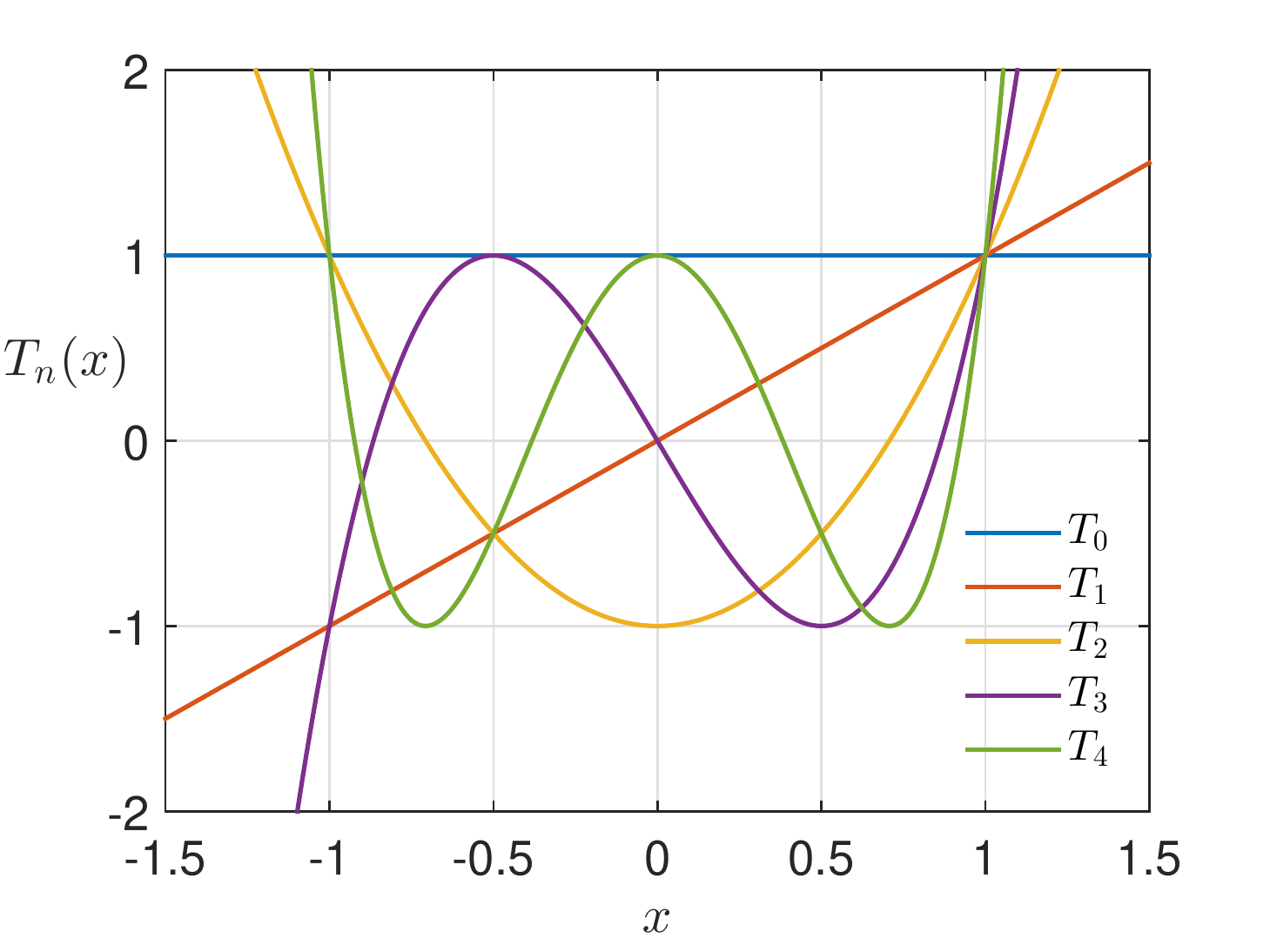}
\caption{\label{fig:chebyshevT} Polyn\^omes de Chebyshev de premi\`ere esp\`ece pour $n=0,1,2,3,4$.}
\end{figure}
\begin{align*}
&U_0(x)=0,\\
&U_1(x)=\sqrt{1-x^2},\\
&U_2(x)=\sqrt{1-x^2}2x,\\
&U_3(x)=\sqrt{1-x^2}(4x^2-1),\\
&U_4(x)=\sqrt{1-x^2}(8x^3-4x),\\
&U_5(x)=\sqrt{1-x^2}(16x^4-12x^2+1).
\end{align*}
\begin{figure}[!ht]
\centering
\includegraphics[scale=0.7]{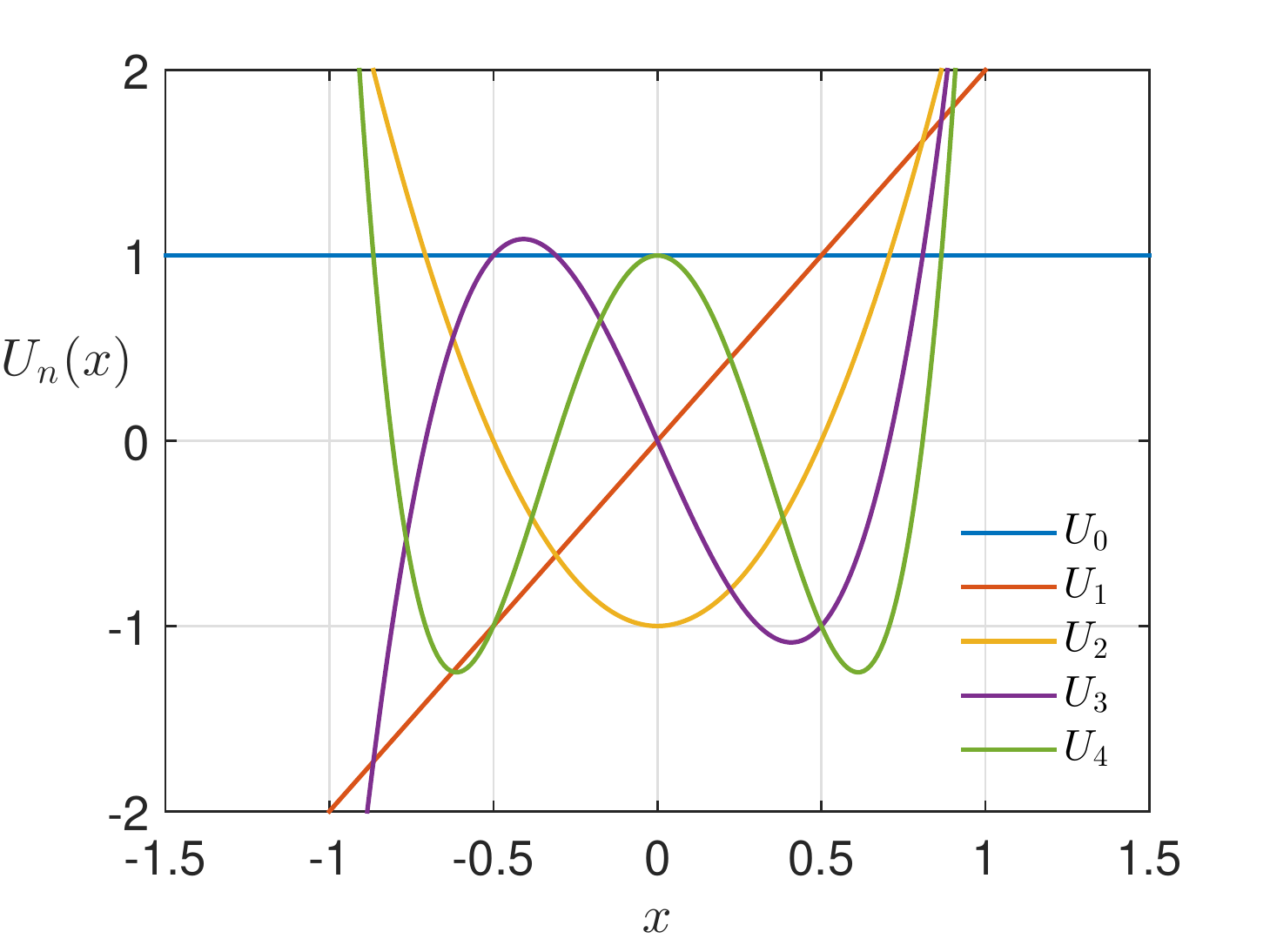}
\caption{\label{fig:chebyshevU} Polyn\^omes de Chebyshev de deuxi\`eme esp\`ece pour $n=0,1,2,3,4$.}
\end{figure}
Elles sont repr\'esent\'ees sur les figures Fig.~(\ref{fig:chebyshevT}) et (\ref{fig:chebyshevU}) respectivement.
\subsection{Fonctions g\'en\'eratrices}
Les fonctions g\'en\'eratrices des polyn\^omes de Chebyshev sont donn\'ees par
\begin{align}
&\frac{1-t^{2}}{1-2tx+t^{2}}=T_0(x)+2\sum_{n=1}^{\infty}t^{n}T_n(x)\label{eq:5_54}\\
&\frac{\sqrt{1-t^{2}}}{1-2tx+t^{2}}=\sum_{n=0}^{\infty}t^{n}U_{n+1}(x)\label{eq:5_55}
\end{align}
\textit{D\'emonstration}\\
Pour d\'emontrer (\ref{eq:5_54}), on fait le changement de variable $x=\cos\theta=(e^{i\theta}+e^{-i\theta})/2$. On trouve
\begin{align*}
\frac{1-t^{2}}{1-2tx+t^{2}}&=\frac{1-t^{2}}{1-(e^{i\theta}+e^{-i\theta})t+t^{2}}\\
&=\frac{1-t^{2}}{(1-e^{i\theta}t)(1-e^{-i\theta}t)}\\
&=(1-t^2)\sum_{r=0}^{\infty}\left(te^{i\theta}\right)^{r}\sum_{s=0}^{\infty}\left(te^{-i\theta}\right)^{s}\\
&=(1-t^2)\sum_{r,s=0}^{\infty}e^{i(r-s)\theta}t^{r+s}\\
&=\sum_{r,s=0}^{\infty}e^{i(r-s)\theta}t^{r+s}-\sum_{r,s=0}^{\infty}e^{i(r-s)\theta}t^{r+s+2}
\end{align*}
En posant $n=r+s$, on obtient le coefficient de $t^{0}$, \i.e. pour $n=0$, en consid\'erant $r=0$ et 
$s=0$ . Donc le coefficient est donn\'ee par
\begin{align*}
e^{i(0-0)\theta}=1=T_0(x)
\end{align*}
On obtient le coefficient de $t^{1}$ pour $n=1$ en prenant soit $r=1$ et $s=0$ o\`u inversement, donc
\begin{align*}
e^{i\theta}+e^{-i\theta}&=2\cos\theta\\
&=2x\\
&=2T_1(x)
\end{align*}
pour $n\geq 2$, le coefficient de $t^{n}$ est obtenue soit en posant $n=s+r$ (\i.e., $s=n-r$) dans la premi\`ere somme ou bien en posant $n=r+s+2$ (\i.e., $s=n-r-2$) dans la deuxi\`eme somme, donc le coefficient de $t^{n}$ est donn\'e par
\begin{align*}
\sum_{r=0}^{n}e^{i(r-(n-r))\theta}-\sum_{r=0}^{n-2}e^{i(r-(n-r-2))\theta}
&=e^{-in\theta}\sum_{r=0}^{n}e^{i2r\theta}-e^{-i(n-2)\theta}\sum_{r=0}^{n-2}e^{i2r\theta}\\
&=e^{-in\theta}\frac{1-\left(e^{i2\theta}\right)^{n+1}}{1-e^{i2\theta}}-e^{-i(n-2)\theta}\frac{1-\left(e^{i2\theta}\right)^{n-1}}{1-e^{i2\theta}}\\
&=\frac{e^{-in\theta}-e^{i(n+2)\theta}-e^{-i(n-2)\theta}+e^{in\theta}}{1-e^{i2\theta}}\\
&=\frac{e^{-in\theta}\left[1-e^{i2\theta}\right]+e^{in\theta}\left[1-e^{i2\theta}\right]}{1-e^{i2\theta}}\\
&=e^{in\theta}+e^{-in\theta}\\
&=2\cos(n\theta)\\
&=2T_n(x)
\end{align*}
On obtient finalement
\begin{align*}
\frac{1-t^{2}}{1-2tx+t^{2}}&=T_0(x)+2tT_1(x)+2\sum_{n=2}^{\infty}t^{n}T_n(x)\\
&=T_0(x)+2\sum_{n=1}^{\infty}t^{n}T_n(x)
\end{align*}
On peut utiliser exactement la m\^eme preuve pour d\'emontrer l'expression (Eq.~(\ref{eq:5_55})) pour $U_n(x)$.
\subsection{Relations d'orthogonalit\'e}
\begin{align}
&\int_{-1}^{1}\frac{T_n(x)T_m(x)}{\sqrt{1-x^{2}}}dx=\left\lbrace\begin{array}{ll}
0 &  m\neq n\\
\pi/2 & m=n\neq0\\
\pi & m=n=0
\end{array}\right.\label{eq:5_56}\\
&\int_{-1}^{1}\frac{U_n(x)U_m(x)}{\sqrt{1-x^{2}}}dx=\left\lbrace\begin{array}{ll}
0 &  m\neq n\\
\pi/2 & m=n\neq0\\
0 & m=n=0
\end{array}\right.\label{eq:5_57}
\end{align}
\textit{D\'emonstration}\\
Si on pose $x=\cos\theta$, on a
\begin{align*}
\int_{-1}^{1}\frac{T_n(x)T_m(x)}{\sqrt{1-x^{2}}}dx&=\int_{\pi}^{0}\frac{T_n(\cos\theta)T_m(\cos\theta)}{\sqrt{1-\cos^{2}\theta}}(-\sin\theta d\theta)\\
&=\int_{0}^{\pi}\cos(n\theta)\cos(m\theta)d\theta\\
&=\int_{0}^{\pi}\frac{1}{2}\left\lbrace\cos(n+m)\theta+\cos(n-m)\theta\right\rbrace d\theta\nonumber\\
&=\frac{1}{2}\Bigl[-\frac{1}{n+m}\sin(n+m)\theta - \frac{1}{n-m}\sin(n-m)\theta \Bigr]_0^{\pi}\nonumber\\
&=0 \hspace{10mm} \mbox{si $n-m\neq 0$}
\end{align*}
l'int\'egrale est \'egale a $0$ pour $n\neq m$ car $n+m$ et $n-m$ sont des nombres entiers et donc le \textit{sin} est toujours nulle. Par contre si $n=m$, on a (\`a partir de la deuxi\`eme ligne du raisonnement pr\'ec\'edent)
\begin{align*}
\int_{-1}^{1}\frac{T_n(x)T_m(x)}{\sqrt{1-x^{2}}}dx&=\int_{0}^{\pi}\cos^{2}(n\theta)d\theta\\
&=\int_{0}^{\pi}\frac{1}{2}\left(1+\cos(2n\theta)\right)d\theta\\
&=\frac{1}{2}\left[\theta-\frac{1}{2n}\sin(2n\theta)\right]_0^{\pi}\\
&=\frac{\pi}{2} \hspace{10mm} \mbox{si $n\neq 0$}
\end{align*}
si $n=m=0$, on a
\begin{align*}
\int_{-1}^{1}\frac{T_0(x)T_0(x)}{\sqrt{1-x^{2}}}dx&=\int_{0}^{\pi}d\theta\nonumber\\
&=\pi
\end{align*}
On peut suivre les m\^eme \'etapes pour d\'emontrer la deuxi\`eme relation d'orthogonalit\'e (pour $U_n(x)$).
\subsection{Relations de r\'ecurrence}
\begin{align}\label{eq:5_58}
&(a)\hspace{1.5mm}T_{n+1}(x)-2xT_n(x)+T_{n-1}(x)=0\\
&(b)\hspace{1.5mm}(1-x^{2})T_n'(x)=-nxT_n(x)+nT_{n-1}(x)\nonumber\\
&(c)\hspace{1.5mm}U_{n+1}(x)-2xU_n(x)+U_{n-1}(x)=0\nonumber\\
&(d)\hspace{1.5mm}(1-x^{2})U_n'(x)=-nxU_n(x)+nU_{n-1}(x)\nonumber
\end{align}
\textit{D\'emonstration}\\
(a) Si on remplace $x$ par $\cos\theta$ dans l'expression (\ref{eq:5_58})(a), on obtient
\begin{align*}
\cos\left((n+1)\theta\right) - 2\cos(\theta)\cos(n\theta)+\cos\left((n-1)\theta\right) &=\cos(n\theta)\cos\theta -\sin(n\theta)\sin\theta \\
-2\cos(\theta)\cos(n\theta) &+\cos(n\theta)\cos\theta + \sin(n\theta)\sin\theta\\
&=0
\end{align*}
(b) Si on remplace maintenant $x$ par $\cos\theta$ dans l'Eq.~(\ref{eq:5_58})(b), on obtient
\begin{align*}
(1-\cos^{2}\theta)\frac{d}{d(\cos\theta)}\cos(n\theta)=-n\cos(\theta)\cos(n\theta)+n\cos(n-1)\theta
\end{align*}
le c\^ot\'e gauche donne
\begin{align*}
(1-\cos^{2}\theta)\frac{d}{d(\cos\theta)}\cos(n\theta)&=\sin^{2}\theta\left(-\frac{1}{\sin\theta}\frac{d}{d\theta}\cos(n\theta)\right)\\
&=-\sin\theta(-n\sin(n\theta))\\
&=n\sin(\theta)\sin(n\theta)
\end{align*}
et le c\^ot\'e droit donne
\begin{align*}
-n\cos(\theta)\cos(n\theta)+n\cos(n-1)\theta &=-n\cos(\theta)\cos(n\theta)+n(\cos(n\theta)\cos\theta+\sin(n\theta)\sin\theta)\\
&=n\sin(\theta)\sin(n\theta)
\end{align*}
Les deux c\^otes sont \'eguax, la propri\'et\'e est ainsi d\'emontr\'e.\\
La preuve des deux relations de r\'ecurrence pour $U_n(x)$ sont similaires \`a celles de $T_n(x)$.
\section{Exercices}
$\mathbf{Exercice\hspace{1mm} 1:}$ \\
1. \hspace{1.5mm} Montrer que
\begin{align*}
\int_{-1}^{1}\frac{P_n(t)}{\sqrt{1+x^2-2xt}}dt=\frac{2x^{n}}{2n+1}
\end{align*}
2. \hspace{1.5mm}En utilisant la formule de Rodrigues, montrer que
\begin{align*}
P_{n+1}'(x)=P_{n-1}'(x)+(2n+1)P_n(x)
\end{align*}
3. \hspace{1.5mm} Utiliser l'identit\'e
\begin{align*}
\int_{-1}^{1}f(x)\frac{d^{n}}{dx^{n}}(x^2-1)^{n}dx=(-1)^{n}\int_{-1}^{1}(x^2-1)^{n}\frac{d^{n}}{dx^{n}}f(x)
\end{align*}
pour d\'evelopper la fonction $f(x)=\ln\left(\frac{1+x}{1-x}\right)$ sous la forme s\'erie de polyn\^omes de Legendre pour $-1<x<1$.\\
$\mathbf{Solution}$\\
1. \hspace{1.5mm} Soit $\phi(t)=(1+x^2-2xt)^{-1/2}$ la fonction g\'en\'eratrice de $P_l^{m}(x)$
\begin{align*}
\phi(t)=\sum_{l=0}^{\infty}x^{l}P_l^{m}(t)
\end{align*}
donc
\begin{align*}
\int_{-1}^{1}\frac{P_n(t)}{\sqrt{1+x^2-2xt}}dt&=\int_{-1}^{1}P_n(t)\phi(t) dt\\
&=\sum_{l=0}^{\infty}x^{l}\int_{-1}^{1}P_n(t)P_l(t) dt\\
&=\sum_{l=0}^{\infty}x^{l}\frac{2\delta_{nl}}{2n+1}\\
&=\frac{2x^{n}}{2n+1}
\end{align*}
2. \hspace{1.5mm} D'apr\`es la formule de Rodrigues on a
\begin{align*}
P_{n+1}'(x)&=\frac{1}{2^{n+1}(n+1)!}\frac{d^{n+2}}{dx^{n+2}}\left[(x^2-1)^{n+1}\right]\\
&=\frac{2(n+1)}{2^{n+1}(n+1)!}\frac{d^{n+1}}{dx^{n+1}}\left[x(x^2-1)^{n}\right]\\
&=\frac{1}{2^{n}n!}\frac{d^{n}}{dx^{n}}\left[(x^2-1)^{n} + 2nx^2(x^2-1)^{n-1}\right]\\
&=\frac{1}{2^{n}n!}\frac{d^{n}}{dx^{n}}\left[(x^2-1)^{n} + 2n(x^2-1+1)(x^2-1)^{n-1}\right]\\
&=\frac{1}{2^{n}n!}\frac{d^{n}}{dx^{n}}\left[(x^2-1)^{n} + 2n(x^2-1)(x^2-1)^{n-1} + 2n(x^2-1)^{n-1}\right]\\
&=\frac{1}{2^{n}n!}\frac{d^{n}}{dx^{n}}\left[(x^2-1)^{n} + 2n(x^2-1)^{n} + 2n(x^2-1)^{n-1}\right]\\
&=\frac{2n+1}{2^{n}n!}\frac{d^{n}}{dx^{n}}\left[(x^2-1)^{n}\right] + \frac{2n}{2^{n}n!}\frac{d}{dx}\frac{d^{n-1}}{dx^{n-1}}\left[(x^2-1)^{n-1}\right]\\
&=(2n+1)P_n(x)+ P_{n-1}'(x)
\end{align*}
3. \hspace{1.5mm} Puisque la fonction $f(x)$ est impair, l'indice $n$ devrait \^etre impair. On a
\begin{align*}
c_n&=\frac{2n+1}{2}\int_{-1}^{1}f(x)P_n(x)dx\\
&=\frac{2n+1}{2}\frac{1}{2^{n}n!}\int_{-1}^{1}f(x)\frac{d^{n}}{dx^{n}}\left(x^2-1\right)^{n} dx\\
&=(-1)^{n}\frac{2n+1}{2}\frac{1}{2^{n}n!}\int_{-1}^{1}\left(x^2-1\right)^{n}\frac{d^{n}}{dx^{n}}\ln\left(\frac{1+x}{1-x}\right) dx
\end{align*}
Dans la derni\`ere ligne nous avons utilis\'e l'identit\'e donn\'ee. En plus on a
\begin{align*}
\frac{d^{n}}{dx^{n}}\left[\ln(1+x) - \ln(1-x)\right]=\frac{(-1)^{n-1}(n-1)!}{(1+x)!}+\frac{(n-1)!}{(1-x)^n}
\end{align*}
et
\begin{align*}
(x^2-1)^n&=(x+1)^n(x-1)^n\\
&=-(1+x)^n(1-x)^n \hspace{20mm} \mbox{si n est impair}
\end{align*}
donc
\begin{align*}
c_n&=\frac{(2n+1)(n-1)!}{2^{n+1}n!}\int_{-1}^{1}\left[(-1)^{n-1}(1-x)^{n}+(1+x)^{n}\right]dx\\
&=\frac{2n+1}{2^{n+1}n}\left\lbrace \left[-\frac{(1-x)^{n+1}}{n+1}\right]_{-1}^{1} + \left[\frac{(1+x)^{n+1}}{n+1}\right]_{-1}^{1}\right\rbrace\\
&=\frac{2n+1}{2^{n+1}n}\left(\frac{2^{n+1}}{n+1}+\frac{2^{n+1}}{n+1}\right)\\
&=\frac{2(2n+1)}{n(n+1)}
\end{align*}
$\mathbf{Exercice\hspace{1mm} 2:}$ Montrer que
\begin{align*}
&1. \hspace{1.5mm} \sqrt{1-x^2}T_n(x)=U^{n+1}(x)-xU_n(x)\\
&2. \hspace{1.5mm} T_{n+m}(x)+T_{n-m}(x)=2T_n(x)T_m(x)
\end{align*}
$\mathbf{Solutions}$\\
1. \hspace{1.5mm} Posons $x=\cos\theta$ et en utilisant les d\'efinitions Eq.~(\ref{eq:5_48}) et (\ref{eq:5_49}) on obtient
\begin{align*}
T_n(\cos\theta)=\cos(n\theta)\\
U_n(\cos\theta)=\sin(n\theta)
\end{align*}
Donc, nous devons montrer
\begin{align*}
\sin\theta \cos(n\theta)=\sin((n+1)\theta) - \cos\theta \sin(n\theta)
\end{align*}
Mais
\begin{align*}
\sin((n+1)\theta) - \cos\theta \sin(n\theta)&= \sin(n\theta)\cos\theta +\cos(n\theta)\sin\theta  - \cos\theta \sin(n\theta)\\
&=\cos(n\theta)\sin\theta
\end{align*}
ce qui prouve le r\'esultat.\\
2. \hspace{1.5mm} Posons $x=\cos\theta$ on a
\begin{align*}
T_{n+m}(\cos\theta)=\cos((n+m)\theta)\\
T_{n-m}(\cos\theta)=\cos((n-m)\theta)
\end{align*}
Donc
\begin{align*}
T_{n+m}(\cos\theta)+T_{n-m}(\cos\theta)&=\cos(n\theta)\cos(m\theta)-\sin(n\theta)\sin(m\theta)\\
&\hspace{10mm}+\cos(n\theta)\cos(m\theta)+\sin(n\theta)\sin(m\theta)\\
&=2\cos(n\theta)\cos(m\theta)\\
&=2T_n(x)T_m(x)
\end{align*}
$\mathbf{Exercice\hspace{1mm} 3:}$ Pour une fonction continue, on peut \'ecrire $f(x)=\sum_{n=0}^{\infty}c_nH_n(x)$.\\ Montrer que
\begin{align*}
c_n=\frac{1}{2^{n}\sqrt{\pi}n!}\int_{-\infty}^{\infty}f(x)H_n(x)e^{-x^2}dx
\end{align*}
$\mathbf{Solution}$\\
On a
\begin{align*}
\int_{-\infty}^{\infty}e^{-x^2}f(x)H_n(x)dx &= \sum_{r=0}^{\infty}c_r\int_{-\infty}^{\infty}e^{-x^2}H_r(x)H_n(x)\\
&=\sum_{r=0}^{\infty}c_r2^{n}n!\sqrt{\pi}\delta_{rn}\\
&=c_n2^{n}n!\sqrt{\pi}
\end{align*}
Donc
\begin{align*}
c_n=\frac{1}{2^{n}\sqrt{\pi}n!}\int_{-\infty}^{\infty}f(x)H_n(x)e^{-x^2}dx
\end{align*}
$\mathbf{Exercice\hspace{1mm} 4:}$ Soit $\psi_n(x)$ tel que
\begin{align*}
\psi_n(x)=\frac{1}{\sqrt{2^{n}\sqrt{\pi}n!}}e^{-x^2/2}H_n(x)
\end{align*}
Montrer que les op\'erateurs $\frac{1}{\sqrt{2}}\left(x+\frac{d}{dx}\right)$ et $\frac{1}{\sqrt{2}}\left(x-\frac{d}{dx}\right)$ sont respectivement des op\'erateurs d'annihilation et de creation pour $\psi_n(x)$\\
$\mathbf{Solution}$\\
On a
\begin{align*}
\frac{1}{\sqrt{2}}\left(x+\frac{d}{dx}\right)\psi_n(x)&=\frac{1}{\sqrt{2^{n+1}\sqrt{\pi}n!}}\left(x+\frac{d}{dx}\right)e^{-x^2/2}H_n(x)\\
&=\frac{1}{\sqrt{2^{n+1}\sqrt{\pi}n!}}\left(xe^{-x^2/2}H_n(x)- xe^{-x^2/2}H_n(x) + e^{-x^2/2}H_n'(x)\right)
\end{align*}
Mais d'apr\`es l'Eq.~(\ref{eq:5_31}-a) on a $H_n'(x)=2nH_{n-1}(x)$, donc
\begin{align*}
\frac{1}{\sqrt{2}}\left(x+\frac{d}{dx}\right)\psi_n(x)&=\frac{\sqrt{n}}{\sqrt{2^{n-1}\sqrt{\pi}(n-1)!}}e^{-x^2/2}2nH_{n-1}(x)\\
&=\sqrt{n}\psi_{n-1}(x)
\end{align*}
Pour $\frac{1}{\sqrt{2}}\left(x-\frac{d}{dx}\right)$ on a
\begin{align*}
\frac{1}{\sqrt{2}}\left(x-\frac{d}{dx}\right)\psi_n(x)&=\frac{1}{\sqrt{2^{n+1}\sqrt{\pi}n!}}\left(x-\frac{d}{dx}\right)e^{-x^2/2}H_n(x)\\
&=\frac{1}{\sqrt{2^{n+1}\sqrt{\pi}n!}}\left(xe^{-x^2/2}H_n(x)+ xe^{-x^2/2}H_n(x) - e^{-x^2/2}H_n'(x)\right)\\
&=\frac{1}{\sqrt{2^{n+1}\sqrt{\pi}n!}}\left(2xe^{-x^2/2}H_n(x)- e^{-x^2/2}2nH_{n-1}(x)\right)
\end{align*}
Mais d'apr\`es l'Eq.~(\ref{eq:5_31}-b) il en r\'esulte que
\begin{align*}
\frac{1}{\sqrt{2}}\left(x-\frac{d}{dx}\right)\psi_n(x)&=\frac{1}{\sqrt{2^{n+1}\sqrt{\pi}n!}}e^{-x^2/2}H_{n+1}(x)\\
&=\sqrt{n+1}\psi_{n+1}(x)
\end{align*}
$\mathbf{Exercice\hspace{1mm} 5:}$ Prouver en utilisant la r\`egle de Leibniz, la repr\'esentation
de Rodrigues pour les polyn\^omes de Laguerre associ\`es
\begin{align*}
L_n^{k}(x)=\frac{1}{n!}x^{-k}e^{x}\frac{d^{n}}{dx^{n}}\left(e^{-x}x^{n+k}\right)
\end{align*}
$\mathbf{Solution}$\\
En utilisant le r\`egle de Leibniz
\begin{align*}
\frac{d^{n}}{dx^{n}}(uv)=\sum_{r=0}^{n}\frac{n!}{r!(n-r)!}\frac{d^{n-r}u}{dx^{n-r}}\frac{d^{r}v}{dx^{r}}
\end{align*}
on a
\begin{align*}
\frac{1}{n!}x^{-k}e^{x}\frac{d^{n}}{dx^{n}}\left(x^{n+k}e^{-x}\right)=\frac{1}{n!}x^{-k}e^{x}\sum_{r=0}^{n}\frac{n!}{r!(n-r)!}\frac{d^{n-r}x^{n+k}}{dx^{n-r}}\frac{d^{r}e^{-x}}{dx^{r}}
\end{align*}
Mais
\begin{align*}
\frac{d^{p}}{dx^{p}}x^{q}&=q(q-1)\ldots(q-p+1)x^{q-p}\\
&=\frac{q!}{(q-p)!}x^{q-p}
\end{align*}
donc on a
\begin{align*}
\frac{1}{n!}x^{-k}e^{x}\frac{d^{n}}{dx^{n}}\left(x^{n+k}e^{-x}\right)&=\frac{1}{n!}x^{-k}e^{x}\sum_{r=0}^{n}\frac{n!}{r!(n-r)!}\frac{(n+k)!}{(k+r)!}x^{k+r}(-1)^{r}e^{-x}\\
&=\sum_{r=0}^{n}(-1)^{r}\frac{(n+k)!}{r!(n-r)!(k+r)!}x^{r}\\
&=L_n^{k}(x)
\end{align*}
\chapter{Les fonctions hyperg\'eom\'etriques}
\label{chap: hyper_fct}
Le terme s\'eries hyperg\'eometriques a \'et\'e introduit par  John Wallis en 1656, et elles ont \'et\'e \'etudi\'ees syst\'ematiquement pour la premi\`ere fois par Gauss en 1812. Les fonctions hyperg\'eom\'etriques sont introduites comme g\'en\'eralisation de la notion de s\'eries g\'eom\'etriques. L'utilit\'e des fonctions hyperg\'eom\'etriques r\'eside dans le fait que presque toutes les fonctions sp\'eciales peuvent
\^etre exprim\'ees en termes de fonctions hyperg\'eom\'etriques. Les fonctions hypergom\'etriques ont \'et\'e utilis\'ees dans une large gamme de probl\`emes en mechanique classique et quantique, en ing\'enierie et en math\'ematiques appliqu\'ees. En physique, ils sont tr\`es utiles dans les probl\`emes de forces centrales, par exemple dans l'\'etude de l'atome d'hydrog\`ene et de l'oscillateur harmonique en m\'ecanique quantique. Leur int\'er\^et pour les math\'ematiques r\'eside dans le fait que de nombreuses \'equations (bien connues) aux d\'eriv\'ees partielles peuvent \^etre r\'eduites à l'\'equation hyperg\'eom\'etrique de Gauss par s\'eparation des variables.

\section{Fonction hyperg\'eometrique de Gauss}
La fonction hyperg\'eometrique de Gauss est la g\'en\'eralization des s\'eries g\'eom\'etriques suivantes
\begin{align*}
(1-z)^{-1}&=\sum_{r=0}^{\infty}z^{r},\\
(1-z)^{-\alpha}&=\sum_{n=0}^{\infty}\frac{\alpha(\alpha+1)(\alpha+2)\ldots (\alpha+r-1)}{r!}z^{r}\\
&=\sum_{n=0}^{\infty}\frac{(\alpha)_r}{r!}z^{r}
\end{align*}
o\`u $(\alpha)_r$ est appell\'e l'indice de Pochhammer et il est d\'efinit par
\begin{align}\label{eq:6_1}
(\alpha)_r&=\alpha(\alpha+1)\ldots(\alpha+r-1)\\
&=\frac{\Gamma(\alpha+r)}{\Gamma(\alpha)};\hspace{10mm} \mbox{$r$ est un entier positif}\nonumber
\end{align}
Une g\'en\'eralization des s\'eries pr\'ec\'edentes est donn\'ee par
\begin{align}\label{eq:6_2}
{_2F_1}(\alpha,\beta;\gamma;x)=\sum_{r=0}^{\infty}\frac{(\alpha)_r(\beta)_r}{(\gamma)_r}\frac{x^{r}}{r!}
\end{align}
o\`u $_2F_1$ est appel\'ee la fonction hyperg\'eom\'etrique de Gauss. D'apr\'es la d\'efinition (\ref{eq:6_2}), cette fonction converge si $|x|<1$ et diverge si $|x|>1$. Pour $x=1$ la serie (\ref{eq:6_2}) converge si $\gamma>\alpha+\beta$ et si $x=-1$ elle converge si $\gamma>\alpha+\beta-1$.\\
Il d\'ecoule directement de la d\'efinition (\ref{eq:6_2}) la relation de sym\'etrie suivante
\begin{align*}
_2F_1(\alpha,\beta;\gamma;x)={_2}F_1(\beta,\alpha;\gamma;x)
\end{align*}
\subsection{\'Equation hyperg\'eom\'etrique de Gauss}
$_2F_1(\alpha,\beta;\gamma;x)$ est une solution de l'\'equation hyperg\'eom\'etrique (appell\'ee aussi \'equation hyperg\'eom\'etrique de Gauss)
\begin{align}\label{eq:6_3}
x(1-x)\frac{d^{2}y}{dx^{2}} + \left\lbrace \gamma -(\alpha + \beta +1)x\right\rbrace \frac{dy}{dx}-\alpha\beta y=0
\end{align}
\textit{D\'emonstration}\\
La m\'ethode de Frobenius utilis\'ee pour r\'esoudre l'\'equation de Bessel 
peut \^etre utiliser ici pour r\'esoure l'\'equation hypg\'eom\'etrique de Gauss et la solution finale est donn\'ee par la fonction $_2F_1(\alpha,\beta;\gamma;x)$.
\subsection{Relation avec d'autres fonctions sp\'eciales}
\begin{align}\label{eq:6_4}
&(a)\hspace{1.5mm} P_n(x)= {_2F_1}\left(-n,n+1;1:\frac{1-x}{2}\right),\nonumber\\
&(b)\hspace{1.5mm}T_n(x)={_2F_1}\left(-n,n;\frac{1}{2};\frac{1-x}{2}\right)\\
&(c)\hspace{1.5mm} U_n(x)=\sqrt{1-x^2}n\hspace{0.2mm}{_2F_1}\left(-n+1,n+1;\frac{3}{2};\frac{1-x}{2}\right)\nonumber
\end{align}
\textit{D\'emonstration}\\
On d\'emontre seulement la relation (\ref{eq:6_4})(a), pour les autres relations la preuve est similaire. A partir de la d\'efinition (\ref{eq:6_2}) on a
\begin{align*}
_2F_1\left(-n,n+1;1:\frac{1-x}{2}\right)=\sum_{r=0}^{\infty}\frac{(-n)_r(n+1)_r}{(1)_r}\frac{\left[(1-x)/2\right]^{r}}{r!}
\end{align*}
On prend $n$ un entier non n\'egatif car le polyn\^ome de Legendre est d\'efini seulement pour ces valeurs. On a donc
\begin{align*}
(-n)_r&=(-n)(-n+1)(-n+2)\ldots(-n+r-1)\\
&=(-1)^{r}n(n-1)(n-2)\ldots(n-r+1)\\
&=\left\lbrace\begin{array}{ll}
(-1)^{r}\frac{n!}{(n-r)!} & si\hspace{1.5mm} r\leq n\\
0 & si\hspace{1.5mm} r\geq n+1
\end{array}\right.
\end{align*}
on a aussi
\begin{align*}
(n+1)_r&=(n+1)(n+2)\ldots(n+r)\\
&=\frac{(n+1)!}{n!}
\end{align*}
et
\begin{align*}
(1)_r! &=1.2.3\ldots(1+r-1)\\
&=r!
\end{align*}
donc
\begin{align*}
_2F_1\left(-n,n+1;1;\frac{1-x}{2}\right)&=\sum_{r=0}^{\infty}(-1)^{r}\frac{n!}{(n-r)!}\frac{(n+r)!}{n!}\frac{1}{r!}\frac{\left(1-x\right)^{r}}{2^{r}r!}\\
&=\sum_{r=0}^{\infty}\frac{(-1)^{r}}{2^r}\frac{(n+r)!}{(n-r)!(r!)^2}\left(1-x\right)^{r}\\
&=\sum_{r=0}^{\infty}\frac{(n+r)!}{2^{r}(n-r)!(r!)^2}\left(x-1\right)^{r}\\
\end{align*}
Montrons que cette s\'erie est \'egale \`a $P_n(x)$, pour cela on va d\'evelopper le polyn\^ome de Legendre $P_n(x)$ en s\'erie au voisinage de $x=1$. D'apr\`es le th\'eor\`eme de Taylor, on a
\begin{align}\label{eq:6_5}
P_n(x)=\sum_{r=0}^{\infty}P_n^{(r)}(1)\frac{(x-1)^{r}}{r!}
\end{align}
o\`u 
\begin{align*}
P_n^{(r)}(1) = \frac{d^{r}}{dx^{r}}P_n(x)\Big|_{x=1}
\end{align*}
Pour calculer $P_n^{(r)}(1)$ on va d\'eriver $r$ fois l'expression (\ref{eq:5_34}) de la fonction g\'en\'eratrice du polyn\^ome de Legendre, on obtient
\begin{align*}
\sum_{r=0}^{\infty}P_n^{(r)}(1)t^{n}&=\frac{d^{r}}{dx^{r}}(1-2tx+t^2)^{-1/2}\\
&=(-2t)^{r}\left(-\frac{1}{2}\right)\left(-\frac{1}{2}-1\right)\left(-\frac{1}{2}-2\right)\ldots\left(-\frac{1}{2}-r+1\right)\left(1-2tx+t^2\right)^{-\frac{1}{2}-r}\\
&=2^{r}t^{r}\frac{1}{2}\left(\frac{1}{2}+1\right)\left(\frac{1}{2}+2\right)\ldots\left(\frac{1}{2}+r-1\right)\left(1-2tx+t^2\right)^{-\frac{1}{2}-r}\\
&=t^{r}\frac{(2r)!}{2^{r}r!}\left(1-2tx+t^2\right)^{-\frac{1}{2}-r}
\end{align*}
posons $x=1$, on obtient
\begin{align*}
\sum_{n=0}^{\infty}P_n^{(r)}(1)t^{n}&=t^{r}\frac{(2r)!}{2^{r}r!}\left(1-2t+t^2\right)^{-\frac{1}{2}-r}\\
&=t^{r}\frac{(2r)!}{2^{r}r!}(1-t)^{-1-2r}
\end{align*}
En utilisant la formule du bin\^ome de Newton pour $(1-t)^{-1-2r}$, on a
\begin{align*}
\frac{1}{(1-t)^{2r+1}}&=\sum_{s=0}^{\infty}\frac{(2r+s)!}{(2r)!s!}t^{s}
\end{align*}
donc
\begin{align*}
\sum_{n=0}^{\infty}P_n^{(r)}(1)t^{n}&=t^{r}\frac{(2r)!}{2^{r}r!}\sum_{s=0}^{\infty}\frac{(2r+s)!}{(2r)!r!}t^{s}\\
&=\frac{1}{2^{r}r!}\sum_{s=0}^{\infty}\frac{(2r+s)!}{s!}t^{r+s}\\
&=\frac{1}{2^{r}r!}\sum_{n=r}^{\infty}\frac{(n+r)!}{(n-r)!}t^{n}
\end{align*}
Par identification on obtient
\begin{align*}
P_n^{(r)}(1)=\left\lbrace\begin{array}{ll}
\frac{1}{2^{r}r!}\frac{(n+r)!}{(n-r)!} & si \hspace{1.5mm} n\geq r\\
0 & si \hspace{1.5mm} n<r
\end{array}\right.
\end{align*}
on ins\`ere l'\'expression de $P_n^{(r)}(1)$ dans (\ref{eq:6_5}), il en r\'esulte que
\begin{align*}
P_n(x)&=\sum_{r=0}^{n}\frac{1}{2^{r}r!}\frac{(n+r)!}{(n-r)!}\frac{(x-1)^{r}}{r!}\\
&=\sum_{r=0}^{n}\frac{1}{2^{r}(r!)^2}\frac{(n+r)!}{(n-r)!}(x-1)^r\\
&=\hspace{0.5mm} _2F_1\left(-n,n+1;1;\frac{1-x}{2}\right)
\end{align*}
o\`u la derni\`ere ligne r\'esulte de la d\'efinition (\ref{eq:6_2}).
\subsection{Repr\'esentation int\'egrale}
\begin{align}\label{eq:6_6}
_2F_1(\alpha,\beta;\gamma;x)=\frac{\Gamma(\gamma)}{\Gamma(\beta)\Gamma(\gamma-\beta)}\int_{0}^{1}t^{\beta-1}(1-t)^{\gamma-\beta-1}(1-xt)^{-\alpha}dt \hspace{5mm} si\hspace{1mm} \gamma>\beta>0
\end{align}
\textit{D\'emonstration}\\
A partir de la d\'efinition de la fonction $_2F_1$ on a
\begin{align*}
_2F_1(\alpha,\beta;\gamma;x)&=\sum_{r=0}^{\infty}\frac{(\alpha)_r(\beta)_r}{(\gamma)_r}\frac{x^r}{r!}\\
&=\sum_{r=0}^{\infty}\frac{\Gamma(\alpha+r)}{\Gamma(\alpha)}\frac{\Gamma(\beta+r)}{\Gamma(\beta)}\frac{\Gamma(\gamma)}{\Gamma(\gamma+r)}\frac{x^r}{r!}\\
&=\frac{\Gamma(\gamma)}{\Gamma(\alpha)\Gamma(\beta)}\frac{1}{\Gamma(\gamma-\beta)}\sum_{r=0}^{\infty}\Gamma(\alpha+r)\frac{\Gamma(\gamma-\beta)\Gamma(\beta+r)}{\Gamma(\gamma+r)}\frac{x^r}{r!}\\
&=\frac{\Gamma(\gamma)}{\Gamma(\alpha)\Gamma(\beta)}\frac{1}{\Gamma(\gamma-\beta)}\sum_{r=0}^{\infty}\Gamma(\alpha+r)B(\gamma-\beta,\beta+r)\frac{x^r}{r!}\\
&=\frac{\Gamma(\gamma)}{\Gamma(\alpha)\Gamma(\beta)}\frac{1}{\Gamma(\gamma-\beta)}\sum_{r=0}^{\infty}\Gamma(\alpha+r)\int_{0}^{1}(1-t)^{\gamma-\beta-1}t^{\beta+r-1}dt\frac{x^r}{r!}\\
&=\frac{\Gamma(\gamma)}{\Gamma(\beta)\Gamma(\gamma-\beta)}\int_{0}^{1}(1-t)^{\gamma-\beta-1}t^{\beta+r-1}\sum_{r=0}^{\infty}\frac{\Gamma(\alpha+r)}{\Gamma(\alpha)}\frac{x^r}{r!} dt\\
&=\frac{\Gamma(\gamma)}{\Gamma(\beta)\Gamma(\gamma-\beta)}\int_{0}^{1}(1-t)^{\gamma-\beta-1}t^{\beta+r-1}(1-xt)^{-\alpha} dt
\end{align*}
Notons que la derni\`ere \'egalit\'e est une cons\'equence de la formule du bin\^ome de Newton.
\section{Fonction hyperg\'eom\'etrique confluente}
La fonction hyperg\'eom\'etrique confluente (ou fonction de Kummer) est d\'efinie par
\begin{align}\label{eq:6_7}
{_1F_1}(\alpha;\beta;x)=\sum_{r=0}^{\infty}\frac{(\alpha)_r}{(\beta)_r}\frac{x^{r}}{r!}
\end{align}
et elle est la solution de l'\'equation hyperg\'eom\'etrique confluente
\begin{align}\label{eq:6_8}
x^2\frac{d^{2}y}{dx^{2}} + (\beta -x)\frac{dy}{dx}-\alpha y=0
\end{align}
qui peut \^etre r\'esolu par la m\'ethode de Frobenius. La fonction hyperg\'eom\'etrique confluente converge pour toutes les valeurs de $x$.
\subsection{Relation avec d'autres fonctions sp\'eciales}
\begin{align}\label{eq:6_9}
&(a)\hspace{1.5mm} P_n^{m}(x)=\frac{(n+m)!}{(n-m)!}\frac{(1-x^2)^{m/2}}{2^{m}m!}{_1F_1}\left(m-n,m+n+1;m+1;\frac{1-x}{2}\right).\nonumber\\
&(b)\hspace{1.5mm} J_n(x)=\frac{e^{-ix}}{n!}\left(\frac{x}{2}\right)^{n}{_1F_1}\left(n+1;2n+1;2ix\right)\nonumber\\
&(c)\hspace{1.5mm} H_{2n}(x)=(-1)^{n}\frac{(2n)!}{n!}{_1F_1}\left(-n;\frac{1}{2};x^2\right)\\
&(d)\hspace{1.5mm} H_{2n+1}(x)=(-1)^{n}\frac{2(2n+1)!}{n!}x\hspace{0.2mm}{_1F_1}\left(-n;\frac{3}{2};x^2\right)\nonumber\\
&(e) \hspace{1.5mm}L_n(x)={_1F_1}\left(-n;1;x \right)\nonumber\\
&(f)\hspace{1.5mm} L_n^{k}(x)=\frac{\Gamma(n+k+1)}{n!\Gamma(k+1)}{_1F_1}\left(-n;k+1;x\right)\nonumber
\end{align}
\textit{D\'emonstration}\\
La preuve des \'equations Eq.~(\ref{eq:6_9})(a)-(f) est similaire \`a celle de l'Eq.~(\ref{eq:6_4})(a).
\subsection{Repr\'esentation int\'egrale}
\begin{align}\label{eq:6_10}
_1F_1(\alpha;\beta;x)=\frac{\Gamma(\beta)}{\Gamma(\alpha)\Gamma(\beta-\alpha)}\int_{0}^{1}(1-t)^{\beta-\alpha-1}t^{\alpha-1}e^{xt}dt \hspace{5mm} si\hspace{0.5mm}\beta>\alpha>0
\end{align}
\textit{D\'emonstration}\\
La preuve est trivial, on refait les m\^emes \'etapes de la d\'emonstration de la repr\'esentation 
int\'egrale de la fonction hyperg\'eom\'etrique de Gauss.
\section{Fonctions hyperg\'eom\'etriques g\'en\'eralis\'ees}
Les fonctions hyperg\'eometriques g\'en\'eralis\'ees not\'ees $_{m}F_n(\alpha_1,\alpha_2,\ldots,\alpha_m;\beta_1,\beta_2,\ldots,\beta_n;x)$ sont d\'efinies par
\begin{align}\label{eq:6_11}
{_{m}F_n}(\alpha_1,\alpha_2,\ldots,\alpha_m;\beta_1,\beta_2,\ldots,\beta_n;x)=\sum_{r=0}^{\infty}\frac{(\alpha_1)_r(\alpha_2)_r\ldots(\alpha_n)_r}{(\beta_1)_r(\beta_2)_r\ldots(\beta_m)_r}\frac{x^{r}}{r!}
\end{align}
o\`u $m$ et $n$ prennent les valeurs $0, 1, 2, \ldots$. Une autre notation utilis\'ee pour les fonctions hyperg\'eom\'etriques est la suivante
\begin{align*}
{_mF_n}\left[\begin{array}{l}
\alpha_1, \alpha_2, \ldots, \alpha_n;\\
\beta_1,\beta_2,\ldots,\beta_m;
\end{array}\hspace{-2mm}\begin{array}{c}
x
\end{array}\right]
\end{align*}
Si $m=n+1$, le rayon de convergence de la s\'erie ${_mF_n}$ est l'unit\'e. Si $m<n+1$, le rayon de convergence est \'egal \`a l'infini.
\section{Exercices}
$\mathbf{Exercice\hspace{1mm}  1:}$ Montrer que
\begin{align*}
&1. \hspace{2mm}{_2F_1}(\alpha,\beta;\beta;x)=(1-x)^{-\alpha}\\
&2. \hspace{2mm}{_2F_1}(1,1;2;x)=-\frac{\ln(1-x)}{x}\\
&3. \hspace{2mm}{_2F_1}\left(\frac{1}{2},1;\frac{3}{2};x^{2}\right)=\frac{1}{2x}\ln\left(\frac{1+x}{1-x}\right)\\
&4. \hspace{2mm}{_2F_1}\left(\frac{1}{2},1;\frac{3}{2};-x^{2}\right)=\frac{\arctan x}{x}\\
&5. \hspace{2mm}{_2F_1}\left(\frac{1}{2},\frac{1}{2};\frac{3}{2};x^{2}\right)=\frac{\arcsin x}{x}
\end{align*}
$\mathbf{Solutions}$\\
1. \hspace{1.5mm} On a
\begin{align*}
{_2F_1}(\alpha,\beta;\beta;x)&=\sum_{r=0}^{\infty}\frac{(\alpha)_r(\beta)_r}{(\beta)_r}\frac{x^{r}}{r!}\\
&=\sum_{r=0}^{\infty}(\alpha)_r\frac{x^{r}}{r!}\\
&= \sum_{r=0}^{\infty}\frac{\alpha(\alpha+1)\ldots(\alpha+r-1)}{r!}x^{r}\\
&= \sum_{r=0}^{\infty}(-1)^{r}\frac{(-\alpha)(-\alpha-1)\ldots(-\alpha-r+1)}{r!}x^{r}\\
&=\sum_{r=0}^{\infty}(-1)^{r}C_r^{-\alpha}x^{r}\\
&=(1-x)^{-\alpha}
\end{align*}
o\`u $C_{r}^{-\alpha}$ est le coefficient binomial. Il est d\'efini par
\begin{align*}
C_r^{\alpha}=\frac{\alpha(\alpha-1)\ldots(\alpha-r+1)}{r!}
\end{align*}
2. \hspace{1.5mm} On a
\begin{align*}
{_2F_1}(1,1;2;x)&=\sum_{r=0}^{\infty}\frac{(1)_r(1)_r}{(2)_r}\frac{x^{r}}{r!}\\
&=\sum_{r=0}^{\infty}\frac{r!r!}{(r+1)!}\frac{x^{r}}{r!}\\
&=\sum_{r=0}^{\infty}\frac{x^{r}}{r+1}\\
&=\frac{1}{x}\sum_{r=0}^{\infty}\frac{x^{r+1}}{r+1}\\
&=\frac{1}{x}\sum_{r=1}^{\infty}\frac{x^{r}}{r}\\
&=\frac{1}{x}\sum_{r=0}^{\infty}\frac{x^{r+1}}{r+1}\\
&=-\frac{\ln(1-x)}{x}
\end{align*}
3. \hspace{1.5mm} On a
\begin{align*}
{_2F_1}\left(\frac{1}{2},1;\frac{3}{2};x^2\right)=\sum_{r=0}^{\infty}\frac{(\frac{1}{2})_r(1)_r}{(\frac{3}{2})_r}\frac{x^{2r}}{r!}
\end{align*}
Mais
\begin{align*}
\left(\frac{1}{2}\right)_r=\frac{\Gamma\left(\frac{1}{2}+r\right)}{\Gamma\left(\frac{1}{2}\right)}
\end{align*}
et
\begin{align*}
\left(\frac{3}{2}\right)_r&=\frac{\Gamma\left(\frac{3}{2}+r\right)}{\Gamma\left(\frac{3}{2}\right)}\\
&=\frac{\left(\frac{1}{2}+r\right)\Gamma\left(\frac{1}{2}+r\right)}{\left(\frac{1}{2}\right)\Gamma\left(\frac{1}{2}\right)}
\end{align*}
donc
\begin{align*}
{_2F_1}\left(\frac{1}{2},1;\frac{3}{2};x\right)&=\sum_{r=0}^{\infty}\frac{x^{2r}}{2r+1}\\
&=\frac{1}{x}\sum_{r=0}^{\infty}\frac{x^{2r+1}}{2r+1}\\
\end{align*}
D'autre part nous avons
\begin{align*}
\ln(1+x)-\ln(1-x)&=\sum_{r=1}^{\infty}(-1)^{r+1}\frac{x^r}{r} + \sum_{r=1}^{\infty}\frac{x^r}{r}\\
&=2\sum_{r=0}^{\infty}\frac{x^{2r+1}}{2r+1}
\end{align*}
et on a finalement
\begin{align*}
{_2F_1}\left(\frac{1}{2},1;\frac{3}{2};x\right)=\frac{1}{2x}\ln\left(\frac{1+x}{1-x}\right)
\end{align*}
4. \hspace{1.5mm} On a
\begin{align*}
{_2F_1}\left(\frac{1}{2},1;\frac{3}{2};-x^2\right)&=\sum_{r=0}^{\infty}(-1)^{r}\frac{(\frac{1}{2})_r(1)_r}{(\frac{3}{2})_r}\frac{x^{2r}}{r!}\\
&=\sum_{r=0}^{\infty}(-1)^{r}\frac{x^{2r}}{2r+1}\\
&=\frac{1}{x}\sum_{r=0}^{\infty}(-1)^{r}\frac{x^{2r+1}}{2r+1}\\
&=\frac{\arctan x}{x}
\end{align*}
5. \hspace{1.5mm} On a
\begin{align*}
{_2F_1}\left(\frac{1}{2},\frac{1}{2};\frac{3}{2};x^2\right)&=\sum_{r=0}^{\infty}(-1)^{r}\frac{(\frac{1}{2})_r(\frac{1}{2})_r}{(\frac{3}{2})_r}\frac{x^{2r}}{r!}\\
&=\sum_{r=0}^{\infty}\frac{\Gamma\left(\frac{1}{2}+r\right)}{\Gamma\left(\frac{1}{2}\right)}\frac{x^{2r}}{r!}\\
&=\frac{1}{x}\sum_{r=0}^{\infty}\frac{(2r)!\sqrt{\pi}}{r!2^{2r}\sqrt{\pi}}\frac{1}{2r+1}\frac{x^{2r+1}}{r!}\\
&=\frac{1}{x}\sum_{r=0}^{\infty}\frac{(2r)!}{(r!)^{2}4^{r}}\frac{x^{2r+1}}{2r+1}\\
&=\frac{\arcsin x}{x}
\end{align*}
$\mathbf{Exercice\hspace{1mm} 2:}$ Montrer que
\begin{align*}
{_2F_1}(\alpha,\beta;\gamma;1)=\frac{\Gamma(\gamma)\Gamma(\gamma-\alpha-\beta)}{\Gamma(\gamma-\alpha)\Gamma(\gamma-\beta)}
\end{align*}
$\mathbf{Solution}$\\
D'apr\`es la repr\'esentation int\'egrale Eq.~(\ref{eq:6_6}) on a
\begin{align*}
{_2F_1}(\alpha,\beta;\gamma;1)&=\frac{\Gamma(\gamma)}{\Gamma(\beta)\Gamma(\gamma-\beta)}\int_{0}^{\infty}t^{\beta-1}(1-t)^{\gamma-\beta-1}(1-t)^{-\alpha}dt\\
&=\frac{\Gamma(\gamma)}{\Gamma(\beta)\Gamma(\gamma-\beta)}\int_{0}^{\infty}t^{\beta-1}(1-t)^{\gamma-\alpha-\beta-1}dt\\
&=\frac{\Gamma(\gamma)}{\Gamma(\beta)\Gamma(\gamma-\beta)}B(\beta,\gamma-\alpha-\beta)\\
&=\frac{\Gamma(\gamma)}{\Gamma(\beta)\Gamma(\gamma-\beta)}\frac{\Gamma(\beta)\Gamma(\gamma-\alpha-\beta)}{\Gamma(\gamma-\alpha)}\\
&=\frac{\Gamma(\gamma)\Gamma(\gamma-\alpha-\beta)}{\Gamma(\gamma-\alpha)\Gamma(\gamma-\beta)}
\end{align*}
$\mathbf{Exercice\hspace{1mm} 3:}$ Montrer que
\begin{align*}
&1. \hspace{2mm}\frac{d^{n}}{dx^{n}}{_1F_1}(\alpha;\beta;x)=\frac{(\alpha)_n}{(\beta)_n}{_1F_1}(\alpha+n;\beta+n;x)\\
&2. \hspace{2mm}\frac{d^{n}}{dx^{n}}{_2F_1}(\alpha,\beta;\gamma;x)=\frac{(\alpha)_n(\beta)_n}{(\gamma)_n}{_2F_1}(\alpha+n,\beta+n;\gamma +n;x)\\
&3. \hspace{2mm}x{_1F_1}(\alpha+1;\beta+1;x)+\beta{_1F_1}(\alpha;\beta;x)-\beta{_1F_1}(\alpha+1;\beta;x)=0
\end{align*}
$\mathbf{Solutions}$\\
1. \hspace{1.5mm} On a
\begin{align*}
\frac{d^{n}}{dx^{n}}{_1F_1}(\alpha;\beta;x)&=\sum_{r=0}^{\infty}\frac{(\alpha)_r}{(\beta)_r}\frac{1}{r!}\frac{d^{n}}{dx^{n}}x^{r}
\end{align*}
Si $r<n$ on a
\begin{align*}
\frac{d^{n}}{dx^{n}}x^{r}=0
\end{align*}
et si $r\geq n$ on a
\begin{align*}
\frac{d^{n}}{dx^{n}}x^{r}&=r(r-1)(r-2)\ldots(r-n+1)x^{r-n}\\
&=\frac{r!}{(r-n)!}x^{r-n}
\end{align*}
donc
\begin{align*}
\frac{d^{n}}{dx^{n}}{_1F_1}(\alpha;\beta;x)&=\sum_{r=n}^{\infty}\frac{(\alpha)_r}{(\beta)_r}\frac{1}{r!}\frac{r!}{(r-n)!}x^{r-n}
\end{align*}
Posons $s=r-n$ on trouve
\begin{align*}
\frac{d^{n}}{dx^{n}}{_1F_1}(\alpha;\beta;x)=\sum_{s=0}^{\infty}\frac{(\alpha)_{s+n}}{(\beta)_{s+n}}\frac{x^s}{s!}
\end{align*}
Mais
\begin{align*}
(\alpha)_{s+n}&=\frac{\Gamma(\alpha+s+n)}{\Gamma(\alpha)}\\
&=\frac{1}{\Gamma(\alpha)}(\alpha+s+n-1)\Gamma(\alpha+s+n-1)\\
&=\frac{1}{\Gamma(\alpha)}(\alpha+s+n-1)(\alpha+s+n-2)\Gamma(\alpha+s+n-2)\\
&\vdots\\
&=\frac{1}{\Gamma(\alpha)}(\alpha+s+n-1)(\alpha+s+n-2)\ldots(\alpha+n+1)(\alpha+n)\Gamma(\alpha+n)\\
&=\frac{\Gamma(\alpha+n)}{\Gamma(\alpha)}(\alpha+n)(\alpha+n+1)\ldots(\alpha+n+s-2)(\alpha+n+s-1)\\
&=(\alpha)_n(\alpha+n)_s
\end{align*}
Donc il en r\'esulte
\begin{align*}
\frac{d^{n}}{dx^{n}}{_1F_1}(\alpha;\beta;x)&=\sum_{s=0}^{\infty}\frac{(\alpha)_n(\alpha+n)_s}{(\beta)_n(\beta+n)_s}\frac{x^s}{s!}\\
&=\frac{(\alpha)_n}{(\beta)_n}{_1F_1}(\alpha+n;\beta+n;x)
\end{align*}
2. \hspace{1.5mm} De la m\^eme mani\`ere on a
\begin{align*}
\frac{d^{n}}{dx^{n}}{_2F_1}(\alpha,\beta;\gamma;x)&=\sum_{r=0}^{\infty}\frac{(\alpha)_r(\beta)_r}{(\gamma)_r}\frac{1}{r!}\frac{d^{n}}{dx^{n}}x^{r}\\
&=\sum_{r=n}^{\infty}\frac{(\alpha)_r(\beta)_r}{(\gamma)_r}\frac{1}{r!}\frac{r!}{(r-n)!}x^{r-n}\\
&=\sum_{s=0}^{\infty}\frac{(\alpha)_{s+n}(\beta)_{s+n}}{(\gamma)_{s+n}}\frac{x^{s}}{s!}\\
&=\sum_{s=0}^{\infty}\frac{(\alpha)_{n}(\alpha+n)_{s}(\beta)_{n}(\beta+n)_{s}}{(\gamma)_{n}(\gamma+n)_{s}}\frac{x^{s}}{s!}\\
&=\frac{(\alpha)_{n}(\beta)_{n}}{(\gamma)_{n}}{_2F_1}(\alpha+n,\beta+n;\gamma+n;x)
\end{align*}
3. \hspace{1.5mm} On a
\begin{align*}
x{_1F_1}(\alpha+1&;\beta+1;x)+\beta{_1F_1}(\alpha;\beta;x)-\beta{_1F_1}(\alpha+1;\beta;x)\\
&=x\sum_{r=0}^{\infty}\frac{(\alpha+1)_r}{(\beta+1)_r}\frac{x^r}{r!}+\beta\sum_{r=0}^{\infty}\frac{(\alpha)_r}{(\beta)_r}\frac{x^r}{r!}-\beta\sum_{r=0}^{\infty}\frac{(\alpha+1)_r}{(\beta)_r}\frac{x^r}{r!}
\end{align*}
En remplaçant
\begin{align*}
(\alpha+1)_r&=\frac{(\alpha)_{r+1}}{(\alpha)_1}\\
&=\frac{(\alpha)_{r+1}}{\alpha},
\end{align*}
et
\begin{align*}
(\beta+1)_r=\frac{(\beta)_{r+1}}{\alpha}
\end{align*}
dans le premier terme et
\begin{align*}
(\alpha+1)_r&=\frac{(\alpha)_{r+1}}{\alpha}\\
&=\frac{(\alpha)_{r}(\alpha+r)}{\alpha}\\
&=(\alpha)_{r}+\frac{r}{\alpha}(\alpha)_{r}
\end{align*}
dans le dernier terme on obtient
\begin{align*}
x{_1F_1}(\alpha+1&;\beta+1;x)+\beta{_1F_1}(\alpha;\beta;x)-\beta{_1F_1}(\alpha+1;\beta;x)\\
&=\frac{\beta}{\alpha}\sum_{r=0}^{\infty}\frac{(\alpha)_{r+1}}{(\beta)_{r+1}}\frac{x^{r+1}}{r!}+\beta\sum_{r=0}^{\infty}\frac{(\alpha)_r}{(\beta)_r}\frac{x^r}{r!}-\beta\sum_{r=0}^{\infty}\frac{(\alpha)_r}{(\beta)_r}\left(1+\frac{r}{\alpha}\right)\frac{x^r}{r!}\\
&=\frac{\beta}{\alpha}\sum_{r=0}^{\infty}\frac{(\alpha)_{r}}{(\beta)_{r}}\frac{rx^r}{r!}+\beta\sum_{r=0}^{\infty}\frac{(\alpha)_r}{(\beta)_r}\frac{x^r}{r!}-\beta\sum_{r=0}^{\infty}\frac{(\alpha)_r}{(\beta)_r}\frac{x^r}{r!}-\beta\sum_{r=0}^{\infty}\frac{(\alpha)_r}{(\beta)_r}\frac{rx^r}{r!}\\
&=0
\end{align*}
$\mathbf{Exercice\hspace{1mm}  4:}$ Montrer que la fonction de Bessel $J_n(x)$ est donn\'ee par
\begin{align*}
J_n(x)=\frac{(x/2)^{n}}{\Gamma(n+1)}e^{-ix}{_1F_1}\left(n+\frac{1}{2};2n+1;2ix\right)
\end{align*}
$\mathbf{Solution}$\\
En utilisant la repr\'esentation int\'egrale de la fonction hyperg\'eom\'etrique confluante Eq.~(\ref{eq:6_10}) on a
\begin{align*}
{_1F_1}\left(n+\frac{1}{2};2n+1;2ix\right)=\frac{\Gamma(2n+1)}{\Gamma(n+1/2)\Gamma(n+1/2)}\int_{0}^{1}(1-t)^{n-1/2}t^{n-1/2}e^{2ixt}dt
\end{align*}
En plus on a
\begin{align*}
\Gamma(2n+1)=\frac{(2n)!}{n!}\Gamma(n+1)
\end{align*}
et
\begin{align*}
\Gamma\left(n+\frac{1}{2}\right)=\frac{(2n)!\sqrt{\pi}}{n!2^{2n}}
\end{align*}
donc
\begin{align*}
\frac{(x/2)^{n}}{\Gamma(n+1)}e^{-ix}{_1F_1}\left(n+\frac{1}{2};2n+1;2ix\right)&=\left(\frac{x}{2}\right)^n\frac{2^{2n}}{\sqrt{\pi}\Gamma(1+1/2)}\int_{0}^{1}(1-t)^{n-1/2}t^{n-1/2}e^{ix(2t-1)}dt
\end{align*}
Posons $u=2t-1$ on obtient
\begin{align*}
\frac{(x/2)^{n}}{\Gamma(n+1)}e^{ix}{_1F_1}\left(n+\frac{1}{2};2n+1;2ix\right)&=\frac{2^{n}x^n}{\sqrt{\pi}\Gamma(1+1/2)}\int_{-1}^{1}\left(\frac{1-u}{2}\right)^{n-1/2}\left(\frac{1+u}{2}\right)^{n-1/2}e^{ixu}\frac{du}{2}\\
&=\frac{1}{\sqrt{\pi}\Gamma(1+1/2)}\left(\frac{x}{2}\right)^n\int_{-1}^{1}\left(1-u^2\right)^{n-1/2}e^{ixu}du\\
&=J_n(x)
\end{align*}

\end{document}